\documentclass[12pt]{book}
\usepackage{here,amsfonts,amssymb,pstricks,pst-node,
pst-plot,pst-coil,epsfig,graphicx,makeidx}
\usepackage{amsmath}
\makeindex

\def\reals{\mathbb{R}}
\def\complex{\mathbb{C}}
\def\integs{\mathbb{Z}}
\def\natnums{\mathbb{N}}
\def\eucreal{\mathbb{E}}
\def\projr#1{{\mathbb{P}}^{#1}}
\def\affs{\s{E}}
\def\mapdef#1#2#3{#1\co #2\rightarrow #3}
\def\chull#1{\s{C}(#1)}
\def\chullb#1{\mathrm{conv}(#1)}
\def\vector#1{\fakoverrightarrow{#1}}
\def\affvec#1{\vector{#1}}
\def\famil#1#2#3{(#1_{#2})_{#2\in #3}}
\def\lincombin#1#2#3#4{\sum_{#3\in #4} #1_{#3} #2_{#3}}
\def\Ker{\mathrm{Ker}\,}
\def\mdeg{m}
\def\ndeg{n}
\def\affreal{\mathbb{A}}
\def\libvecbo#1#2{{\bf #1#2}}
\def\novect#1{#1}
\def\interio#1{\buildrel \circ\over #1}
\def\dimm{\mathrm{dim}}
\def\norme#1{\left\|#1\right\|}
\def\smnorme#1{\|#1\|}
\def\adher#1{\overline{#1}}
\def\fr#1{\mathrm{Fr}\>#1}
\def\transpos#1{#1^{\top}}
\def\rprospac#1{{\mathbb{RP}}^{#1}}
\def\dInt{\mathrm{Int}\,}

\def\id{\mathrm{id}}

\def\buildrell#1\over#2{\mathrel{\mathop{\null#1}\limits_{#2}}}

\def\dBd{\partial}
\def\convclo{\mathrm{conv}}

\newbox\sample

\newtheorem{thm}{Theorem}[chapter]
\newtheorem{lemma}[thm]{Lemma}
\newtheorem{cor}[thm]{Corollary}
\newtheorem{prop}[thm]{Proposition}

\newtheorem{defin}{Definition}[chapter]

\newbox\sample
\newif\ifproofmode
\newif\ifsymindex
\global\symindexfalse
\newwrite\inx
\def\indsyma#1#2{\ifproofmode\marginpar{$\scriptstyle#1$}\fi%
\ifx#2\empty\write\inx{$\noexpand#1$,\space\thepage}%
\write\inx{\string\newline}\else%
\write\inx{$\noexpand#1$,\space#2,\space\thepage}%
\write\inx{\string\newline}\fi\ignorespaces}%
\def\indsym#1#2{\ifsymindex%
\ifproofmode\marginpar{$\scriptstyle#1$}\fi%
\ifx#2\empty\write\inx{\string\item \space$\noexpand#1$,\space\thepage}%
\else%
\write\inx{\string\item \space$\noexpand#1$,\space#2,\space\thepage}%
\fi\ignorespaces\fi}%
\def\nsindex#1{\index{#1}\ignorespaces}

\newskip\dangerskipb
\newskip\dangerskip
\def\hang{\hangindent\dangerskip}

\def\s#1{{\cal #1}}

\def\lag{\left\langle}
\def\rag{\right\rangle}
\def\blag{\bigl\langle}
\def\brag{\bigr\rangle}

\def\proof{\noindent{\it Proof\/}.\enspace}

\def\endproof{\bigskip}

\def\remark{\bigskip\noindent{\bf Remark:}\enspace}
\def\endremark{\bigskip}
\def\remarks{\bigskip\noindent{\bf Remarks:}\enspace}

\font\manual=manfnt at 12pt
\def\danbend{{\manual\char127}}
\def\datanger{\medbreak\begingroup\clubpenalty=10000
 \def\par{\endgraf\endgroup\medbreak} \noindent\hang\hangafter=-2
 \hbox to0pt{\hskip-3.5pc\danbend\hfill}}
\outer\def\danger{\datanger}%
\def\ddatanger{\medbreak\begingroup\clubpenalty=10000
 \def\par{\endgraf\endgroup\medbreak} \noindent\hang\hangafter=-2
 \hbox to0pt{\hskip-3.5pc\danbend\kern1pt%
\danbend\hfill}}

\def\dobdownarrow{\mathop{\vbox{\kern2pt \hbox{$\Big\downarrow$}\kern-16.5pt
                          \nointerlineskip\hbox{$\Big\downarrow$}}}}

\def\lrightarrow{\hbox to 25pt{\rightarrowfill}}

\def\supexp{exp(m,n,p)=m^{m^{m^{\cdot^{\cdot^{\cdot^{m^{p}}}}}}}
\vbox{\hbox{$\Big\}\scriptstyle n$}\kern0pt}}

\def\supexpo#1#2#3{#1^{#1^{\cdot^{\cdot^{\cdot^{#1^{#2}}}}}}
\vbox{\hbox{$\Big\}\scriptstyle #3$}\kern0pt}}

\def\sqr#1#2{{\vcenter{\hrule height .#2pt
         \hbox{\vrule width.#2pt height#1pt \kern#1pt
             \vrule width.#2pt}
         \hrule height.#2pt}}}

\def\bigsquare{\mathchoice\sqr76\sqr76\sqr{2.1}3\sqr{1.5}3}

\def\lag{\langle}
\def\rag{\rangle}

\def\co{\colon}

\newskip\bogcentering \bogcentering= 0pt plus 1000pt minus 1000pt 

\def\matth{\mathsurround=0pt}
\def\fakrightarrowfill{$\matth \mathord- \mkern-6mu
  \cleaders\hbox{$\mkern-2mu \mathord- \mkern-2mu$}\hfill
 \mkern-6mu \mathord\rightarrow$}

\def\fakoverrightarrow#1{\vbox{\ialign{##\crcr
  \fakrightarrowfill\crcr\noalign{\kern-1pt\nointerlineskip}
 $\hfil\displaystyle{#1}\hfil$\crcr}}}

\def\cases#1{\left\{\,\vcenter{\normalbaselines\matth
  \ialign{$##\hfil$&\quad##\hfil\crcr#1\crcr}}\right.}

\newif\ifdtatp

\def\displaty{%
\global \dtatptrue \openup \jot \matth \everycr{\noalign{\ifdtatp \global 
\dtatpfalse \vskip -\lineskiplimit \vskip \normallineskiplimit \else 
\penalty \interdisplaylinepenalty \fi }}}

\def\displaylignes#1{\displaty
   \halign{\hbox to\displaywidth{$\displaystyle##$}\crcr
   #1\crcr}}
\def\eqaligneno#1{\displaty \tabskip=\bogcentering
 \halign to\displaywidth{\hfil$\displaystyle{##}$\tabskip=0pt
 &$\displaystyle{{}##}$\hfil\tabskip=\bogcentering
 &\llap{$##$}\tabskip=0pt\crcr
 #1\crcr}}
\def\leqaligneno#1{\displaty \tabskip=\bogcentering
 \halign to\displaywidth{\hfil$\displaystyle{##}$\tabskip=0pt
 &$\displaystyle{{}##}$\hfil\tabskip=\bogcentering
 &\kern-\displaywidth\rlap{$##$}\tabskip=\displaywidthpt\crcr
 #1\crcr}}
\def\ligne{\hbox to\hsize}
\newdimen\nouvpagewidth
\newdimen\offwidth
\newdimen\lawidthoui
\nouvpagewidth=\lawidthoui
\def\kboxit#1{\vbox{\hrule\hbox{\vrule\kern3pt
              \vbox{\kern3pt#1\kern3pt}\kern3pt\vrule}\hrule}}

\def\kboxitb#1{\vbox{\hrule\hbox{\vrule\kern3pt
              \vbox{\kern3pt#1\kern3pt}\kern3pt\vrule}\hrule}}

\def\laboxaround#1{
\aboxaround{\hbox to\hsize{\hfill\box2\hfill}}{#1}
}

\def\boxar#1#2{
\aboxaround{\hbox to\hsize{\hfill#1\hfill}}{#2}
}

\def\aboxaround#1#2{
\setbox4=\vbox{\hsize #2\noindent\strut#1\strut}
\kboxitb{\box4}}

\def\kframeit#1{\vbox{\hrule\hbox{\vrule\kern5pt
              \vbox{\kern5pt#1\kern5pt}\kern5pt\vrule}\hrule}}

\newskip\savnormalbaselineskip
\newskip\savnormallineskip
\newdimen\savnormallineskiplimit

\def\res{\upharpoonright}
\input xy
\xyoption{all}
\title{Notes on Convex Sets, Polytopes, Polyhedra\\
Combinatorial Topology, Voronoi Diagrams and
Delaunay Triangulations}
\author{Jean Gallier\\
Department of Computer and Information Science\\
University of Pennsylvania\\
Philadelphia, PA 19104, USA\\
e-mail: {\tt jean@cis.upenn.edu}
\ \\
}
\begin{document}
\maketitle
\ \vfill\eject
\begin{center}
{\large \bf
Notes on Convex Sets, Polytopes, Polyhedra
Combinatorial Topology, Voronoi Diagrams and
Delaunay Triangulations}

\vspace{1cm}
Jean Gallier
\end{center}

\vspace{2cm}

\noindent
{\bf Abstract:}
Some basic mathematical tools such as
convex sets, polytopes and combinatorial topology,
are used quite heavily in applied fields such as
geometric modeling, meshing, computer vision, medical imaging and robotics.
This report may be viewed as a tutorial and a set of notes
on convex sets, polytopes, polyhedra, combinatorial topology,
Voronoi Diagrams and Delaunay Triangulations. 
It is intended for a broad audience
of mathematically inclined readers.

One of my (selfish!) motivations in writing these notes
was to understand the concept of {\it shelling\/}
and how it is used to prove the famous Euler-Poincar\'e
formula  (Poincar\'e, 1899) and the more recent 
{\it Upper Bound Theorem\/} (McMullen, 1970) for polytopes.
Another of my motivations was to give a ``correct''
account of Delaunay triangulations and Voronoi diagrams in terms of
(direct and inverse) stereographic projections
onto a sphere and prove rigorously that
the projective map that sends the (projective) sphere
to the (projective) paraboloid works
correctly, that is, maps the 
Delaunay triangulation and Voronoi diagram w.r.t. the lifting onto the sphere
to the Delaunay diagram and Voronoi diagrams w.r.t. the traditional lifting onto 
the paraboloid. 
Here, the problem is that this map
is only well defined (total) in projective space
and we are forced to define the notion of convex polyhedron
in projective space.

\medskip
It turns out that in order to achieve (even partially)
the above goals, I found that it was necessary to include
quite a bit of background material on convex sets,
polytopes, polyhedra and projective spaces.
I have included
a rather thorough treatment of the equivalence of
$\s{V}$-polytopes and $\s{H}$-polytopes and also of
the equivalence of $\s{V}$-polyhedra and $\s{H}$-polyhedra,
which is a bit harder. 
In particular, the {\it Fourier-Motzkin elimination\/} method
(a version of Gaussian elimination for inequalities) is
discussed in some detail.
I also had to include some material
on projective spaces, projective maps  and polar duality w.r.t.
a nondegenerate quadric in order to define
a suitable notion of ``projective polyhedron''
based on cones. To the best of our knowledge, this 
notion of projective polyhedron is new.
We also believe that some of our proofs establishing
the equivalence of $\s{V}$-polyhedra and $\s{H}$-polyhedra
are new.


\tableofcontents
\vfill\eject
\chapter[Introduction]
{Introduction}
\label{chap0}
\section{Motivations and Goals}
\label{sec0a}
For the past eight years or so I have been teaching
a graduate course whose main goal is to expose students to 
some fundamental concepts of geometry, keeping in
mind  their applications
to geometric modeling, meshing, computer vision, medical imaging,
robotics, {\it etc.} 
The audience has been primarily computer science students
but  a fair number of mathematics students and also students
from other engineering disciplines (such as Electrical, Systems, 
Mechanical and Bioengineering) have  been attending
my classes. In the  past three years, I have been focusing
more on convexity, polytopes and
combinatorial topology, as concepts and tools from these areas
have been used increasingly in meshing and also
in computational biology and medical imaging.
One of my (selfish!) motivations
was to understand the concept of {\it shelling\/}
and how it is used to prove the famous Euler-Poincar\'e
formula  (Poincar\'e, 1899) and the more recent 
{\it Upper Bound Theorem\/} (McMullen, 1970) for polytopes.
Another of my motivations was to give a ``correct''
account of Delaunay triangulations and Voronoi diagrams in terms of
(direct and inverse) stereographic projections
onto a sphere and prove rigorously that
the projective map that sends the (projective) sphere
to the (projective) paraboloid works
correctly, that is, maps the 
Delaunay triangulation and  Voronoi diagram w.r.t. the lifting onto the sphere
to the Delaunay triangulation and Voronoi diagram w.r.t. the lifting onto 
the paraboloid. Moreover, the projections of these polyhedra
onto the hyperplane $x_{d+1} = 0$, 
from the sphere or from the paraboloid, are identical.
Here, the problem is that this map
is only well defined (total) in projective space
and we are forced to define the notion of convex polyhedron
in projective space.

\medskip
It turns out that in order to achieve (even partially)
the above goals, I found that it was necessary to include
quite a bit of background material on convex sets,
polytopes, polyhedra and projective spaces.
I have included
a rather thorough treatment of the equivalence of
$\s{V}$-polytopes and $\s{H}$-polytopes and also of
the equivalence of $\s{V}$-polyhedra and $\s{H}$-polyhedra,
which is a bit harder. 
In particular, the {\it Fourier-Motzkin elimination\/} method
(a version of Gaussian elimination for inequalities) is
discussed in some detail.
I also had to include some material
on projective spaces, projective maps  and polar duality w.r.t.
a nondegenerate quadric, in order to define
a suitable notion of ``projective polyhedron''
based on cones. This notion turned out to be indispensible
to give a correct treatment of the Delaunay and Voronoi
complexes using  inverse stereographic projection onto
a sphere and to prove rigorously  that the
well known projective map between the sphere and the paraboloid
maps the  Delaunay triangulation and the Voronoi diagram w.r.t. the sphere to
the more traditional  Delaunay triangulation and Voronoi diagram w.r.t. the 
paraboloid. To the best of our knowledge, this notion of projective polyhedron is new.
We also believe that some of our proofs establishing
the equivalence of $\s{V}$-polyhedra and $\s{H}$-polyhedra
are new.

\medskip
Chapter \ref{chapter3} on combinatorial topology
is hardly original. However, most texts covering this
material are either old fashion or too advanced.
Yet, this material is used extensively
in meshing and geometric modeling.
We tried to give a rather intuitive yet rigorous
exposition. We decided to introduce the terminology
{\it combinatorial manifold\/}, a notion usually
referred to as {\it triangulated manifold\/}.

\medskip
A recurring theme in these notes is the process of
``conification'' (algebraically, ``homogenization''),
that is, forming a cone from some geometric
object. Indeed, ``conification'' turns an object into
a set of lines, and since lines play the role of points in 
projective geometry, ``conification'' (``homogenization'')
is the way to
``projectivize'' geometric affine objects. Then, these (affine) 
objects appear as ``conic sections'' of cones by  hyperplanes,
just the  way the classical conics (ellipse, hyperbola,
parabola) appear as conic sections. 

\medskip
It is  worth warning our readers that
convexity and polytope theory is deceptively simple.
This is a subject where most intuitive propositions
fail as soon as the dimension of the space is greater than $3$
(definitely $4$), because our human intuition is not very
good in dimension greater than $3$. Furthermore, rigorous proofs
of seemingly very simple facts are often quite complicated
and may require sophisticated tools (for example, shellings,
for a correct proof of the Euler-Poincar\'e formula).
Nevertheless, readers are urged to strenghten their
geometric intuition; they should just be very vigilant!
This is another case where Tate's famous saying is more
than pertinent:
``Reason geometrically,  prove algebraically.''

\medskip
At first, 
these notes were meant as a complement to Chapter 3 
(Properties of Convex Sets: A Glimpse) of
my book ({\sl Geometric Methods and Applications\/},
\cite{Gallbook2}). However, they turn out to cover much more material.
For the reader's convenience, I have included  Chapter 3
of my book as part of Chapter
\ref{chap1b} of these notes. I also assume some familiarity
with affine geometry. The reader may wish to review
the basics of affine geometry. These can be found  in any standard 
geometry text
(Chapter 2 of Gallier \cite{Gallbook2} covers more than
needed for these notes).

\medskip
Most of the material on convex sets is 
taken  from
Berger \cite{Berger90b} ({\sl Geometry II\/}).
Other relevant sources include Ziegler \cite{Ziegler97}, 
Gr\"unbaum \cite{Grunbaum}
Valentine \cite{Valentine}, Barvinok \cite{Barvinok},
Rockafellar \cite{Rockafellar},
Bourbaki (Topological Vector Spaces)
\cite{BourbakiEVT} and Lax \cite{Lax},
the last four dealing with affine spaces
of infinite dimension.
As to  polytopes and polyhedra, ``the'' classic reference is
Gr\"unbaum \cite{Grunbaum}. Other good references include
Ziegler \cite{Ziegler97}, Ewald  \cite{Ewald},  Cromwell
\cite{Cromwell} and Thomas \cite{Thomas}.

\medskip
The recent book by Thomas contains an excellent
and easy going presentation of polytope theory.
This book also gives an introduction to the theory of
triangulations of point configurations, including the definition
of secondary polytopes and state polytopes, which
happen to play a role in certain areas of biology.
For this, a quick but very efficient presentation of
Gr\"obner bases  is provided. We highly recommend
Thomas's book \cite{Thomas} as further reading.
It is also an excellent preparation for
the more advanced book by Sturmfels \cite{Sturmfels96}.
However, in our opinion, the ``bible'' on polytope theory is
without any contest,
Ziegler \cite{Ziegler97}, a masterly and beautiful piece of mathematics.
In fact, our Chapter \ref{chap4}
is heavily inspired by Chapter 8 of Ziegler.
However, the pace of Ziegler's book is quite brisk and we 
hope that our more pedestrian account will inspire readers
to go back and read the masters.

\medskip
In a not too distant future, I would like to
write about constrained Delaunay triangulations,
a formidable topic, please be patient! 

\medskip
I wish to thank Marcelo Siqueira for catching many typos and mistakes 
and for his many helpful suggestions
regarding the presentation. At least a third of this manuscript
was written while I was on sabbatical at INRIA, Sophia Antipolis,
in the Asclepios Project. My deepest thanks to Nicholas Ayache
and his colleagues (especially Xavier Pennec and Herv\'e Delingette)
for inviting me to spend a wonderful and very productive year
and for making me feel perfectly at home within the Asclepios Project.

\chapter{Basic Properties of Convex Sets} 
\label{chap1b}
\section{Convex Sets}
\label{secconvex1}
Convex sets play a very important role in geometry.
In this chapter we state and prove some of the ``classics''
of convex affine geometry: Carath\'eodory's theorem, Radon's theorem,
and Helly's theorem. These theorems share the property
that they are easy to state, but they are  deep, and their
proof, although rather short, requires a lot of creativity.

\medskip
Given an affine space $E$, recall that a 
subset $V$ of $E$ is {\it convex\/} if \nsindex{convex!set}
for any two points $a, b\in V$, we have $c\in V$ for every point
$c = (1 -\lambda) a + \lambda b$, with $0\leq \lambda \leq 1$
($\lambda \in\reals$). Given any two points $a, b$,
the notation $[a, b]$ is often used to denote
\indsym{[a, b]}{line segment from $a$ to $b$}
the line segment between $a$ and $b$, that is,
\[[a, b] = \{c\in E\ |\ c = (1 - \lambda)a + \lambda b,\, 0\leq \lambda\leq 1\},\]
and thus a set $V$ is convex if $[a, b]\subseteq V$
for any two points $a, b\in V$ ($a = b$ is allowed).
The empty set is trivially convex, every one-point set $\{a\}$ is convex,
and the entire affine space $E$ is of course convex.

\medskip
It is obvious that the intersection of 
any family (finite or infinite) of convex sets is convex.
Then, given any (nonempty) subset $S$ of $E$, there is a smallest
convex set containing $S$ denoted by $\chull{S}$  or $\mathrm{conv}(S)$
and called
the {\it convex hull of $S$\/}
(namely, the intersection of all convex sets containing $S$).
\index{convex hull!definition}\indsym{\chull{S}}{convex hull of $S$}%
\index{convex hull}%
The {\it affine hull\/} of a subset, $S$, of $E$ is 
the smallest affine set containing $S$ and it will be denoted
by $\lag S\rag$ or $\mathrm{aff}(S)$.

\medskip
A good understanding of what $\chull{S}$ is, and good
methods for computing it, are essential.
First, we have the following simple but crucial lemma:

\begin{lemma}
\label{convexlem1}
Given an affine space $\blag E, \affvec{E}, +\brag$, for any family
$\famil{a}{i}{I}$ of points in $E$, the set $V$ of convex combinations
$\lincombin{\lambda}{a}{i}{I}$ (where 
$\sum_{i\in I} \lambda_i = 1$ and $\lambda_i \geq 0$) is the convex
hull of $\famil{a}{i}{I}$.
\end{lemma}

\medskip
\proof If $\famil{a}{i}{I}$ is empty, then $V = \emptyset$, because
of the condition $\sum_{i\in I} \lambda_i = 1$.
As in the case of affine combinations, it is easily
shown by induction that any convex combination
can be obtained by computing convex combinations
of two points at a time. As a consequence,
if $\famil{a}{i}{I}$ is nonempty, then 
the smallest convex subspace containing $\famil{a}{i}{I}$
must contain the set $V$ of all convex combinations
$\lincombin{\lambda}{a}{i}{I}$. 
Thus, it is enough to show
that $V$ is closed under convex  combinations, which
is immediately verified.
$\bigsquare$
\endproof

\medskip
In view of Lemma \ref{convexlem1},
it is obvious that any affine subspace
of $E$ is convex. Convex sets also arise
in terms of hyperplanes.
\nsindex{hyperplane}
Given a hyperplane $H$, if $\mapdef{f}{E}{\reals}$ is any nonconstant
affine form defining $H$ (i.e., $H = \Ker{f}$),
we can define the two subsets 
\[H_{+}(f) = \{a\in E\ | \ f(a) \geq 0\}\quad
\hbox{and}\quad H_{-}(f) = \{a\in E\ | \ f(a) \leq 0\},\]
called {\it (closed) half-spaces associated with $f$\/}.
\nsindex{closed half-space!definition}%
\nsindex{closed half-space!associated with $f$}%
\indsym{H_{+}(f)}{closed half-space associated with $f$}%
\indsym{H_{-}(f)}{closed half-space associated with $f$}

\medskip
Observe that if $\lambda > 0$, then $H_{+}(\lambda f) = H_{+}(f)$,
but if $\lambda < 0$, then  $H_{+}(\lambda f) = H_{-}(f)$,
and similarly for  $H_{-}(\lambda f)$. However, the set
\[\{H_{+}(f), H_{-}(f)\}\] 
depends only  on the hyperplane $H$,
and the choice of a specific $f$ defining $H$ amounts
to the choice of one of the two half-spaces.
For this reason, we will also say that $H_{+}(f)$ and $H_{-}(f)$
are the {\it closed half-spaces associated with $H$\/}.
\nsindex{closed half-space!associated with $H$}
Clearly, $H_{+}(f) \cup H_{-}(f) = E$ and
$H_{+}(f) \cap H_{-}(f) = H$. It is immediately verified
that  $H_{+}(f)$ and  $H_{-}(f)$ are convex.
Bounded convex sets arising as the intersection of a finite
family of half-spaces associated with hyperplanes
play a major role in convex geometry 
and topology (they are called
{\it convex polytopes\/}).

\medskip
It  is natural to wonder whether Lemma \ref{convexlem1}
can be sharpened
in two directions: (1) Is it possible to have a fixed  bound 
on the number of points involved in  the convex combinations?
(2) Is it necessary to consider convex combinations
of all points, or is it possible to 
consider only a subset with special properties?

\medskip
The answer is yes in both cases.
In case 1, assuming that the affine space $E$ has dimension
$\mdeg$, Carath\'eodory's theorem asserts that it
is enough to consider convex combinations of
$\mdeg + 1$ points. For example, in the plane
$\affreal^2$, the convex hull of a set $S$ of points is
the union of all triangles (interior points included)
with vertices in $S$. In case 2, 
the theorem of Krein and Milman asserts that a convex set 
\nsindex{Krein and Milman's theorem}
that is also compact is the convex hull of its extremal points
(given a convex set $S$, a point $a\in S$ is extremal if
$S - \{a\}$ is also convex,
see Berger \cite{Berger90b} or Lang \cite{Lang96}). 
Next, we prove Carath\'eodory's theorem.
\section{Carath\'eodory's Theorem}
\label{secconvex2}
The proof of  Carath\'eodory's theorem
is really beautiful. It proceeds by contradiction
and uses a minimality argument.
\begin{thm}
\label{caratheo}
Given any affine space $E$ of dimension $\mdeg$,
for any (nonvoid) family 
$S = (a_i)_{i\in L}$ in $E$,
the convex hull $\chull{S}$ of $S$ is equal to 
the set of convex combinations of families
of $\mdeg + 1$ points of $S$.
\end{thm}
\nsindex{Carath\'eodory's theorem}

\medskip
\proof By Lemma \ref{convexlem1}, 
\[\chull{S} = \biggl\{\sum_{i\in I}\lambda_i a_i\ |\ 
a_i \in S,\, \sum_{i\in I}\lambda_i = 1,\,  \lambda_i\geq 0,\,
I\subseteq L,\, \hbox{$I$ finite}\biggr\}.\]
We would like to prove that
\[\chull{S} = \biggl\{\sum_{i\in I}\lambda_i a_i\ |\ 
a_i \in S,\, \sum_{i\in I}\lambda_i = 1,\,  \lambda_i\geq 0,\,
I\subseteq L,\, |I| = \mdeg + 1\biggr\}.\]
We proceed by contradiction.
If the theorem is false, there is some point $b\in \chull{S}$
such that $b$ can be expressed as a convex combination
$b = \sum_{i\in I} \lambda_i a_i$, 
where $I\subseteq L$ is a finite set of cardinality 
$|I| = q$ with $q \geq \mdeg + 2$, and $b$ cannot be expressed
as any convex combination $b = \sum_{j\in J} \mu_j a_j$
of strictly fewer than $q$ points in $S$, that is,
where $|J| <  q$. Such a  point $b\in\chull{S}$ 
is a convex combination
\[b = \lambda_1 a_1 + \cdots + \lambda_q a_q,\]
where $\lambda_1 + \cdots + \lambda_q = 1$
and $\lambda_i > 0$ $(1 \leq i \leq q$). 
We shall prove that $b$ can be written
as a convex combination of $q - 1$ of the $a_i$.
Pick any origin $O$ in $E$. Since
there are $q > \mdeg + 1$ points
$a_1, \ldots, a_q$, these points
are affinely dependent, and by Lemma 2.6.5 
from Gallier \cite{Gallbook2},
there is a family $(\mu_1, \ldots,\mu_q)$
all scalars not all null, such that
$\mu_1 + \cdots + \mu_q = 0$ and
\[\sum_{i = 1}^{q} \mu_i\libvecbo{O}{a_i} = \novect{0}.\]
Consider the set $T\subseteq \reals$ defined by
\[T = \{t\in \reals\ |\ \lambda_i + t\mu_i  \geq 0,\, \mu_i \not= 0,\, 
1\leq i\leq q\}.\]
The set $T$ is nonempty, since it contains $0$.
Since $\sum_{i = 1}^{q}\mu_i = 0$ and the $\mu_i$ are not all null,
there are some $\mu_h, \mu_k$ such that $\mu_h < 0$ and $\mu_k > 0$,
which implies that $T = [\alpha, \beta]$, where
\[\alpha = \max_{1\leq i \leq q}\{-\lambda_i/\mu_i\ |\ \mu_i > 0\}\quad
\hbox{and}\quad
\beta = \min_{1\leq i \leq q}\{-\lambda_i/\mu_i\ |\ \mu_i < 0\}\]
($T$ is the intersection of the closed half-spaces 
$\{t\in \reals\ |\ \lambda_i + t\mu_i \geq 0,\, \mu_i\not= 0\}$). 
Observe that $\alpha < 0 < \beta$, since  $\lambda_i > 0$ for all 
$i = 1, \ldots, q$.

\medskip
We claim that there is some
$j$  $(1\leq j \leq q)$ such that
\[\lambda_j + \alpha\mu_j = 0.\]
Indeed, since
\[
\alpha = \max_{1\leq i \leq q}\{-\lambda_i/\mu_i\ |\ \mu_i > 0\},
\]
as the set on the right hand side is finite, the maximum is achieved
and there is some index $j$ so that
$\alpha = -\lambda_j/\mu_j$.
If  $j$ is some index such that 
$\lambda_j + \alpha\mu_j = 0$, since
$\sum_{i = 1}^{q} \mu_i\libvecbo{O}{a_i} = \novect{0}$, we have
\[\eqaligneno{
b &= \sum_{i = 1}^{q} \lambda_i a_i =
O + \sum_{i = 1}^{q} \lambda_i\libvecbo{O}{a_i} + \novect{0},\cr
&= O + \sum_{i = 1}^{q} \lambda_i\libvecbo{O}{a_i} + 
\alpha\biggl(\sum_{i = 1}^{q} \mu_i \libvecbo{O}{a_i}\biggr),\cr
&= O + \sum_{i = 1}^{q} (\lambda_i + \alpha\mu_i) \libvecbo{O}{a_i},\cr
&= \sum_{i = 1}^{q} (\lambda_i + \alpha\mu_i) a_i,\cr
&= \sum_{i = 1,\, i\not= j}^{q} (\lambda_i + \alpha\mu_i) a_i,\cr
}\]
since $\lambda_j + \alpha\mu_j = 0$.
Since $\sum_{i = 1}^{q} \mu_i = 0$,
$\sum_{i = 1}^{q} \lambda_i = 1$, 
and $\lambda_j + \alpha\mu_j = 0$, we have
\[\sum_{i = 1,\, i\not= j}^{q} \lambda_i + \alpha\mu_i = 1,\]
and since $\lambda_i + \alpha\mu_i \geq  0$ for $i = 1, \ldots, q$, 
the above shows that $b$ can be expressed as a 
convex combination of $q - 1$ points from $S$.
However, this
contradicts the assumption that $b$ cannot be expressed
as a convex combination of strictly fewer than $q$ points from $S$,
and the theorem is proved.
$\bigsquare$
\endproof

\medskip
If $S$ is a finite (of infinite) set of points in
the affine plane $\affreal^2$, Theorem \ref{caratheo}
confirms our intuition that $\chull{S}$ is the union
of triangles (including interior points)
whose vertices belong to $S$. Similarly,
the convex hull of a set $S$ of points in
$\affreal^3$ is the union of tetrahedra
(including interior points)
whose vertices belong to $S$.
We get the feeling that triangulations
play a crucial role, which is of course true!

\medskip
Now that we have given an answer to the first
question posed
at the end of Section \ref{secconvex1}
we give an answer to the second question.

\section{Vertices,  Extremal Points and Krein and Milman's Theorem}
\label{sec3}
First, we define the notions of separation
and of separating hyperplanes.
For this, recall the definition of the 
closed (or open) half--spaces determined by a hyperplane.

\medskip
Given a hyperplane $H$, if $\mapdef{f}{E}{\reals}$ is any nonconstant
affine form defining $H$ (i.e., $H = \Ker{f}$),
we define the {\it closed half-spaces associated with $f$\/} by
$$\eqaligneno{
H_{+}(f) &= \{a\in E\ | \ f(a) \geq 0\},\cr
H_{-}(f) &= \{a\in E\ | \ f(a) \leq 0\}.\cr
}$$

Observe that if $\lambda > 0$, then $H_{+}(\lambda f) = H_{+}(f)$,
but if $\lambda < 0$, then  $H_{+}(\lambda f) = H_{-}(f)$,
and similarly for  $H_{-}(\lambda f)$. 

\medskip
Thus, the set
$\{H_{+}(f), H_{-}(f)\}$ depends only  on the hyperplane, $H$,
and the choice of a specific $f$ defining $H$ amounts
to the choice of one of the two half-spaces.

\medskip
We also define the {\it open half--spaces associated with $f$\/}
as the two sets
$$\eqaligneno{
\interio{H}_{+}\!(f) &= \{a\in E\ | \ f(a) > 0\},\cr
\interio{H}_{-}\!(f) &= \{a\in E\ | \ f(a) < 0\}.\cr
}$$

The set 
$\{\interio{H}_{+}\!(f), \interio{H}_{-}\!(f)\}$ 
only depends on the hyperplane $H$.
Clearly, we have $\interio{H}_{+}\!(f) = H_{+}(f) - H$ and
$\interio{H}_{-}\!(f) = H_{-}(f) - H$.

\begin{defin}
\label{separdef}
{\em
Given an affine space, $X$, and two nonempty subsets,
$A$ and $B$, of $X$, we say that a hyperplane $H$
{\it separates (resp. strictly separates) $A$ and $B$\/}
if $A$ is in one and $B$ is in the other
of the two half--spaces (resp.
open half--spaces) determined by $H$.
}
\end{defin}

\medskip
The special case of separation where
$A$ is convex and $B = \{a\}$, for some point, $a$, in $A$,
is of particular importance.

\begin{defin}
\label{suporthyp}
{\em
Let $X$ be an affine space and let $A$ be any nonempty subset
of $X$. A {\it supporting hyperplane of $A$\/} is any
hyperplane, $H$, containing some point, $a$, of $A$, and separating
$\{a\}$ and $A$. We say that $H$ is a {\it supporting hyperplane of $A$
at $a$\/}.
}
\end{defin}

\medskip
Observe that if $H$ is a supporting hyperplane of $A$ at $a$,
then we must have $a\in \partial A$. Otherwise,
there would be some open ball $B(a, \epsilon)$  of center $a$ 
contained in $A$ and so there would be points of 
$A$ (in $B(a, \epsilon)$) in both half-spaces determined by $H$,
contradicting the fact that
$H$ is a supporting hyperplane of $A$ at $a$.
Furthermore, $H\>\cap \interio{A}\> = \> \emptyset$.

\medskip
One should experiment with various pictures and realize that
supporting hyperplanes at a point may not exist (for example,
if $A$ is not convex),
may not be unique, and may have several distinct supporting
points!

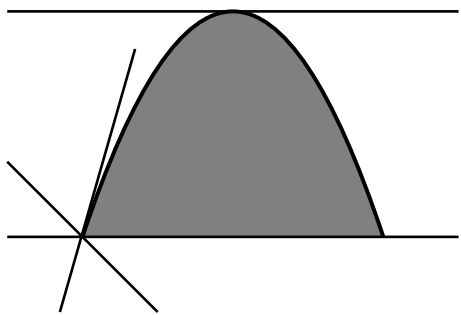
\begin{figure}
  \begin{center}
    \begin{pspicture}(0,0)(6,4.5)
    \parabola[fillstyle=solid,fillcolor=gray,linewidth=1.5pt](1,1)(3,4)
    \pnode(0,4){v1}
    \pnode(6,4){v2}
    \ncline[linewidth=1pt]{v1}{v2}
    \pnode(0,1){v3}
    \pnode(6,1){v4}
    \ncline[linewidth=1pt]{v3}{v4}
    \pnode(0,2){v5}
    \pnode(2,0){v6}
    \ncline[linewidth=1pt]{v5}{v6}
    \pnode(0.7,0){v7}
    \pnode(1.7,3.5){v8}
    \ncline[linewidth=1pt]{v7}{v8}
    \end{pspicture}
  \end{center}
  \caption{Examples of supporting hyperplanes}
  \label{supportfig1}
\end{figure}

\medskip
Next, we need to define various types of boundary
points of closed convex sets.

\begin{defin}
\label{vertexpt}
{\em
Let $X$ be an affine space of dimension $d$.
For any nonempty closed and convex subset, $A$,
of dimension $d$, a point $a\in \partial A$
has {\it order $k(a)$\/} if the intersection
of all the supporting hyperplanes of $A$ at $a$
is an affine subspace of dimension $k(a)$.
We say that $a\in \partial A$ is a {\it vertex\/}
if $k(a) = 0$; we say that $a$ is {\it smooth\/}
if $k(a) = d - 1$, i.e., if the supporting
hyperplane at $a$ is unique.
}
\end{defin}

\medskip
A vertex is a boundary point, $a$, such that
there are $d$ independent supporting hyperplanes
at $a$. 
A $d$-simplex has boundary points of order
$0, 1, \ldots, d - 1$.
The following proposition is shown in
Berger \cite{Berger90b} (Proposition 11.6.2):

\begin{prop}
\label{vertexlem}
The set of vertices of 
a closed and convex subset is countable.
\end{prop}

\medskip
Another important concept is that of an extremal point.

\begin{defin}
\label{extremept}
{\em
Let $X$ be an affine space.
For any nonempty convex subset, $A$,
a point $a\in \partial A$ is {\it extremal\/} (or {\it extreme\/})
if $A \{a\}$ is still convex.
}
\end{defin}

\medskip
It is fairly obvious that a point $a\in \partial A$
is extremal if it does not belong to any closed
nontrivial line segment $[x, y]\subseteq A$ ($x \not= y$).

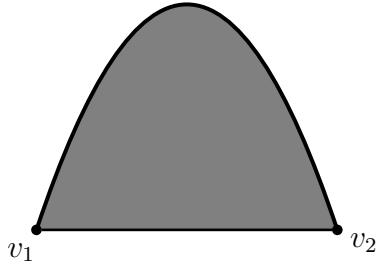
\begin{figure}
  \begin{center}
    \begin{pspicture}(0,0.8)(6,4.5)
    \parabola[fillstyle=solid,fillcolor=gray,linewidth=1.5pt](1,1)(3,4)
    \cnode[fillstyle=solid,fillcolor=black](1,1){2pt}{v1}
    \cnode[fillstyle=solid,fillcolor=black](5,1){2pt}{v2}
    \ncline[linewidth=1pt]{v1}{v2}
    \uput[-120](1,1){$v_1$}
    \uput[-30](5,1){$v_2$}
    \end{pspicture}
  \end{center}
  \caption{Examples of vertices and extreme points}
  \label{vertexfig1}
\end{figure}

\medskip
Observe that a vertex is extremal, but the converse
is false. For example, in Figure \ref{vertexfig1}, all the points on
the arc of parabola, including $v_1$ and $v_2$, are extreme points.
However, only $v_1$ and $v_2$ are vertices.
Also, if $\dimm\> X \geq 3$, the set of extremal points
of a compact convex may not be closed.

\medskip
Actually, it is not at all obvious that a 
nonempty compact convex set possesses extremal points.
In fact, a stronger results holds
(Krein and Milman's theorem). In preparation for the
proof of this important theorem, observe that any
compact (nontrivial) interval of $\affreal^1$ has two 
extremal points, its two endpoints.
We need the following lemma:

\begin{lemma}
\label{comphull}
Let $X$ be an affine space of dimension $n$, and
let $A$ be a nonempty compact and convex set. Then,
$A = \s{C}(\partial A)$, i.e., $A$ is equal to the convex hull of its
boundary.
\end{lemma}

\proof
Pick any $a$ in $A$, and consider any line, $D$, through
$a$. Then, $D\cap A$ is closed and convex. However,
since $A$ is compact, it follows that 
$D\cap A$ is a closed interval $[u, v]$ containing
$a$, and $u, v\in \partial A$. Therefore,
$a\in \s{C}(\partial A)$, as desired.
$\bigsquare$

\medskip
The following important theorem
shows that only
extremal points matter as far as determining
a compact and convex subset from its boundary.
The proof of Theorem \ref{KreinMilman}
makes use of a proposition due to Minkowski
(Proposition \ref{suportprop1}) which will be proved
in Section \ref{sec2}.

\begin{thm} (Krein and Milman, 1940)
\label{KreinMilman}
Let $X$ be an affine space of dimension $n$.
Every compact and convex nonempty subset, $A$, is
equal to the convex hull of its set of extremal points.
\end{thm}

\proof
Denote the set of extremal points of $A$ by
$\mathrm{Extrem}(A)$. We proceed by induction on $d = \dimm\> X$.
When $d = 1$, the convex and compact subset $A$
must be a closed interval $[u, v]$, or a single point.
In either cases, the theorem holds trivially.
Now, assume $d \geq 2$, and assume that the theorem
holds for $d - 1$. It is easily verified that
\[\mathrm{Extrem}(A\cap H) = (\mathrm{Extrem}(A))\cap H,\]
for every supporting hyperplane $H$ of $A$ (such hyperplanes
exist, by Minkowski's proposition (Proposition \ref{suportprop1})). 
Observe that Lemma \ref{comphull} implies that
if we can prove that
\[\partial A \subseteq \s{C}(\mathrm{Extrem}(A)),\]
then, since $A = \s{C}(\partial A)$,  
we will have established that
\[A = \s{C}(\mathrm{Extrem}(A)).\]
Let $a\in \partial A$, and let $H$ be a supporting
hyperplane of $A$ at $a$ (which exists, by Minkowski's proposition).
Now, $A\cap H$ is convex and $H$ has dimension $d - 1$, and by the
induction hypothesis, we have
\[A\cap H = \s{C}(\mathrm{Extrem}(A\cap H)).\]
However, 
\begin{eqnarray*}
 \s{C}(\mathrm{Extrem}(A\cap H)) & = &  \s{C}((\mathrm{Extrem}(A)) \cap H)\\
 & = &  \s{C}(\mathrm{Extrem}(A)) \cap H \subseteq \s{C}(\mathrm{Extrem}(A)),
\end{eqnarray*}
and so, $a\in A\cap H \subseteq \s{C}(\mathrm{Extrem}(A))$.
Therefore, we proved that
\[\partial A \subseteq \s{C}(\mathrm{Extrem}(A)),\]
from which we deduce that
$A = \s{C}(\mathrm{Extrem}(A))$,
as explained earlier.
$\bigsquare$

\remark
Observe that Krein and Milman's theorem implies
that any nonempty compact and convex set has
a nonempty subset of extremal points.
This is intuitively obvious, but hard to prove!
Krein and Milman's theorem also applies to
infinite dimensional affine spaces, provided that they
are locally convex, see Valentine \cite{Valentine}, Chapter 11,
Bourbaki \cite{BourbakiEVT}, Chapter II, Barvinok \cite{Barvinok},
Chapter 3, or Lax \cite{Lax}, Chapter 13.

\medskip
We conclude this  chapter
with three other classics of convex geometry.

\section{Radon's and Helly's Theorems and Centerpoints}
\label{secconvex3}
We begin with  {\it Radon's theorem\/}.
\begin{thm}
\label{radon}
Given any affine space $E$ of dimension $\mdeg$,
for every subset $X$ of $E$, if $X$ has  at least
$\mdeg + 2$ points, then there is a partition of $X$ into
two nonempty disjoint subsets $X_1$ and $X_2$ such that the
convex hulls of $X_1$ and $X_2$ have a nonempty intersection.
\end{thm}
\nsindex{Radon's theorem}

\medskip
\proof Pick some origin $O$ in $E$.
Write $X = (x_i)_{i\in L}$ for some index set $L$
(we can let $L = X$).
Since by assumption $|X| \geq m + 2$
where $m = \dimm{(E)}$, $X$ is affinely dependent, and
by Lemma 2.6.5
from Gallier \cite{Gallbook2},
there is a family $(\mu_k)_{k\in L}$ (of finite support)
of scalars, not all null,   such that
\[\sum_{k \in L}\mu_k = 0\quad\hbox{and}\quad
\sum_{k \in L} \mu_k\libvecbo{O}{x_k} = \novect{0}.\]
Since $\sum_{k \in L}\mu_k = 0$, the $\mu_k$ are not all null,
and $(\mu_k)_{k\in L}$ has finite support, the sets 
\[I = \{i\in L\ |\ \mu_i > 0\}\quad\hbox{and}\quad
J = \{j\in L\ |\ \mu_j < 0\}\]
are nonempty, finite,  and obviously disjoint.
Let
\[X_1 = \{x_i\in X\ |\ \mu_i > 0\}\quad\hbox{and}\quad
X_2 = \{x_i\in X\ |\ \mu_i \leq 0\}.\]
Again, since the $\mu_k$ are not all null and
$\sum_{k \in L}\mu_k = 0$, the sets $X_1$ and $X_2$
are nonempty, and obviously
\[X_1\cap X_2 = \emptyset\quad\hbox{and}\quad
X_1\cup X_2 = X.\]
Furthermore, the definition of 
$I$ and $J$ implies  that
$(x_i)_{i\in I}\subseteq X_1$ and $(x_j)_{j\in J}\subseteq X_2$.
It remains to prove that $\chull{X_1}\cap \chull{X_2} \not=\emptyset$.
The definition of $I$ and $J$ implies that
\[\sum_{k \in L} \mu_k\libvecbo{O}{x_k} = \novect{0}\]
can be written as
\[
\sum_{i \in I} \mu_i\libvecbo{O}{x_i} +
\sum_{j \in J} \mu_j\libvecbo{O}{x_j} =  \novect{0},
\]
that is, as
\[
\sum_{i \in I} \mu_i\libvecbo{O}{x_i} =
\sum_{j \in J} -\mu_j\libvecbo{O}{x_j},
\]
where 
\[
\sum_{i \in I} \mu_i = \sum_{j \in J} -\mu_j = \mu,
\]
with $\mu > 0$.
Thus, we have
\[\sum_{i \in I} \frac{\mu_i}{\mu}\,\libvecbo{O}{x_i} =
\sum_{j \in J} -\frac{\mu_j}{\mu}\,\libvecbo{O}{x_j},\]
with 
\[\sum_{i \in I} \frac{\mu_i}{\mu} =
\sum_{j \in J} -\frac{\mu_j}{\mu}= 1,\]
proving that
$\sum_{i \in I} ({\mu_i}/{\mu}) x_i\in \chull{X_1}$ and
$\sum_{j \in J} -({\mu_j}/{\mu}) x_j\in \chull{X_2}$
are identical, and thus that
$\chull{X_1}\cap \chull{X_2} \not=\emptyset$.
$\bigsquare$
\endproof

\medskip
Next, we prove 
a version of {\it Helly's theorem\/}.

\begin{thm}
\label{helly}
Given any affine space $E$ of dimension $\mdeg$,
for every family $\{K_1, \ldots, K_n\}$ of $n$ convex subsets of $E$, 
if $n\geq \mdeg + 2$ and the intersection 
$\bigcap_{i \in I} K_i$ of any $\mdeg + 1$ of the
$K_i$ is nonempty (where $I\subseteq \{1, \ldots, n\}$, $|I| = \mdeg + 1$), 
then $\bigcap_{i = 1}^{n} K_i$ is nonempty.
\end{thm}
\nsindex{Helly's theorem}

\medskip
\proof 
The proof is by induction on $n \geq m + 1$
and uses Radon's theorem in the induction step.
For $n = \mdeg + 1$, the assumption of the theorem is that the intersection
of any family of $\mdeg + 1$ of the $K_i$'s is nonempty, and the
theorem holds trivially. Next,
let $L = \{1, 2, \ldots, n + 1\}$, where $n + 1 \geq \mdeg + 2$.
By the induction hypothesis,
$C_i = \bigcap_{j \in (L - \{i\})} K_j$
is nonempty for every $i\in L$. 

\medskip
We claim that
$C_i\cap C_j \not= \emptyset$ for some
$i \not= j$. If so, as $C_i\cap C_j = \bigcap_{k = 1}^{n+1} K_k$, 
we are done.
So, let us assume that the $C_i$'s are pairwise disjoint.
Then, we can pick a set $X = \{a_1, \ldots, a_{n+1}\}$ such that
$a_i\in C_i$, for every $i\in L$.
By Radon's Theorem, there are
two nonempty disjoint sets $X_1, X_2\subseteq X$ such that $X = X_1\cup X_2$
and $\chull{X_1}\cap \chull{X_2} \not=\emptyset$.
However, $X_1 \subseteq K_j$ for every $j$  with
$a_j\notin X_1$. This is because $a_j\notin K_j$ for every $j$,
and so, we get 
\[
X_1 \subseteq \bigcap_{a_j\notin X_1} K_j.
\]
Symetrically, we also have
\[
X_2 \subseteq \bigcap_{a_j\notin X_2} K_j.
\]
Since the $K_j$'s are convex and 
\[
\left(\bigcap_{a_j\notin X_1} K_j\right)\cap \left(\bigcap_{a_j\notin X_2} K_j\right)
= \bigcap_{i = 1}^{n+1} K_i,
\]
it follows that
$\chull{X_1}\cap \chull{X_2} \subseteq
\bigcap_{i = 1}^{n+1} K_i$,  so that
$\bigcap_{i = 1}^{n+1} K_i$ is nonempty, contradicting the fact
that $C_i\cap C_j = \emptyset$ for all $i \not= j$.
$\bigsquare$
\endproof

\medskip
A more general version of Helly's theorem is proved
in Berger \cite{Berger90b}.
An amusing corollary of Helly's theorem is 
the following result:
Consider $n\geq 4$ parallel line segments in the affine plane
$\affreal^2$. If every three of these line segments meet
a line, then all of these line segments meet a common line.

\medskip
We conclude this chapter
with a nice application of Helly's Theorem to
the existence of centerpoints. Centerpoints generalize the notion
of median to higher dimensions. Recall that if we have a set of
$n$ data points, $S = \{a_1, \ldots, a_n\}$, on the real line, 
a {\it median\/} for $S$ is a point, $x$, such that at least $n/2$
of the points in $S$ belong to both  intervals $[x, \infty)$
and $(-\infty, x]$. 

\medskip
Given any hyperplane, $H$, recall that the closed half-spaces
determined by $H$ are denoted $H_+$ and $H_-$ and that
$H \subseteq H_+$ and  $H \subseteq H_-$.
We let $\interio{{H_+}} = H_+ - H$ and
$\interio{{H_-}} = H_- - H$ be the {\it open half-spaces\/}
determined by $H$. 

\begin{defin}
\label{centerpt}
{\em
Let $S = \{a_1, \ldots, a_n\}$ be a set of $n$ points in $\affreal^d$.
A point, $c\in \affreal^d$, is a {\it centerpoint of $S$\/}
iff for every hyperplane, $H$, whenever the closed half-space
$H_+$ (resp. $H_-$) contains $c$, then
$H_+$ (resp. $H_-$) contains at least $\frac{n}{d+1}$ points from $S$.
}
\end{defin}

\medskip
So, for $d = 2$, for each line, $D$, if the closed half-plane
$D_+$ (resp. $D_-$) contains $c$, then $D_+$ (resp. $D_-$)
contains at least a third of the points from $S$.
For $d = 3$, for each plane, $H$, if the closed half-space
$H_+$ (resp. $H_-$) contains $c$, 
then $H_+$ (resp. $H_-$) contains at least a fourth of the points from $S$,
{\it etc.}

\medskip
Observe that a point, $c\in \affreal^d$, is a centerpoint of $S$
iff $c$ belongs to every open half-space, $\interio{{H_+}}$
(resp.  $\interio{{H_-}}$) containing at least $\frac{d n}{d+1} + 1$
points from $S$. 

\medskip
Indeed, if $c$ is a centerpoint of $S$ and $H$ is any hyperplane
such that  $\interio{{H_+}}$ (resp.  $\interio{{H_-}}$)
contains at least  $\frac{d n}{d+1} + 1$ points from $S$, then  
$\interio{{H_+}}$ (resp. $\interio{{H_-}}$) must
contain $c$ as otherwise, the closed  half-space, $H_-$ (resp. $H_+$)
would contain $c$ and at most $n - \frac{d n}{d+1} -1 = \frac{n}{d+1} - 1$
points from $S$, a contradiction. Conversely, assume that 
$c$ belongs to every open half-space, $\interio{{H_+}}$
(resp.  $\interio{{H_-}}$) containing at least $\frac{d n}{d+1} + 1$
points from $S$. Then, for any hyperplane, $H$, if
$c\in H_+$ (resp. $c\in H_-$) but $H_+$ contains at most 
$\frac{n}{d+1} - 1$ points from $S$, then the open half-space,
$\interio{{H_-}}$ (resp.  $\interio{{H_+}}$) would contain
at least $n - \frac{n}{d+1} + 1 = \frac{d n}{d +1} + 1$ points from $S$
but not $c$, a contradiction.

\medskip
We are now ready to prove the existence of centerpoints.

\begin{thm}
\label{centerpt1}
Every finite set,
$S = \{a_1, \ldots, a_n\}$, of $n$ points in $\affreal^d$ has some
centerpoint.
\end{thm}

\proof
We will use the second characterization of centerpoints involving
open half-spaces containing at least  $\frac{d n}{d+1} + 1$ points.

\medskip
Consider the family of sets, 
\begin{eqnarray*}
\s{C} & = &\left\{\mathrm{conv}(S\>\cap \interio{{H_+}})  \mid (\exists H)
\left(|S\>\cap \interio{{H_+}}| > \frac{d n }{d + 1}\right)\right\} \\
& &  \cup  \left\{\mathrm{conv}(S\>\cap \interio{{H_-}}) \mid (\exists H)
\left(|S\>\cap \interio{{H_-}}| > \frac{d n }{d + 1}\right)\right\},
\end{eqnarray*}
where $H$ is a hyperplane.

\medskip
As $S$ is finite, $\s{C}$ consists of a finite number of convex sets,
say $\{C_1, \ldots, C_m\}$. If we prove that
$\bigcap_{i = 1}^m C_i \not= \emptyset$ we are done, because
$\bigcap_{i = 1}^m C_i$ is the set of  centerpoints of $S$.

\medskip
First, we prove by induction on $k$ (with $1\leq k \leq d + 1$),
that any intersection of $k$ of the $C_i$'s has at least
$\frac{(d + 1 - k)n }{d + 1} + k$ elements from $S$.
For $k = 1$, this holds by definition of the $C_i$'s.

\medskip
Next, consider the intersection of $k + 1\leq d + 1$ of the $C_i$'s,
say $C_{i_1} \cap \cdots \cap C_{i_k} \cap C_{i_{k+1}}$.
Let
\begin{eqnarray*}
A & = & S\cap (C_{i_1} \cap \cdots \cap C_{i_k} \cap C_{i_{k+1}}) \\
B & = & S\cap (C_{i_1} \cap \cdots \cap C_{i_k}) \\
C & = & S \cap C_{i_{k+1}}.
\end{eqnarray*}
Note that $A = B\cap C$.
By the induction hypothesis, $B$ contains at least
$\frac{(d + 1 - k)n}{d + 1} + k$ elements from $S$. As
$C$ contains at least $\frac{d n}{d + 1} + 1$ points from $S$, 
and as
\[
|B \cup C| =  |B| + |C| - |B\cap C| = |B| + |C| - |A|
\]
and $|B\cup C| \leq n$, we get $n \geq |B| + |C| - |A|$, that is,
\[
|A| \geq |B| + |C| - n.
\]
It follows that
\[
|A| \geq \frac{(d + 1 - k)n}{d + 1} + k + \frac{d n}{d + 1} + 1 - n
\]
that is,
\[
|A| \geq \frac{(d + 1 - k)n + dn - (d + 1)n}{d + 1} + k + 1 =
 \frac{(d + 1 - (k + 1))n}{d+1} + k + 1,
\]
establishing the induction hypothesis. 

\medskip
Now, if $m \leq d + 1$, the above claim for $k = m$
shows that $\bigcap_{i = 1}^m C_i \not= \emptyset$ and we are done.
If $m \geq d + 2$, the above claim
for $k = d + 1$ shows that any intersection
of $d + 1$ of the $C_i$'s is nonempty. Consequently,
the  conditions for applying Helly's Theorem are satisfied
and therefore,
\[
\bigcap_{i = 1}^m C_i \not= \emptyset.
\]
However, $\bigcap_{i = 1}^m C_i$ is the set of centerpoints
of $S$ and we are done.
$\bigsquare$

\remark
The above proof actually shows that the set of centerpoints of $S$
is a convex set. In fact, it is a finite intersection of
convex hulls of finitely many points, so it is the convex
hull of finitely many points, in other words, a polytope.

\medskip
Jadhav and Mukhopadhyay have given a linear-time  algorithm  for computing
a centerpoint of a finite set of points in the plane.
For $d\geq 3$, it appears that the best that can be done (using
linear programming) is $O(n^d)$. However, there are good
approximation algorithms (Clarkson, Eppstein, Miller, Sturtivant
and Teng) and in $\mathbb{E}^3$  there is a near quadratic algorithm
(Agarwal, Sharir and Welzl).

\chapter[Separation and Supporting  Hyperplanes]
{Separation and Supporting  Hyperplanes}
\label{chap1}
\section{Separation Theorems and Farkas Lemma}
\label{sec1}
It seems intuitively rather obvious that if $A$ and $B$ are two
nonempty disjoint convex sets in $\affreal^2$, then
there is a line, $H$, separating them, in the sense that
$A$ and $B$ belong to the two (disjoint) open half--planes
determined by $H$. However, this is not always true!
For example, this fails if both $A$ and $B$ are closed
and unbounded (find an example).
Nevertheless, the result is true if both $A$ and $B$ are open,
or if the notion of separation is weakened a little bit.
The key result, from which most separation results follow,
is a geometric version of the {\it Hahn-Banach theorem\/}.
In the sequel, we restrict our attention to real
affine spaces of finite dimension. Then, if $X$ is an
affine space of dimension $d$, there is an affine
bijection $f$ between $X$ and $\affreal^d$.

\medskip
Now, $\affreal^d$ is a topological space, under the usual
topology on $\reals^d$ (in fact, $\affreal^d$ is a metric space).
Recall that if $a = (a_1, \ldots, a_d)$ and
$b = (b_1, \ldots, b_d)$ are any two points in $\affreal^d$,
their {\bf Euclidean distance\/}, $d(a, b)$, is given by
\[
d(a, b)= \sqrt{(b_1 - a_1)^2 + \cdots + (b_d - a_d)^2},
\]
which is also the {\it norm\/}, $\norme{\libvecbo{a}{b}}$,
of the vector $\libvecbo{a}{b}$
and that for any $\epsilon > 0$, the {\it open ball of
center $a$ and radius $\epsilon$\/}, $B(a, \epsilon)$, is given by
\[
B(a, \epsilon) = \{b\in \affreal^d \mid d(a, b) < \epsilon\}.
\]
A subset $U \subseteq \affreal^d$ is {\it open\/} (in the
{\it norm topology\/}) if either $U$ is empty or for
every point, $a\in U$, there is some (small) open ball,
$B(a, \epsilon)$, contained in $U$. A subset $C\subseteq \affreal^d$
is {\it closed\/} iff $\affreal^d - C$ is open.
For example, the {\it closed balls\/},
$\overline{B(a, \epsilon)}$, where
\[
\overline{B(a, \epsilon)} = \{b\in \affreal^d \mid d(a, b) \leq \epsilon\},
\]
are closed.
A subset $W\subseteq \affreal^d$ 
is {\it bounded\/} iff there is some ball (open or closed), $B$,
so that  $W \subseteq B$.
A subset $W\subseteq \affreal^d$ is {\it compact\/}
iff  every family, $\{U_i\}_{i\in I}$, that is an open cover of
$W$ (which means that $W = \bigcup_{i\in I} (W\cap U_i)$,
with each $U_i$ an open set) possesses a finite subcover
(which means that there is a finite subset, $F \subseteq I$,
so that  $W = \bigcup_{i\in F} (W\cap U_i)$).
In $\affreal^d$, it can be shown that a subset $W$ is
compact iff $W$ is closed and bounded.
Given a function, $\mapdef{f}{\affreal^m}{\affreal^n}$,
we say that $f$ is {\it continuous\/} if $f^{-1}(V)$ is open
in $\affreal^m$ whenever $V$ is open in $\affreal^n$.
If  $\mapdef{f}{\affreal^m}{\affreal^n}$ is a continuous
function, although it is generally {\bf false\/} that
$f(U)$ is open if $U \subseteq \affreal^m$ is open, it is
easily checked that
$f(K)$ is compact if  $K \subseteq \affreal^m$ is compact.

\medskip
An affine space $X$ of dimension $d$ becomes a topological space
if we give it the topology for which the open subsets are 
of the form $f^{-1}(U)$, where $U$ is any open subset
of $\affreal^d$ and $\mapdef{f}{X}{\affreal^d}$ is an affine
bijection.

\medskip
Given any subset, $A$, of a topological space, $X$, 
the smallest closed set containing  $A$ is denoted
by $\adher{A}$, and is called the {\it closure\/} or {\it adherence of $A$\/}.
A subset, $A$, of $X$, is {\it dense in $X$\/} if $\adher{A} = X$.  
The largest open set contained in $A$ is denoted by $\interio{A}$,
and is called the {\it interior of $A$\/}.
The set, $\fr{A} = \adher{A}\>\cap \adher{X-A}$, is called the 
{\it boundary\/} (or {\it frontier\/}) of $A$. We also denote the boundary
of $A$ by $\dBd A$.

\medskip
In order to prove the Hahn-Banach theorem, we 
will need two lemmas.  Given any two distinct points
$x, y\in X$, we let
\[\>]x, y[\> = \{(1 - \lambda)x + \lambda y \in X \mid
0 < \lambda < 1\}.\]
Our first lemma (Lemma \ref{lem1}) is intuitively quite
obvious so the reader might be puzzled by the length of its proof.
However, after proposing several wrong proofs, we realized
that its proof is more subtle than it might appear.
The proof below is due to Valentine \cite{Valentine}.
See if you can find a shorter (and correct) proof! 

\begin{lemma}
\label{lem1}
Let $S$ be a nonempty convex set and let $x\in \>\interio{S}$
and $y\in \adher{S}$. Then, we have
$\>]x, y[\>\subseteq\> \interio{S}$.
\end{lemma}

\proof
Let $z\in \>]x, y[\>$, that is, $z = (1 - \lambda)x + \lambda y$,
with $0 < \lambda < 1$. Since $x\in \>\interio{S}$,
we can find some open subset, $U$, contained in $S$ so that $x\in U$.
It is easy to check that 
the central magnification of center $z$,
$H_{z, \frac{\lambda - 1}{\lambda}}$, 
maps $x$ to $y$. Then, $H_{z, \frac{\lambda - 1}{\lambda}}(U)$
is an open subset containing $y$ and as $y\in \adher{S}$, we have
$H_{z, \frac{\lambda - 1}{\lambda}}(U) \cap S \not= \emptyset$.
Let $v \in H_{z, \frac{\lambda - 1}{\lambda}}(U) \cap S$
be a point of $S$ in this intersection. Now, there is a unique
point, $u\in U \subseteq S$, such that
$H_{z, \frac{\lambda - 1}{\lambda}}(u) = v$ and, as $S$ is convex,
we deduce that $z = (1 - \lambda)u + \lambda v\in S$. Since $U$ is open,
the set
\[  
(1 - \lambda)U + \lambda v = \{(1 - \lambda) w + \lambda v \mid w\in U\}
\subseteq S 
\]
is also open and $z\in (1 - \lambda)U + \lambda v$, which shows
that $z\in \interio{S}$.
$\bigsquare$

\begin{cor}
\label{cor1}
If $S$ is convex, then $\interio{S}$ is also convex,
and we have  $\interio{S}\> = \>\interio{{\overline{S}}}$.
Furthermore, if $\interio{S}\>\not=\emptyset$, then
$\adher{S} = \overline{{\interio{S}}}$.  
\end{cor}

\danger
Beware that if $S$ is a closed set, then the convex hull,
$\mathrm{conv}(S)$, of $S$ is not necessarily closed!
(Find a counter-example.) 
However, it can be shown that if $S$ is compact, then
$\mathrm{conv}(S)$ is also compact and thus, closed.

\medskip
There is a simple criterion to test whether
a convex set has an empty interior, based on the
notion of dimension of a convex set.

\begin{defin}
\label{dimconv}
{\em
The {\it dimension\/} of a nonempty convex subset,
$S$, of $X$, denoted by $\dimm\> S$, 
is the dimension of the smallest
affine subset, $\lag S\rag$, containing $S$.
}
\end{defin}

\begin{prop}
\label{dimprop}
A nonempty convex set $S$ has a nonempty interior iff
$\dimm\> S = \dimm\> X$.
\end{prop}

\proof
Let $d = \dimm\> X$.
First, assume that $\interio{S}\>\not=\emptyset$.
Then, $S$ contains some open ball of center $a_0$,
and in it, we can find a frame
$(a_0, a_1, \ldots, a_d)$ for $X$.
Thus, $\dimm\> S = \dimm\> X$.
Conversely, let $(a_0, a_1, \ldots, a_d)$ be a frame
of $X$, with $a_i \in S$, for $i = 0, \ldots, d$.
Then, we have
\[\frac{a_0 + \cdots + a_d}{d+1} \in \>\interio{S},\]
and $\interio{S}$ is nonempty.
$\bigsquare$

\danger
Proposition \ref{dimprop} is false in infinite dimension.

\medskip
We leave the following property as an exercise:

\begin{prop}
\label{proplco}
If $S$ is convex, then $\adher{S}$ is also convex.
\end{prop}

\medskip
One can also easily prove that convexity is preserved
under direct image and inverse image by an affine map.

\medskip
The next lemma, which seems intuitively obvious,
is the core of the proof of the
Hahn-Banach theorem. This is the case where
the affine space has dimension two.
First, we need to define what is a convex cone.

\begin{defin}
\label{affcone}
{\em
A convex set, $C$, is a {\it convex cone with
vertex $x$\/} if $C$ is invariant under all central
magnifications, $H_{x, \lambda}$, of center $x$ and 
ratio $\lambda$, with $\lambda > 0$
(i.e., $H_{x, \lambda}(C) = C$).
}
\end{defin}

\medskip
Given a convex set, $S$, and a point, $x\notin S$, 
we can define
\[\mathrm{cone}_x(S) = \bigcup_{\lambda > 0} H_{x, \lambda}(S).\]
It is easy to check that this is a convex cone.

\begin{lemma}
\label{lem2}
Let $B$ be a nonempty open and convex subset of $\affreal^2$,
and let $O$ be a point of $\affreal^2$ so that $O\notin B$. 
Then, there is some line, $L$, through $O$, so that $L\cap B = \emptyset$.
\end{lemma}

\begin{figure}
  \begin{center}
    \begin{pspicture}(-2,-0.3)(6.5,5)
\pspolygon*[fillstyle=solid,linecolor=lightgray]
(0,0)(0,5)(5,5)(5,0)
    \pscircle*[fillstyle=solid,linecolor=darkgray](2,2){2}
    \pnode(0,0){v1}
    \pnode(6,0){v2}
    \pnode(0,5){v3}
    \pnode(-2,0){v4}
    \cnode[fillstyle=solid,fillcolor=black](4,0){2pt}{v5}
    \cnode[fillstyle=solid,fillcolor=black](0,0){2pt}{v6}
    \ncline[linewidth=2pt]{v4}{v2}
    \ncline[linewidth=1pt]{v1}{v3}
    \uput[0](1.8,1.8){$B$}
    \uput[-120](0,0){$O$}
    \uput[45](4.1,4.1){$C$}
    \uput[-90](4,0){$x$}
    \uput[30](6,0){$L$}
    \end{pspicture}
  \end{center}
  \caption{Hahn-Banach Theorem in the plane (Lemma \ref{lem2})}
  \label{HanhBanachpf1}
\end{figure}
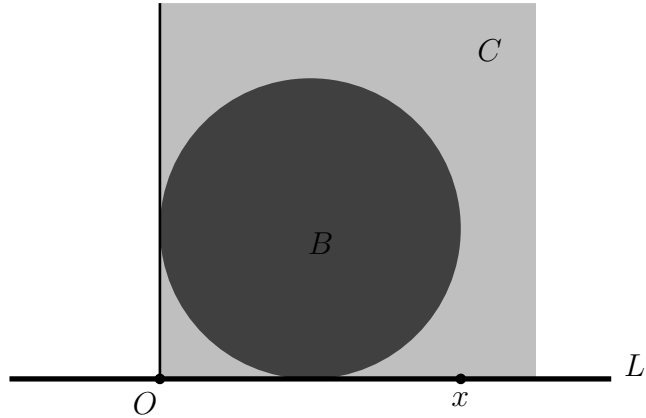

\proof
Define the convex cone $C = \mathrm{cone}_O(B)$.
As $B$ is open, it is easy to check that each $H_{O, \lambda}(B)$
is open and since $C$ is the union of the 
$H_{O, \lambda}(B)$ (for $\lambda > 0$),
which are open, $C$ itself is open. Also,
$O\notin C$. We claim that a least
one point, $x$, of the boundary, $\partial C$, of $C$, is distinct
from $O$. Otherwise, $\partial C = \{O\}$ and
we claim that $C = \affreal^2 - \{O\}$, 
which is not convex, a contradiction. 
Indeed, as $C$ is convex it is connected,
$\affreal^2 - \{O\}$ itself
is connected and $C \subseteq \affreal^2 - \{O\}$.
If $C \not= \affreal^2 - \{O\}$, pick some
point $a\not= O$ in $\affreal^2 - C$ and some point
$c\in C$. Now, a basic property of connectivity
asserts that every continuous
path from $a$ (in the exterior of $C$) to $c$ (in the interior of $C$)
must intersect the boundary of $C$, namely, $\{O\}$.
However, there are plenty of paths from $a$ to $c$ that avoid $O$,
a contradiction. Therefore,  $C = \affreal^2 - \{O\}$.

\medskip
Since $C$ is open and $x\in \partial C$,
we have $x\notin C$. Furthermore, we claim that
$y = 2O - x$ (the symmetric of $x$ w.r.t. $O$)
does not belong to $C$ either. Otherwise,
we would have $y\in \>\interio{C}\> = C$
and $x\in \adher{C}$, and by Lemma \ref{lem1},
we would get $O\in C$, a contradiction.
Therefore, the line through $O$ and $x$ misses
$C$ entirely (since $C$ is a cone), and thus, $B\subseteq C$.
$\bigsquare$

\medskip
Finally, we come to the 
Hahn-Banach theorem.

\begin{thm}
\label{hahnBanach} (Hahn-Banach Theorem, geometric form)
Let $X$ be a (finite-dimensional) affine space,
$A$ be a nonempty open and convex subset of $X$ and $L$
be an affine subspace of $X$ so that $A\cap L = \emptyset$.
Then, there is some hyperplane, $H$, containing $L$,  that is
disjoint from $A$. 
\end{thm}

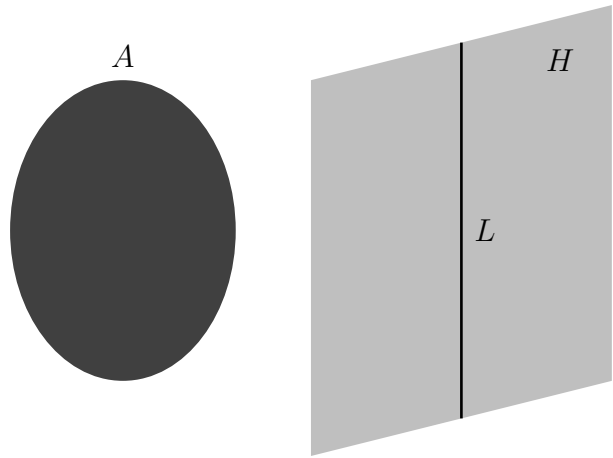
\begin{figure}
  \begin{center}
    \begin{pspicture}(0,0)(8,6)
\pspolygon*[fillstyle=solid,linecolor=lightgray]
(4,0)(4,5)(8,6)(8,1)
    \psellipse*[fillstyle=solid,linecolor=darkgray](1.5,3)(1.5,2)
    \pnode(6,0.5){v1}
    \pnode(6,5.5){v2}
    \ncline[linewidth=1pt]{v1}{v2}
    \uput[90](1.5,5){$A$}
    \uput[0](6,3){$L$}
    \uput[45](7,5){$H$}
    \end{pspicture}
  \end{center}
  \caption{Hahn-Banach Theorem, geometric form (Theorem \ref{hahnBanach})}
  \label{HanhBanachf1}
\end{figure}

\proof
The case where $\dimm\> X = 1$ is trivial. Thus, we
may assume that  $\dimm\> X \geq 2$.
We reduce the proof to the case where $\dimm\> X = 2$.
Let $V$ be an affine subspace of $X$ of
maximal dimension containing $L$ and so that
$V\cap A = \emptyset$. Pick an origin $O\in L$ in $X$,
and consider the vector space $X_O$.
We would like to prove that $V$ is a hyperplane, i.e.,
$\dimm\> V = \dimm\> X - 1$.
We proceed by contradiction. Thus, assume that
$\dimm\> V \leq \dimm\> X - 2$.
In this case, the quotient space $X/V$ has
dimension at least $2$. We also know that $X/V$ is isomorphic to 
the orthogonal complement, $V^{\perp}$, of $V$ so we may identify
$X/V$ and $V^{\perp}$. The (orthogonal) projection map,
$\mapdef{\pi}{X}{V^{\perp}}$, is linear, continuous, and we can show 
that $\pi$ maps the open subset $A$ to an open subset $\pi(A)$,
which is also convex (one way to prove that $\pi(A)$ is open
is to observe that for any point, $a\in A$, a small open ball 
of center $a$ contained in $A$ is projected by $\pi$ to an open ball 
contained in $\pi(A)$ and as $\pi$ is surjective, $\pi(A)$ is open).
Furthermore, $0\notin \pi(A)$.
Since  $V^{\perp}$ has dimension at least $2$, there is some
plane $P$ (a subspace of dimension $2$) intersecting
$\pi(A)$, and thus, we obtain a nonempty open and convex 
subset $B = \pi(A)\cap P$ in the plane
$P \cong \affreal^2$. So, we can apply Lemma \ref{lem2}
to $B$ and the point $O = 0$ in $P\cong \affreal^2$ to find
a line, $l$,  (in $P$) through $O$ with $l\cap B = \emptyset$. 
But then, $l\cap \pi(A) = \emptyset$ and
$W = \pi^{-1}(l)$ is an affine subspace such that
$W\cap A = \emptyset$ and $W$ properly contains $V$,
contradicting the maximality of $V$.
$\bigsquare$

\remark
The geometric form of the Hahn-Banach theorem also 
holds when the dimension of $X$ is infinite but
a slightly more sophisticated proof is required.
Actually, all that is needed is to prove that
a maximal affine subspace containing $L$ and disjoint 
from $A$ exists. This can be done using Zorn's lemma.
For other proofs,
see Bourbaki \cite{BourbakiEVT}, Chapter 2,
Valentine \cite{Valentine}, Chapter 2,
Barvinok \cite{Barvinok}, Chapter 2, or Lax \cite{Lax}, Chapter 3.

\danger
Theorem \ref{hahnBanach} is false if we omit the assumption that
$A$ is open. For a counter-example, let $A\subseteq \affreal^2$
be the union of the half space $y < 0$ with the closed
segment $[0, 1]$ on the $x$-axis and let $L$ be the point
$(2, 0)$ on the boundary of $A$.
It is also false if $A$ is closed! (Find a counter-example).

\medskip
Theorem \ref{hahnBanach} has many important corollaries.
For example, we will eventually prove that for any two nonempty
disjoint convex sets, $A$ and $B$, there is a hyperplane
separating $A$ and $B$, but this will take some work
(recall the definition of a separating hyperplane given in
Definition \ref{separdef}).
We begin with the following version of the Hahn-Banach theorem:

\begin{thm}
\label{hahnBanach2} (Hahn-Banach, second version)
Let $X$ be a (finite-dimensional) affine space,
$A$ be a nonempty convex subset of $X$ with nonempty interior and $L$
be an affine subspace of $X$ so that $A\cap L = \emptyset$.
Then, there is some hyperplane, $H$, containing $L$ and
separating $L$ and $A$.
\end{thm}

\proof
Since $A$ is convex,  by Corollary \ref{cor1}, 
$\interio{A}$ is also convex. By hypothesis,
$\interio{A}$ is nonempty.
So, we can apply Theorem \ref{hahnBanach} to the nonempty
open and convex $\interio{A}$ and to the affine subspace $L$.
We get a hyperplane $H$ containing $L$
such that $\interio{A}\cap\> H = \emptyset$. However, 
$A \subseteq \adher{A} = \adher{\interio{A}}$ and
$\adher{\interio{A}}$ is contained in the closed half space
($H_+$ or $H_-$) containing $\interio{A}$, so
$H$ separates $A$ and $L$.
$\bigsquare$

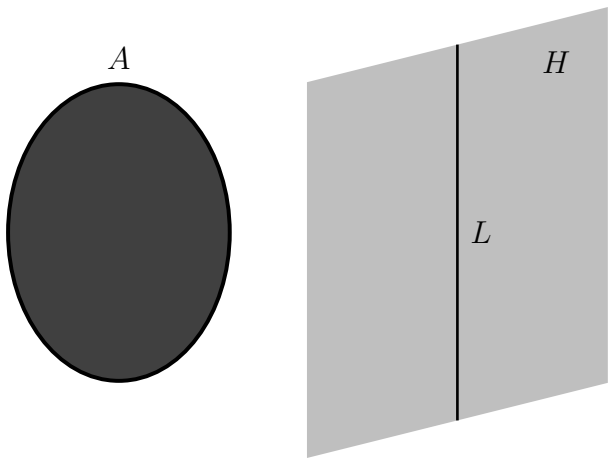
\begin{figure}
  \begin{center}
    \begin{pspicture}(0,0)(8,6)
\pspolygon*[fillstyle=solid,linecolor=lightgray]
(4,0)(4,5)(8,6)(8,1)
    \psellipse[fillstyle=solid,fillcolor=darkgray,linewidth=1.5pt](1.5,3)(1.5,2)
    \pnode(6,0.5){v1}
    \pnode(6,5.5){v2}
    \ncline[linewidth=1pt]{v1}{v2}
    \uput[90](1.5,5){$A$}
    \uput[0](6,3){$L$}
    \uput[45](7,5){$H$}
    \end{pspicture}
  \end{center}
  \caption{Hahn-Banach Theorem, second version (Theorem \ref{hahnBanach2})}
  \label{HanhBanachf2}
\end{figure}

\begin{cor}
\label{separcor1}
Given an affine space, $X$, let $A$ and $B$ be two nonempty 
disjoint convex subsets and assume that $A$ has nonempty
interior ($\interio{A}\> \not= \emptyset$).
Then, there is a hyperplane separating $A$ and $B$.
\end{cor}

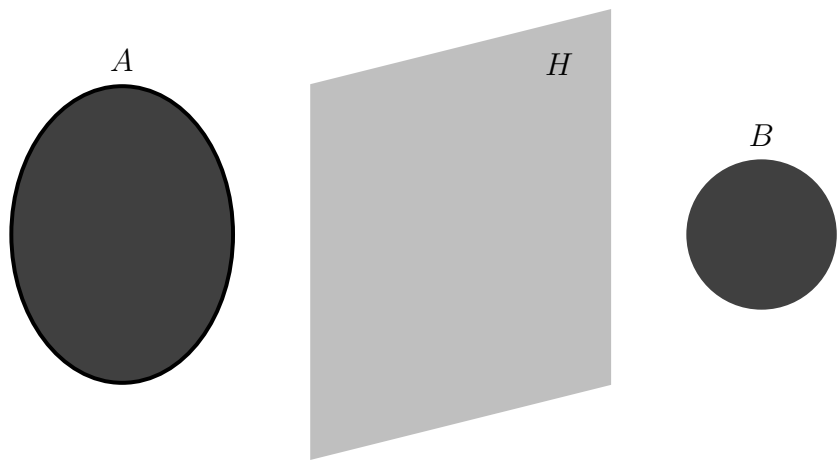
\begin{figure}
  \begin{center}
    \begin{pspicture}(0,0)(11,6)
\pspolygon*[fillstyle=solid,linecolor=lightgray]
(4,0)(4,5)(8,6)(8,1)
    \psellipse[fillstyle=solid,fillcolor=darkgray,linewidth=1.5pt](1.5,3)(1.5,2)
    \pscircle*[fillstyle=solid,linecolor=darkgray](10,3){1}
    \uput[90](1.5,5){$A$}
    \uput[90](10,4){$B$}
    \uput[45](7,5){$H$}
    \end{pspicture}
  \end{center}
  \caption{Separation Theorem, version 1 (Corollary \ref{separcor1})}
  \label{HanhBanachf3}
\end{figure}

\proof
Pick some origin $O$ and consider the vector space $X_O$.
Define $C = A - B$ (a special case of the Minkowski sum) as 
follows:
\[A - B = \{a - b \mid a\in A,\> b\in B\}
= \bigcup_{b\in B} (A - b).
\]
It is easily verified that $C = A - B$ is convex
and has nonempty interior (as a union of subsets having
a nonempty interior). 
Furthermore $O\notin C$, since $A\cap B = \emptyset$.%
\footnote{
Readers who prefer a purely affine argument may define
$C = A - B$ as the {\it affine\/} subset 
\[
A - B = \{O + a - b \mid a\in A,\> b\in B\}.
\]
Again, $O\notin C$ and $C$ is convex.
By adjusting $O$ we can pick the affine form, $f$,
defining a separating hyperplane,  $H$, of $C$ and $\{O\}$, so that
$f(O + a - b) \leq f(O)$, for all $a\in A$ and all $b\in B$, {\it i.e.\/},
$f(a) \leq f(b)$.
}
(Note that the definition depends on the choice
of $O$, but this has no effect on the proof.)
Since $\interio{C}$ is nonempty,  
we can apply Theorem \ref{hahnBanach2} to 
$C$ and to the affine subspace $\{O\}$ and we get a hyperplane, $H$,
separating $C$ and $\{O\}$.
Let $f$ be any linear form defining the hyperplane $H$.
We may assume 
that $f(a - b) \leq 0$, for all $a\in A$ and all $b\in B$, {\it i.e.\/},
$f(a) \leq f(b)$. Consequently, if we let
$\alpha = \sup\{f(a) \mid a\in A\}$ (which makes sense, since
the set $\{f(a) \mid a\in A\}$ is bounded), we have
$f(a) \leq \alpha$ for all $a\in A$ and $f(b) \geq \alpha$
for all $b\in B$, which shows that the affine hyperplane defined
by $f - \alpha$ separates $A$ and $B$. 
$\bigsquare$

\medskip
\remark
Theorem \ref{hahnBanach2} and Corollary \ref{separcor1} 
also hold in the infinite dimensional case,
see Lax \cite{Lax}, Chapter 3, or Barvinok, Chapter 3.

\medskip
Since a hyperplane, $H$, separating $A$ and $B$ as in
Corollary \ref{separcor1} is the boundary of 
each of the two half--spaces that it determines,
we also obtain the following
corollary:

\begin{cor}
\label{separcor2}
Given an affine space, $X$, let $A$ and $B$ be two nonempty 
disjoint open and convex subsets. Then, there is
a hyperplane strictly separating $A$ and $B$.
\end{cor}

\danger
Beware that Corollary \ref{separcor2} {\it fails\/} for
{\it closed\/} convex sets. 
However,  Corollary \ref{separcor2} holds if we also
assume that $A$ (or $B$) is compact.

\medskip
We need to review the notion of distance 
from a point to a subset.
Let $X$ be a metric space with distance
function, $d$. Given any point, $a\in X$, and any nonempty subset, 
$B$, of $X$, we let
\[d(a, B) = \inf_{b\in B} d(a, b)\]
(where $\inf$ is the notation for least upper bound).

\medskip
Now, if $X$ is an affine space of dimension $d$,
it can be given a metric structure by
giving the corresponding vector space a metric structure,
for instance, the metric induced by a Euclidean
structure. We have the following
important property:
For any nonempty closed  subset,
$S \subseteq X$ (not necessarily convex),
and any point, $a\in X$,
there is some point $s\in S$
``achieving the distance from $a$ to $S$,'' i.e.,
so that
\[d(a, S) = d(a, s).\]
The proof uses the fact that the distance
function is continuous and that
a continuous function attains its minimum
on a compact  set, and is left as
an exercise.

\begin{cor}
\label{separcor3}
Given an affine space, $X$, let $A$ and $B$ be two nonempty 
disjoint closed and convex subsets, with $A$ compact.
Then, there is
a hyperplane strictly separating $A$ and $B$.
\end{cor}

\medskip\noindent
{\it Proof sketch\/}.
First, we pick an origin $O$ and 
we give $X_O\cong \affreal^n$ a Euclidean 
structure. Let $d$ denote the associated distance.
Given any  subsets $A$ of $X$,  let
\[A + B(O, \epsilon) = \{x\in X \mid d(x, A) < \epsilon\},\]
where $B(a, \epsilon)$ denotes the open ball, 
$B(a, \epsilon) = \{x\in X \mid d(a, x) < \epsilon\}$,
of center $a$ and radius $\epsilon > 0$.
Note that
\[
A + B(O, \epsilon) 
= \bigcup_{a\in A} B(a, \epsilon),
\]
which shows that $A + B(O, \epsilon)$ is open;
furthermore it is easy to see that
if $A$ is convex, then $A + B(O, \epsilon)$ is also
convex. Now, the function $a \mapsto d(a, B)$ (where
$a\in A$) is continuous and since $A$ is compact, it
achieves its minimum, $d(A, B) =  \min_{a\in A} d(a, B)$,
at some point, $a$, of $A$. Say, $d(A, B) = \delta$.
Since $B$ is closed, there is some $b \in B$ so that 
$d(A, B) = d(a, B) = d(a, b)$
and since $A\cap B = \emptyset$, we must have $\delta > 0$.
Thus, if we pick $\epsilon < \delta/2$, we see that
\[(A + B(O, \epsilon)) \cap (B + B(O, \epsilon)) = \emptyset.\]
Now, $A + B(O, \epsilon)$ and $B + B(O, \epsilon)$ are open, convex
and disjoint and  we
conclude by applying Corollary \ref{separcor2}.
$\bigsquare$

\medskip
A ``cute'' application of Corollary \ref{separcor3} is
one of the many versions of ``Farkas Lemma'' (1893-1894, 1902), a basic result
in the theory of linear programming.
For  any vector, $x = (x_1, \ldots, x_n) \in \reals^n$, and any real,
$\alpha\in \reals$, 
write $x \geq \alpha$ iff $x_i \geq \alpha$, for $i = 1, \ldots, n$.

\begin{lemma}
\label{FarkasI} (Farkas Lemma, Version I)
Given any $d \times n$ real matrix, $A$, and any vector, $z\in \reals^d$,
exactly one of the
following alternatives occurs:
\begin{enumerate}
\item[(a)]
The linear system, $Ax = z$, has a solution,  
$x = (x_1, \ldots, x_n)$,  such that
$x\geq 0$ and $x_1 + \cdots + x_n = 1$, or
\item[(b)]
There is some $c \in \reals^d$ and some $\alpha\in \reals$
such that $\transpos{c} z < \alpha$ and $\transpos{c} A \geq \alpha$.
\end{enumerate}
\end{lemma}

\proof
Let $A_1, \ldots, A_n \in\reals^d$ be the $n$ points corresponding to
the columns of $A$. Then, either 
$z\in \mathrm{conv}(\{A_1, \ldots, A_n\})$ or
$z\notin \mathrm{conv}(\{A_1, \ldots, A_n\})$. In the first case, 
we have a convex combination
\[
z = x_1 A_1 + \cdots +  x_n A_n
\]
where $x_i \geq 0$ and $x_1 + \cdots + x_n = 1$,
so $x = (x_1, \ldots, x_n)$ is a solution satisfying (a).

\medskip
In the second case, by Corollary \ref{separcor3}, there is a hyperplane,
$H$, strictly separating $\{z\}$ and $\mathrm{conv}(\{A_1, \ldots, A_n\})$,
which is obviously closed. In fact, observe that
$z\notin \mathrm{conv}(\{A_1, \ldots, A_n\})$ iff 
there is a hyperplane,
$H$, such that $z\in \interio{H}_-$ and $A_i \in H_+$,
for $i = 1, \ldots, n$.
As the affine hyperplane, $H$, is the zero locus
of an equation of the form
\[
c_1y_1 + \cdots + c_dy_d  = \alpha,
\]
either $\transpos{c} z < \alpha$
and 
$\transpos{c} A_i \geq \alpha$ for $i = 1, \ldots, n$,
that is, $\transpos{c} A \geq \alpha$,
or 
$\transpos{c} z > \alpha$ and $\transpos{c} A \leq \alpha$.
In the second case, 
$\transpos{(-c)} z < -\alpha$
and 
$\transpos{(-c)} A \geq -\alpha$,
so  (b) is satisfied by either $c$ and $\alpha$  or by $-c$ and $-\alpha$. 
$\bigsquare$ 

\remark
If we relax the requirements on solutions of $Ax = z$ and only require
$x \geq 0$ ($x_1 + \cdots + x_n = 1$ is no longer required)
then, in condition (b),
we can take $\alpha = 0$. This is another version of Farkas Lemma.
In this case, instead of 
considering the convex hull of $\{A_1, \ldots, A_n\}$
we are considering the convex cone,
\[
\mathrm{cone}(A_1, \ldots, A_n) = \{\lambda A_1 + \cdots + \lambda_n A_n \mid
\lambda_i \geq 0,\, 1\leq i \leq n\},
\]
that is, we are dropping the condition $\lambda_1 + \cdots + \lambda_n = 1$.
For this version of Farkas Lemma we need the following separation
lemma:

\begin{prop}
\label{separcor3b}
Let $C \subseteq \eucreal^d$
be any closed convex cone with vertex $O$. Then, for every point, $a$, not
in $C$, there is a hyperplane, $H$, passing through $O$
separating $a$ and $C$ with $a\notin H$.
\end{prop}

\proof
Since $C$ is closed and convex and $\{a\}$ is compact and convex,
by Corollary \ref{separcor3}, there is a hyperplane, $H'$,
strictly separating $a$ and $C$.
Let $H$ be the hyperplane through $O$ parallel to $H'$.
Since $C$ and $a$ lie in the two disjoint open half-spaces
determined by $H'$, the point $a$ cannot belong to $H$.
Suppose that some point, $b\in C$, lies in the  open half-space
determined by $H$ and $a$. Then, the line, $L$, through $O$ and $b$
intersects $H'$ in some point, $c$, and as $C$ is a cone, the half line
determined by $O$ and $b$ is contained in $C$. So, $c\in C$
would belong to $H'$, a contradiction. Therefore, $C$
is contained in the closed half-space determined by $H$
that does not contain $a$, as claimed.
$\bigsquare$

\begin{lemma}
\label{FarkasII} (Farkas Lemma, Version II)
Given any $d \times n$ real matrix, $A$, and any vector, $z\in \reals^d$,
exactly one of the
following alternatives occurs:
\begin{enumerate}
\item[(a)]
The linear system, $Ax = z$, has a solution,  
$x$, such that $x\geq 0$, or
\item[(b)]
There is some $c \in \reals^d$
such that $\transpos{c} z < 0$ and $\transpos{c} A \geq 0$.
\end{enumerate}
\end{lemma}

\proof
The proof is analogous to the proof of Lemma \ref{FarkasI} except
that it uses Proposition \ref{separcor3b} instead of  Corollary \ref{separcor3}
and either $z\in \mathrm{cone}(A_1, \ldots, A_n)$ or
$z\notin \mathrm{cone}(A_1, \ldots, A_n)$. 
$\bigsquare$

\medskip
One can show that Farkas II implies Farkas I.
Here is another version of Farkas Lemma having to do with
a system of inequalities, $Ax \leq z$. Although, this version
may seem weaker that Farkas II, it is actually equivalent to it!

\begin{lemma}
\label{FarkasIII} (Farkas Lemma, Version III)
Given any $d \times n$ real matrix, $A$, and any vector, $z\in \reals^d$,
exactly one of the
following alternatives occurs:
\begin{enumerate}
\item[(a)]
The system of inequalities, $Ax \leq z$, has a solution,  
$x$, or
\item[(b)]
There is some $c \in \reals^d$
such that $c\geq 0$, $\transpos{c} z < 0$ and $\transpos{c} A = 0$.
\end{enumerate}
\end{lemma}

\proof
We use two tricks from linear programming:
\begin{enumerate}
\item
We convert the system of inequalities, $Ax \leq z$, into a system of
equations by introducing
a vector of ``slack variables'', $\gamma = (\gamma_1, \ldots, \gamma_d)$,
where the system of equations is
\[
(A, I)\binom{x}{\gamma} = z,
\]
with $\gamma \geq 0$.
\item
We replace each ``unconstrained variable'', $x_i$, by
$x_i = X_i - Y_i$, with $X_i, Y_i \geq 0$.
\end{enumerate}
Then, the original system $Ax \leq z$ has a solution, $x$
(unconstrained), iff the  system of equations
\[
(A, -A, I)
\begin{pmatrix}
X \\
Y \\
\gamma
\end{pmatrix}
 = z
\]
has a solution with $X, Y, \gamma \geq 0$.
By Farkas II, this system has no solution iff there exists some
$c\in \reals^d$ with $\transpos{c} z < 0$ and
\[
\transpos{c} (A, -A, I) \geq 0,
\]
that is, $\transpos{c} A \geq 0$,  $-\transpos{c} A \geq 0$, and
$c \geq 0$. However, these four conditions reduce to
$\transpos{c} z < 0$, $\transpos{c} A  = 0$ and $c\geq 0$.
$\bigsquare$

\medskip
Finally, we have the separation theorem
announced earlier for arbitrary
nonempty  convex subsets.

\begin{thm}
\label{separcor4} (Separation of disjoint convex sets)
Given an affine space, $X$, let $A$ and $B$ be two nonempty 
disjoint convex subsets.
Then, there is a hyperplane separating $A$ and $B$.
\end{thm}

\begin{figure}
  \begin{center}
    \begin{pspicture}(0,0)(10,8)
\pspolygon*[fillstyle=solid,linecolor=gray]
(2,4)(2.5,5)(8.5,5)(8,4)
\pspolygon*[fillstyle=solid,linecolor=lightgray]
(1.5,1)(4.5,7)(8.5,7)(5.5,1)
\pspolygon*[fillstyle=solid,linecolor=gray]
(1.5,3)(2.02,4.05)(8.02,4.05)(7.5,3)
    \psellipse[fillstyle=solid,fillcolor=darkgray,linewidth=1.5pt](5,4)(2,0.8)
    \psellipse[fillstyle=solid,fillcolor=gray,linewidth=1.5pt](3.5,1)(2,0.8)
    \psellipse[fillstyle=solid,fillcolor=gray,linewidth=1.5pt](6.5,7)(2,0.8)
    \pscircle*[fillstyle=solid,linecolor=darkgray](6.9,2.3){1}
    \cnode[fillstyle=solid,fillcolor=black](6.5,7){2pt}{v1}
    \cnode[fillstyle=solid,fillcolor=black](3.5,1){2pt}{v2}
    \cnode[fillstyle=solid,fillcolor=black](5,4){2pt}{v3}
    \uput[0](6.5,7){$x$}
    \uput[0](3.5,1){$-x$}
    \uput[0](7,4){$A$}
    \uput[0](7.9,2.3){$B$}
    \uput[30](8.5,7){$A + x$}
    \uput[0](5.5,5.5){$C$}
    \uput[-30](5.5,1){$A - x$}
    \uput[40](4,2){$D$}
    \uput[45](1.7,3){$H$}
    \uput[0](5,4){$O$}
    \end{pspicture}
  \end{center}
  \caption{Separation Theorem, final version  (Theorem \ref{separcor4})}
  \label{HanhBanachf4}
\end{figure}
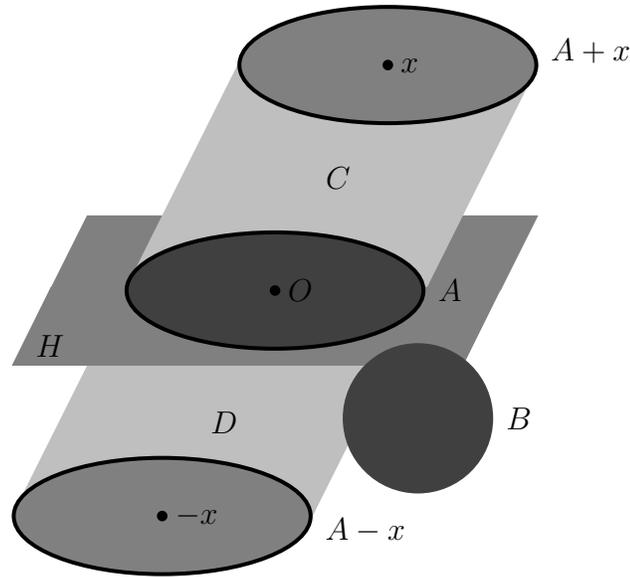

\proof
The proof is by descending induction on $n = \dimm\> A$.
If $\dimm\> A = \dimm\> X$, we know from Proposition \ref{dimprop}
that $A$ has nonempty interior
and we conclude using Corollary \ref{separcor1}.
Next, asssume that the induction hypothesis holds if
$\dimm\> A \geq n$ and assume $\dimm\> A = n - 1$.
Pick an origin $O\in A$ and let $H$ be a hyperplane
containing $A$.
Pick $x\in X$ outside $H$ and define 
$C = \convclo(A \cup \{A + x\})$ 
where $A + x = \{a + x \mid a\in A\}$ and
$D = \convclo(A\cup \{A -x\})$ 
where $A - x = \{a - x \mid a\in A\}$.
Note that $C\cup D$ is convex.
If $B\cap C\not= \emptyset$
and $B\cap D\not=\emptyset$, then the convexity of $B$ and
$C \cup D$ implies that $A\cap B \not= \emptyset$, a contradiction.
Without loss of generality, assume that $B\cap C = \emptyset$.
Since $x$ is outside $H$, we have
$\dimm\> C = n$ and by the induction hypothesis, there is a hyperplane,
$H_1$ separating $C$ and $B$. As $A \subseteq C$, we see that
$H_1$ also separates $A$ and $B$.
$\bigsquare$

\remarks
\begin{enumerate}
\item[(1)]
The reader should compare this proof (from Valentine
\cite{Valentine}, Chapter II) with 
Berger's proof using compactness of the projective space
$\projr{d}$ \cite{Berger90b} (Corollary 11.4.7).
\item[(2)]
Rather than using the Hahn-Banach theorem 
to deduce separation results, one may proceed differently
and use the following intuitively obvious lemma,
as in Valentine \cite{Valentine} (Theorem 2.4):
\begin{lemma}
\label{complem1}
If $A$ and $B$ are two nonempty convex sets such that $A \cup B = X$
and $A\cap B = \emptyset$, then $V = \adher{A}\cap \adher{B}$
is a hyperplane.
\end{lemma}
One can then deduce Corollaries \ref{separcor1} and
Theorem \ref{separcor4}.
Yet another approach is followed in Barvinok
\cite{Barvinok}.
\item[(3)]
How can some of the above results be generalized to
infinite dimensional affine spaces, especially
Theorem \ref{hahnBanach} and Corollary \ref{separcor1}?
One approach is to simultaneously relax the notion of interior
and tighten a little the notion of closure, in a more ``linear and 
less topological'' fashion, as in Valentine \cite{Valentine}.

\medskip
Given any subset $A \subseteq X$ (where $X$ may be infinite
dimensional, but is a Hausdorff topological vector space),
say that a point $x\in X$ is {\it linearly accessible from $A$\/}
iff there is some $a\in A$ with $a\not= x$ and
$]a, x[\> \subseteq A$. We let $\mathrm{lina}\> A$ be the
set of all points linearly accessible from $A$ and
$\mathrm{lin}\> A = A \cup \mathrm{lina}\> A$.

\medskip
A point $a\in A$ is a {\it core point of $A$\/} iff for every $y\in X$,
with $y\not= a$, there is some $z\in ]a, y[\>$,  such that
$[a, z] \subseteq A$. The set of all core points is denoted
$\mathrm{core}\> A$.

\medskip
It is not difficult to prove that $\mathrm{lin}\> A \subseteq \adher{A}$
and $\interio{A}\> \subseteq \mathrm{core}\> A$.
If $A$ has nonempty interior, then 
$\mathrm{lin}\> A = \adher{A}$ and $\interio{A}\> = \mathrm{core}\> A$.
Also, if $A$ is convex, then $\mathrm{core}\> A$ and $\mathrm{lin}\> A$
are  convex. Then, Lemma \ref{complem1} still holds
(where $X$ is not necessarily finite dimensional)
if we redefine $V$ as $V = \mathrm{lin}\> A \cap \mathrm{lin}\> B$
and allow the possibility that $V$ could be $X$ itself.
Corollary  \ref{separcor1} also holds in the general case
if we assume that $\mathrm{core}\> A$ is nonempty.
For details, see Valentine \cite{Valentine}, Chapter I and II.
\item[(4)]
Yet another approach is to define the notion of an
algebraically open convex set, as in Barvinok \cite{Barvinok}.
A convex set, $A$, is {\it algebraically open\/} iff
the intersection of $A$ with every line, $L$,  is an open interval,
possibly empty or infinite at either end (or all of $L$).
An open convex set is algebraically open. Then, the Hahn-Banach
theorem holds provided that $A$ is an algebraically open convex
set and similarly, Corollary  \ref{separcor1} also holds 
provided  $A$ is algebraically open.
For details, see Barvinok \cite{Barvinok}, Chapter 2 and 3.
We do not know how the notion ``algebraically open'' relates
to the concept of $\mathrm{core}$.
\item[(5)]
Theorems \ref{hahnBanach}, \ref{hahnBanach2} and
Corollary \ref{separcor1}  are proved in Lax \cite{Lax} using the
notion of {\it gauge function\/} in the more general case
where $A$ has some core point (but beware that Lax uses the
terminology {\it interior point\/} instead of core point!).
\end{enumerate}

\medskip
An important special case of separation is the case where
$A$ is convex and $B = \{a\}$, for some point, $a$, in $A$.

\section{Supporting Hyperplanes and Minkowski's Proposition}
\label{sec2}
Recall the definition of a supporting hyperplane given in Definition
\ref{suporthyp}.
We have the following important proposition
first proved by Minkowski (1896):

\begin{prop} (Minkowski)
\label{suportprop1}
Let $A$ be a nonempty closed and convex subset. Then, for every
point $a\in \partial A$, there is a supporting hyperplane 
to $A$ through $a$.
\end{prop}

\proof
Let $d = \dimm\> A$. If $d < \dimm\> X$
(i.e., $A$ has empty interior), then
$A$ is contained in some affine subspace $V$ of dimension
$d < \dimm\> X$, and any hyperplane containing $V$ is
a supporting hyperplane for every $a\in A$.
Now, assume $d = \dimm\> X$, so that
$\interio{A}\>\not= \> \emptyset$.
If $a\in \partial A$, then $\{a\}\>\cap \interio{A}\>= \> \emptyset$.
By Theorem \ref{hahnBanach}, there is a hyperplane $H$ separating
$\interio{A}$ and $L = \{a\}$. However, by Corollary \ref{cor1},
since $\interio{A}\>\not= \> \emptyset$ and $A$ is closed, we have
\[A = \adher{A} = \overline{{\interio{A}}}.\]
Now, the half--space containing $\interio{A}$ is closed, and thus,
it contains $\overline{{\interio{A}}} = A$. Therefore, $H$ separates
$A$ and $\{a\}$.
$\bigsquare$

\danger
Beware that Proposition \ref{suportprop1} is false
when the dimension of $X$ is infinite and
when $\interio{A}\> = \emptyset$.

\medskip
The proposition below gives a sufficient condition
for a closed subset to be convex.

\begin{prop}
\label{suportprop2}
Let $A$ be a closed subset with nonempty interior. If 
there is a supporting hyperplane for every point $a\in \partial A$, 
then $A$ is convex.
\end{prop}

\proof
We leave it as an exercise (see Berger \cite{Berger90b},
Proposition 11.5.4).
$\bigsquare$

\danger
The condition that $A$ has nonempty interior is crucial!

\medskip
The proposition below characterizes closed convex sets
in terms of (closed) half--spaces. It is another intuitive fact whose
rigorous proof is nontrivial.

\begin{prop}
\label{suportprop3}
Let $A$ be a nonempty closed and convex subset.
Then, $A$ is the intersection of all the closed half--spaces
containing it. 
\end{prop}

\proof
Let $A'$ be the intersection of all the
closed half--spaces containing $A$. It is immediately checked that
$A'$ is closed and convex and that $A \subseteq A'$. 
Assume that $A' \not= A$, 
and pick $a\in A' - A$. Then, we can apply Corollary \ref{separcor3}
to $\{a\}$ and $A$ and we find a hyperplane, $H$, strictly separating
$A$ and $\{a\}$; this shows that $A$ belongs to one of the two half-spaces
determined by $H$, yet $a$ does not belong to the same half-space,
contradicting the definition of $A'$.
$\bigsquare$

\section{Polarity and Duality}
\label{sec4}
Let $E = \eucreal^n$ be a Euclidean space of dimension $n$.
Pick any origin, $O$, in $\eucreal^n$
(we may assume $O = (0, \ldots, 0)$).
We know that the inner product on $E = \eucreal^n$ induces
a duality between $E$ and its dual $E^*$
(for example, see Chapter 6, Section 2 of Gallier
\cite{Gallbook2}), namely,
$u\mapsto \varphi_u$, where $\varphi_u$ is the linear
form defined by $\varphi_u(v) = u\cdot v$, for all $v\in E$.
For geometric purposes, it is more convenient to recast
this duality as a correspondence between points and hyperplanes,
using the notion of polarity with respect to the unit sphere,
$S^{n-1} = \{a\in \eucreal^n \mid \smnorme{\libvecbo{O}{a}} = 1\}$.

\medskip
First, we need the following simple fact: For every hyperplane, $H$,
not passing through $O$, there is a {\it unique\/} point, $h$, so that
\[
H = \{a\in \eucreal^n \mid \libvecbo{O}{h}\cdot \libvecbo{O}{a} = 1\}.
\]
Indeed, any hyperplane, $H$, in $\eucreal^n$
is the null set of some equation of the form
\[
\alpha_1 x_1 + \cdots + \alpha_n x_n = \beta,
\] 
and if $O\notin H$, then $\beta \not= 0$. Thus, any hyperplane, $H$,
not passing through $O$  is defined by an equation of the form
\[
h_1 x_1 + \cdots + h_n x_n = 1,
\] 
if we set $h_i = \alpha_i/\beta$. So, if we let $h = (h_1, \ldots, h_n)$,
we see that
\[
H = \{a\in \eucreal^n \mid \libvecbo{O}{h}\cdot \libvecbo{O}{a} = 1\},
\]
as claimed. 
%
\medskip
Now, assume that 
\[
H = \{a\in \eucreal^n \mid \libvecbo{O}{h_1}\cdot \libvecbo{O}{a} = 1\}
=  \{a\in \eucreal^n \mid \libvecbo{O}{h_2}\cdot \libvecbo{O}{a} = 1\}.
\]
The functions 
$a \mapsto \libvecbo{O}{h_1}\cdot \libvecbo{O}{a} - 1$ and
$a \mapsto \libvecbo{O}{h_2}\cdot \libvecbo{O}{a} - 1$ are two affine
forms defining the same hyperplane, so there is a nonzero
scalar, $\lambda$, so that
\[
\libvecbo{O}{h_1}\cdot \libvecbo{O}{a} - 1 = 
\lambda(\libvecbo{O}{h_2}\cdot \libvecbo{O}{a} - 1)
\quad\hbox{for all $a\in \eucreal^n$}
\]
(see Gallier \cite{Gallbook2}, Chapter 2, Section 2.10).
In particular, for $a = O$, we find that
$\lambda = 1$, and so,
\[
\libvecbo{O}{h_1}\cdot \libvecbo{O}{a}  = 
\libvecbo{O}{h_2}\cdot \libvecbo{O}{a}
\quad\hbox{for all $a$}, 
\]
which implies $h_1 = h_2$. This proves the uniqueness of $h$.
  
\medskip
Using the above, we make the following definition:

\begin{defin}
\label{polarpt}
{\em
Given any point, $a\not= O$, the {\it polar hyperplane
of $a$ (w.r.t. $S^{n-1}$)\/} or {\it dual of $a$\/} is the hyperplane, 
$a^{\dagger}$,
given by
\[
a^{\dagger} = \{b\in \eucreal^n \mid  \libvecbo{O}{a}\cdot 
\libvecbo{O}{b}  = 1\}.
\]
Given a hyperplane, $H$, not containing $O$, the {\it pole
of $H$ (w.r.t $S^{n-1}$\/}) or {\it dual of $H$\/}
is the (unique) point, $H^{\dagger}$, so that
\[
H = \{a\in \eucreal^n \mid 
\libvecbo{O}{H^{\dagger}}\cdot \libvecbo{O}{a} = 1\}.
\]
}
\end{defin}

\medskip
We often abbreviate polar hyperplane to polar.
We immediately check that $a^{{\dagger}{\dagger}} = a$ and $H^{{\dagger}{\dagger}} = H$, so,
we obtain a bijective correspondence between $\eucreal^n - \{O\}$ and
the set of hyperplanes not passing through $O$.

\medskip
When $a$ is outside the sphere $S^{n-1}$, there is a nice geometric
interpetation for the polar hyperplane, $H = a^{\dagger}$.
Indeed, in this case, since
\[
H = a^{\dagger} = \{b\in \eucreal^n \mid  
\libvecbo{O}{a}\cdot \libvecbo{O}{b}  = 1\}
\]
and $\norme{\libvecbo{O}{a}} > 1$, the hyperplane $H$ intersects
$S^{n-1}$ (along an $(n-2)$-dimensional sphere) and if $b$ is any point
on $H\cap S^{n-1}$, we claim that $\libvecbo{O}{b}$ and
$\libvecbo{b}{a}$ are orthogonal. This means that
$H\cap S^{n-1}$ is the set of points on  $S^{n-1}$ where the
lines through $a$ and tangent to $S^{n- 1}$ touch $S^{n-1}$
(they form a cone tangent to $S^{n-1}$ with apex $a$).
Indeed, as $\libvecbo{O}{a} = \libvecbo{O}{b} + \libvecbo{b}{a}$
and $b\in H\cap S^{n - 1}$ i.e.,
$\libvecbo{O}{a}\cdot \libvecbo{O}{b}  = 1$
and $\norme{\libvecbo{O}{b}}^2 = 1$,  we get
\[
1 = \libvecbo{O}{a}\cdot \libvecbo{O}{b}  =
(\libvecbo{O}{b} + \libvecbo{b}{a}) \cdot \libvecbo{O}{b} =
\norme{\libvecbo{O}{b}}^2 + \libvecbo{b}{a}\cdot \libvecbo{O}{b} = 
1 + \libvecbo{b}{a}\cdot \libvecbo{O}{b},
\]
which implies $\libvecbo{b}{a}\cdot \libvecbo{O}{b} = 0$.
When $a\in S^{n - 1}$, the hyperplane $a^{\dagger}$ is tangent to $S^{n -1}$
at $a$.

\begin{figure}
  \begin{center}
    \begin{pspicture}(0,0)(5,6)
    \pscircle[linewidth=1pt](2,3){1.5}
    \cnode[fillstyle=solid,fillcolor=black](5,3){2pt}{v1}
    \cnode[fillstyle=solid,fillcolor=black](2,3){2pt}{c1}
    \cnode[fillstyle=solid,fillcolor=black](2.75,4.29){2pt}{b1}
    \pnode(2.75,4.29){v2}
    \pnode(2.75,1.71){v3}
    \pnode(2,3){v4}
    \pnode(2.75,0){v5}
    \pnode(2.75,6){v6}
    \pnode(0,3){v7}
    \pnode(6,3){v8}
    \ncline[linewidth=1pt]{v1}{v2}
    \ncline[linewidth=1pt]{v1}{v3}
    \ncline[linewidth=1pt]{v4}{v2}
    \ncline[linewidth=1pt]{v4}{v3}
    \ncline[linewidth=2pt]{v5}{v6}
    \ncline[linewidth=1pt,linestyle=dashed]{v7}{v8}
    \uput[60](5,3){$a$}
    \uput[30](2.75,0.8){$a^{\dagger}$}
    \uput[120](2,3){$O$}
    \uput[45](2.75,4.29){$b$}
    \end{pspicture}
  \end{center}
  \caption{The polar, $a^{\dagger}$, of a point, $a$, 
outside the sphere  $S^{n-1}$}
  \label{polarfig1}
\end{figure}
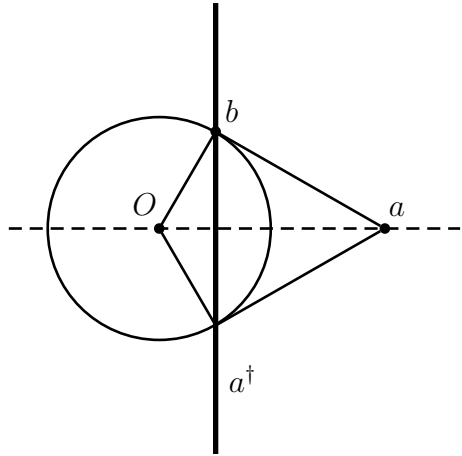

\medskip
Also, observe that for any point, $a\not= O$, and 
any hyperplane, $H$, not passing through $O$, if $a\in H$, 
then, $H^{\dagger}\in a^{\dagger}$, i.e, the pole, $H^{\dagger}$, of
$H$ belongs to the polar, $a^{\dagger}$, of $a$.
Indeed, $H^{\dagger}$ is the unique point so that
\[
H = \{b\in \eucreal^n \mid \libvecbo{O}{H^{\dagger}}\cdot \libvecbo{O}{b} = 1\}
\]
and 
\[
a^{\dagger} = \{b\in \eucreal^n \mid  \libvecbo{O}{a}\cdot \libvecbo{O}{b}  = 1\};
\]
since $a\in H$, we have
$\libvecbo{O}{H^{\dagger}}\cdot \libvecbo{O}{a} = 1$, which shows that
$H^{\dagger}\in a^{\dagger}$. 

\medskip
If $a = (a_1, \ldots, a_n)$, the equation of the polar hyperplane,
$a^{\dagger}$, is
\[
a_1 X_1 + \cdots + a_n X_n = 1.
\]

\remark
As we noted, polarity in a Euclidean space suffers from the minor
defect that the polar of the origin is 
undefined and, similarly, the pole of a hyperplane through the
origin does not make sense. If we embed $\eucreal^n$ into
the projective space, $\projr{n}$, by adding a ``hyperplane
at infinity''  (a copy of $\projr{n - 1}$), 
thereby  viewing $\projr{n}$ as the disjoint union
$\projr{n} = \eucreal^n \cup \projr{n - 1}$, 
then the polarity correspondence can be defined everywhere.
Indeed, the polar of the origin is the hyperplane at 
infinity ($\projr{n - 1}$) and since $\projr{n - 1}$ can 
be viewed as the set of hyperplanes through the origin in $\eucreal^n$,
the pole of a hyperplane through the origin is the corresponding
``point at infinity'' in $\projr{n - 1}$.

\medskip
Now, we would like to extend this correspondence to subsets of 
$\eucreal^n$, in particular, to convex sets.
Given a hyperplane, $H$, not containing $O$, we denote by $H_-$
the closed half-space containing $O$.

\begin{defin}
\label{polar}
{\em
Given any subset, $A$, of $\eucreal^n$, the set
\[
A^* = \{b\in \eucreal^n \mid  \libvecbo{O}{a}\cdot \libvecbo{O}{b}  \leq 1,
\quad\hbox{for all $a\in A$}\} = 
\bigcap_{
\begin{subarray}{l}
a\in A \\
a\not= O
\end{subarray}
} (a^{\dagger})_-,
\]
is called the {\it polar dual\/} or {\it reciprocal\/} of $A$.
}
\end{defin}

\medskip
For simplicity of notation, we write $a^{\dagger}_-$ for 
$(a^{\dagger})_-$. Observe that $\{O\}^* = \eucreal^n$, so it is convenient to
set $O^{\dagger}_- = \eucreal^n$, even though $O^{\dagger}$
is undefined.
By definition, $A^*$ is convex even if $A$
is not. Furthermore, note that
\begin{enumerate}
\item[(1)]
$A\subseteq A^{**}$.
\item[(2)]
If $A \subseteq B$, then
$B^* \subseteq A^*$.
\item[(3)]
If $A$ is convex and closed, then
$A^* = (\partial A)^*$.
\end{enumerate}

\medskip
It follows immediately from (1) and (2) that
$A^{***} = A^*$.
Also, if $B^n(r)$ is the (closed) ball of radius $r > 0$ and center $O$,
it is obvious by definition that
$B^n(r)^* = B^n(1/r)$.

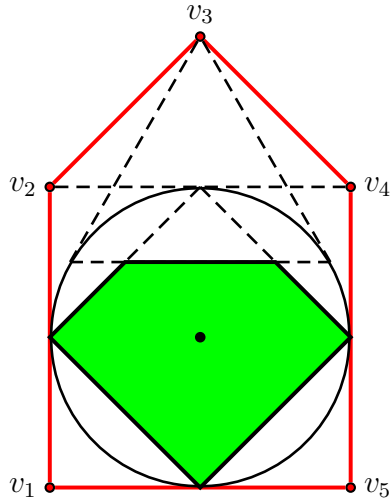
\begin{figure}
  \begin{center}
    \begin{pspicture}(0,0)(4,6.5)
    \cnode[fillstyle=solid,fillcolor=red](0,0){2pt}{v1}
    \cnode[fillstyle=solid,fillcolor=red](0,4){2pt}{v2}
    \cnode[fillstyle=solid,fillcolor=red](2,6){2pt}{v3}
    \cnode[fillstyle=solid,fillcolor=red](4,4){2pt}{v4}
    \cnode[fillstyle=solid,fillcolor=red](4,0){2pt}{v5}
    \pnode(0.27,3){u1}
    \pnode(3.73,3){u2}
    \pnode(2,4){u3}
    \pnode(0,2){u4}
    \pnode(4,2){u5}
    \ncline[linewidth=1.5pt,linecolor=red]{v1}{v2}
    \ncline[linewidth=1.5pt,linecolor=red]{v2}{v3}
    \ncline[linewidth=1.5pt,linecolor=red]{v3}{v4}
    \ncline[linewidth=1.5pt,linecolor=red]{v4}{v5}
    \ncline[linewidth=1.5pt,linecolor=red]{v5}{v1}
    \ncline[linewidth=1pt,linestyle=dashed]{u1}{u2}
    \ncline[linewidth=1pt,linestyle=dashed]{v3}{u1}
    \ncline[linewidth=1pt,linestyle=dashed]{v3}{u2}
    \ncline[linewidth=1pt,linestyle=dashed]{u3}{u4}
    \ncline[linewidth=1pt,linestyle=dashed]{u3}{u5}
    \ncline[linewidth=1pt,linestyle=dashed]{v2}{v4}
\pspolygon[fillstyle=solid,fillcolor=green,linewidth=1.5pt]
(2,0)(0,2)(1,3)(3,3)(4,2)
    \pscircle[linewidth=1pt](2,2){2}
    \cnode[fillstyle=solid,fillcolor=black](2,2){2pt}{c1}
    \uput[180](0,0){$v_1$}
    \uput[180](0,4){$v_2$}
    \uput[90](2,6){$v_3$}
    \uput[0](4,4){$v_4$}
    \uput[0](4,0){$v_5$}
    \end{pspicture}
  \end{center}
  \caption{The polar dual of a polygon}
  \label{polarfig2}
\end{figure}

\medskip
In Figure \ref{polarfig2}, the polar dual of the polygon
$(v_1, v_2, v_3, v_4, v_5)$ is the polygon shown in green.
This polygon is cut out by the half-planes
determined by the polars of the vertices
$(v_1, v_2, v_3, v_4, v_5)$ and  containing the center of the circle.
These polar lines are all easy to determine
by drawing for each vertex, $v_i$,
the tangent lines to the circle and joining the contact
points. The construction of the polar of $v_3$
is shown in detail.

\remark
We chose a different notation for polar hyperplanes
and polars ($a^{\dagger}$ and $H^{\dagger}$) and polar
duals ($A^*$), to avoid the potential confusion between
$H^{\dagger}$ and $H^*$, where $H$ is a hyperplane (or 
$a^{\dagger}$ and $\{a\}^*$, where $a$ is a point).
Indeed, they are completely different! For example,
the polar dual of a hyperplane is either a line
orthogonal to $H$ through $O$, if $O\in H$, or a semi-infinite
line through $O$ and orthogonal to $H$ whose endpoint is
the pole, $H^{\dagger}$, of $H$, whereas, $H^{\dagger}$ is
a single point! Ziegler (\cite{Ziegler97}, Chapter 2)
use the notation $A^{\bigtriangleup}$
instead of $A^*$ for the polar dual of $A$.

\medskip
We would like to investigate the duality induced by 
the operation $A \mapsto A^*$.
Unfortunately, 
it is not always the case that $A^{**} = A$, but
this is true when $A$ is closed and convex, as shown in the following
proposition:

\begin{prop}
\label{propdual1}
Let $A$ be any subset of $\eucreal^n$ (with origin $O$). 
\begin{enumerate}
\item[(i)]
If $A$ is bounded, then $O \in\> \interio{{A^*}}$; if
$O\in\> \interio{A}$, then $A^*$ is bounded.
\item[(ii)]
If $A$ is a closed and convex subset containing $O$, then
$A^{**} = A$.
\end{enumerate}
\end{prop}

\proof
(i)
If $A$ is bounded, then $A \subseteq B^n(r)$ for some $r > 0$ large
enough. Then, \\
$B^n(r)^* = B^n(1/r) \subseteq A^*$, so that
$O \in\> \interio{{A^*}}$.
If $O\in\> \interio{A}$, then $B^n(r) \subseteq A$ for some
$r$ small enough, so $A^* \subseteq B^n(r)^* = B^r(1/r)$
and $A^*$ is bounded.

\medskip
(ii)
We always have $A \subseteq A^{**}$. We prove that if
$b\notin A$, then $b\notin A^{**}$; this shows that
$A^{**} \subseteq A$ and thus, $A = A^{**}$.
Since $A$ is closed and convex and $\{b\}$ is compact (and convex!),
by Corollary \ref{separcor3}, there is a hyperplane, $H$, 
strictly separating $A$ and $b$ and, in particular, $O\notin H$,
as $O\in A$. If $h = H^{\dagger}$ is the pole of $H$, we have
\[
\libvecbo{O}{h}\cdot \libvecbo{O}{b} > 1
\quad\hbox{and}\quad
\libvecbo{O}{h}\cdot \libvecbo{O}{a} < 1,
\quad\hbox{for all $a\in A$} 
\]
since
$H_- = \{a\in \eucreal^n \mid \libvecbo{O}{h}\cdot \libvecbo{O}{a} \leq 1\}$.
This shows that $b\notin A^{**}$, 
since
\begin{eqnarray*}
A^{**} & = & \{c\in \eucreal^n \mid  
\libvecbo{O}{d}\cdot \libvecbo{O}{c}  \leq 1\quad\hbox{for all $d\in A^*$}\} \\
& = &
 \{c\in \eucreal^n \mid  (\forall d\in \eucreal^n)(\hbox{if}\quad
\libvecbo{O}{d}\cdot \libvecbo{O}{a}  \leq 1
\quad\hbox{for all $a\in A$},
\quad\hbox{then}\quad
\libvecbo{O}{d}\cdot \libvecbo{O}{c}  \leq 1)\},
\end{eqnarray*}
just let $c = b$ and $d = h$.
$\bigsquare$

\medskip
\remark
For an arbitrary subset, $A\subseteq \eucreal^n$, it can be shown
that $A^{**} = \overline{\chullb{A \cup \{O\}}}$,
the topological closure of the convex hull of $A \cup \{O\}$.

\medskip
Proposition \ref{propdual1} will play a key role in studying
polytopes, but before doing this, we need one more proposition.

\begin{prop}
\label{propdual2}
Let $A$ be any closed convex subset  of $\eucreal^n$
such that $O\in\> \interio{A}$. The polar hyperplanes of the points
of the boundary of $A$
constitute the set of supporting hyperplanes of $A^*$.
Furthermore, for any $a\in \partial A$, the points
of $A^*$ where  $H = a^{\dagger}$ is a supporting hyperplane of $A^*$
are the poles of supporting hyperplanes of $A$ at $a$.
\end{prop}

\proof
Since $O\in\> \interio{A}$, we have $O\notin \partial A$,
and so, for every  $a\in \partial A$, the polar hyperplane $a^{\dagger}$
is well-defined. 
Pick any $a\in \partial A$ 
and let $H = a^{\dagger}$ be its polar hyperplane.
By definition, $A^* \subseteq H_-$,
the closed half-space determined by $H$ and containing $O$.
If $T$ is any  supporting hyperplane to $A$ at $a$, as 
$a\in T$, we have $t = T^{\dagger} \in a^{\dagger} = H$. Furthermore,
it is a simple exercise to prove that 
$t\in (T_-)^*$ (in fact, $(T_-)^*$ is the interval with endpoints
$O$ and $t$).
Since $A \subseteq T_-$ (because $T$ is a supporting hyperplane
to $A$ at $a$), we deduce that $t\in A^*$, and thus,
$H$ is a supporting hyperplane to $A^*$ at $t$.
By Proposition \ref{propdual1}, as $A$ is closed and convex,
$A^{**} = A$; it follows that all supporting hyperplanes to $A^*$ are indeed
obtained this way.
$\bigsquare$

\chapter[Polyhedra and Polytopes]
{Polyhedra and Polytopes}
\label{chap2}
\section{Polyhedra, $\s{H}$-Polytopes and $\s{V}$-Polytopes}
\label{sec5}
There are two natural ways to define a convex polyhedron, $A$:
\begin{enumerate}
\item[(1)]
As the convex hull of a finite set of points.
\item[(2)]
As a subset of $\eucreal^n$ cut out by a finite number
of hyperplanes, more precisely, as the intersection of a finite
number of (closed) half-spaces.
\end{enumerate}

\medskip
As stated, these two definitions are not equivalent
because (1) implies that a polyhedron is bounded,
whereas (2) allows unbounded subsets. Now, if we require in
(2) that the convex set $A$ is bounded, it is quite
clear for $n = 2$ that the two definitions  (1) and (2)
are equivalent; for $n = 3$, it is intuitively clear that definitions 
(1) and (2) are still equivalent, but  proving  this equivalence
rigorously does not appear to be that easy.
What about the equivalence when $n \geq 4$?

\medskip
It turns out that definitions (1) and (2) are equivalent for all $n$,
but this is a nontrivial theorem and 
a rigorous proof does not come by so cheaply. Fortunately,
since we have Krein and Milman's theorem at our disposal
and polar duality, we
can give a rather short proof. The hard direction 
of the equivalence consists in
proving that definition (1) implies definition (2). This is where
the duality induced by polarity becomes handy, especially,
the fact that $A^{**}= A$! (under the right hypotheses).
First, we give precise definitions (following
Ziegler \cite{Ziegler97}).

\begin{defin}
\label{polytopedef}
{\em
Let $\affs$ be any affine Euclidean space of finite dimension, $n$.%
\footnote{This means that the vector space, $\vector{\affs}$,
associated with $\affs$ is a Euclidean space.
}
An {\it $\s{H}$-polyhedron\/} in $\affs$, 
for short, a {\it polyhedron\/}, 
is any subset, $P = \bigcap_{i = 1}^p C_i$,
of $\affs$ defined as the intersection of a finite number
of closed half-spaces, $C_i$; an {\it $\s{H}$-polytope\/} in $\affs$ 
is a bounded polyhedron and
a {\it $\s{V}$-polytope\/} is the convex hull, $P = \convclo(S)$, 
of a finite set of points, $S \subseteq \affs$.
}
\end{defin}

\medskip
Obviously, polyhedra and polytopes are convex and closed 
(in $\affs$). Since the notions of $\s{H}$-polytope
and $\s{V}$-polytope are equivalent (see Theorem \ref{equivpoly}), we 
often use the simpler locution polytope.
Examples of an $\s{H}$-polyhedron and of a $\s{V}$-polytope
are shown in Figure \ref{polyfig1a}.

\begin{figure}
  \begin{center}
    \begin{pspicture}(0,-1.5)(5,4.5)
\pspolygon*[fillstyle=solid,linecolor=green]
(5,4)(2,1)(3,0)(5,0)
    \pnode(5,4){u1}
    \pnode(0,-1){u2}
    \pnode(0,3){u3}
    \pnode(4,-1){u4}
    \pnode(0,0){u5}
    \pnode(5,0){u6}
    \cnode[fillstyle=solid,fillcolor=red](2,1){2pt}{w1}
    \cnode[fillstyle=solid,fillcolor=red](3,0){2pt}{w2}
    \ncline[linewidth=1pt]{u1}{u2}
    \ncline[linewidth=1pt]{u3}{u4}
    \ncline[linewidth=1pt]{u5}{u6}
    \ncline[linewidth=2pt,linecolor=blue]{w1}{w2}
    \ncline[linewidth=2pt,linecolor=blue]{w1}{u1}
    \ncline[linewidth=2pt,linecolor=blue]{w2}{u6}
    \uput[-90](2.5,-1){(a)}    
    \end{pspicture}
\hskip 2cm
    \begin{pspicture}(0,-1.5)(5,4.5)
\pspolygon*[fillstyle=solid,linecolor=green]
(3,4)(0,1)(1,0)(3,0)(5,2)
    \cnode[fillstyle=solid,fillcolor=red](0,1){2pt}{w1}
    \cnode[fillstyle=solid,fillcolor=red](1,0){2pt}{w2}
    \cnode[fillstyle=solid,fillcolor=red](3,0){2pt}{w3}
    \cnode[fillstyle=solid,fillcolor=red](5,2){2pt}{w4}
    \cnode[fillstyle=solid,fillcolor=red](3,4){2pt}{w5}
    \ncline[linewidth=2pt,linecolor=blue]{w1}{w2}
    \ncline[linewidth=2pt,linecolor=blue]{w2}{w3}
    \ncline[linewidth=2pt,linecolor=blue]{w3}{w4}
    \ncline[linewidth=2pt,linecolor=blue]{w4}{w5}
    \ncline[linewidth=2pt,linecolor=blue]{w5}{w1}
    \uput[-90](2.5,-1){(b)}    
    \end{pspicture}
  \end{center}
  \caption{(a) An $\s{H}$-polyhedron. (b) A $\s{V}$-polytope}
  \label{polyfig1a}
\end{figure}

\medskip
Note that Definition \ref{polytopedef} allows $\s{H}$-polytopes
and $\s{V}$-polytopes to have an empty interior,
which is somewhat of an inconvenience. This is not a problem,
since we may always restrict ourselves to the affine
hull of $P$ (some affine space, $E$, of dimension $d\leq n$, 
where $d = \dimm(P)$, as in Definition \ref{dimconv}) as we now show.

\begin{prop}
\label{polydef1}
Let $A\subseteq \affs$ be a $\s{V}$-polytope or
an $\s{H}$-polyhedron, let $E = \mathrm{aff}(A)$ 
be the affine hull of $A$ in $\affs$ (with the Euclidean
structure on $E$ induced by the Euclidean structure on $\affs$)
and write $d = \mathrm{dim}(E)$.
Then, the following assertions hold:
\begin{enumerate}
\item[(1)]
The set, $A$, is a $\s{V}$-polytope in $E$ (i.e., viewed as a subset of $E$)
iff $A$ is a $\s{V}$-polytope in $\affs$.
\item[(2)]
The set, $A$, is an $\s{H}$-polyhedron in $E$ 
(i.e., viewed as a subset of $E$)
iff $A$ is an $\s{H}$-polyhedron in $\affs$.
\end{enumerate}
\end{prop}

\proof
(1) This follows immediately because $E$ is an affine subspace
of $\affs$ and every affine subspace of $\affs$
is closed under affine combinations and so, {\it a fortiori\/},
under convex combinations. We leave the details as an easy
exercise.

\medskip
(2) 
Assume $A$ is an $\s{H}$-polyhedron in $\affs$
and that $d < n$. By definition, $A = \bigcap_{i = 1}^p C_i$, 
where the $C_i$ are closed half-spaces determined by some hyperplanes,
$H_1, \ldots, H_p$, in $\affs$.  (Observe that
the hyperplanes, $H_i$'s, associated with the closed half-spaces,
$C_i$, may not be distinct.  For example, we may have 
$C_i = (H_i)_+$ and $C_j = (H_i)_-$, for the two closed
half-spaces determined  by $H_i$.)
As $A \subseteq E$, we have
\[
A = A\cap E = \bigcap_{i = 1}^p (C_i \cap E),
\]
where $C_i \cap E$ is one of the closed half-spaces
determined by the hyperplane, $H_i' = H_i\cap E$,  in $E$.
Thus, $A$ is also an $\s{H}$-polyhedron in $E$.

\medskip
Conversely, assume that $A$ is an $\s{H}$-polyhedron in $E$
and that $d < n$. As any hyperplane,
$H$, in $\affs$ can be written as the intersection,
$H = H_- \cap H_+$, of the two closed half-spaces that it bounds,
$E$ itself can be written as the intersection,
\[
E = \bigcap_{i = 1}^{p} E_i =
\bigcap_{i = 1}^{p} (E_i)_+ \cap (E_i)_-,
\]
of finitely many half-spaces in $\affs$.
Now, as $A$ is an $\s{H}$-polyhedron in $E$, we have
\[
A = \bigcap_{j = 1}^q C_j,
\]
where the $C_j$ are closed half-spaces in $E$ determined
by some hyperplanes, $H_j$, in $E$. However, each $H_j$ can be extended
to a hyperplane, $H_j'$, in $\affs$, and so,
each $C_j$ can be extended to a closed half-space, $C_j'$, in $\affs$
and we still have
\[
A = \bigcap_{j = 1}^q C_j'.
\]
Consequently, we get
\[
A = A\cap E = \bigcap_{i = 1}^{p} ((E_i)_+ \cap (E_i)_-)
\cap\bigcap_{j = 1}^q C_j',
\]
which proves that $A$ is also an $\s{H}$-polyhedron in $\affs$.
$\bigsquare$

\medskip
The following simple proposition shows that we may assume
that $\affs = \eucreal^n$:

\begin{prop}
\label{polydef2}
Given any two affine Euclidean spaces, $E$ and $F$,
if $\mapdef{h}{E}{F}$ is any affine map then:
\begin{enumerate}
\item[(1)]
If $A$ is any $\s{V}$-polytope in $E$, then $h(E)$ is
a $\s{V}$-polytope in $F$.
\item[(2)]
If $h$ is bijective and $A$ is any $\s{H}$-polyhedron in $E$, 
then $h(E)$ is an $\s{H}$-polyhedron in $F$.
\end{enumerate}
\end{prop}

\proof
(1) 
As any affine map preserves affine combinations it also
preserves convex combination. Thus,
$h(\mathrm{conv}(S)) = \mathrm{conv}(h(S))$, for any $S\subseteq E$.

\medskip
(2)
Say $A = \bigcap_{i = 1}^p C_i$ in $E$.
Consider any half-space, $C$, in $E$ and assume that
\[
C = \{x\in E \mid \varphi(x) \leq 0\},
\]
for some affine form, $\varphi$, defining the hyperplane,
$H = \{x\in E\mid \varphi(x) = 0\}$. 
Then, as $h$ is bijective, we get 
\begin{eqnarray*}
h(C) & = & \{h(x)\in F \mid \varphi(x)  \leq 0\} \\
& = & \{y\in F \mid \varphi(h^{-1}(y))  \leq 0\} \\
& = & \{y\in F \mid (\varphi\circ h^{-1})(y)  \leq 0\}.
\end{eqnarray*}
This shows that $h(C)$ is one of the closed half-spaces in $F$ determined
by the hyperplane, $H' = \{y\in F\mid (\varphi\circ h^{-1})(y) = 0\}$. 
Furthermore, as $h$ is bijective, it preserves intersections so
\[
h(A) = 
h\left(\bigcap_{i = 1}^p C_i\right) = \bigcap_{i = 1}^p h(C_i),
\]
a finite intersection of closed half-spaces.
Therefore, $h(A)$ is an $\s{H}$-polyhedron in $F$.
$\bigsquare$

\medskip
By Proposition  \ref{polydef2} we may assume that $\affs = \eucreal^d$
and by Proposition \ref{polydef1} we may assume that
$\mathrm{dim}(A) = d$. These propositions justify the type of argument
beginning with: ``We may assume that $A\subseteq \eucreal^d$
has dimension $d$, that is, that $A$ has nonempty interior''.
This kind of reasonning will occur many times.

\medskip
Since the boundary of a closed half-space, $C_i$, is 
a hyperplane, $H_i$, and since hyperplanes are defined
by affine forms, a closed half-space is defined by
the locus of points satisfying a  ``linear'' inequality
of the form $a_i\cdot x \leq b_i$ or  $a_i\cdot x \geq b_i$,
for some vector $a_i\in \reals^n$ and some $b_i\in \reals$. 
Since  $a_i\cdot x \geq b_i$ is equivalent to
$(-a_i)\cdot x \leq -b_i$, we may restrict our attention
to inequalities with a $\leq$ sign.
Thus, if $A$ is the $d\times p$ matrix whose $i^{\mathrm{th}}$ row
is $a_i$, we see that the $\s{H}$-polyhedron, $P$, is defined by 
the system of linear inequalities, $Ax \leq b$, where
$b = (b_1, \ldots, b_p)\in \reals^p$. We write
\[
P = P(A, b),\quad{with}\quad
P(A, b) = \{x\in \reals^n \mid Ax \leq b\}.
\]
An equation, $a_i\cdot x = b_i$, may be handled as the conjunction of 
the two inequalities $a_i\cdot x \leq  b_i$
and $(-a_i)\cdot x \leq -b_i$.
Also, if $0\in P$, observe that we must have $b_i \geq 0$ for
$i = 1, \ldots, p$. In this case, every inequality for which
$b_i > 0$ can be normalized by dividing both sides by $b_i$,
so we may assume that $b_i = 1$ or $b_i = 0$.
This observation will be useful to show that the polar dual
of an $\s{H}$-polyhedron is a $\s{V}$-polyhedron.

\remark
Some authors call  ``convex''  polyhedra
and ``convex'' polytopes what we have simply called
polyhedra and polytopes. Since Definition \ref{polytopedef} implies that
these objects are convex and since we are not going to 
consider non-convex polyhedra in this chapter, we stick to the
simpler terminology.

\begin{figure}
  \begin{center}
    \begin{pspicture}(0,0)(6.5,6.2)
\pspolygon[fillstyle=solid,fillcolor=darkgray,linewidth=1.5pt]
(0,2.5)(1.5,0.75)(3,0)(2.5,1.7)(0.7,3)
\pspolygon[fillstyle=solid,fillcolor=darkgray,linewidth=1.5pt]
(3,0)(4.9,0.6)(5.8,2)(4.3,3)(2.5,1.7)
\pspolygon[fillstyle=solid,fillcolor=lightgray,linewidth=1.5pt]
(4.3,3)(5.8,2)(6.5,4)(5,5.7)(3.7,5)
\pspolygon[fillstyle=solid,fillcolor=gray,linewidth=1.5pt]
(2.5,1.7)(4.3,3)(3.7,5)(1.5,5)(0.7,3)
\pspolygon[fillstyle=solid,fillcolor=gray,linewidth=1.5pt]
(0,2.5)(0.7,3)(1.5,5)(1,5.5)(0.4,4.2)
\pspolygon[fillstyle=solid,fillcolor=white,linewidth=1.5pt]
(1.5,5)(3.7,5)(5,5.7)(3.2,6)(1,5.5)
    \cnode[fillstyle=solid,fillcolor=black](0,2.5){2pt}{v1}
    \cnode[fillstyle=solid,fillcolor=black](1.5,0.75){2pt}{v2}
    \cnode[fillstyle=solid,fillcolor=black](3,0){2pt}{v3}
    \cnode[fillstyle=solid,fillcolor=black](4.9,0.6){2pt}{v4}
    \cnode[fillstyle=solid,fillcolor=black](5.8,2){2pt}{v5}
    \cnode[fillstyle=solid,fillcolor=black](6.5,4){2pt}{v6}    
    \cnode[fillstyle=solid,fillcolor=black](5,5.7){2pt}{v7}    
    \cnode[fillstyle=solid,fillcolor=black](3.2,6){2pt}{v8}
    \cnode[fillstyle=solid,fillcolor=black](1,5.5){2pt}{v9}
    \cnode[fillstyle=solid,fillcolor=black](0.4,4.2){2pt}{v10}
    \cnode[fillstyle=solid,fillcolor=black](1.5,5){2pt}{v11}
    \cnode[fillstyle=solid,fillcolor=black](0.7,3){2pt}{v12}
    \cnode[fillstyle=solid,fillcolor=black](2.5,1.7){2pt}{v13}
    \cnode[fillstyle=solid,fillcolor=black](4.3,3){2pt}{v14}
    \cnode[fillstyle=solid,fillcolor=black](3.7,5){2pt}{v15}
    \end{pspicture}
  \end{center}
  \caption{Example of a polytope (a dodecahedron)}
  \label{polyfig1}
\end{figure}
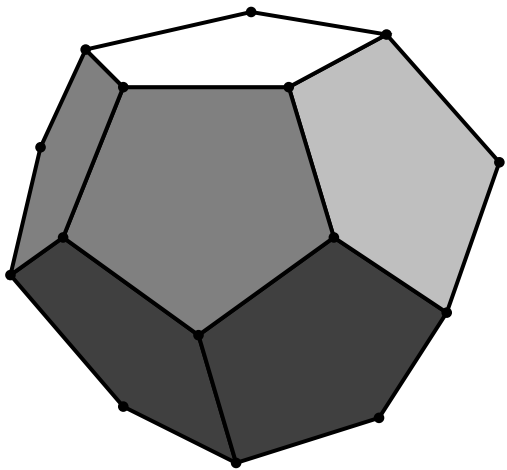

\medskip
One should consult Ziegler \cite{Ziegler97}, Berger \cite{Berger90b},
Grunbaum \cite{Grunbaum} and especially Cromwell \cite{Cromwell},
for pictures of polyhedra and polytopes.
Figure \ref{polyfig1} shows the picture a polytope whose
faces are all pentagons. This polytope is called a
{\it dodecahedron\/}. The dodecahedron has $12$ faces, $30$ edges and 
$20$ vertices.

\medskip
Even better and a lot more entertaining,
take a look at the spectacular web sites
of George Hart,

\medskip\noindent
{\it Virtual Polyedra\/}:
http://www.georgehart.com/virtual-polyhedra/vp.html, 

\medskip\noindent
{\it George Hart\/}'s web site:
http://www.georgehart.com/

\medskip\noindent
and also

\medskip\noindent
{\it Zvi Har'El\/}'s web site:
http://www.math.technion.ac.il/~rl/

\medskip\noindent
The {\it Uniform Polyhedra\/} web site:
http://www.mathconsult.ch/showroom/unipoly/

\medskip\noindent
{\it Paper Models of Polyhedra\/}: 
http://www.korthalsaltes.com/

\medskip\noindent
Bulatov's {\it Polyhedra Collection\/}:
http://www.physics.orst.edu/~bulatov/polyhedra/

\medskip\noindent
Paul Getty's {\it Polyhedral Solids\/}:
http://home.teleport.com/~tpgettys/poly.shtml

\medskip\noindent
Jill Britton's {\it Polyhedra Pastimes\/}:
http://ccins.camosun.bc.ca/~jbritton/jbpolyhedra.htm

\medskip\noindent
and many other web sites dealing with polyhedra in one way or another
by searching for ``polyhedra'' on {\it Google\/}!

\medskip
Obviously, an $n$-simplex is a $\s{V}$-polytope.
The {\it standard $n$-cube\/} is the set
\[
\{(x_1, \ldots, x_n)\in \eucreal^n \mid |x_i| \leq 1,
\quad 1 \leq i \leq n\}.
\]
The standard cube is a $\s{V}$-polytope.
The {\it standard $n$-cross-polytope\/} (or {\it $n$-co-cube\/}) is the set
\[
\{(x_1, \ldots, x_n)\in \eucreal^n \mid \sum_{i = 1}^n |x_i| \leq 1\}.
\]
It is also a $\s{V}$-polytope.

\medskip
What happens if we take the dual of a $\s{V}$-polytope
(resp. an $\s{H}$-polytope)?
The following proposition, although very simple, is an important
step in answering the above question: 

\begin{prop} 
\label{polytopdual}
Let $S = \{a_{i}\}_{i = 1}^p$ be a finite set of points in 
$\eucreal^n$ and let $A = \convclo(S)$ be its convex hull.
If $S\not= \{O\}$, 
then, the dual, $A^*$, of $A$ w.r.t. the center $O$
is an  $\s{H}$-polyhedron; furthermore,
if $O\in\> \interio{A}$, then $A^*$ is an $\s{H}$-polytope, i.e.,
the dual of a $\s{V}$-polytope with nonempty interior
is an $\s{H}$-polytope. If $A = S = \{O\}$, then $A^* = \eucreal^d$.
\end{prop} 

\proof
By definition,  we have
\[
A^* = \{b\in \eucreal^n \mid \libvecbo{O}{b} \cdot 
(\sum_{j = 1}^{p} \lambda_j \libvecbo{O}{a_j}) \leq 1,
\quad \lambda_j \geq 0,\> \sum_{j = 1}^{p} \lambda_j = 1\},
\]
and the right hand side is clearly equal to
$\bigcap_{i = 1}^p \{b\in \eucreal^n \mid \libvecbo{O}{b}\cdot 
\libvecbo{O}{a_i} \leq 1\} = 
\bigcap_{i = 1}^p\, (a_i^{\dagger})_-$,
which is a polyhedron. (Recall that $(a_i^{\dagger})_- = \eucreal^n$
if $a_i = O$.)
If  $O\in\> \interio{A}$, then
$A^*$ is bounded 
(by Proposition \ref{propdual1}) and so,
$A^*$ is an $\s{H}$-polytope.
$\bigsquare$

\medskip
Thus, the dual of the convex hull of a finite set of 
points, $\{a_1, \ldots, a_p\}$, is the intersection of the
half-spaces containing $O$ determined by the polar hyperplanes
of the points $a_i$.

\medskip
It is convenient to restate Proposition \ref{polytopdual} using
matrices. First, observe that the proof of Proposition \ref{polytopdual} 
shows that
\[
\mathrm{conv}(\{a_1, \ldots, a_p\})^* = 
\mathrm{conv}(\{a_1, \ldots, a_p\}\cup \{O\})^*.
\] 
Therefore, we may assume that not all $a_i = 0$ ($1\leq i \leq p$).
If we pick  $O$ as an origin, then every point $a_j$
can be identified with a vector in $\eucreal^n$ and $O$
corresponds to the zero vector, $0$.
Observe that any set of $p$ points, $a_j\in \eucreal^n$, 
corresponds to the $n\times p$ matrix, $A$, whose $j^{\mathrm{th}}$
column is $a_j$. 
Then, the equation of the the polar hyperplane, 
$a_j^{\dagger}$, of any $a_j \> (\not= 0)$ is $a_j\cdot x = 1$, that is
\[
\transpos{a_j}x = 1.
\]
Consequently, the system of inequalities defining 
$\mathrm{conv}(\{a_1, \ldots, a_p\})^*$  
can be written in matrix form as
\[
\mathrm{conv}(\{a_1, \ldots, a_p\})^* = 
\{x\in \reals^n \mid \transpos{A} x \leq \mathbf{1}\},
\]
where $\mathbf{1}$ denotes the vector of $\reals^p$
with all coordinates equal to $1$. We write \\ 
$P(\transpos{A}, \mathbf{1}) = 
\{x\in \reals^n \mid \transpos{A} x \leq \mathbf{1}\}$.
There is a useful converse of this property as proved in the next 
proposition.

\begin{prop}
\label{polytopdual2}
Given any set of $p$ points, $\{a_1, \ldots, a_p\}$, in $\reals^n$
with $\{a_1, \ldots, a_p\} \not= \{0\}$, if $A$ is the
$n\times p$ matrix whose $j^{\mathrm{th}}$ column is $a_j$, then
\[
\mathrm{conv}(\{a_1, \ldots, a_p\})^* = 
P(\transpos{A}, \mathbf{1}),
\]
with 
$P(\transpos{A}, \mathbf{1}) = 
\{x\in \reals^n \mid \transpos{A} x \leq \mathbf{1}\}$.

\medskip
Conversely, given any $p\times n$ matrix, $A$, not equal to the
zero matrix, we have
\[
P(A, \mathbf{1})^* = \mathrm{conv}(\{a_1, \ldots, a_p\}\cup \{0\}),
\]
where $a_i\in\reals^n$ is the  $i^{\mathrm{th}}$ row of $A$ or, equivalently,
\[
P(A, \mathbf{1})^* = \{x\in \reals^n \mid x = \transpos{A}t, \>
t\in \reals^p, \> t \geq 0,\>
\mathbb{I}t = 1\}, 
\]
where $\mathbb{I}$ is the row vector of length $p$ whose
coordinates are all equal to $1$.  
\end{prop}

\proof
Only the second part needs a proof. Let
$B =  \mathrm{conv}(\{a_1, \ldots, a_p\}\cup \{0\})$,
where $a_i\in\reals^n$ is the  $i^{\mathrm{th}}$ {\it row\/} of $A$.
Then, by the first part,
\[
B^* = P(A, \mathbf{1}).
\]
As $0\in B$, by Proposition \ref{propdual1}, we have
$B = B^{**} =  P(A, \mathbf{1})^*$, as claimed.
$\bigsquare$

\remark
Proposition \ref{polytopdual2}
still holds if $A$ is the zero matrix because then,
the inequalities  $\transpos{A} x \leq \mathbf{1}$ 
(or  $A x \leq \mathbf{1}$) are trivially satisfied.
In the first case, $P(\transpos{A}, \mathbf{1}) = \eucreal^d$
and in the second case, $P(A, \mathbf{1})  = \eucreal^d$.

\medskip
Using the above, the reader should check that the dual of a 
simplex is a simplex
and that the dual of an $n$-cube is an $n$-cross polytope.

\medskip
Observe that not every $\s{H}$-polyhedron is of the form
$P(A, \mathbf{1})$. Firstly,  $0$ belongs 
to the interior of $P(A, \mathbf{1})$ and, secondly
cones with apex $0$ can't be described in this form.
However, we will see in Section \ref{sec5c} that
the full class of polyhedra can be captured
is we allow inequalities of the form $\transpos{a}x \leq 0$.
In order to find the corresponding ``$\s{V}$-definition''
we will need to add positive combinations of vectors
to convex combinations of points. Intuitively, these
vectors correspond to ``points at infinity''.

\medskip
We will see shortly that if $A$ is an 
$\s{H}$-polytope and if $O\in\> \interio{A}$,
then $A^*$ is also an $\s{H}$-polytope.
For this, we will prove first that an $\s{H}$-polytope is
a $\s{V}$-polytope. This requires taking a closer look at
polyhedra. 

\medskip
Note that some of the hyperplanes cutting out a polyhedron
may be redundant. If \\
$A = \bigcap_{i = 1}^t C_i$ is a polyhedron
(where each closed half-space, $C_i$, is associated with a
hyperplane, $H_i$, so that $\partial C_i = H_i$),
we say that  {\it $\bigcap_{i = 1}^t C_i$ is an irredundant
decomposition of $A$\/} if $A$ cannot be expressed as 
$A = \bigcap_{i = 1}^m C_i'$ with $m < t$ (for some closed
half-spaces, $C_i'$). The following proposition shows that
the $C_i$ in an irredundant decomposition of $A$ are 
uniquely determined by $A$.

\begin{prop}
\label{irredund}
Let $A$ be a polyhedron with nonempty
interior and assume that \\
$A  = \bigcap_{i = 1}^t C_i$ is
an irredundant decomposition of $A$. Then,
\begin{enumerate}
\item[(i)] 
Up to order, the $C_i$'s are uniquely determined by $A$.
\item[(ii)] 
If $H_i = \partial C_i$ is the boundary of $C_i$, then
$H_i \cap A$ is a polyhedron with nonempty interior in $H_i$,
denoted $\mathrm{Facet}_i\, A$, and called a facet of $A$.
\item[(iii)] 
We have $\partial A = \bigcup_{i = 1}^t\> \mathrm{Facet}_i\, A$,
where the union is irredundant, i.e., $\mathrm{Facet}_i\, A$
is not a subset of $\mathrm{Facet}_j\, A$, for all $i \not= j$.
\end{enumerate}
\end{prop}

\proof
(ii) Fix any $i$ and consider $A_i = \bigcap_{j \not= i} C_j$.
As $A = \bigcap_{i = 1}^t C_i$ is an irredundant decomposition, 
there is some $x\in A_i - C_i$. Pick any $a\in\ \interio{A}$.
By Lemma \ref{lem1}, we get $b = [a, x]\cap H_i \in\> \interio{{A_i}}$,
so $b$ belongs to the interior of $H_i \cap A_i$ in $H_i$.

\medskip
(iii)
As $\partial A = A\> -  \interio{A}\> = A\cap (\interio{A})^c$
(where $B^c$ denotes the complement of a subset $B$ of $\eucreal^n$)
and $\partial C_i = H_i$,
we get
\begin{eqnarray*}
\partial A & = & \Biggl(\bigcap_{i = 1}^t C_i\Biggr)\> - 
\interio{{\Biggl(\bigcap_{j = 1}^t C_j\Biggr)}} \\
& = & \Biggl(\bigcap_{i = 1}^t C_i\Biggr) - 
\Biggl(\bigcap_{j = 1}^t \interio{C_j}\Biggr) \\
& = & \Biggl(\bigcap_{i = 1}^t C_i\Biggr) \cap 
\Biggl(\bigcap_{j = 1}^t \interio{C_j}\Biggr)^c \\
& = & \Biggl(\bigcap_{i = 1}^t C_i\Biggr) \cap 
\Biggl(\bigcup_{j = 1}^t\>  (\interio{C_j})^c\Biggr) \\
& = & \bigcup_{j = 1}^t  \Biggl( 
\Bigl(\bigcap_{i = 1}^t C_i\Bigr)\>\cap (\interio{C_j})^c\Biggr) \\
& = & \bigcup_{j = 1}^t  \Biggl(\partial C_j 
\cap \Bigl(\bigcap_{i \not= j} C_i\Bigr)\Biggr) \\
& = & \bigcup_{j = 1}^t \, (H_j\cap A) = \bigcup_{j = 1}^t\> 
\mathrm{Facet}_j\, A.
\end{eqnarray*}
If we had $\mathrm{Facet}_i\, A \subseteq \mathrm{Facet}_j\, A$,
for some $i \not= j$, then, by (ii), as the affine hull of
$\mathrm{Facet}_i\, A$ is $H_i$ and the affine hull of
$\mathrm{Facet}_j\, A$  is $H_j$, we would have
$H_i \subseteq  H_j$, a contradiction.

\medskip
(i)
As the decomposition is irredundant, the $H_i$ are pairwise distinct.
Also, by (ii), each facet, $\mathrm{Facet}_i\, A$, has dimension
$d - 1$ (where $d = \mathrm{dim}\, A$). Then, 
in  (iii), we can show that the decomposition of
$\partial A$ as a union of potytopes of dimension $d - 1$
whose pairwise nonempty intersections have dimension at most $d - 2$
(since they are contained in pairwise distinct hyperplanes)
is unique up to permutation. Indeed, assume that
\[
\partial A = F_1\cup \cdots \cup F_m = G_1\cup \cdots \cup G_n,
\]
where the $F_i$'s and $G_j'$ are polyhedra of dimension $d - 1$
and each of the unions is irredundant.
Then, we claim that for each $F_i$, there is some $G_{\varphi(i)}$
such that $F_i \subseteq G_{\varphi(i)}$. If not, $F_i$ would be
expressed as a union
\[
F_i = (F_i\cap G_{i_1}) \cup \cdots \cup (F_i\cap G_{i_k}) 
\]
where $\mathrm{dim}(F_i\cap G_{i_j}) \leq d - 2$, since the
hyperplanes containing $F_i$ and the $G_j$'s are pairwise distinct,
which is absurd, since $\mathrm{dim}(F_i) = d - 1$.
By symmetry, for each $G_j$, there is some $F_{\psi(j)}$ such that
$G_j \subseteq F_{\psi(j)}$. But then,
$F_i \subseteq F_{\psi(\varphi(i))}$ for all $i$ and
$G_j \subseteq G_{\varphi(\psi(j))}$ for all $j$
which implies $\psi(\varphi(i)) = i$ for all $i$
and $\varphi(\psi(j)) = j$ for all $j$ since
the unions are irredundant. Thus,
$\varphi$ and $\psi$ are mutual inverses and
the $B_j$'s are just a permutation of the $A_i$'s,
as claimed. Therefore, the facets, $\mathrm{Facet}_i\, A$,
are uniquely determined by $A$ and 
so are the hyperplanes, $H_i = \mathrm{aff}(\mathrm{Facet}_i\, A)$, and
the half-spaces, $C_i$, that they determine.
$\bigsquare$

\medskip
As a consequence, if $A$ is a polyhedron,
then so are its facets and the same holds for $\s{H}$-polytopes.
If $A$ is an $\s{H}$-polytope and
$H$ is a hyperplane with $H\>\cap \interio{A}\> \not= \emptyset$,
then $H\cap A$ is an $\s{H}$-polytope whose facets are of the form $H\cap F$,
where $F$ is a facet of $A$.

\medskip
We can use induction and define $k$-faces, for
$0 \leq k \leq n - 1$.

\begin{defin}
\label{kfaces}
{\em
Let $A\subseteq \eucreal^n$ be a polyhedron with nonempty interior.
We define a  {\it $k$-face of $A$\/} to be a facet of a $(k+1)$-face of 
$A$, for $k = 0, \ldots, n - 2$,
where an {\it $(n - 1)$-face\/} is just a facet of $A$.
The $1$-faces are called {\it edges\/}.
Two $k$-faces are {\it adjacent\/} if their
intersection is a $(k - 1)$-face.
}
\end{defin}

\medskip
The polyhedron $A$ itself is also called 
a {\it face\/} (of itself) or {\it $n$-face\/} and the $k$-faces of $A$
with $k \leq n - 1$ are called {\it proper faces\/} of $A$. 
If $A = \bigcap_{i = 1}^t C_i$ is an irredundant decomposition of
$A$ and $H_i$ is the boundary of $C_i$, then
the hyperplane, $H_i$, is called the {\it supporting hyperplane\/}
of the facet $H_i\cap A$.
We suspect that the $0$-faces of a polyhedron are vertices
in the sense of Definition \ref{vertexpt}. This is true and, in fact,
the vertices of a polyhedron coincide with its extreme points
(see Definition \ref{extremept}).

\begin{prop}
\label{vertexp1}
Let $A\subseteq \eucreal^n$ be a polyhedron with nonempty interior.
\begin{enumerate}
\item[(1)]
For any point, $a\in \partial A$, on the boundary of $A$,
the intersection of all the supporting hyperplanes to $A$ at $a$
coincides with the intersection of all the faces that contain $a$.
In particular, points of order $k$ of $A$ are those points
in the relative interior of the $k$-faces of $A$%
\footnote{
Given a convex set, $S$, in $\affreal^n$, its {\it relative interior\/}
is its interior in the affine hull of $S$ (which might be of dimension
strictly less than $n$).}; thus, 
$0$-faces coincide with the vertices of $A$.
\item[(2)]
The vertices of $A$ coincide with the extreme points of $A$.
\end{enumerate}
\end{prop}

\proof
(1)
If $H$ is a supporting hyperplane to $A$ at $a$, then, one of the
half-spaces, $C$, determined by $H$, satisfies $A = A\cap C$.
It follows from Proposition \ref{irredund} that if $H\not= H_i$
(where the hyperplanes $H_i$ are the supporting hyperplanes of the 
facets of $A$), then $C$ is redundant, from which (1) follows.
 
\medskip
(2)
If $a\in \partial A$ is not extreme, then $a\in [y, z]$,
where $y, z\in \partial A$. However, this implies that
$a$ has order $k\geq 1$, i.e, $a$ is not a vertex.
$\bigsquare$

\section{The Equivalence of  $\s{H}$-Polytopes and $\s{V}$-Polytopes}
\label{sec5b}
We are now ready for the theorem showing the equivalence
of $\s{V}$-polytopes and $\s{H}$-polytopes. This is a nontrivial
theorem usually attributed to Weyl and Minkowski
(for example, see Barvinok \cite{Barvinok}).

\begin{thm}
\label{equivpoly} (Weyl-Minkowski)
If $A$ is an $\s{H}$-polytope, then $A$ has a finite number
of extreme points (equal to its vertices) and $A$ is the convex hull
of its set of vertices; thus, an $\s{H}$-polytope is a $\s{V}$-polytope.
Moreover, $A$ has a finite number of $k$-faces
(for $k = 0, \ldots, d - 2$, where $d = \dimm(A)$).
Conversely, the convex hull of a finite set of points is an
$\s{H}$-polytope. As a consequence, a $\s{V}$-polytope
is an $\s{H}$-polytope.
\end{thm}

\proof
By restricting ourselves to the affine hull of $A$
(some $\eucreal^d$, with $d \leq n$) we may assume
that $A$ has nonempty interior.
Since an $\s{H}$-polytope has finitely many facets,
we deduce by induction that an $\s{H}$-polytope
has a finite number of $k$-faces, for
$k = 0, \ldots, d - 2$. In particular, an $\s{H}$-polytope has 
finitely many vertices. By proposition \ref{vertexp1},
these vertices are the extreme points of $A$
and since an  $\s{H}$-polytope is compact and convex, by the theorem
of Krein and Milman (Theorem \ref{KreinMilman}),
$A$ is the convex hull of its set of vertices.

\medskip
Conversely, again, we may assume that $A$ has nonempty interior
by restricting ourselves to the affine hull of $A$.
Then, pick an origin, $O$,
in the interior of $A$ and consider the dual, $A^*$, of $A$.
By Proposition \ref{polytopdual}, the convex set $A^*$ is
an $\s{H}$-polytope. By the first part of the proof 
of Theorem \ref{equivpoly}, the $\s{H}$-polytope, $A^*$, is the convex
hull of its vertices. Finally, as the hypotheses of
Proposition \ref{propdual1} and Proposition \ref{polytopdual} (again)
hold, we deduce that $A = A^{**}$ is an $\s{H}$-polytope.
$\bigsquare$

\medskip
In view of Theorem \ref{equivpoly}, we are justified
in dropping the $\s{V}$ or $\s{H}$ in front of polytope,
and will do so from now on.
Theorem \ref{equivpoly} has some interesting corollaries
regarding  the dual of a polytope.

\begin{cor}
\label{dualpoly1}
If $A$ is any polytope in $\eucreal^n$ such that
the interior of $A$ contains the origin, $O$, then the dual,
$A^*$, of $A$ is also a polytope whose interior contains $O$ and
$A^{**} = A$.
\end{cor}

\begin{cor}
\label{dualpoly2}
If $A$ is any polytope in $\eucreal^n$ 
whose interior contains the origin, $O$, then the $k$-faces
of $A$ are in bijection with the $(n - k - 1)$-faces of the dual polytope, 
$A^*$. This correspondence is as follows:
If $Y = \mathrm{aff}(F)$  is the 
$k$-dimensional subspace determining the $k$-face, $F$, of $A$
then the subspace, $Y^* = \mathrm{aff}(F^*)$, 
determining the corresponding face, $F^*$, of
$A^*$, is the intersection of the polar hyperplanes of points in $Y$.
\end{cor} 

\proof
Immediate from Proposition \ref{vertexp1} and Proposition
\ref{propdual2}.
$\bigsquare$

\medskip
We also have the following proposition whose proof
would not be that simple if we only had the notion
of an $\s{H}$-polytope (as a matter of fact, there is a way
of proving Theorem \ref{equivpoly}
using Proposition \ref{mappoly})

\begin{prop}
\label{mappoly}
If $A\subseteq \eucreal^n$ is a polytope
and $\mapdef{f}{\eucreal^n}{\eucreal^m}$
is an affine map, then $f(A)$ is a polytope in $\eucreal^m$.
\end{prop}

\proof
Immediate, since an $\s{H}$-polytope is a $\s{V}$-polytope and since
affine maps send convex sets to convex sets.
$\bigsquare$

\medskip
The reader should check that the Minkowski sum of polytopes
is a polytope.

\medskip
We were able to give a short proof of Theorem \ref{equivpoly}
because we relied on a powerful theorem, namely,
Krein and Milman. A drawback of 
this approach is that it bypasses the interesting 
and important problem of designing algorithms  for finding
the vertices of an $\s{H}$-polyhedron from the sets of
inequalities defining it. A method for doing this is
Fourier-Motzkin elimination, see Ziegler \cite{Ziegler97}
(Chapter 1) and Section \ref{sec5c}.
This is also a special case of
{\it linear programming\/}.

\medskip
It is also possible to generalize the notion of
$\s{V}$-polytope to polyhedra using the notion of cone
and to generalize the equivalence theorem to $\s{H}$-polyhedra
and $\s{V}$-polyhedra. 

\section{The Equivalence of  $\s{H}$-Polyhedra and $\s{V}$-Polyhedra}
\label{sec5c}
The equivalence of $\s{H}$-polytopes  and $\s{V}$-polytopes
can be generalized to polyhedral sets, {\it i.e.\/}, finite
intersections of closed half-spaces that are not necessarily bounded.
This equivalence was first proved by Motzkin in the early 1930's.
It can be proved in several ways, some involving cones.

\begin{defin}
\label{conedef}
{\em
Let $\affs$ be any affine Euclidean space of finite dimension, $n$
(with associated vector space, $\vector{\affs}$). 
A subset, $C \subseteq \vector{\affs}$, is a {\it cone\/}
if $C$ is closed under linear combinations involving only
nonnnegative scalars called {\it positive combinations\/}.
Given a subset, 
$V \subseteq \vector{\affs}$, the {\it conical hull\/} or {\it positive
hull\/} of $V$ is the set
\[
\mathrm{cone}(V) = \Bigl\{\sum_{I} \lambda_i v_i \mid \{v_i\}_{i\in I} 
\subseteq V,\>
\lambda_i \geq 0\quad\hbox{for all $i\in I$}\Bigr\}.
\]
A {\it $\s{V}$-polyhedron\/} or {\it polyhedral set\/}
is a subset, $A \subseteq \affs$,
such that
\[
A = \convclo(Y) + \mathrm{cone}(V) 
= \{a + v \mid a\in \convclo(Y),\> v\in \mathrm{cone}(V)\}, 
\]
where $V\subseteq \vector{\affs}$ is a finite set of vectors and
$Y\subseteq \affs$ is a finite set of points.

\medskip
A set, $C\subseteq \vector{\affs}$, is a {\it $\s{V}$-cone\/} or
{\it polyhedral cone\/}
if $C$ is the positive hull of a finite set of vectors, that is,
\[
C = \mathrm{cone}(\{u_1, \ldots, u_p\}),
\]
for some vectors, $u_1, \ldots, u_p\in \vector{\affs}$.
An {\it $\s{H}$-cone\/} is any subset of $\vector{\affs}$ given by a finite 
intersection of closed half-spaces
cut out by hyperplanes through $0$.
}
\end{defin}

\medskip
The positive hull, $\mathrm{cone}(V)$, of $V$ is also denoted
$\mathrm{pos}(V)$. Observe that a $\s{V}$-cone can be viewed as
a polyhedral set for which $Y = \{O\}$, a single point. 
However, if we take the point $O$ as the origin, we may view 
a $\s{V}$-polyhedron, $A$, for which $Y = \{O\}$, as a $\s{V}$-cone.
We will switch back and forth between these two views of cones
as we find it convenient, this should not cause any confusion.
In this section, we favor
the view that $\s{V}$-cones are special kinds of $\s{V}$-polyhedra.
As a consequence, a ($\s{V}$ or $\s{H}$)-cone
always contains $0$, sometimes called an {\it apex\/} of the cone.

\medskip
A set of the form $\{a + tu \mid t\geq 0\}$, where $a\in \affs$
is a point and $u\in \vector{\affs}$ is a nonzero vector, is called a
{\it half-line\/} or {\it ray\/}. Then, we see that
a $\s{V}$-polyhedron, $A = \mathrm{conv}(Y) + \mathrm{cone}(V)$,
is the convex hull of the union of a finite set of points
with a finite set of rays. In the case of a $\s{V}$-cone,
all these rays meet in a common point, an apex of the cone.

\medskip
Propositions \ref{polydef1} and \ref{polydef2}
generalize easily to $\s{V}$-polyhedra and cones.

\begin{prop}
\label{polydef1b}
Let $A\subseteq \affs$ be a $\s{V}$-polyhedron or
an $\s{H}$-polyhedron, let $E = \mathrm{aff}(A)$ 
be the affine hull of $A$ in $\affs$ (with the Euclidean
structure on $E$ induced by the Euclidean structure on $\affs$)
and write $d = \mathrm{dim}(E)$.
Then, the following assertions hold:
\begin{enumerate}
\item[(1)]
The set, $A$, is a $\s{V}$-polyhedron in $E$ (i.e., viewed as a subset of $E$)
iff $A$ is a $\s{V}$-polyhedron in $\affs$.
\item[(2)]
The set, $A$, is an $\s{H}$-polyhedron in $E$ 
(i.e., viewed as a subset of $E$)
iff $A$ is an $\s{H}$-polyhedron in $\affs$.
\end{enumerate}
\end{prop}

\proof
We already proved (2) in Proposition \ref{polydef1}.
For (1), observe that the direction, $\vector{E}$, of $E$ is a linear subspace
of $\vector{\affs}$. Consequently, $E$ is closed under
affine combinations and $\vector{E}$ is closed under linear combinations
and the result follows immediately. 
$\bigsquare$

\begin{prop}
\label{polydef2b}
Given any two affine Euclidean spaces, $E$ and $F$,
if $\mapdef{h}{E}{F}$ is any affine map then:
\begin{enumerate}
\item[(1)]
If $A$ is any $\s{V}$-polyhedron in $E$, then $h(E)$ is
a $\s{V}$-polyhedron in $F$.
\item[(2)]
If $\mapdef{g}{\vector{E}}{\vector{F}}$ is any linear map
and if $C$ is any $\s{V}$-cone in $\vector{E}$, then $g(C)$ is
a $\s{V}$-cone in $\vector{F}$.
\item[(3)]
If $h$ is bijective and $A$ is any $\s{H}$-polyhedron in $E$, 
then $h(E)$ is an $\s{H}$-polyhedron in $F$.
\end{enumerate}
\end{prop}

\proof
We already proved (3) in Proposition \ref{polydef2}.
For (1), using the fact that 
$h(a + u) = h(a) + \vector{h}(u)$ for any affine map, $h$,
where $\vector{h}$ is
the linear map associated with $h$, we get
\[
h(\mathrm{conv}(Y) + \mathrm{cone}(V)) =
\mathrm{conv}(h(Y)) + \mathrm{cone}(\vector{h}(V)).
\]
For (1), as $g$ is linear, we get
\[
g(\mathrm{cone}(V)) = \mathrm{cone}(g(V)),
\]
establishing the proposition.
$\bigsquare$

\medskip
Propositions \ref{polydef1b}  and \ref{polydef2b}
allow us to assume that $\affs = \eucreal^d$ and
that our ($\s{V}$ or $\s{H}$)-polyhedra, $A\subseteq \eucreal^d$,
have nonempty interior ({\it i.e.\/} $\mathrm{dim}(A) = d$).

\medskip
The generalization of Theorem \ref{equivpoly} is that
every $\s{V}$-polyhedron, $A$, is an $\s{H}$-polyhedron and
conversely. At first glance, it may seem that
there is a small problem  when $A = \eucreal^d$.
Indeed,  Definition \ref{conedef} allows the possibility
that $\mathrm{cone}(V) = \eucreal^d$ for some finite subset,
$V \subseteq \reals^d$. This is because it is possible to generate
a basis of $\reals^d$ using finitely many positive combinations. 
On the other hand the definition of an $\s{H}$-polyhedron, $A$,
(Definition \ref{polytopedef}) assumes that $A\subseteq \eucreal^n$ 
is cut out by {\it at least one\/} hyperplane.
So, $A$ is always contained in some half-space of $\eucreal^n$ and
strictly speaking, $\eucreal^n$ is not an $\s{H}$-polyhedron!
The simplest way to circumvent this difficulty is to decree
that $\eucreal^n$ itself is a polyhedron, which we do.


\medskip
Yet another solution is to assume that, unless stated otherwise,
every finite set of vectors, $V$, that we consider
when defining a polyhedron
has the property that there is some hyperplane, $H$, through the
origin so that all the vectors in $V$ lie in only one of the
two closed half-spaces determined by $H$.
But then,  the polar dual
of a polyhedron can't be a single point! Therefore,
we stick to our decision that $\eucreal^n$ itself is a polyhedron.

\medskip
To prove the equivalence of $\s{H}$-polyhedra and $\s{V}$-polyhedra,
Ziegler proceeds as follows: First, he shows that
the equivalence of   $\s{V}$-polyhedra and  $\s{H}$-polyhedra
reduces to the equivalence of $\s{V}$-cones and $\s{H}$-cones
using an ``old trick'' of projective geometry, namely, ``homogenizing''
\cite{Ziegler97} (Chapter 1).
Then, he uses two dual versions of Fourier-Motzkin elimination
to pass from $\s{V}$-cones to $\s{H}$-cones and conversely.

\medskip
Since the homogenization method is an important technique
we will describe it in some detail. However, 
it turns out that the double dualization technique
used in the proof of Theorem \ref{equivpoly} can be easily
adapted to prove that every $\s{V}$-polyhedron is an
$\s{H}$-polyhedron. Moreover, it can also be used to prove
that  every $\s{H}$-polyhedron is a $\s{V}$-polyhedron!
So, we will not describe the
version of Fourier-Motzkin elimination used to go from
$\s{V}$-cones to $\s{H}$-cones. However, we will present
the Fourier-Motzkin elimination method used to go from
$\s{H}$-cones to $\s{V}$-cones.

\medskip
Here is the generalization of Proposition \ref{polytopdual}
to polyhedral sets. In order to avoid confusion between 
the origin of $\eucreal^d$ and the center of
polar duality
we will denote the origin by $O$ and the center of
our polar duality by $\Omega$.
Given any nonzero vector,
$u\in \reals^d$, let $u^{\dagger}_-$ be the closed half-space
\[
u^{\dagger}_- = \{x\in \reals^d \mid x\cdot u \leq 0\}.
\]
In other words, $u^{\dagger}_-$ is the closed half-space bounded
by the hyperplane through $\Omega$ normal to $u$ and on the
``opposite side'' of $u$.

\begin{prop}
\label{polyhedraldual}
Let $A = \convclo(Y) + \mathrm{cone}(V)\subseteq \eucreal^d$ 
be a $\s{V}$-polyhedron
with $Y = \{y_1, \ldots$, $y_p\}$ and
$V = \{v_1, \ldots, v_q\}$.  Then, for any point,
$\Omega$, if $A\not= \{\Omega\}$, then
the polar dual, $A^*$, of $A$ w.r.t. $\Omega$ is the
$\s{H}$-polyhedron given by
\[
A^* = \bigcap_{i = 1}^p (y_i^{\dagger})_- \cap
\bigcap_{j = 1}^q (v_j^{\dagger})_-. 
\]
Furthermore, if $A$ has nonempty interior
and $\Omega$ belongs to the interior of $A$, then
$A^*$ is bounded, that is, $A^*$ is an $\s{H}$-polytope.
If $A = \{\Omega\}$, then $A^*$ is the special polyhedron, $A^* = \eucreal^d$.
\end{prop}

\proof
By definition of $A^*$ w.r.t. $\Omega$, we have
\begin{eqnarray*}
A^* & = & \left\{x\in \eucreal^d\> \left| \>\libvecbo{\Omega}{x}\cdot
\libvecbo{\Omega}{\left(
\sum_{i = 1}^p \lambda_i y_i + \sum_{j = 1}^q \mu_j v_j 
\right)} \leq 1,\> \lambda_i \geq 0,\> 
\sum_{i = 1}^p \lambda_i = 1,\> \mu_j \geq 0
\right. \right\} \\
& = & \left\{x\in \eucreal^d\> \left| \>
\sum_{i = 1}^p \lambda_i \libvecbo{\Omega}{x}\cdot\libvecbo{\Omega}{y_i} +
 \sum_{j = 1}^q \mu_j \libvecbo{\Omega}{x}\cdot  v_j \leq 1,\> 
\lambda_i \geq 0,\> 
\sum_{i = 1}^p \lambda_i = 1,\> \mu_j \geq 0
\right. \right\}.
\end{eqnarray*}

\medskip
When $\mu_j = 0$ for $j = 1, \ldots, q$, we get
\[
\sum_{i = 1}^p \lambda_i \libvecbo{\Omega}{x}\cdot\libvecbo{\Omega}{y_i}
\leq 1, \quad \lambda_i \geq 0,\> \sum_{i = 1}^p \lambda_i = 1
\]
and we check that
\begin{eqnarray*}
 \left\{x\in \eucreal^d\> \left| \>
\sum_{i = 1}^p \lambda_i \libvecbo{\Omega}{x}\cdot\libvecbo{\Omega}{y_i}
\leq 1, \> \lambda_i \geq 0,\> \sum_{i = 1}^p \lambda_i = 1
\right. \right\} & = &
\bigcap_{i = 1}^p \{x\in \eucreal^d \mid 
\libvecbo{\Omega}{x}\cdot\libvecbo{\Omega}{y_i} \leq 1\} \\
& = & \bigcap_{i = 1}^p (y_i^{\dagger})_- .
\end{eqnarray*}

\medskip
The points in $A^*$ must also satisfy the conditions
\[
 \sum_{j = 1}^q \mu_j \libvecbo{\Omega}{x}\cdot  v_j \leq 1 - \alpha,\quad 
 \mu_j \geq 0, \> \mu_j > 0\> \hbox{for some $j$, $1\leq j \leq q$},
\] 
with $\alpha \leq 1$ (here 
$\alpha = \sum_{i = 1}^p \lambda_i 
\libvecbo{\Omega}{x}\cdot\libvecbo{\Omega}{y_i}$).
In particular, for every $j\in \{1, \ldots, q\}$,
if we set $\mu_k = 0$ for 
$k \in \{1, \ldots, q\} - \{j\}$, we should have
\[
\mu_j \libvecbo{\Omega}{x}\cdot  v_j \leq 1 - \alpha
\quad\hbox{for all}\quad \mu_j > 0,
\]
that is,
\[
 \libvecbo{\Omega}{x}\cdot  v_j \leq \frac{1 - \alpha}{\mu_j}
\quad\hbox{for all}\quad \mu_j > 0,
\]
which is equivalent to
\[
 \libvecbo{\Omega}{x}\cdot  v_j \leq 0.
\]
Consequently, if $x\in A^*$, we must also have
\[
x \in \bigcap_{j = 1}^q  \{x\in \eucreal^d \mid 
 \libvecbo{\Omega}{x}\cdot  v_j \leq 0\} =
\bigcap_{j = 1}^q   (v_j^{\dagger})_-.
\]
Therefore,
\[
A^* \subseteq  \bigcap_{i = 1}^p (y_i^{\dagger})_- \cap
\bigcap_{j = 1}^q (v_j^{\dagger})_-. 
\]
However, the reverse inclusion is obvious and thus,
we have equality. If $\Omega$ belongs to the interior
of $A$, we know from Proposition \ref{propdual1} that $A^*$ is bounded.
Therefore, $A^*$ is indeed an $\s{H}$-polytope of the above form.
$\bigsquare$

\medskip
It is  fruitful to restate Proposition \ref{polyhedraldual}
in terms of matrices (as we did for Proposition \ref{polytopdual}).
First, observe that
\[
(\mathrm{conv}(Y) + \mathrm{cone}(V))^* = 
(\mathrm{conv}(Y\cup \{\Omega\}) + \mathrm{cone}(V))^*. 
\]
If we pick $\Omega$ as an origin then we can represent the points
in $Y$ as vectors.
The old origin is still denoted $O$
and $\Omega$ is now denoted  $0$. The zero vector is denoted
$\mathbf{0}$. 

\medskip
If $A = \mathrm{conv}(Y) + \mathrm{cone}(V) \not= \{0\}$,
let $Y$ be the $d\times p$ matrix whose $i^{\mathrm{th}}$ column
is $y_i$ and let $V$ is the $d\times q$ matrix whose  $j^{\mathrm{th}}$ column
is $v_j$. Then Proposition  \ref{polyhedraldual} says that
\[
(\mathrm{conv}(Y) + \mathrm{cone}(V))^* = 
\{x\in \reals^d \mid \transpos{Y}x \leq \mathbf{1},\>
\transpos{V}x \leq \mathbf{0}\}.
\]
We write $P(\transpos{Y}, \mathbf{1}; \transpos{V}, \mathbf{0}) =
\{x\in \reals^d \mid \transpos{Y}x \leq \mathbf{1},\>
\transpos{V}x \leq \mathbf{0}\}$.

\medskip
If  $A = \mathrm{conv}(Y) + \mathrm{cone}(V) = \{0\}$,
then both $Y$ and $V$ must be zero matrices but then,
the inequalities $\transpos{Y}x \leq \mathbf{1}$ and
$\transpos{V}x \leq \mathbf{0}$ are trivially satisfied
by all $x\in \eucreal^d$, so even in this case we  have
\[
\eucreal^d =
(\mathrm{conv}(Y) + \mathrm{cone}(V))^* =
P(\transpos{Y}, \mathbf{1}; \transpos{V}, \mathbf{0}). 
\]
The converse of Proposition  \ref{polyhedraldual}
also holds as shown below.

\begin{prop}
\label{polytopdual3}
Let  $\{y_1, \ldots, y_p\}$  be any set of points in
$\eucreal^d$ and 
let  $\{v_1, \ldots, v_q\}$ be any set of nonzero vectors in $\reals^d$.
If $Y$ is the $d\times p$ matrix whose $i^{\mathrm{th}}$ column
is $y_i$ and $V$ is the $d\times q$ matrix whose  $j^{\mathrm{th}}$ column
is $v_j$, then
\[
 (\mathrm{conv}(\{y_1, \ldots, y_p\}) 
\cup \mathrm{cone}(\{v_1, \ldots, v_q\}))^*
=  P(\transpos{Y}, \mathbf{1}; \transpos{V}, \mathbf{0}),
\]
with 
$P(\transpos{Y}, \mathbf{1}; \transpos{V}, \mathbf{0}) = 
\{x\in \reals^d \mid \transpos{Y}x \leq \mathbf{1},\>
\transpos{V}x \leq \mathbf{0}\}$. 

\medskip
Conversely, given any $p\times d$ matrix, $Y$, and any
$q\times d$ matrix, $V$, we have
\[
P(Y, \mathbf{1}; V, \mathbf{0})^* = 
 \mathrm{conv}(\{y_1, \ldots, y_p\}\cup \{0\}) 
\cup \mathrm{cone}(\{v_1, \ldots, v_q\}),
\]
where $y_i\in\reals^n$ is the  $i^{\mathrm{th}}$ row of $Y$ and
$v_j\in\reals^n$ is the  $j^{\mathrm{th}}$ row of $V$ 
or, equivalently,
\[
P(Y, \mathbf{1}; V, \mathbf{0})^* = 
\{x\in \reals^d \mid x = \transpos{Y}u + \transpos{V}t, \>
u \in \reals^p, \> t\in \reals^q,\> u, t \geq 0,\>
\mathbb{I}u = 1\}, 
\]
where $\mathbb{I}$ is the row vector of length $p$ whose
coordinates are all equal to $1$. 
\end{prop}

\proof
Only the second part needs a proof. Let
\[
B =   \mathrm{conv}(\{y_1, \ldots, y_p\}\cup \{0\}) 
\cup \mathrm{cone}(\{v_1, \ldots, v_q\}),
\]
where $y_i\in\reals^p$ is the  $i^{\mathrm{th}}$ {\it row\/} of $Y$
and  $v_j\in\reals^q$ is the  $j^{\mathrm{th}}$ {\it row\/} of $V$.
Then, by the first part,
\[
B^* = P(Y, \mathbf{1}; V, \mathbf{0}).
\]
As $0\in B$, by Proposition \ref{propdual1}, we have
$B = B^{**} =   P(Y, \mathbf{1}; V, \mathbf{0})$,
as claimed.
$\bigsquare$

\medskip
Proposition \ref{polytopdual3} has the following
important Corollary:

\begin{prop}
\label{polytopdual4}
The following assertions hold:
\begin{enumerate}
\item[(1)]
The polar dual, $A^*$, of every  $\s{H}$-polyhedron, is a $\s{V}$-polyhedron.
\item[(2)]
The polar dual, $A^*$, of every  $\s{V}$-polyhedron, is an $\s{H}$-polyhedron.
\end{enumerate}
\end{prop}

\proof
(1)
We may assume that $0\in A$, in which case,
$A$ is of the form $A = P(Y, \mathbf{1}; V, \mathbf{0})$.
By the second part of Proposition \ref{polytopdual3}, 
$A^*$ is a $\s{V}$-polyhedron.

\medskip
(2)
This is the first part of Proposition \ref{polytopdual3}.
$\bigsquare$

\medskip
We can now use  Proposition \ref{polyhedraldual}, Proposition \ref{propdual1}
and Krein and Millman's Theorem
to prove that every $\s{V}$-polyhedron is an $\s{H}$-polyhedron.

\begin{prop}
\label{VtoHpolyhedral}
Every  $\s{V}$-polyhedron, $A$,
is an $\s{H}$-polyhedron.
Furthermore, if $A\not= \eucreal^d$, then $A$ is of the form
$A = P(Y, \mathbf{1})$.
\end{prop}

\proof
Let $A$ be a $\s{V}$-polyhedron of dimension $d$.
Thus, $A\subseteq \eucreal^d$ has nonempty interior so we can
pick some point,  $\Omega$, in the interior of $A$.
If  $d = 0$, then $A = \{0\} = \eucreal^{0}$ and we are done.
Otherwise, by Proposition \ref{polyhedraldual},
the polar dual, $A^*$, of $A$ w.r.t. $\Omega$ is an $\s{H}$-polytope.
Then, as in the proof of Theorem \ref{equivpoly}, using
Krein and Millman's Theorem we deduce that
$A^*$ is a $\s{V}$-polytope. Now, as $\Omega$ belongs to 
$A$, by Proposition \ref{propdual1} (even if $A$ is not bounded)
we have $A = A^{**}$ and 
by Proposition \ref{polytopdual} (or
Proposition \ref{polyhedraldual})
we conclude that $A = A^{**}$ is an $\s{H}$-polyhedron
of the form $A =  P(Y, \mathbf{1})$.
$\bigsquare$

\medskip
Interestingly, we can now prove easily that every
$\s{H}$-polyhedron is a $\s{V}$-polyhedron.

\begin{prop}
\label{HtoVpolyhedral}
Every $\s{H}$-polyhedron is a $\s{V}$-polyhedron.
\end{prop}

\proof
Let $A$ be an $\s{H}$-polyhedron of dimension $d$.
By Proposition \ref{polytopdual4}, the polar dual, $A^*$, of $A$
is a $\s{V}$-polyhedron. By Proposition \ref{VtoHpolyhedral},
$A^*$ is an $\s{H}$-polyhedron and again,
by  Proposition \ref{polytopdual4},
we deduce that  $A^{**} = A$ is a $\s{V}$-polyhedron
($A = A^{**}$ because $0\in A$).
$\bigsquare$

\medskip
Putting together Propositions \ref{VtoHpolyhedral} and
\ref{HtoVpolyhedral} we obtain our main theorem:

\begin{thm} (Equivalence of $\s{H}$-polyhedra and $\s{V}$-polyhedra)
\label{equivpolyhedra}
Every $\s{H}$-polyhedron is a $\s{V}$-polyhedron and
conversely.
\end{thm}

\medskip
Even though we proved the main result of this section,
it is instructive to consider a more computational proof
making use of cones and an elimination method
known as {\it Fourier-Motzkin elimination\/}.

\section[Fourier-Motzkin Elimination and Cones]
{Fourier-Motzkin Elimination and the Polyhedron-Cone
Correspondence}
\label{secFourMotz}
The problem with the converse of Proposition \ref{VtoHpolyhedral}
when $A$ is unbounded ({\it i.e.\/}, not compact) is that  Krein and Millman's
Theorem does not apply. We need to take into account ``points
at infinity'' corresponding to certain vectors. The trick we used
in Proposition \ref{VtoHpolyhedral} is that the polar dual of
a $\s{V}$-polyhedron with nonempty interior is an
{\it $\s{H}$-polytope\/}. This reduction to polytopes
allowed us to use Krein and Millman to convert
an $\s{H}$-polytope to a $\s{V}$-polytope and then again we took
the polar dual.

\medskip
Another trick  is to switch to cones by ``homogenizing''.
Given any subset, $S\subseteq \eucreal^d$, we can form the 
cone, $C(S) \subseteq \eucreal^{d+1}$, by ``placing'' a copy of $S$ 
in the hyperplane, $H_{d+1} \subseteq \eucreal^{d+1}$, of equation
$x_{d+1} = 1$, and drawing all the half-lines
from the origin through any point of $S$. If $S$ is given by $m$
polynomial inequalities of the form
\[
P_i(x_1, \ldots, x_d) \leq b_i,
\]
where $P_i(x_1, \ldots, x_d)$ is a polynomial of total degree $n_i$,
this amounts to forming the new homogeneous inequalities
\[
x_{d+1}^{n_i}P_i\left(\frac{x_1}{x_{d+1}}, \ldots, \frac{x_d}{x_{d+1}}\right) -
b_ix_{d+1}^{n_i} \leq 0
\]
together with $x_{d+1} \geq 0$. In particular, if the $P_i$'s
are linear forms (which means that $n_i = 1$),
then our inequalities are of the form
\[
a_i\cdot x \leq b_i,
\]
where $a_i\in \reals^d$ is some vector
and the homogenized inequalities are
\[
a_i\cdot x - b_ix_{d+1} \leq 0.
\]

\medskip
It will be convenient to formalize the process of
putting a copy of a subset, $S\subseteq \eucreal^d$, 
in the hyperplane, $H_{d+1} \subseteq \eucreal^{d+1}$, of equation
$x_{d+1} = 1$, as follows: For every point, $a\in \eucreal^d$, let
\[
\widehat{a} = \binom{a}{1} \in \eucreal^{d+1}
\]
and let $\widehat{S} = \{\widehat{a} \mid a\in S\}$.
Obviously, the map $S \mapsto \widehat{S}$ is a bijection
between the subsets of $\eucreal^d$   and the subsets of
$H_{d+1}$. We will use this bijection to identify
$S$ and $\widehat{S}$ and use the simpler notation, $S$,
unless confusion arises.
Let's apply this to polyhedra.

\medskip
Let $P\subseteq \eucreal^d$ be an $\s{H}$-polyhedron.
Then, $P$ is cut out by 
$m$ hyperplanes, $H_i$, and for each $H_i$, there is a nonzero
vector, $a_i$, and some $b_i\in \reals$ so that
\[
H_i = \{x\in \eucreal^d \mid a_i\cdot x = b_i\}
\]
and $P$ is given by
\[
P = \bigcap_{i = 1}^m\,  \{x\in \eucreal^d \mid a_i\cdot x \leq  b_i\}.
\]
If $A$ denotes the $m\times d$ matrix whose $i$-th row is
$a_i$ and $b$ is the vector $b = (b_1, \ldots, b_m)$, 
then we can write
\[
P = P(A, b) = \{x\in \eucreal^d \mid Ax \leq b\}.
\]
 
\medskip
We ``homogenize'' $P(A, b)$ as follows:
Let $C(P)$ be the subset of $\eucreal^{d+1}$ defined by
\begin{eqnarray*}
C(P) & = & \left\{
\binom{x}{x_{d+1}}\in \reals^{d+1} \mid Ax \leq x_{d+1}b,\> x_{d+1} \geq 0 
\right\} \\
& = & \left\{
\binom{x}{x_{d+1}}\in \reals^{d+1} \mid Ax -  x_{d+1}b\leq 0,\> -x_{d+1} \leq 0 
\right\}.
\end{eqnarray*}
Thus, we see that $C(P)$ is the $\s{H}$-cone given by the
system of inequalities
\[
\begin{pmatrix}
A & -b \\
0 & -1
\end{pmatrix}\binom{x}{x_{d+1}} \leq \binom{0}{0}
\]
and that
\[
\widehat{P} = C(P)\cap H_{d+1}.
\]
Conversely, if $Q$ is any $\s{H}$-cone in $\eucreal^{d+1}$
(in fact, any $\s{H}$-polyhedron), it is clear that \\ 
$P = Q\cap H_{d+1}$ is an $\s{H}$-polyhedron in $H_{d+1} \approx \eucreal^d$.

\medskip
Let us now assume that 
$P\subseteq \eucreal^d$ is a $\s{V}$-polyhedron,
$P = \mathrm{conv}(Y) + \mathrm{cone}(V)$, where \\
$Y = \{y_1, \ldots, y_p\}$ and $V = \{v_1, \ldots, v_q\}$.
Define $\widehat{Y} = \{\widehat{y}_1, \ldots, \widehat{y}_p\}
\subseteq  \eucreal^{d+1}$, 
and \\
$\widehat{V} =  \{\widehat{v}_1, \ldots, \widehat{v}_q\}
\subseteq  \eucreal^{d+1}$,  by 
\[
\widehat{y}_i = \binom{y_i}{1},
\qquad
\widehat{v}_j = \binom{v_j}{0}.
\] 

\medskip
We check immediately that
\[
C(P) = \mathrm{cone}(\{\widehat{Y}\cup \widehat{V}\})
\]
is a $\s{V}$-cone in $\eucreal^{d+1}$ such that
\[
\widehat{P} = C(P) \cap H_{d+1},
\]
where $H_{d+1}$ is the hyperplane of equation $x_{d+1} = 1$.

\medskip
Conversely, if $C = \mathrm{cone}(W)$ is a
$\s{V}$-cone in $\eucreal^{d+1}$, with $w_{i d+1} \geq 0$
for every $w_i\in W$, we prove next that $P = C\cap H_{d+1}$ is
a $\s{V}$-polyhedron. 

\begin{prop} (Polyhedron--Cone Correspondence)
\label{equivcone}
We have the following correspondence between polyhedra in
$\eucreal^d$ and cones in $\eucreal^{d+1}$:
\begin{enumerate}
\item[(1)]
For any $\s{H}$-polyhedron,  $P\subseteq \eucreal^d$, if
$P = P(A, b) = \{x\in \eucreal^d \mid Ax \leq b\}$, where
$A$ is an $m\times d$-matrix  and $b\in \reals^m$,
then $C(P)$ given by
\[
\begin{pmatrix}
A & -b \\
0 & -1
\end{pmatrix}\binom{x}{x_{d+1}} \leq \binom{0}{0}
\]
is an $\s{H}$-cone in $\eucreal^{d+1}$ and $\widehat{P} = C(P)\cap H_{d+1}$,
where $H_{d+1}$ is the hyperplane of equation $x_{d+1} = 1$.
Conversely, if $Q$ is any $\s{H}$-cone in $\eucreal^{d+1}$
(in fact, any $\s{H}$-polyhedron), then
$P = Q\cap H_{d+1}$ is an $\s{H}$-polyhedron in $H_{d+1} \approx \eucreal^d$.
\item[(2)]
Let $P\subseteq \eucreal^d$ be any $\s{V}$-polyhedron, 
where $P = \mathrm{conv}(Y) + \mathrm{cone}(V)$ with
$Y = \{y_1, \ldots, y_p\}$ and $V = \{v_1, \ldots, v_q\}$.
Define $\widehat{Y} = \{\widehat{y}_1, \ldots, \widehat{y}_p\}
\subseteq  \eucreal^{d+1}$, 
and
$\widehat{V} =  \{\widehat{v}_1, \ldots, \widehat{v}_q\}
\subseteq  \eucreal^{d+1}$,  by 
\[
\widehat{y}_i = \binom{y_i}{1},
\qquad
\widehat{v}_j = \binom{v_j}{0}.
\] 
Then, 
\[
C(P) = \mathrm{cone}(\{\widehat{Y}\cup \widehat{V}\})
\]
is a $\s{V}$-cone in $\eucreal^{d+1}$ such that
\[
\widehat{P} = C(P) \cap H_{d+1},
\]
Conversely, if $C = \mathrm{cone}(W)$ is a
$\s{V}$-cone in $\eucreal^{d+1}$, with $w_{i\, d+1} \geq 0$
for every $w_i\in W$, then $P = C\cap H_{d+1}$ is
a $\s{V}$-polyhedron in $H_{d+1} \approx \eucreal^d$.
\end{enumerate}
\end{prop}

\proof
We already proved  everything
except the last part of the proposition.
Let
\[
\widehat{Y} = \left\{\left. \frac{w_i}{w_{i\, d+1}} \>  
\right| \> w_i \in W,\> w_{i\, d + 1} > 0\right\}
\]
and 
\[
\widehat{V} = \{w_j \in W \mid w_{j\, d + 1} = 0\}.
\]
We claim that
\[
 P = C\cap H_{d+1} = \mathrm{conv}(\widehat{Y}) + \mathrm{cone}(\widehat{V}),
\]
and thus, $P$ is a $\s{V}$-polyhedron.

\medskip
Recall that any element, $z\in C$, can be written as
\[
z = \sum_{k = 1}^s \mu_k w_k, \quad \mu_k \geq 0. 
\]
Thus, we have
\begin{eqnarray*}
z & = & \sum_{k = 1}^s \mu_k w_k \quad \mu_k \geq 0 \\
 & = &  \sum_{w_{i\, d+1} > 0} \mu_i w_i 
+ \sum_{w_{j\, d+1} = 0} \mu_j w_j  \quad \mu_i, \mu_j \geq 0 \\
 & = &  \sum_{w_{i\, d+1} > 0} w_{i\, d+1}\mu_i \frac{w_i}{w_{i\, d+1}} 
+ \sum_{w_{j\, d+1} = 0} \mu_j w_j, \quad \mu_i, \mu_j \geq 0 \\
 & = &  \sum_{w_{i\, d+1} > 0} \lambda_i \frac{w_i}{w_{i\, d+1}} 
+ \sum_{w_{j\, d+1} = 0} \mu_j w_j, \quad \lambda_i, \mu_j \geq 0.
\end{eqnarray*}
Now, $z\in C\cap H_{d+1}$ iff $z_{d+1} = 1$ iff
$\sum_{i = 1}^p \lambda_i = 1$ (where $p$ is the number of
elements in $\widehat{Y}$),
since the $(d+1)^{\mathrm{th}}$ 
coordinate of $\frac{w_i}{w_{i\, d+1}}$ is equal to $1$,
and the above shows that 
\[
 P = C\cap H_{d+1} = \mathrm{conv}(\widehat{Y}) + \mathrm{cone}(\widehat{V}),
\]
as claimed.
$\bigsquare$

\medskip
By Proposition \ref{equivcone}, if $P$ is an $\s{H}$-polyhedron, then
$C(P)$ is an $\s{H}$-cone. If we can prove that every $\s{H}$-cone
is a $\s{V}$-cone, then again, Proposition \ref{equivcone}
shows that $\widehat{P} = C(P)\cap H_{d+1}$ is a $\s{V}$-polyhedron
and so, $P$ is a $\s{V}$-polyhedron.
Therefore, in order to prove that every  $\s{H}$-polyhedron is
a $\s{V}$-polyhedron it suffices to show that
every $\s{H}$-cone is a $\s{V}$-cone.

\medskip
By a similar argument,  Proposition \ref{equivcone} shows that
in order to prove that every  $\s{V}$-polyhedron is
an $\s{H}$-polyhedron it suffices to show that
every $\s{V}$-cone is an $\s{H}$-cone. We will not prove this
direction again since we already have it by Proposition \ref{VtoHpolyhedral}.

\medskip
It remains to prove that every $\s{H}$-cone is a $\s{V}$-cone.
Let $C\subseteq \eucreal^d$ be an $\s{H}$-cone. Then, $C$ is cut out by 
$m$ hyperplanes, $H_i$, through $0$. For each $H_i$, there is a nonzero
vector, $u_i$, so that
\[
H_i = \{x\in \eucreal^d \mid u_i\cdot x = 0\}
\]
and $C$ is given by
\[
C = \bigcap_{i = 1}^m\,  \{x\in \eucreal^d \mid u_i\cdot x \leq  0\}.
\]
If $A$ denotes the $m\times d$ matrix whose $i$-th row is
$u_i$, then we can write
\[
C = P(A, 0) = \{x\in \eucreal^d \mid Ax \leq 0\}.
\]
Observe that $C = C_0(A) \cap H_w$, where
\[
C_0(A) = \left\{\binom{x}{w}\in \reals^{d + m} \> 
\left| \> Ax \leq w \right. \right\}
\]
is an $\s{H}$-cone in $\eucreal^{d + m}$ and 
\[
H_w =  \left\{\binom{x}{w}\in \reals^{d + m} \> 
\left| \> w = 0 \right. \right\}
\]
is an affine subspace in $\eucreal^{d + m}$.

\medskip
We claim that $C_0(A)$ is a $\s{V}$-cone.
This follows by observing that for every $\binom{x}{w}$
satisfying $A x \leq w$, we can write
\[
\binom{x}{w} = \sum_{i = 1}^d |x_i|(\mathrm{sign}(x_i))\binom{e_i}{Ae_i}
+ \sum_{j = 1}^m (w_j - (Ax)_j)\binom{0}{e_j},
\]
and then
\[
C_0(A) = \mathrm{cone}\left(
\left\{
\pm\binom{e_i}{Ae_i} \> \left| \> 1\leq i \leq d  \right. \right\}
\cup
\left\{
\binom{0}{e_j} \> \left| \> 1\leq j \leq m  \right. \right\}
\right).
\]

\medskip
Since $C = C_0(A) \cap H_w$ is now the intersection of a $\s{V}$-cone
with an affine subspace, to prove that $C$ is a $\s{V}$-cone it is enough to
prove that the intersection of a $\s{V}$-cone with a hyperplane
is also a $\s{V}$-cone. For this, we use {\it Fourier-Motzkin
elimination\/}. It suffices to prove the result for a hyperplane, $H_k$,
in $\eucreal^{d + m}$ of equation $y_k = 0$ ($1\leq k \leq d + m$).

\begin{prop} (Fourier-Motzkin Elimination)
\label{FourierMotzkin1}
Say $C = \mathrm{cone}(Y)\subseteq \eucreal^d$ is a $\s{V}$-cone.
Then, the intersection $C\cap H_k$ (where $H_k$ is the hyperplane
of equation $y_k = 0$) is a $\s{V}$-cone,
$C\cap H_k= \mathrm{cone}(Y^{/k})$, with
\[
Y^{/k} = \{y_i \mid y_{i k} = 0\} \cup \{y_{i k}y_j - y_{j k}y_i
\mid y_{i k} >0,\, \> y_{j k} < 0\},
\]
the set of vectors obtained from $Y$ by ``eliminating
the $k$-th coordinate''. Here, each $y_i$ is a vector in $\reals^d$.
\end{prop}

\proof
The only nontrivial direction is to prove that
$C\cap H_k \subseteq \mathrm{cone}(Y^{/k}) $. For this,
consider any $v = \sum_{i = 1}^d  t_i y_i\in C\cap H_k$, with $t_i \geq 0$
and $v_k = 0$. Such a $v$ can be written
\[
v = \sum_{i \mid y_{i k} = 0} t_i y_i +
 \sum_{i \mid y_{i k} > 0} t_i y_i + \sum_{j \mid y_{j k} < 0} t_j y_j 
\]
and as $v_k = 0$, we have 
\[
 \sum_{i \mid y_{i k} > 0} t_i y_{i k} + \sum_{j \mid y_{j k} < 0} t_j y_{j k}
= 0.
\]
If $t_i y_{i k} = 0$  for $i = 1,\ldots, d$, we are done. 
Otherwise, we can write
\[
\Lambda = \sum_{i \mid y_{i k} > 0} t_i y_{i k} =
\sum_{j \mid y_{j k} < 0} -t_j y_{j k} > 0.
\]
Then, 
\begin{eqnarray*}
v  & =  & \sum_{i \mid y_{i k} = 0} t_i y_i + 
 \frac{1}{\Lambda} \sum_{i \mid y_{i k} > 0}\left(
\sum_{j \mid y_{j k} < 0} -t_j y_{j k}
\right) t_iy_i 
 + \frac{1}{\Lambda}  \sum_{j \mid y_{j k} < 0}\left(
\sum_{i \mid y_{i k} > 0} t_i y_{i k}
\right) t_j y_{j} \\
& = &
\sum_{i \mid y_{i k} = 0} t_i y_i + 
  \sum_{
\begin{subarray}{l}
i \mid y_{i k} > 0 \\
j \mid y_{j k} < 0
\end{subarray}
} 
\frac{t_it_j}{\Lambda}
\left(y_{i k} y_j   - y_{j k} y_i\right) .
\end{eqnarray*}
Since the $k^{\mathrm{th}}$ coordinate of  $y_{i k} y_j   - y_{j k} y_i$
is $0$, the above shows that any $v\in C\cap H_{k}$
can be written as a positive combination of vectors in
$Y^{/k}$.
$\bigsquare$

\medskip
As discussed above, Proposition
\ref{FourierMotzkin1} implies (again!)

\begin{cor}
\label{HtoVpolyhedral2}
Every $\s{H}$-polyhedron is a $\s{V}$-polyhedron.
\end{cor}

\medskip
Another way of proving that every  $\s{V}$-polyhedron is an 
$\s{H}$-polyhedron is to use cones. This method is interesting
in its own right so we discuss it briefly.

\medskip
Let $P = \mathrm{conv}(Y) + \mathrm{cone}(V) \subseteq \eucreal^d$
be a $\s{V}$-polyhedron.
We can view $Y$ as a $d\times p$ matrix whose $i$th column
is the $i$th vector in $Y$ and $V$ as $d\times q$ matrix
whose $j$th column is the $j$th vector in $V$. Then, we can write 
\[
P = \{x\in \reals^d \mid (\exists u\in \reals^p)(\exists t\in \reals^d)
(x = Y u + V t,\> u \geq 0,\> \mathbb{I}u = 1,\> t \geq 0)\}, 
\]
where $\mathbb{I}$ is the row vector 
\[
\mathbb{I} = \underbrace{(1, \ldots, 1)}_{p}.
\]
Now, observe that $P$ can be interpreted as the projection of
the $\s{H}$-polyhedron, $\widetilde{P} \subseteq \eucreal^{d + p + q}$,
given by
\[
\widetilde{P} = \{(x, u, t) \in  \reals^{d + p + q} \mid
x = Y u + V t,\> u \geq 0,\> \mathbb{I}u = 1,\> t \geq 0\}
\]
onto $\reals^d$.
Consequently, if we can prove that the projection of
an $\s{H}$-polyhedron is also an  $\s{H}$-polyhedron,
then we will have proved that every  $\s{V}$-polyhedron is an 
$\s{H}$-polyhedron. Here again, it is possible that 
$P = \eucreal^d$, but that's fine since $\eucreal^d$
has been declared to be an $\s{H}$-polyhedron.

\medskip
In view of Proposition \ref{equivcone} and the
discussion that followed, it is enough to prove
that the projection of any $\s{H}$-cone is an $\s{H}$-cone.
This can be done by using a type of  Fourier-Motzkin 
elimination dual to the method used in Proposition 
\ref{FourierMotzkin1}.
We state the result without proof and refer the
interested reader to Ziegler \cite{Ziegler97},
Section 1.2--1.3, for full details.

\begin{prop}
\label{FourierMotzkin2}
If $C = P(A, 0)\subseteq \eucreal^d$ is an $\s{H}$-cone, then
the projection, $\mathrm{proj}_k(C)$, onto the hyperplane, $H_k$,
of equation $y_k = 0$ is given by 
$\mathrm{proj}_k(C) = \mathrm{elim}_k(C) \cap H_k$, with \\
$\mathrm{elim}_k(C) = 
\{x \in \reals^d \mid (\exists t\in \reals)(x + te_k \in P)\} = 
\{z - te_k \mid z\in P,\> t\in \reals\}
= P(A^{/k}, 0)$ 
and where the rows of $A^{/k}$ are given by
\[
A^{/k} = \{a_i \mid a_{i\, k} = 0\} \cup
\{a_{i\, k}a_j - a_{j\, k}a_i \mid a_{i\, k} > 0,\> a_{j\, k} < 0\}.
\]
\end{prop}

\medskip
It should be noted that both Fourier-Motzkin elimination
methods generate a quadratic number of new vectors or
inequalities at each step and thus they lead to a combinatorial explosion.
Therefore, these methods become intractable rather quickly.
The problem is that many of the new vectors or inequalities
are redundant. Thereore, it is important to find ways
of detecting redundancies and there are various methods for 
doing so. Again, the interested reader should consult
Ziegler \cite{Ziegler97}, Chapter 1.

\medskip
There is yet another way of proving that an $\s{H}$-polyhedron
is a $\s{V}$-polyhedron without using Fourier-Motzkin elimination.
As we already observed, Krein and Millman's theorem does not apply
if our polyhedron is unbounded. Actually, the full power of
Krein and Millman's theorem is not needed to show that an $\s{H}$-polytope
is a $\s{V}$-polytope. The crucial point is that if $P$
is an $\s{H}$-polytope with nonempty interior, then {\it every\/} line, $\ell$,
through any point, $a$, in the interior of $P$ intersects $P$ in
a line segment. This is because $P$ is compact and $\ell$ is closed,
so $P \cap \ell$ is a compact subset of a line thus,
a closed interval $[b, c]$ with $b < a < c$, as
$a$ is in the interior of $P$.
Therefore, we can use induction on the dimension
of $P$ to show that every point in $P$ is a convex combination
of vertices of  the facets of $P$. Now, if $P$ is unbounded
and cut out by at least two half-spaces (so, $P$ is not a half-space), then
we claim that for every point, $a$, in the interior of $P$,
there is {\it some\/} line through $a$ that intersects two facets of $P$.
This is because if we pick the origin in the interior of $P$, we may assume
that $P$ is given by an irredundant intersection,
$P = \bigcap_{i = 1}^t (H_i)_-$, and for any point, $a$, in the
interior of $P$, there is a line, $\ell$,  through $P$ {\it in general position
w.r.t. $P$\/}, which means that $\ell$ is not parallel to any of the
hyperplanes $H_i$ and intersects all of them in distinct points
(see Definition \ref{genpos}). Fortunately, lines in general
position always exist, as shown in Proposition \ref{perturb1}.
Using this fact, we can prove the following result:

\begin{prop}
\label{HtoVagain}
Let $P \subseteq \eucreal^d$ be an $\s{H}$-polyhedron,
$P = \bigcap_{i = 1}^t (H_i)_-$ (an irredundant decomposition), 
with nonempty interior. If $t = 1$, that is, $P = (H_1)_-$ is a half-space,
then 
\[
P = a + \mathrm{cone}(u_1, \ldots, u_{d-1}, -u_1, \ldots, -u_{d-1}, u_{d}), 
\]
where $a$ is any point in $H_1$, the vectors
$u_1, \ldots, u_{d-1}$ form a basis of the direction of $H_1$ and
$u_{d}$ is normal to (the direction of) $H_1$. (When $d = 1$,
$P$ is the half-line, $P = \{a + t u_1 \mid t\geq 0\}$.)
If $t \geq 2$, then every point, $a\in P$, can be written
as a convex combination, $a = (1 - \alpha)b + \alpha c\>$ 
($0 \leq \alpha \leq 1$),
where $b$ and $c$ belong to two distinct facets, $F$ and $G$, of $P$ and 
where 
\[
F = \mathrm{conv}(Y_F) + \mathrm{cone}(V_F)
\quad\hbox{and}\quad 
G = \mathrm{conv}(Y_G) + \mathrm{cone}(V_G),
\]
for some finite sets of points, $Y_F$ and $Y_G$ and some finite
sets of vectors, $V_F$ and $V_G$.  (Note: $\alpha = 0$ or $\alpha = 1$
is allowed.) 
Consequently, every $\s{H}$-polyhedron is a $\s{V}$-polyhedron.
\end{prop}

\proof
We proceed by induction on the dimension, $d$, of $P$.
If $d = 1$, then $P$ is either a closed interval, $[b, c]$, or
a half-line, $\{a + tu \mid t\geq 0\}$, where $u\not= 0$.
In either case, the proposition is clear.

\medskip
For the induction step, assume $d > 1$. Since every facet, $F$, of
$P$ has dimension $d - 1$,  the induction hypothesis holds for $F$, that
is, there exist a finite set of points, $Y_F$, and a finite set of
vectors, $V_F$, so that
\[
F = \mathrm{conv}(Y_F) + \mathrm{cone}(V_F).
\]
Thus, every point on the boundary of $P$ is of the desired form.
Next, pick any point, $a$, in the interior of $P$. Then, from our previous 
discussion, there is a line, $\ell$, through $a$
in general position w.r.t. $P$. Consequently, the intersection points,
$z_i = \ell \cap H_i$, of the line $\ell$ 
with the hyperplanes supporting the facets 
of $P$ exist and are all distinct. If we give $\ell$ an orientation,
the $z_i$'s can be sorted and there is a unique $k$ such that
$a$ lies between $b = z_k$ and $c = z_{k+1}$. But then, 
$b\in F_k = F$ and $c\in F_{k+1} = G$, where $F$ and $G$
the facets of $P$ supported by $H_k$ and $H_{k+1}$, and
$a = (1 - \alpha) b + \alpha c$, a convex combination.
Consequently, every point in $P$ is a mixed convex $+$ positive
combination of finitely many points associated with the
facets of $P$ and finitely many vectors associated with
the directions of the supporting hyperplanes of the facets $P$.
Conversely, it is easy to see that any such mixed combination 
must belong to $P$ and therefore, $P$ is a $\s{V}$-polyhedron.
$\bigsquare$

\medskip
We conclude this section with a version of Farkas Lemma
for polyhedral sets.

\begin{lemma}
\label{FarkasIV} (Farkas Lemma, Version IV)
Let $Y$ be any $d\times p$ matrix and $V$ be any
$d\times q$ matrix. For every $z\in \reals^d$,
exactly one of the
following alternatives occurs:
\begin{enumerate}
\item[(a)]
There exist $u\in \reals^p$ and $t\in \reals^q$, with
$u \geq 0$, $t\geq 0$, $\mathbb{I}u = 1$ and $z = Yu + Vt$.
\item[(b)]
There is some vector,  $(\alpha, c) \in \reals^{d+1}$,
such that $\transpos{c} y_i  \geq  \alpha$ for all $i$ with
$1 \leq i \leq p$, $\transpos{c} v_j  \geq  0$
for all $j$ with $1 \leq j \leq q$,  and $\transpos{c} z < \alpha$.
\end{enumerate}
\end{lemma}

\proof
We use Farkas Lemma, Version II (Lemma \ref{FarkasII}).
Observe that (a) is equivalent to the problem:
Find $(u, t)\in \reals^{p + q}$, so that
\[
\binom{u}{t} \geq \binom{0}{0}
\quad\hbox{and}\quad 
\begin{pmatrix}
\mathbb{I} & \mathbb{O} \\
Y & V
\end{pmatrix}
\binom{u}{t} = \binom{1}{z},
\]
which is exactly Farkas II(a). Now, the second alternative of
Farkas II says that there is no solution as above if there
is some $(-\alpha, c)\in \reals^{d+1}$ so that
\[
(-\alpha, \transpos{c})\binom{1}{z} < 0
\quad\hbox{and}\quad 
(-\alpha, \transpos{c})
\begin{pmatrix}
\mathbb{I} & 0 \\
Y & V
\end{pmatrix} 
\geq (\mathbb{O}, \mathbb{O}).
\]
These are equivalent to
\[
-\alpha + \transpos{c} z < 0,\quad
-\alpha \mathbb{I} +  \transpos{c}Y \geq \mathbb{O}, \quad
\transpos{c}V \geq \mathbb{O},
\]
namely,
$\transpos{c} z < \alpha$, $\>\transpos{c}Y \geq \alpha\mathbb{I}$
and $\transpos{c}V \geq \mathbb{O}$, which are indeed
the conditions of Farkas IV(b), in matrix form.
$\bigsquare$

\medskip
Observe that Farkas IV can be viewed as a separation criterion
for polyhedral sets. This version subsumes Farkas I and Farkas II.

\chapter[Projective Spaces and Polyhedra, Polar Duality]
{Projective Spaces, Projective Polyhedra, Polar Duality w.r.t.
a Nondegenerate Quadric}
\label{chap2b}
\section{Projective Spaces}
\label{sec5d}
The fact that not just points but also vectors are needed
to deal with unbounded polyhedra is a hint that perhaps
the notions of polytope and polyhedra can be unified by
``going projective''. However, we have to be careful because
projective geometry does not accomodate well the notion  of
convexity. This is because convexity has to do with
convex combinations, but the essense of
projective geometry is that everything is defined up to
{\it non-zero\/} scalars, without any requirement that
these scalars be positive. 

\medskip
It is possible to develop a theory
of {\it oriented projective geometry\/} (due to J. Stolfi 
\cite{Stolfi})
in wich convexity is nicely accomodated. However, in this approach,
every point comes as a pair, $($positive point, negative point$)$,
and although it is a very elegant theory, we find it a bit
unwieldy. However, since all we really want is to ``embed''
$\eucreal^d$ into its {\it projective completion\/}, $\projr{d}$,
so that we can deal with ``points at infinity'' and ``normal point''
in a uniform manner in particular, with respect to projective
transformations, we will content ourselves with
a definition of the notion of a projective polyhedron using
the notion of polyhedral cone. Thus, we will not attempt to define
a general notion of convexity. 

\medskip
We begin with a ``crash course'' on (real) projective spaces.
There are many texts on projective geometry.
We suggest starting with Gallier \cite{Gallbook2}
and then move on to far more comprehensive treatments such as
Berger (Geometry II) \cite{Berger90b} or Samuel \cite{Samuel}.

\begin{defin}
\label{rpdef}
{\em
The {\it (real) projective space\/}, $\rprospac{n}$, is the set of all lines
through the origin in $\reals^{n + 1}$, i.e., the set of one-dimensional
subspaces of $\reals^{n+1}$  (where $n \geq 0$). Since a one-dimensional
subspace, $L \subseteq \reals^{n+1}$, is spanned by any nonzero vector,
$u\in L$, we can view $\rprospac{n}$ as the set of equivalence classes
of nonzero vectors in $\reals^{n+1} - \{0\}$ modulo the equivalence relation,
\[
u \sim v
\quad\hbox{iff}\quad v = \lambda u,
\quad\hbox{for some}\quad \lambda \in \reals,\> \lambda\not= 0.
\]
We have the projection,
$\mapdef{p}{(\reals^{n+1} - \{0\})}{\rprospac{n}}$, given by
$p(u) = [u]_{\sim}$, the equivalence class of $u$ modulo $\sim$.
Write $[u]$ (or $\lag u\rag$) for the line, 
\[
[u] = \{\lambda u \mid \lambda\in \reals\},
\]
defined by the nonzero vector, $u$.
Note that $[u]_{\sim} = [u] - \{0\}$, for every $u\not= 0$, so
the map $[u]_{\sim} \mapsto [u]$ is a bijection which allows
us to identify $[u]_{\sim}$ and $[u]$. Thus, we will use
both notations interchangeably  as  convenient.
}
\end{defin}

\medskip
The projective space,  $\rprospac{n}$, is sometimes denoted
$\mathbb{P}(\reals^{n+1})$.
Since every line, $L$, in $\reals^{n+1}$ intersects the sphere $S^{n}$
in two antipodal points, we can view $\rprospac{n}$ as the quotient of
the sphere $S^{n}$ by identification of antipodal points. We call this the
{\it spherical model\/} of $\rprospac{n}$.

\medskip
A more subtle construction consists in considering the (upper) half-sphere
instead of the sphere, where the upper half-sphere $S^n_+$
is set of points on the sphere $S^n$ such that $x_{n+1} \geq 0$.
This time, every line through the center
intersects the (upper) half-sphere in a single point, except
\index{half-sphere $S^n_+$}%
\indsym{S^n_+}{upper half $n$-sphere}%
on the boundary of the half-sphere, where it intersects
in two  antipodal points $a_+$ and $a_-$.
Thus, the projective space $\rprospac{n}$ is the quotient space obtained
from the (upper) half-sphere $S^n_+$ 
by identifying antipodal points $a_+$ and $a_-$
on the boundary of the half-sphere. 
We call this model of $\rprospac{n}$  the {\it half-spherical model\/}.
\index{half-spherical model of projective geometry}%

\medskip
When $n = 2$, we get a circle. 
When $n = 3$, the upper half-sphere is homeomorphic to a closed disk
(say, by orthogonal projection onto the $xy$-plane), and
$\rprospac{2}$ is in bijection with a closed disk in which
antipodal points on its boundary (a unit circle) 
have been identified. This is hard to visualize!
In this model of the real projective space, projective lines are great
semicircles on the upper half-sphere, with 
antipodal points on the boundary identified.
Boundary points correspond to  points at infinity.
By orthogonal projection, these great semicircles correspond to 
semiellipses, with 
antipodal points on the boundary identified.
Traveling along such a  projective ``line,'' when we reach
a boundary point, we ``wrap around''!
In general, the upper half-sphere $S^n_+$ is homeomorphic
to the closed unit ball in $\reals^n$, whose boundary is the
$(n-1)$-sphere $S^{n-1}$. \nsindex{closed unit ball}
For example, the projective space $\rprospac{3}$ is in bijection with the
closed unit ball in $\reals^3$, with 
antipodal points on its boundary (the sphere $S^2$)
identified!

\medskip
Another useful way of ``visualizing'' $\rprospac{n}$ is 
to use the hyperplane, $H_{n+1}\subseteq \reals^{n+1}$, of equation 
$x_{n + 1} = 1$.
Observe that for every line, $[u]$, through the origin in $\reals^{n + 1}$,
if $u$ does not belong to the hyperplane, $H_{n+1}(0) \cong \reals^n$, of 
equation $x_{n + 1} = 0$, then $[u]$ intersects $H_{n+1}$ is a unique
point, namely, 
\[
\left(\frac{u_1}{u_{n+1}}, \ldots, \frac{u_n}{u_{n+1}}, 1\right),
\]
where $u = (u_1, \ldots, u_{n+1})$. The lines, $[u]$, for which
$u_{n+1} = 0$ are ``points at infinity''. Observe that 
the set of lines in $H_{n+1}(0) \cong \reals^n$ is the set of points
of the projective space, $\rprospac{n-1}$, and so, 
$\rprospac{n}$ can be written as the disjoint union
\[
\rprospac{n} = \reals^n \amalg \rprospac{n-1}.
\]

We can repeat the above analysis on $\rprospac{n-1}$ and
so we can think of $\rprospac{n}$ as the disjoint union
\[
\rprospac{n} = \reals^n \amalg \reals^{n-1} \amalg \cdots \amalg
\reals^1 \amalg \reals^{0},
\]
where $\reals^{0} = \{0\}$ consist of a single point.
The above shows that there is an embedding, 
$\reals^n \hookrightarrow \rprospac{n}$, given by 
$(u_1, \ldots, u_n) \mapsto  (u_1, \ldots, u_n, 1)$.

\medskip
It will also be very useful to use  homogeneous coordinates.
Given any point, \\
$p = [u]_{\sim}\in \rprospac{n}$, the set
\[
\{(\lambda u_1, \ldots, \lambda u_{n+1}) \mid \lambda \not= 0\}
\]
is called the set of  {\it homogeneous coordinates\/} of $p$.
Since $u \not= 0$, observe that for all homogeneous coordinates,
$(x_1, \ldots, x_{n+1})$, for $p$, some $x_i$ must be non-zero.
The traditional notation for the homogeneous coordinates of a point,
$p = [u]_{\sim}$, is 
\[
(u_1\co \cdots \co u_n \co u_{n+1}).
\]

\medskip
There is a useful bijection between certain kinds of subsets of
$\reals^{d+1}$ and subsets of $\rprospac{d}$.
For  any subset, $S$, of $\reals^{d+1}$,  let
\[
-S = \{-u \mid u \in S\}.
\]
Geometrically, $-S$ is the reflexion of $S$ about $0$.
Note that for {\it any\/} nonempty subset, \\
$S\subseteq \reals^{d+1}$, with $S\not= \{0\}$, the sets
$S$, $-S$  and $S \cup -S$ all induce
the {\it same\/} set of points in projective space, $\rprospac{d}$, since
\begin{eqnarray*}
p(S - \{0\} ) & = & \{[u]_{\sim} \mid u \in S - \{0\}\} \\
& = &  \{[-u]_{\sim} \mid u \in S- \{0\} \} \\ 
& = &  \{[u]_{\sim} \mid u \in -S - \{0\}\} = p((-S) - \{0\} ) \\
& = & \{[u]_{\sim} \mid u\in S - \{0\}\} 
\cup \{[u]_{\sim} \mid u\in (-S) - \{0\}\} =
p((S \cup -S) - \{0\}),
\end{eqnarray*}
because $[u]_{\sim} = [-u]_{\sim}$.
Using these facts we  obtain a bijection between
subsets of $\rprospac{d}$ and certain subsets of $\reals^{d+1}$.

\medskip
We say that a set, $S\subseteq \reals^{d+1}$, is {\it symmetric\/} iff
$S = -S$. Obviously, $S\cup -S$ is symmetric for any set, $S$.
Say that a subset, $C\subseteq \reals^{d+1}$, is a {\it double cone\/}
iff for every $u\in C - \{0\}$, the entire line, $[u]$,
spanned  by $u$ is contained in $C$. 
Again, we exclude the trivial double cone, $C = \{0\}$.
Thus, every double cone can be viewed as a set of lines through $0$.
Note that a double cone is symmetric.
Given any nonempty subset, $S\subseteq \rprospac{d}$, let
$v(S) \subseteq \reals^{d+1}$ be the set of vectors,
\[
v(S) = \bigcup_{[u]_{\sim}\in S} [u]_{\sim} \cup \{0\}.
\]
Note that $v(S)$ is a double cone.

\begin{prop}
\label{vSprop}
The map, $v\co S \mapsto v(S)$, from the set of nonempty subsets of
$\rprospac{d}$ to the set of nonempty, nontrivial double cones in $\reals^{d+1}$
is a bijection.
\end{prop}

\proof
We already noted that $v(S)$ is nontrivial double cone. Consider the map, 
\[
\mathrm{ps}\co S \mapsto p(S) 
= \{[u]_{\sim}\in \rprospac{d} \mid u\in S - \{0\}\}.
\]
We leave it as an easy exercise to check that
$\mathrm{ps} \circ v = \id$ and $v\circ \mathrm{ps} = \id$,
which shows that $v$ and $\mathrm{ps}$ are mutual inverses.
$\bigsquare$


\medskip
Given any subspace, $X\subseteq \reals^{n+1}$, with
$\mathrm{dim}\, X = k+1 \geq 1$ and $0 \leq k \leq n$,
a {\it $k$-dimensional projective
subspace\/} of  $\rprospac{n}$ is the image,
$Y = p(X - \{0\})$, of $X - \{0\}$ under the projection $p$.
We often write $Y = \mathbb{P}(X)$.
When $k = n-1$, we say that $Y$ is a {\it projective hyperplane\/}
or simply a {\it hyperplane\/}. When $k = 1$, we say that
$Y$ is a {\it projective line\/} or simply a {\it line\/}.
It is easy to see that every (projective) hyperplane is 
the kernel (zero set) of some linear equation of the form
\[
a_1x_1 + \cdots + a_{n+1} x_{n+1} = 0,
\]
where one of the $a_i$ is nonzero.
Conversely, the kernel of any such linear equation is a hyperplane.
Furthermore, given a (projective) hyperplane, 
$H\subseteq \rprospac{n}$, the linear equation defining
$H$ is unique up to a nonzero scalar.

\medskip
For any $i$, with $1\leq i \leq n + 1$, the set 
\[
U_i = \{(x_1\co \cdots\co x_{n+1}) \in \rprospac{n} \mid x_i \not= 0\}
\]
is a subset of  $\rprospac{n}$ called an {\it affine patch\/}
of $\rprospac{n}$. We have a bijection, $\mapdef{\varphi_i}{U_i}{\reals^n}$,  
between $U_i$ and $\reals^n$ given by
\[
\varphi_i \co (x_1\co \cdots\co x_{n+1}) \mapsto 
\left(\frac{x_1}{x_{i}}, \ldots, \frac{x_{i-1}}{x_{i}}, 
\frac{x_{i+1}}{x_{i}}, \ldots,  \frac{x_{n+1}}{x_{i}}\right).
\]
This map is well defined because if
$(y_1, \ldots, y_{n+1}) \sim (x_1, \ldots, x_{n+1})$,
that is, \\
$(y_1, \ldots, y_{n+1}) = \lambda  (x_1, \ldots, x_{n+1})$,
with $\lambda\not= 0$, then
\[
\frac{y_j}{y_{i}} = \frac{\lambda x_j}{\lambda x_{i}} =
\frac{x_j}{x_{i}} \qquad (1 \leq i \leq n + 1),
\]
since $\lambda \not = 0$ and $x_{i}, y_{i} \not= 0$.
The inverse, $\mapdef{\psi_i}{\reals^n}{U_i \subseteq \rprospac{n}}$,
of $\varphi_i$ is given by 
\[
\psi_i\co (x_1, \cdots, x_{n}) \mapsto (x_1\co \cdots x_{i-1}\co 1
\co x_{i}\co \cdots \co x_n).
\]

\medskip
Observe that the  bijection, $\varphi_i$, between $U_i$ and $\reals^n$
can also be viewed as the bijection 
\[
(x_1\co \cdots\co x_{n+1}) \mapsto 
\left(\frac{x_1}{x_{i}}, \ldots, \frac{x_{i-1}}{x_{i}}, 1, 
\frac{x_{i+1}}{x_{i}}, \ldots,  \frac{x_{n+1}}{x_{i}}\right),
\]
between $U_i$ and the
hyperplane, $H_i\subseteq \reals^{n+1}$, of equation $x_i = 1$.
We will make heavy use of these bijections. For example,
for any subset,  $S\subseteq \rprospac{n}$, the
``view of $S$ from the patch $U_i$'', $S\res U_{i}$, is in bijection
with $v(S)\cap H_i$, where $v(S)$ is the double cone
associated with $S$ (see Proposition \ref{vSprop}).

\medskip
The affine patches, $U_1, \ldots, U_{n+1}$, cover the projective space
$\rprospac{n}$, in the sense that every
$(x_1\co \cdots\co x_{n+1}) \in \rprospac{n}$ belongs
to one of the $U_i$'s, as not all $x_i = 0$.
The $U_i$'s turn out to be open subsets of $\rprospac{n}$
and they have nonempty overlaps. When we restrict ourselves to one of the
$U_i$, we have an ``affine view of $\rprospac{n}$ from $U_i$''.
In particular, on the affine patch $U_{n + 1}$,
we have the ``standard view'' of $\reals^n$  embedded
into $\rprospac{n}$  as $H_{n+1}$, the hyperplane of equation
$x_{n+1} = 1$. The complement, $H_i(0)$, of $U_i$ in  $\rprospac{n}$
is the (projective) hyperplane of equation $x_i = 0$
(a copy of $\rprospac{n - 1}$).
With respect to the affine patch, $U_i$, the hyperplane,
$H_i(0)$, plays the role of {\it hyperplane (of points) at infinity\/}.

\medskip
From now on, for simplicity of notation, we will write $\projr{n}$
for  $\rprospac{n}$. We need to define projective maps. Such maps
are induced  by linear maps. 

\begin{defin}
\label{projmapdef}
{\em
Any injective linear map,
$\mapdef{h}{\reals^{m+1}}{\reals^{n+1}}$, induces a map,
$\mapdef{\mathbb{P}(h)}{\projr{m}}{\projr{n}}$, defined by
\[
\mathbb{P}(h)([u]_{\sim}) = [h(u)]_{\sim}
\]
and called a {\it projective map\/}.
When $m = n$ and $h$ is bijective, the map $\mathbb{P}(h)$ is also
bijective and it is called a {\it projectivity\/}.
}
\end{defin}

We have to check that this definition makes sense, that is, 
it is compatible with the equivalence relation, $\sim$. 
For this, assume that $u \sim v$, that is
\[
v = \lambda u,
\]
with $\lambda\not= 0$ (of course, $u, v \not= 0$).
As $h$ is linear, we get
\[
h(v) = h(\lambda u) = \lambda h(u),
\]
that is, $h(u) \sim h(v)$, which shows that $[h(u)]_{\sim}$
does on depend on the representative chosen in the
equivalence class of $[u]_{\sim}$.
It is also easy to check that whenever two linear maps, $h_1$ and $h_2$,
induce the same projective map, {\it i.e.\/}, if
$\mathbb{P}(h_1) = \mathbb{P}(h_2)$, then there is a nonzero scalar,
$\lambda$, so that $h_2 = \lambda h_1$.

\medskip
Why did we require $h$ to be injective? Because if $h$ has a nontrivial kernel,
then, any nonzero vector, $u\in \Ker(h)$, is mapped to
$0$, but as $0$ does {\bf not} correspond to any point of $\projr{n}$,
the map $\mathbb{P}(h)$  is undefined on $\mathbb{P}(\Ker(h))$.

\medskip
In some case, we allow projective maps induced by 
non-injective linear maps $h$. In this case, $\mathbb{P}(h)$
is a map whose domain is  $\projr{n} - \mathbb{P}(\Ker(h))$.
An example is  the map, $\mapdef{\sigma_N}{\projr{3}}{\projr{2}}$,
given by
\[
(x_1\co x_2\co x_3\co x_4) \mapsto (x_1\co x_2\co x_4 - x_3),
\]
which is undefined at the point $(0\co 0\co 1\co 1)$.
This map is the ``homogenization'' of the central
projection (from the north pole, $N = (0, 0, 1)$) from 
$\eucreal^3$ onto $\eucreal^2$.

\danger
Although a projective map, $\mapdef{f}{\projr{n}}{\projr{n}}$,
is induced by some linear map, $h$, the map $f$ is
{\bf not} linear! This is because linear combinations of
points in $\projr{n}$ {\it do not make any sense!\/}

\medskip
Another way of defining functions (possibly partial) between
projective spaces involves using homogeneous polynomials.
If $p_1(x_1, \ldots, x_{m+1})$, $\ldots, p_{n+1}(x_1, \ldots, x_{m+1})$
are $n + 1$ homogeneous polynomials all of the same degree, $d$,
and if these $n+1$ polynomials do not vanish simultaneously, then
we claim that the function, $f$, given by
\[
f(x_1\co \cdots \co x_{m+1}) = 
(p_1(x_1, \ldots, x_{m+1})\co \cdots \co p_{n+1}(x_1, \ldots, x_{m+1}))
\]
is indeed a well-defined function from $\projr{m}$ to $\projr{n}$.
Indeed, if $(y_1, \ldots, y_{m+1}) \sim (x_1, \ldots, x_{m+1})$,
that is,
$(y_1, \ldots, y_{m+1}) = \lambda  (x_1, \ldots, x_{m+1})$,
with $\lambda\not= 0$, as the $p_i$ are homogeneous of degree $d$,
\[
p_i(y_1, \ldots, y_{m+1}) = 
p_i(\lambda x_1, \ldots, \lambda  x_{m+1}) = \lambda^d p_i(x_1, \ldots, x_{m+1}),
\]
and so,
\begin{eqnarray*}
f(y_1\co \cdots \co y_{m+1}) & = &
(p_1(y_1, \ldots, y_{m+1})\co \cdots \co p_{n+1}(y_1, \ldots, y_{m+1})) \\
& = & 
(\lambda^d p_1(x_1, \ldots, x_{m+1})\co \cdots \co 
\lambda^d p_{n+1}(x_1, \ldots, x_{m+1})) \\
& = & 
\lambda^d(p_1(x_1, \ldots, x_{m+1})\co \cdots \co 
p_{n+1}(x_1, \ldots, x_{m+1})) \\
& = & 
\lambda^d f(x_1\co \cdots \co x_{m+1}),
\end{eqnarray*}
which shows that 
$f(y_1\co \cdots \co y_{m+1}) \sim f(x_1\co \cdots \co x_{m+1})$, as
required.

\medskip
For example, the map, $\mapdef{\tau_N}{\projr{2}}{\projr{3}}$,
given by
\[
(x_1\co x_2,\co x_3) \mapsto 
(2x_1x_3\co 2x_2x_3\co x_1^2 + x_2^2 - x_3^2\co x_1^2 + x_2^2 + x_3^2), 
\]
is well-defined. It turns out to be the ``homogenization''
of the inverse stereographic map from $\eucreal^2$ to $S^2$
(see Section \ref{sec12}). Observe that
\[
\tau_N(x_1\co x_2\co 0) = (0\co 0 \co x_1^2 + x_2^2 \co x_1^2 + x_2^2)
= (0\co 0 \co 1\co 1),
\] 
that is, $\tau_N$ maps all the points
at infinity (in $H_3(0)$)
to the ``north pole'', $(0\co 0 \co 1\co 1)$. However,
when $x_3 \not= 0$, we can prove that $\tau_N$ is injective
(in fact, its inverse is $\sigma_N$, defined earlier).

\medskip
Most interesting subsets of projective space arise as the
collection of zeros  of a (finite) set of
homogeneous polynomials.
Let us begin with a single homogeneous polynomial,
$p(x_1, \ldots, x_{n+1})$, of degree $d$ and set
\[
V(p) = \{(x_1\co \cdots \co x_{n+1})\in \projr{n} \mid
p(x_1, \ldots, x_{n+1}) = 0\}.
\]
As usual, we need to check that this definition does not depend
on the  specific representative chosen in the equivalence class
of $[(x_1, \ldots, x_{n+1})]_{\sim}$. If 
 $(y_1, \ldots, y_{n+1}) \sim (x_1, \ldots, x_{n+1})$, that is,
$(y_1, \ldots, y_{n+1}) = \lambda  (x_1, \ldots, x_{n+1})$,
with $\lambda\not= 0$, as  $p$ is homogeneous of degree $d$,
\[
p(y_1, \ldots, y_{n+1}) = 
p(\lambda  x_1, \ldots, \lambda x_{n+1}) = \lambda^d p(x_1, \ldots, x_{n+1}),
\]
and as $\lambda\not= 0$,
\[
p(y_1, \ldots, y_{n+1}) = 0 
\quad\hbox{iff}\quad
p(x_1, \ldots, x_{n+1}) = 0, 
\]
which shows that $V(p)$ is well defined.
For a set of homogeneous polynomials (not necessarily of the same degree), 
$\s{E} = \{p_1(x_1, \ldots, x_{n+1}), \ldots, p_s(x_1, \ldots, x_{n+1})\}$,
we set
\[
V(\s{E}) = \bigcap_{i = 1}^s V(p_i) = 
 \{(x_1\co \cdots \co x_{n+1})\in \projr{n} \mid
p_i(x_1, \ldots, x_{n+1}) = 0,\> i = 1\ldots, s\}.
\]
The set, $V(\s{E})$, is usually called the {\it projective variety\/}
defined by $\s{E}$ (or {\it cut out by $\s{E}$\/}).
When $\s{E}$ consists of a single polynomial, $p$, the set $V(p)$ is
called a (projective) {\it hypersurface\/}. 
For example, if 
\[
p(x_1, x_2, x_3, x_4) = x_1^2 + x_2^2 + x_3^2 - x_4^2,
\]
then $V(p)$ is the {\it projective sphere\/} in $\projr{3}$, also
denoted $S^2$. Indeed, if we ``look'' at $V(p)$ on the affine patch
$U_4$, where $x_4\not= 0$, we know that this amounts to setting
$x_4 = 1$, and we do get the set of points $(x_1, x_2, x_3, 1)\in U_4$
satisfying $x_1^2 + x_2^2 + x_3^2 - 1 = 0$, our usual $2$-sphere!
However, if we look at $V(p)$ on the patch $U_1$, where $x_1\not= 0$, 
we see the quadric of equation $1 + x_2^2 + x_3^2 = x_4^2$,
which is not a sphere but a hyperboloid of two sheets!
Nevertheless, if we pick $x_4 = 0$ as the plane at infinity, note that
the projective sphere does not have points at infinity
since the only {\it real\/} solution of $x_1^2 + x_2^2 + x_3^2 = 0$
is $(0, 0, 0)$, but $(0, 0, 0, 0)$ does not correspond to any
point of $\projr{3}$. 

\medskip
Another example is given by 
\[
q = (x_1, x_2, x_3, x_4) = x_1^2 + x_2^2  - x_3x_4,
\]
for which $V(q)$ corresponds to a paraboloid in the patch
$U_4$. Indeed, if we set $x_4 = 1$, we get the set of points in
$U_4$ satisfying $x_3 = x_1^2 + x_2^2$. For this reason, we
denote $V(q)$ by $\s{P}$ and called it a {\it (projective) paraboloid\/}.

\medskip
Given any homogeneous polynomial, $F(x_1, \ldots, x_{d+1})$, 
we will also make use of the {\it hypersurface cone\/}, 
$C(F)\subseteq \reals^{d+1}$,
defined by
\[
C(F) = \{(x_1, \ldots, x_{d+1}) \in \reals^{d+1} \mid 
F(x_1, \ldots, x_{d+1}) = 0\}.
\]
Observe that $V(F) = \mathbb{P}(C(F))$.

\remark
Every variety, $V(\s{E})$, defined by a {\it set\/}
of polynomials, \\
$\s{E} = \{p_1(x_1, \ldots, x_{n+1}), \ldots, p_s(x_1, \ldots, x_{n+1})\}$,
is also the hypersurface defined by the {\it single\/} polynomial
equation,
\[
p_1^2 + \cdots + p_s^2 = 0.
\]
This fact, peculiar to the real field, $\reals$,
is a mixed blessing. On the one-hand, the study of varieties
is reduced to the study of hypersurfaces. On the other-hand,
this is a hint that we should expect that such a study will be hard.

\medskip
Perhaps to the surprise of the novice, there is a bijective
projective map (a projectivity) sending $S^2$ to $\s{P}$.
This map, $\theta$, is given by
\[
\theta(x_1\co x_2\co x_3\co x_4) = (x_1\co x_2 \co x_3 + x_4\co
x_4 - x_3),
\]
whose inverse is given by
\[
\theta^{-1}(x_1\co x_2\co x_3\co x_4) = \left(x_1\co x_2 \co 
\frac{x_3 - x_4}{2}\co \frac{x_3 + x_4}{2}\right).
\]
Indeed, if $(x_1\co x_2\co x_3\co x_4)$ satisfies
\[
 x_1^2 + x_2^2 + x_3^2 - x_4^2 = 0,
\]
and if  $(z_1\co z_2\co z_3\co z_4) = \theta(x_1\co x_2\co x_3\co x_4)$,
then from above,
\[
(x_1\co x_2\co x_3\co x_4) = \left(z_1\co z_2 \co 
\frac{z_3 - z_4}{2}\co \frac{z_3 + z_4}{2}\right),
\]
and by plugging the right-hand sides in the equation of the sphere, we get
\begin{eqnarray*}
z_1^2 + z_2^2 + \left(\frac{z_3 - z_4}{2}\right)^2 
- \left(\frac{z_3 + z_4}{2}\right)^2  
& = &
z_1^2 + z_2^2 + \frac{1}{4}(z_3^2 + z_4^2 - 2z_3z_4
- (z_3^2 + z_4^2 + 2z_3z_4)) \\
& = &
z_1^2 + z_2^2 - z_3z_4 = 0,
\end{eqnarray*}
which is the equation of the paraboloid, $\s{P}$.

\section{Projective Polyhedra}
\label{sec5e}
Following the ``projective doctrine'' which consists
in replacing points by lines through the origin, that is, to
``conify'' everything, we will define a projective polyhedron
as any set of points in $\projr{d}$ induced by a polyhedral cone
in $\reals^{d+1}$. To do so, it is preferable to
consider cones as sets of positive combinations of vectors
(see Definition \ref{conedef}). Just to refresh our memory,
a set, $C\subseteq \reals^d$, is a {\it $\s{V}$-cone\/} or
{\it polyhedral cone\/}
if $C$ is the positive hull of a finite set of vectors, that is,
\[
C = \mathrm{cone}(\{u_1, \ldots, u_p\}),
\]
for some vectors, $u_1, \ldots, u_p\in \reals^d$.
An {\it $\s{H}$-cone\/} is any subset of $\reals^d$ given by a finite 
intersection of closed half-spaces
cut out by hyperplanes through $0$.

\medskip 
A good place to learn about cones (and much more) is
Fulton \cite{Fultontoric}. See also Ewald \cite{Ewald}.

\medskip
By Theorem \ref{equivpolyhedra},
$\s{V}$-cones and $\s{H}$-cones form the same collection
of convex sets (for every $d\geq 0$). Naturally, we can think of
these cones as sets of rays (half-lines) of the form
\[
\lag u\rag_+ =  \{\lambda u \mid \lambda \in \reals,\> \lambda \geq 0\},
\]
where $u\in \reals^d$ is any {\it nonzero\/} vector.
We exclude the trivial cone, $\{0\}$, since $0$
does not define any point in projective space.
When  we ``go projective'', each ray corresponds to the full
line, $\lag u\rag$, spanned by $u$ which can be expressed as
\[
\lag u \rag = \lag u\rag_+ \cup -\lag u\rag_+,
\]
where $-\lag u\rag_+ = \lag u\rag_- = 
\{\lambda u \mid \lambda \in \reals,\> \lambda \leq 0\}$.
Now, if $C\subseteq \reals^d$ is a polyhedral cone, 
obviously $-C$ is also a polyhedral cone and
the set $C \cup -C$ consists of the union of the
two polyhedral cones $C$ and $-C$. Note that $C \cup -C$ 
can be viewed as the
set of all lines determined by the nonzero vectors
in $C$ (and $-C$). It is a double cone.
Unless $C$ is a closed half-space,
$C \cup -C$ is not convex. 
It seems perfectly natural to define
a projective polyhedron as any set of lines induced by a set of the
form $C \cup -C$, where $C$ is a polyhedral cone.

\begin{defin}
\label{projpolyhedon}
{\em
A {\it projective polyhedron\/}
is any subset,  $P \subseteq \projr{d}$, of the form
\[
P = p((C \cup  -C) -\{0\}) = p(C - \{0\}),
\] 
where $C$ is any polyhedral cone ($\s{V}$ or $\s{H}$ cone)  in $\reals^{d+1}$
(with $C \not= \{0\}$). We write \\
$P = \mathbb{P}(C \cup -C)$ or $P = \mathbb{P}(C)$.
}
\end{defin}

\medskip
It is important to observe that because $C \cup -C$ is a double cone
there is a bijection between nontrivial double polyhedral cones 
and projective polyhedra. So,  projective polyhedra are equivalent
to double polyhedral cones. However, the projective interpretation
of the lines induced by
$C\cup -C$ as points in $\projr{d}$ makes
the study of projective polyhedra geometrically more interesting.

\medskip
Projective polyhedra inherit many of the properties of cones
but we have to be careful because we are really dealing with
double cones, $C\cup -C$, and not cones. As a consequence,
there are a few unpleasant surprises, for example, the fact
that the collection of projective polyhedra is {\bf not} closed under
intersection!

\medskip
Before dealing with these issues, let us show that
every ``standard'' polyhedron, $P\subseteq \eucreal^d$, has
a natural projective completion, $\widetilde{P}\subseteq \projr{d}$,
such that on the affine patch $U_{d+1}$ (where $x_{d+1} \not= 0$),
$\widetilde{P}\res U_{d+1} = P$. For this, we use our theorem on
the Polyhedron--Cone Correspondence
(Theorem \ref{equivcone}, part (2)).

\medskip
Let $A = X + U$, where $X$ is a set of points in $\eucreal^d$ and $U$
is a cone in $\reals^d$.
For every point, $x\in X$, and every vector,
$u\in U$, let
\[
\widehat{x} = \binom{x}{1},
\qquad
\widehat{u} = \binom{u}{0},
\] 
and let $\widehat{X} = \{\widehat{x} \mid x\in X\}$
and  $\widehat{U} = \{\widehat{u} \mid u\in U\}$.
Then, 
\[
C(A) = \mathrm{cone}(\{\widehat{X}\cup \widehat{U}\})
\]
is a cone in $\reals^{d+1}$ such that
\[
\widehat{A} = C(A) \cap H_{d+1},
\]
where $H_{d+1}$ is the hyperplane of equation $x_{d+1} = 1$.
If we set $\widetilde{A} = \mathbb{P}(C(A))$, then we get a 
subset of $\projr{d}$ and in the patch $U_{d+1}$, 
the  set $\widetilde{A}\res U_{d+1}$
is in bijection with 
the intersection $(C(A) \cup -C(A)) \cap H_{d+1} = \widehat{A}$,
and thus, in bijection with $A$.
We call $\widetilde{A}$ the {\it projective completion of $A$\/}.
We have an injection, $A  \longrightarrow \widetilde{A}$,
given by
\[
(a_1, \ldots, a_d) \mapsto (a_1\co \cdots \co a_d\co 1),
\]
which is just the map, $\mapdef{\psi_{d+1}}{\reals^d}{U_{d+1}}$.
What the projective completion does is to add to $A$ the
``points at infinity'' corresponding to the vectors in $U$,
that is, the points of $\projr{d}$ corresponding to the
lines in the cone, $U$.
In particular, if $X = \mathrm{conv}(Y)$ and $U =  \mathrm{cone}(V)$, 
for some finite sets $Y = \{y_1, \ldots, y_p\}$ 
and $V = \{v_1, \ldots, v_q\}$, then
$P = \mathrm{conv}(Y) + \mathrm{cone}(V)$ 
is a $\s{V}$-polyhedron and 
$\widetilde{P} = \mathbb{P}(C(P))$ is a
projective polyhedron.
The projective polyhedron, $\widetilde{P} = \mathbb{P}(C(P))$,
is called the {\it projective completion\/} of $P$.

\medskip
Observe that if $C$ is a closed half-space in $\reals^{d+1}$, then
$P = \mathbb{P}(C\cup -C) = \projr{d}$.
Now, if $C\subseteq \reals^{d+1}$ is a polyhedral cone
and $C$ is contained in a closed half-space, it is still possible 
that $C$ contains some nontrivial linear subspace and we would like 
to understand this situation. 

\medskip
The first thing to observe is that
$U = C \cap (-C)$ is the largest linear subspace contained in $C$.
If $C \cap (-C) = \{0\}$, we say that $C$ is a {\it pointed\/}
or  {\it strongly convex\/} cone.
In this case, one immediately realizes that $0$ is an extreme point
of $C$ and so, there is a hyperplane, $H$, through $0$ so that
$C\cap H = \{0\}$, that is, except for its apex, $C$ lies
in one of the open half-spaces determined by $H$.
As a consequence, by a linear change of coordinates, we may assume that
this hyperplane is $H_{d+1}$  and so, for every projective
polyhedron, $P = \mathbb{P}(C)$, if $C$ is pointed then
there is an affine patch (say, $U_{d+1}$) where $P$ has no points
at infinity, that is, $P$ is a polytope!
On the other hand, from another patch, $U_i$,
as $P\res U_i$ is in bijection with $(C\cup -C) \cap H_i$,
the projective polyhedron
$P$ viewed on $U_i$ may consist of {\it two\/} disjoint polyhedra.

\medskip
The situation is very similar to the classical theory
of projective conics or quadrics (for example, see
Brannan, Esplen and Gray, \cite{Brannan}). 
The case where $C$ is a pointed cone corresponds
to the nondegenerate conics or quadrics.
In the case of the conics, depending how we slice a cone, we see 
an ellipse, a parabola or a hyperbola.
For projective polyhedra, when we slice a polyhedral double cone, $C \cup -C$,
we may see a polytope  ({\it elliptic type\/})
a single unbounded polyhedron ({\it parabolic type\/})
or two unbounded polyhedra ({\it hyperbolic type\/}).

\medskip
Now, when $U = C\cap (-C) \not= \{0\}$, the polyhedral cone,
$C$, contains the linear subspace, $U$, and if $C \not= \reals^{d+1}$, then 
for every hyperplane, $H$,  such that $C$ is contained in one of the two closed
half-spaces determined by $H$, the subspace $U\cap H$ is nontrivial.
An  example is the cone, $C\subseteq \reals^3$,
determined by the intersection of two planes through $0$
(a wedge). In this case, $U$ is equal to the line of
intersection of these two planes. 
Also observe that $C\cap (-C) = C$ iff $C = -C$, that is,
iff $C$ is a linear subspace.

\medskip
The situation where $C\cap (-C)\not= \{0\}$ is reminiscent
of the case of cylinders in the theory of quadric surfaces
(see \cite{Brannan} or Berger \cite{Berger90b}).
Now, every cylinder can be viewed as the ruled surface
defined as the family of lines orthogonal to a plane
and touching some  nondegenerate conic. 
A similar decomposition holds for polyhedral cones
as shown below in a proposition borrowed from
Ewald \cite{Ewald} (Chapter V, Lemma 1.6).
We should warn the reader that we have some doubts about
the proof given there and so, we offer a different proof
adapted from the proof of Lemma 16.2 in Barvinok
\cite{Barvinok}.
Given any two subsets, $V, W \subseteq \reals^d$,
as usual, we write $V + W = \{v + w\mid v\in V,\> w\in W\}$
and $v + W = \{v + w\mid w\in W\}$, for any $v\in\reals^d$.

\begin{prop}
\label{cospan1}
For every polyhedral cone, $C\subseteq \reals^d$, if 
$U = C \cap (-C)$, then there is some pointed cone, $C_0$,
so that $U$ and $C_0$ are orthogonal and
\[
C = U + C_0,
\]
with $\mathrm{dim}(U) + \mathrm{dim}(C_0) = \mathrm{dim}(C)$. 
\end{prop} 

\proof
We already know that $U = C\cap (-C)$ is the largest linear
subspace of $C$. Let $U^{\perp}$ be the orthogonal complement of
$U$ in $\reals^d$ and let $\mapdef{\pi}{\reals^d}{U^{\perp}}$
be the orthogonal projection onto $U^{\perp}$. By Proposition \ref{polydef2b},
the projection, $C_0 = \pi(C)$, of $C$ onto $U^{\perp}$ is
a polyhedral cone. We claim that $C_0$ is pointed and that
\[
C = U + C_0.
\]
Since $\pi^{-1}(v) = v + U$ for every $v\in C_0$, we
have $U + C_0 \subseteq C$. On the other hand, by definition of
$C_0$, we also have $C \subseteq U + C_0$, so
$C = U + C_0$. If $C_0$ was not pointed, then it would contain
a linear subspace, $V$, of dimension at least $1$ but then,
$U + V$ would be a linear subspace of $C$ of dimension strictly
greater than $U$, which is impossible. Finally,
$\mathrm{dim}(U) + \mathrm{dim}(C_0) = \mathrm{dim}(C)$
is obvious by orthogonality.
$\bigsquare$

\medskip
The linear subspace, $U = C\cap (-C)$,
is called the {\it cospan\/} of $C$. Both $U$ and $C_0$
are uniquely determined by $C$.
To a great extent, Proposition reduces the study of
non-pointed cones to the study of pointed cones.
We propose to call the projective polyhedra of the form
$P = \mathbb{P}(C)$, where $C$ is a cone with a non-trivial cospan
(a non-pointed cone) a
{\it projective polyhedral cylinder\/}, by analogy with the
quadric surfaces. We also propose to call the projective polyhedra of the form
$P = \mathbb{P}(C)$, where $C$ is a pointed cone, a
{\it projective polytope\/} (or {\it nondegenerate projective
polyhedron\/}).

\medskip
The following propositions show that projective polyhedra behave
well under projective maps and intersection with a hyperplane:

\begin{prop}
\label{projpolyp1}
Given any projective map, $\mapdef{h}{\projr{m}}{\projr{n}}$,
for any projective polyhedron, $P\subseteq \projr{m}$, the image, $h(P)$,
of $P$ is a projective polyhedron in $\projr{n}$.
Even if $\mapdef{h}{\projr{m}}{\projr{n}}$ is a partial map
but $h$ is defined on $P$, then $h(P)$ is a projective polyhedron.
\end{prop}

\proof
The projective map,  $\mapdef{h}{\projr{m}}{\projr{n}}$,
is of the form $h = \mathbb{P}(\widehat{h})$,  for some injective linear map,
$\mapdef{\widehat{h}}{\reals^{m+1}}{\reals^{n+1}}$. Moreover, the projective
polyhedron, $P$, is of the form $P = \mathbb{P}(C)$, for some polyhedral
cone, $C\subseteq \reals^{n+1}$, with
$C = \mathrm{cone}(\{u_1, \ldots, u_p\})$, for some nonzero vector
$u_i\in \reals^{d+1}$. By definition,
\[
\mathbb{P}(h)(P) = \mathbb{P}(h)(\mathbb{P}(C)) =  \mathbb{P}(\widehat{h}(C)).
\]
As $\widehat{h}$ is linear, 
\[
\widehat{h}(C) = \widehat{h}(\mathrm{cone}(\{u_1, \ldots, u_p\})) =
\mathrm{cone}(\{\widehat{h}(u_1), \ldots, \widehat{h}(u_p)\}).
\]
If we let $\widehat{C} = 
\mathrm{cone}(\{\widehat{h}(u_1), \ldots, \widehat{h}(u_p)\})$, then
$\widehat{h}(C) = \widehat{C}$ is a polyhedral cone and so,
\[
\mathbb{P}(h)(P) = \mathbb{P}(\widehat{h}(C)) = \mathbb{P}(\widehat{C})
\]
is a projective cone. This argument does not depend on the injectivity
of $\widehat{h}$, as long as $C \cap \Ker(\widehat{h}) = \{0\}$.
$\bigsquare$

\medskip
Proposition \ref{projpolyp1} together with earlier arguments
shows that every projective polytope,
$P\subseteq \projr{d}$,
is equivalent under some suitable projectivity
to another projective polytope, $P'$, which is a polytope when
viewed in the affine patch, $U_{d+1}$. This property is similar to the fact
that every (non-degenerate) projective conic is projectively
equivalent to an ellipse.

\medskip
Since the notion of a face is defined for arbitrary polyhedra it
is also defined for cones. Consequently, we can define
the notion of a face for projective polyhedra. \label{faceref}
Given a projective polyhedron, $P \subseteq \projr{d}$,
where $P = \mathbb{P}(C)$ for some polyhedral cone 
(uniquely determined by $P$), 
$C \subseteq \reals^{d+1}$, a {\it face of $P$\/} is any
subset of $P$ of the form $\mathbb{P}(F) = p(F - \{0\})$,
for any nontrivial face, $F\subseteq C$, of $C$ ($F\not= \{0\}$).
Consequently, we say that $\mathbb{P}(F)$ is a {\it vertex\/} 
iff $\mathrm{dim}(F) = 1$, an {\it edge\/} iff $\mathrm{dim}(F) =2$
and a {\it facet\/} iff $\mathrm{dim}(F) = \mathrm{dim}(C) - 1$.
The projective polyhedron, $P$, and the empty set are the
{\it improper\/} faces of $P$.
If $C$ is strongly convex, then it is easy to prove that
$C$ is generated by its edges (its one-dimensional faces, these
are rays) in the sense that any set of nonzero vector spanning these
edges generate $C$ (using positive linear combinations).
As a consequence, if $C$ is strongly convex,
we may say that $P$ is ``spanned''
by its vertices, since $P$ is equal to 
$\mathbb{P}($all positive combinations
of vectors representing its vertices$)$.

\remark
Even though we did not define the notion of convex
combination of points in $\projr{d}$, the notion
of projective polyhedron gives us a way to mimic
certain properties of convex sets in the framework of
projective geometry. That's because every projective
polyhedron corresponds to a unique polyhedral cone.

\medskip
If our projective polyhedron is the completion,
$\widetilde{P} = \mathbb{P}(C(P)) \subseteq \projr{d}$, of some polyhedron, 
$P \subseteq \reals^d$, then each face of the cone, $C(P)$,
is of the form $C(F)$, where $F$ is a face of $P$ and so,
each face of $\widetilde{P}$ is of the form $\mathbb{P}(C(F))$,
for some face, $F$, of $P$. In particular, in the affine patch,
$U_{d+1}$, the face,  $\mathbb{P}(C(F))$, is in bijection with the
face, $F$, of $P$. We will usually identify  $\mathbb{P}(C(F))$ and $F$.

\medskip
We now consider the intersection of projective polyhedra
but first, let us make some general remarks about
the intersection of subsets of $\projr{d}$. 
Given any two nonempty subsets, $\mathbb{P}(S)$ and
$\mathbb{P}(S')$, of $\projr{d}$
what is $\mathbb{P}(S)\cap \mathbb{P}(S')$?
It is tempting to say that
\[
\mathbb{P}(S)\cap \mathbb{P}(S') = \mathbb{P}(S\cap S'),
\]
but unfortunately this is generally false!
The problem is that $\mathbb{P}(S)\cap \mathbb{P}(S')$ is the
set of {\it all lines\/} determined by vectors both in $S$ and $S'$
but  there may be some line spanned by some vector
$u \in (-S)\cap S'$ or $u\in S\cap (-S')$ such that
$u$ does not belong to
$S\cap S'$ or $-(S\cap S')$. 

\medskip
Observe that 
\begin{eqnarray*}
-(-S) & = & S \\
-(S\cap S') & = & (-S)\cap (-S').
\end{eqnarray*}
Then, the correct intersection is given by
\begin{eqnarray*}
(S\cup -S) \cap (S' \cup -S') 
& = & (S\cap S') \cup ((-S)\cap (-S')) \cup
(S\cap (-S')) \cup ((-S)\cap S') \\ 
& = & (S\cap S') \cup -(S\cap S') \cup
(S\cap (-S')) \cup -(S\cap (-S')),
\end{eqnarray*}
which is the union of two  double cones
(except for $0$, which belongs to both).
Therefore,
\[
\mathbb{P}(S) \cap \mathbb{P}(S') =  
\mathbb{P}(S\cap S') \cup \mathbb{P}(S\cap (-S')) =
\mathbb{P}(S\cap S') \cup \mathbb{P}((-S)\cap S'),
\]
since $\mathbb{P}(S\cap (-S')) = \mathbb{P}((-S)\cap S')$.

\medskip
Furthermore, if  $S'$ is symmetric ({\it i.e.\/}, $S' = -S'$), then
\begin{eqnarray*}
(S\cup -S) \cap (S' \cup -S') & = & (S\cup -S) \cap S' \\
& = & (S\cap S') \cup ((-S)\cap S') \\  
& = & (S\cap S') \cup -(S\cap (-S')) \\  
& = & (S\cap S') \cup -(S\cap S').
\end{eqnarray*}
Thus, if either $S$ or $S'$ is symmetric, it is true that
\[
\mathbb{P}(S)\cap \mathbb{P}(S') = \mathbb{P}(S\cap S').
\]
Now, if $C$ is a pointed polyhedral cone
then $C\cap (-C) = \{0\}$. Consequently, for any other
polyhedral cone, $C'$, we have
$(C\cap C') \cap ((-C)\cap C') = \{0\}$.
Using these facts we obtain the following result:

\begin{prop}
\label{projpolyp2}
Let $P = \mathbb{P}(C)$  and $P'=\mathbb{P}(C')$  be any two projective 
polyhedra in  $\projr{d}$. If 
$\mathbb{P}(C) \cap \mathbb{P}(C') \not= \emptyset$, then
the following properties hold:
\begin{enumerate}
\item[(1)]
\[
\mathbb{P}(C) \cap \mathbb{P}(C') =  \mathbb{P}(C\cap C') \cup 
\mathbb{P}(C\cap (-C')),
\]
the union of two projective polyhedra. If  $C$ or $C'$ is a pointed cone
i.e., $P$ or $P'$ is a projective polytope, then
$\mathbb{P}(C\cap C')$ and $\mathbb{P}(C\cap (-C'))$ are disjoint. 
\item[(2)]
If $P' = H$, for some hyperplane,  $H \subseteq \projr{d}$,
then $P\cap H$ is a projective polyhedron.
\end{enumerate}
\end{prop}

\proof
We already proved (1) so only (2) remains to be proved.
Of course, we may assume that $P \not= \projr{d}$.
This time, using the equivalence theorem  of $\s{V}$-cones and $\s{H}$-cones
(Theorem \ref{equivpolyhedra}), 
we know that $P$ is of the form $P = \mathbb{P}(C)$,
with $C = \bigcap_{i = 1}^p C_i$, where the $C_i$ are closed
half-spaces in $\reals^{d+1}$. Moreover,
$H = \mathbb{P}(\widehat{H})$, for some hyperplane, 
$\widehat{H}\subseteq \reals^{d+1}$,
through $0$. Now,  as $\widehat{H}$ is symmetric, 
\[
P \cap H  =   \mathbb{P}(C)\cap \mathbb{P}(\widehat{H}) = 
\mathbb{P}(C\cap \widehat{H}).
\]
Consequently,
\begin{eqnarray*}
P \cap H & = &  \mathbb{P}(C \cap \widehat{H}) \\
& = & \mathbb{P}\left(\left(\bigcap_{i = 1}^p C_i\right)\cap \widehat{H} \right).
\end{eqnarray*}
However, $\widehat{H} = \widehat{H}_+ \cap \widehat{H}_-$, where 
$\widehat{H}_+$ and $\widehat{H}_-$ are the two closed
half-spaces determined by $\widehat{H}$  and so,
\[
\widehat{C} = \left(\bigcap_{i = 1}^p C_i\right)\cap \widehat{H} =
\left(\bigcap_{i = 1}^p C_i \right) \cap   
\widehat{H}_+ \cap \widehat{H}_-
\]
is a polyhedral cone. Therefore,
$P\cap H = \mathbb{P}(\widehat{C})$ is a projective polyhedron.
$\bigsquare$

\medskip
We leave it as an instructive exercise to find  explicit
examples where $P\cap P'$ consists of two disjoint projective
polyhedra in $\projr{1}$ (or  $\projr{2}$).

\medskip
Proposition \ref{projpolyp2} can be sharpened a little.

\begin{prop}
\label{projpolyp3}
Let $P = \mathbb{P}(C)$  and $P'=\mathbb{P}(C')$  be any two projective 
polyhedra in  $\projr{d}$. If 
$\mathbb{P}(C) \cap \mathbb{P}(C') \not= \emptyset$, then
\[
\mathbb{P}(C) \cap \mathbb{P}(C') =  \mathbb{P}(C\cap C') \cup 
\mathbb{P}(C\cap (-C')),
\]
the union of two projective polyhedra. If
$C = -C$, i.e., $C$ is a linear subspace (or if $C'$ is a linear subspace),
then 
\[
\mathbb{P}(C) \cap \mathbb{P}(C') =  \mathbb{P}(C\cap C').
\]
Furthermore, if either $C$ or $C'$ is pointed, the above projective
polyhedra are disjoint, else
if $C$ and $C'$ both have nontrivial cospan and 
$\mathbb{P}(C\cap C')$ and $\mathbb{P}(C\cap (-C'))$ 
intersect then
\[
 \mathbb{P}(C\cap C') \cap \mathbb{P}(C\cap (-C')) =
 \mathbb{P}(C\cap (C'\cap (-C'))) \cup 
\mathbb{P}(C'\cap (C\cap (-C))).
\]
Finally, if the two projective polyhedra on the right-hand side intersect, 
then 
\[
\mathbb{P}(C\cap (C'\cap (-C'))) \cap 
\mathbb{P}(C'\cap (C\cap (-C))) =
\mathbb{P}((C\cap (-C))\cap (C'\cap (-C'))).
\]
\end{prop}

\proof
Left as a simple exercise in boolean algebra.
$\bigsquare$

\medskip
In  preparation for  Section \ref{sec13}, we also need the notion
of tangent space at a point of a variety.

\section[Tangent Spaces of Hypersurfaces]
{Tangent Spaces of Hypersurfaces and Projective Hypersurfaces}
\label{sec5f}
Since we only need to consider the case of hypersurfaces
we restrict attention to this case (but the general case
is a straightforward generalization).
Let us begin with a hypersurface of equation
$p(x_1, \ldots, x_d) = 0$ in $\reals^d$, that is, the set
\[
S = V(p) = \{(x_1, \ldots, x_d) \in \reals^d \mid p(x_1, \ldots, x_d) = 0\},
\]
where $p(x_1, \ldots, x_d)$ is a polynomial of total degree, $m$. 

\medskip
Pick any point $a = (a_1, \ldots, a_d)\in \reals^d$.
Recall that there is a version of the Taylor expansion formula for
polynomials such that, for any polynomial, $p(x_1, \ldots, x_d)$,
of total degree $m$, for every
$h = (h_1, \ldots, h_d)\in \reals^d$, we have
\begin{eqnarray*}
p(a + h) & = & p(a) + \sum_{1 \leq |\alpha|\leq m}
\frac{D^{\alpha} p(a)}{\alpha !}  h^{\alpha} \\
& = & 
p(a) + \sum_{i = 1}^d  p_{x_i}(a) h_i + \sum_{2 \leq |\alpha|\leq m}
\frac{D^{\alpha} p(a)}{\alpha !}  h^{\alpha},
\end{eqnarray*}
where we use the {\it multi-index notation\/}, with
$\alpha = (i_1, \ldots, i_d)\in \natnums^d$,
$|\alpha| = i_1 + \cdots + i_d$,
$\alpha ! = i_1 ! \cdots i_d !$, 
$h^{\alpha} = h_1^{i_1}\cdots h_d^{i_d}$,
\[
D^{\alpha} p(a) = 
\frac{\partial^{i_1}}{\partial x_1^{i_1}} \cdots  
\frac{\partial^{i_d} p}{\partial x_d^{i_d}}\, (a),
\]
and
\[
p_{x_i}(a) = \frac{\partial p}{\partial x_i}\, (a).
\]

\medskip
Consider any line, $\ell$, through $a$,  given parametrically by
\[
\ell= \{a + t h \mid t\in \reals\}, 
\]
with $h\not= 0$ and
say $a\in S$ is a point on the hypersurface, $S = V(p)$, which
means that $p(a) = 0$.
The intuitive idea behind the notion of the tangent space to $S$ at $a$
is that it is the set of lines that intersect $S$ at $a$
in a point of {\it multiplicity at least two\/}, which means that the
equation giving the intersection, $S\cap \ell$, namely
\[
p(a + t h) =  p(a_1 + t h_1, \ldots, a_d + t h_d) = 0,
\]
is of the form
\[
t^2 q(a, h)(t) = 0,
\]
where $q(a, h)(t)$ is some polynomial in $t$.
Using Taylor's formula, as $p(a) = 0$, we have
\[
p(a + t h) =   t \sum_{i = 1}^d  p_{x_i}(a) h_i + t^2 q(a, h)(t),
\]
for some polynomial, $q(a, h)(t)$. From this, we see that 
$a$ is an intersection point of multiplicity at least $2$ iff
\begin{equation}
\sum_{i = 1}^d  p_{x_i}(a) h_i = 0.
\tag{$\dagger$}
\end{equation}
Consequently, if $\nabla p(a) = (p_{x_1}(a), \ldots, p_{x_d}(a)) \not = 0$
(that is, if the gradient of $p$ at $a$ is nonzero), 
we see that $\ell$ intersects $S$ at $a$ in a point of multiplicity
at least $2$ iff $h$ belongs to the hyperplane of equation
$(\dagger)$. 

\begin{defin}
\label{afftangspac}
{\em
Let $S = V(p)$ be a hypersurface in $\reals^d$.
For any point, $a\in S$, 
if $\nabla p(a) \not= 0$, then we say
that $a$ is a {\it non-singular\/} point of $S$. 
When $a$ is nonsingular, 
the ({\it affine\/}) {\it tangent space\/}, $T_a (S)$ (or simply,
$T_a S$), to $S$ at $a$ is the hyperplane
through $a$ of equation
\[
\sum_{i = 1}^d  p_{x_i}(a) (x_i - a_i) = 0.
\]
}
\end{defin}

\medskip
Observe that the hyperplane of the direction of $T_a S$ is
the hyperplane through $0$ and parallel to $T_a S$ given by
\[
\sum_{i = 1}^d  p_{x_i}(a) x_i = 0.
\]
When $\nabla p(a) = 0$, we either say that $T_a S$ is undefined
or we set $T_a S = \reals^d$.

\medskip
We now extend the notion of tangent space to projective varieties.
As we will see, this amounts to homogenizing and the result turns out
to be simpler than the affine case!

\medskip
So, let $S = V(F) \subseteq \projr{d}$ be a projective hypersurface,
which means that
\[
S = V(F) = \{(x_1\co \cdots\co x_{d+1}) \in \projr{d} 
\mid F(x_1, \ldots, x_{d+1}) = 0\},
\]
where $F(x_1, \ldots, x_{d+1})$ is a homogeneous
polynomial of total degree, $m$. Again, we say that
a point,  $a\in S$, is {\it non-singular\/} iff
$\nabla F(a) = (F_{x_1}(a), \ldots, F_{x_{d+1}}(a)) \not= 0$.  
For every $i = 1,\ldots, d+1$,  let
\[
z_j = \frac{x_j}{x_i},
\] 
where $j = 1, \ldots, d+1$ and $j \not= i$, and let
$f^{\res i}$ be the result of
``dehomogenizing'' $F$ at $i$, that is,
\[
f^{\res i}(z_1, \ldots, z_{i - 1}, z_{i+1}, \ldots, z_{d+1}) = 
F(z_1, \ldots, z_{i - 1}, 1, z_{i+1}, \ldots, z_{d+1}). 
\]
We define the {\it (projective) tangent space\/}, $T_a S$, to
$a$ at $S$ as the hyperplane, $H$, such that for each affine patch,
$U_i$ where $a_i \not = 0$, if we let 
\[
a^{\res i}_j = \frac{a_j}{a_i},
\]
where $j = 1, \ldots, d+1$ and $j \not= i$, 
then the restriction, $H\res U_i$,
of $H$ to $U_i$ is the affine hyperplane tangent to $S\res U_i$
given by 
\[
\sum_{
\begin{subarray}{c}
j = 1 \\
j \not = i
\end{subarray}
}^{d+1} f^{\res i}_{z_j}(a^{\res i})(z_j - a^{\res i}_j) = 0.
\]
Thus, on the affine patch, $U_i$, the tangent space,
$T_a S$, is given by the homogeneous
equation
\[
\sum_{
\begin{subarray}{c}
j = 1 \\
j \not = i
\end{subarray}
}^{d+1} f^{\res i}_{z_j}(a^{\res i})(x_j - a_j^{\res i} x_i) = 0.
\]
This looks awful but we can make it pretty if we remember
that $F$ is a homogeneous polynomial of degree $m$ and that
we have the {\it Euler relation\/}:
\[
\sum_{j = 1}^{d+1} F_{x_j}(x) x_j = m F,
\]
for every $x = (x_1, \ldots, x_{d+1})\in \reals^{d+1}$.
Using this, we can come up with a clean equation for our
projective tangent hyperplane. It is enough to carry out
the  computations for $i = d+1$. Our tangent hyperplane has the equation
\[
\sum_{j = 1}^d F_{x_j}(a_1^{\res d+1}, \ldots, a_d^{\res d+1}, 1)
(x_j - a_j^{\res d+1} x_{d+1}) = 0,
\]
that is,
\[
\sum_{j = 1}^d F_{x_j}(a_1^{\res d+1}, \ldots, a_d^{\res d+1}, 1)x_j 
+ \sum_{j = 1}^d F_{x_j}(a_1^{\res d+1}, \ldots, a_d^{\res d+1}, 1) 
(- a_j^{\res d+1} x_{d+1}) = 0.
\]
As $F(x_1, \ldots, x_{d+1})$ is homogeneous of degree $m$,
and as $a_{d+1}\not= 0$ on $U_{d+1}$, we have
\[
a_{d+1}^m F(a_1^{\res d+1}, \ldots, a_d^{\res d+1}, 1) = 
F(a_1, \ldots, a_d, a_{d+1}),
\]
so from the above equation we get
\begin{equation}
\sum_{j = 1}^d F_{x_j}(a_1, \ldots, a_{d+1})x_j 
+ \sum_{j = 1}^d F_{x_j}(a_1, \ldots, a_{d+1}) 
(- a_j^{\res d+1} x_{d+1}) = 0.
\tag{$*$}
\end{equation}
Since $a\in S$, we have $F(a) = 0$, so the Euler relation yields
\[
\sum_{j = 1}^{d} F_{x_j}(a_1, \ldots, a_{d+1}) a_j +
F_{x_{d+1}}(a_1, \ldots, a_{d+1}) a_{d+1}  = 0,
\]
which, by dividing by $a_{d+1}$ and
multiplying by $x_{d+1}$, yields
\[
\sum_{j = 1}^{d} F_{x_j}(a_1, \ldots, a_{d+1}) 
(- a_j^{\res d+1} x_{d+1}) = 
F_{x_{d+1}}(a_1, \ldots, a_{d+1}) x_{d+1},
\]
and by plugging this in $(*)$, we get
\[
\sum_{j = 1}^d F_{x_j}(a_1, \ldots, a_{d+1})x_j +
F_{x_{d+1}}(a_1, \ldots, a_{d+1}) x_{d+1} = 0.
\]
Consequently, 
the tangent hyperplane to $S$ at $a$ is given by the equation
\[
\sum_{j = 1}^{d+1} F_{x_j}(a)x_j = 0. 
\]

\begin{defin}
\label{projtangspac}
{\em
Let $S = V(F)$ be a hypersurface in $\projr{d}$, where
$F(x_1, \ldots, x_{d+1})$ is a homogeneous polynomial.
For any point, $a\in S$, 
if $\nabla F(a) \not= 0$, then we say
that $a$ is a {\it non-singular\/} point of $S$.
When $a$ is nonsingular,  
the  ({\it projective\/}) {\it tangent space\/}, $T_a (S)$ (or simply,
$T_a S$), to $S$ at $a$ is the hyperplane
through $a$ of equation
\[
\sum_{i = 1}^{d+1}  F_{x_i}(a)  x_i= 0.
\]
}
\end{defin}

\medskip
For example, if we consider the sphere, $S^2 \subseteq \projr{3}$, of
equation
\[
x^2 + y^2 + z^2 - w^2 = 0,
\]
the tangent plane to $S^2$ at $a = (a_1, a_2, a_3, a_4)$ is given by
\[
a_1 x + a_2 y + a_3 z - a_4 w = 0.
\]

\remark
If $a\in S = V(F)$, as $F(a) = \sum_{i = 1}^{d+1}  F_{x_i}(a)  a_i= 0$ 
(by Euler), the
equation of the tangent plane, $T_a S$, to $S$ at $a$ can also be written as
\[
\sum_{i = 1}^{d+1}  F_{x_i}(a)  (x_i - a_i)= 0.
\]
Now, if $a = (a_1\co \cdots \co a_d\co 1)$ is a point in the affine
patch $U_{d+1}$, then the equation of the intersection 
of $T_a S$ with $U_{d+1}$ is obtained by setting $a_{d+1} = x_{d+1} = 1$,
that is
\[
\sum_{i = 1}^{d}  F_{x_i}(a_1, \ldots, a_d, 1)  (x_i - a_i)= 0,
\]
which is just the equation of the affine hyperplane to 
$S\cap U_{d+1}$ at $a\in U_{d+1}$.

\medskip
It will be convenient to adopt  
the following notational convention: Given any point,
$x = (x_1, \ldots, x_d)\in \reals^d$, written as a row vector,
we let $\mathbf{x}$ denote the corresponding column vector
such that $\transpos{\mathbf{x}} = x$.

\medskip
Projectivities behave well with respect to hypersurfaces and
their tangent spaces. 
Let $S = V(F) \subseteq \projr{d}$
be a projective hypersurface, where $F$ is a homogeneous polynomial
of degree $m$ and let $\mapdef{h}{\projr{d}}{\projr{d}}$
be a projectivity (a bijective projective map).
Assume that $h$ is induced by the invertible $(d + 1)\times (d + 1)$
matrix, $A = (a_{i\, j})$, and write $A^{-1} =  (a^{-1}_{i\, j})$.
For any hyperplane, $H\subseteq \reals^{d+1}$,
if  $\varphi$ is any linear from defining $\varphi$, {\it i.e.\/},
$H = \Ker(\varphi$), then
\begin{eqnarray*}
h(H) & = & \{h(x) \in \reals^{d+1} \mid \varphi(x) = 0\} \\
& = & \{y\in \reals^{d+1} \mid (\exists x\in \reals^{d+1})(y = h(x),\>
 \varphi(x) = 0)\} \\
& = & \{y\in \reals^{d+1} \mid  (\varphi\circ h^{-1})(y) = 0\}.
\end{eqnarray*}
Consequently, if $H$ is given by
\[
\alpha_1x_1 + \cdots + \alpha_{d+1} x_{d+1} = 0
\]
and if we write $\alpha = (\alpha_1, \ldots, \alpha_{d+1})$, then
$h(H)$ is the hyperplane given by the equation
\[
\alpha A^{-1} \mathbf{y} = 0.
\]
Similarly, 
\begin{eqnarray*}
h(S) & = & \{h(x) \in \reals^{d+1} \mid F(x) = 0\} \\
& = & \{y\in \reals^{d+1} \mid (\exists x\in \reals^{d+1})(y = h(x),\>
 F(x) = 0)\} \\
& = & \{y\in \reals^{d+1} \mid  F(\transpos{(A^{-1} \mathbf{y})})  = 0\}
\end{eqnarray*}
is the hypersurface defined by the polynomial 
\[
G(x_1, \ldots, x_{d+1}) = 
F\left(\sum_{j = 1}^{d+1} a^{-1}_{1\, j} x_j, \ldots, 
\sum_{j = 1}^{d+1} a^{-1}_{d+1\, j} x_j\right).
\]
Furthermore, using the chain rule, we get
\[
(G_{x_1}, \ldots, G_{x_{d+1}}) =  (F_{x_1}, \ldots, F_{x_{d+1}}) A^{-1}, 
\]
which shows that a point, $a\in S$, is non-singular iff its image,
$h(a) \in h(S)$, is non-singular on $h(S)$. This also shows that
\[
h(T_a S) = T_{h(a)} h(S),
\]
that is, the projectivity, $h$, preserves tangent spaces.
In summary, we have

\begin{prop}
\label{projectivityp1}
Let $S = V(F) \subseteq \projr{d}$
be a projective hypersurface, where $F$ is a homogeneous polynomial
of degree $m$ and let $\mapdef{h}{\projr{d}}{\projr{d}}$
be a projectivity (a bijective projective map).
Then, $h(S)$ is a hypersurface in $\projr{d}$ and
a point, $a\in S$, is nonsingular for $S$ iff
$h(a)$ is nonsingular for $h(S)$. Furthermore,
\[
h(T_a S) = T_{h(a)} h(S),
\]
that is, the projectivity, $h$, preserves tangent spaces.
\end{prop}

\remark
If $\mapdef{h}{\projr{m}}{\projr{n}}$ is a projective map,
say induced by an injective linear map given by the 
$(n+1)\times (m+1)$ matrix, $A = (a_{i\, j})$, given any hypersurface,
$S = V(F)\subseteq \projr{n}$, we can define the {\it pull-back\/},
$h^*(S)\subseteq \projr{m}$, of $S$, by
\[
h^*(S) = \{x\in\projr{m} \mid F(h(x)) = 0\}.
\]
This is indeed a hypersurface because $F(x_1, \ldots, x_{n+1})$ is a
homogenous polynomial and $h^*(S)$ is the zero locus of the
homogeneous polynomial
\[
G(x_1, \ldots, x_{m+1}) =
F\left(\sum_{j = 1}^{m+1} a_{1\, j} x_j, \ldots, 
\sum_{j = 1}^{m+1} a_{n+1\, j} x_j\right).
\]
If $m = n$ and $h$ is a projectivity, then we have 
\[
h(S) = (h^{-1})^*(S).
\]

\section{Quadrics (Affine, Projective) and Polar Duality}
\label{sec5g}
The case where $S = V(\Phi)\subseteq \projr{d}$ is a hypersurface
given by a homogeneous polynomial, $\Phi(x_1, \ldots, x_{d+1})$,
of degree $2$ will come up a lot and deserves a little more attention.
In this case, if we write $x = (x_1, \ldots, x_{d+1})$, then 
$\Phi(x) = \Phi(x_1, \ldots, x_{d+1})$
is completely determined by a $(d+1)\times (d+1)$ symmetric
matrix, say $F = (f_{i\, j})$, and we have
\[
 \Phi(x) = \transpos{\mathbf{x}} F \mathbf{x} = 
\sum_{i, j = 1}^{d+1} f_{i\, j} x_i x_j.
\]
Since $F$ is symmetric, we can write
\[
\Phi(x) = \sum_{i, j = 1}^{d+1} f_{i\, j} x_i x_j =
\sum_{i = 1}^{d+1} f_{i\, i} x_i^2 + 
2\sum_{
\begin{subarray}{c}
i, j = 1\\
i < j 
\end{subarray}
} ^{d+1}
f_{i\, j} x_i x_j.
\]
The {\it polar form\/}, $\varphi(x, y)$, 
of $\Phi(x)$, is given by
\[
\varphi(x, y) =  
\transpos{\mathbf{x}} F \mathbf{y} = \sum_{i, j = 1}^{d+1} f_{i\, j} x_i y_j,
\]
where $x = (x_1, \ldots, x_{d+1})$  and 
$y = (y_1, \ldots, y_{d+1})$.
Of course,
\[
2\varphi(x, y) = 
\Phi(x + y) - \Phi(x) - 
\Phi(y).
\]
We also check immediately that
\[
2\varphi(x, y) = 2\transpos{\mathbf{x}} F \mathbf{y} =
\sum_{j = 1}^{d+1} \frac{\partial \Phi(x)}{\partial x_j}\, y_j,
\]
and so,
\[
\left(
\frac{\partial \Phi(x)}{\partial x_1}, \ldots, 
\frac{\partial \Phi(x)}{\partial x_{d+1}}
\right) = 2\transpos{\mathbf{x}} F.
\]
The hypersurface, $S = V(\Phi) \subseteq \projr{d}$, is called
a {\it (projective) (hyper-)quadric surface\/}. We say that a quadric
surface, $S = V(\Phi)$, is {\it nondegenerate\/} iff
the matrix, $F$, defining $\Phi$ is invertible. 

\medskip
For example, the sphere, $S^{d} \subseteq \projr{d+1}$, 
is the nondegenerate quadric given by
\[
\transpos{\mathbf{x}}
\begin{pmatrix}
I_{d+1} & \mathbf{0} \\
\mathbb{O} & -1
\end{pmatrix} \mathbf{x} = 0
\] 
and the paraboloid, $\s{P} \subseteq \projr{d+1}$, 
is the nongenerate quadric given by
\[
\transpos{\mathbf{x}}
\begin{pmatrix}
I_{d} & \mathbf{0} & \mathbf{0} \\
\mathbb{O}  & 0  & -\frac{1}{2} \\
\mathbb{O}  & -\frac{1}{2} & 0 \\
\end{pmatrix} \mathbf{x} = 0.
\] 

\medskip
If $\mapdef{h}{\projr{d}}{\projr{d}}$ is a projectivity induced by
some invertible matrix, $A = (a_{i\, j})$, and if
$S = V(\Phi)$ is  a quadric defined by the matrix $F$, 
we immediately check that $h(S)$ is the quadric defined by the matrix
$\transpos{(A^{-1})} F A^{-1}$.
Furthermore, as $A$ is invertible, we see that
$S$ is nondegenerate iff $h(S)$ is nondegenerate.

\medskip
Observe that polar duality w.r.t. the sphere, $S^{d-1}$,
can be expressed by
\[
X^* = \left\{x\in \reals^d \mid (\forall y\in X)\left(
(\transpos{\mathbf{x}}, 1) 
\begin{pmatrix}
I_{d} & \mathbf{0} \\
\mathbb{O} & -1
\end{pmatrix} 
\begin{pmatrix}
\mathbf{y} \\
1
\end{pmatrix} 
  \leq 0 \right)
\right\},
\]
where $X$ is any subset of $\reals^d$.
The above suggests generalizing polar duality with respect to
any nondegenerate quadric.

\medskip
Let $Q = V(\Phi(x_1, \ldots, x_{d + 1}))$ be a nondegenerate
quadric with corresponding polar form $\varphi$
and matrix $F = (f_{i\, j})$. Then, we know that
$\Phi$ induces a natural duality between $\reals^{d+1}$ and
$(\reals^{d+1})^*$, namely, for every $u\in \reals^{d+1}$,
if $\varphi_u \in (\reals^{d+1})^*$ is the linear form given by
\[
\varphi_u(v) = \varphi(u, v)
\]
for every $v\in \reals^{d+1}$,
then the map $u\mapsto \varphi_u$, from $\reals^{d+1}$
to  $(\reals^{d+1})^*$, is a linear isomorphism.

\begin{defin}
\label{polarquad1}
{\em
Let $Q = V(\Phi(x_1, \ldots, x_{d + 1}))$ be a nondegenerate quadric
with corresponding polar form, $\varphi$.
For any $u\in \reals^{d+1}$, with $u\not= 0$, the set
\[
u^{\dagger} = \{v\in \reals^{d+1} \mid \varphi(u, v) = 0\}
= \{v\in \reals^{d+1} \mid \varphi_u(v) = 0\} = \Ker\, \varphi_u
\]
is a hyperplane called the {\it polar of $u$ (w.r.t. $Q$)\/}.
}
\end{defin}

\medskip
In terms of the matrix representation of $Q$, the polar of
$u$ is given by the equation
\[
\transpos{\mathbf{u}} F \mathbf{x} = 0,
\]
or
\[
\sum_{j = 1}^{d+1} \frac{\partial \Phi(u)}{\partial x_j}\, x_j = 0.
\]
Going over to $\projr{d}$, we say that $\mathbb{P}(u^{\dagger})$ is
the {\it polar (hyperplane)\/} of the point $a = [u]\in \projr{d}$
and we write $a^{\dagger}$ for  $\mathbb{P}(u^{\dagger})$.

\medskip
Note that the equation of the polar hyperplane, $a^{\dagger}$, of a point, 
$a \in \projr{d}$, is identical to the equation of the tangent plane
to $Q$ at $a$, except that $a$ is not necessarily on $Q$.
However, if $a\in Q$, then the polar of $a$ is indeed
the tangent hyperplane, $T_a Q$, to $Q$ at $a$.

\begin{prop}
\label{quadtang1}
Let $Q = V(\Phi(x_1, \ldots, x_{d + 1}))\subseteq \projr{d}$ 
be a nondegenerate quadric
with corresponding polar form, $\varphi$, and matrix, $F$.
Then, every point, $a\in Q$, is nonsingular.
\end{prop}

\proof
Since 
\[
\left(
\frac{\partial \Phi(a)}{\partial x_1}, \ldots, 
\frac{\partial \Phi(a)}{\partial x_{d+1}}
\right) = 2\transpos{\mathbf{a}} F,
\]
if $a\in Q$ is singular, then $\transpos{\mathbf{a}} F = 0$
with $a\not= 0$, 
contradicting the fact that $F$ is invertible.
$\bigsquare$

\medskip
The reader should prove the following simple proposition:

\begin{prop}
\label{quadricdual1}
Let $Q = V(\Phi(x_1, \ldots, x_{d + 1}))$ be a nondegenerate quadric
with corresponding polar form, $\varphi$.
Then, the following properties hold:
For any two points, $a, b\in \projr{d}$, 
\begin{enumerate}
\item[(1)]
$a\in b^{\dagger}$ iff $b\in a^{\dagger}$;
\item[(2)]
$a\in a^{\dagger}$ iff $a\in Q$;
\item[(3)]
$Q$ does not contain any hyperplane.
\end{enumerate}
\end{prop}

\remark
As in the case of the sphere, if $Q$ is a nondegenerate quadric
and $a\in \projr{d}$ is any point such that the polar
hyperplane, $a^{\dagger}$,  intersects $Q$, then
there is a nice geometric interpretation for  $a^{\dagger}$.
Observe that for every $b\in Q\cap a^{\dagger}$,
the polar hyperplane, $b^{\dagger}$, is the
tangent hyperplane, $T_b Q$, to $Q$ at $b$ and that
$a\in T_b Q$. Also, if $a\in T_b Q$ for any $b\in Q$, as
$b^{\dagger} = T_b Q$, then
$b \in a^{\dagger}$.
Therefore, $Q\cap a^{\dagger}$ is the
set of contact points of all the tangent hyperplanes to
$Q$ passing through $a$.

\medskip
Every hyperplane, $H\subseteq \projr{d}$, is the polar of a single
point, $a\in \projr{d}$. Indeed, if $H$ is defined by a nonzero
linear form, $f\in (\reals^{d+1})^*$, as $\Phi$ is nondegenerate, 
there is a unique $u\in \reals^{d+1}$, with $u \not= 0$,
so that $f = \varphi_u$,
and as $\varphi_u$ vanishes on $H$, we see that
$H$ is the polar of the point $a = [u]$. If $H$ is also 
the polar of another point, $b = [v]$, then
$\varphi_v$ vanishes on $H$, which means that
\[
\varphi_v = \lambda \varphi_u = \varphi_{\lambda u},
\]
with $\lambda\not= 0$ and this implies $v = \lambda u$,
that is, $a = [u] = [v] = b$
and the pole of $H$ is indeed unique.

\begin{defin}
\label{polarquad2}
{\em
Let $Q = V(\Phi(x_1, \ldots, x_{d + 1}))$ be a nondegenerate quadric
with corresponding polar form, $\varphi$.
The {\it polar dual (w.r.t. $Q$)\/}, $X^*$,  of a subset, 
$X \subseteq \reals^{d+1}$, is given by
\[
X^* = \{v\in \reals^{d+1} \mid (\forall u\in X)
(\varphi(u, v) \leq 0)\}. 
\]
For every subset, $X\subseteq \projr{d}$, we let
\[
X^* = \mathbb{P}((v(X))^*),
\]
where $v(X)$ is the unique double cone associated with $X$ 
as in Proposition \ref{vSprop}.
}
\end{defin}

\medskip
Observe that $X^*$ is always a double cone, even if 
$X\subseteq \reals^{d+1}$ is not.
By analogy with the Euclidean case, for any nonzero vector,
$u\in \reals^{d+1}$, let
\[
(u^{\dagger})_- = \{v\in \reals^{d+1} \mid
\varphi(u, v) \leq 0\}.
\]
Now, we have the following version of Proposition \ref{polytopdual}:

\begin{prop}
\label{conedualp1}
Let $Q = V(\Phi(x_1, \ldots, x_{d + 1}))$ be a nondegenerate quadric
with corresponding polar form, $\varphi$, 
and matrix,  $F = (f_{i\, j})$.
For any nontrivial polyhedral cone, 
$C = \mathrm{cone}(u_1, \ldots, u_p)$,
where $u_i\in \reals^{d+1}$, $u_i \not = 0$, we have
\[
C^* = \bigcap_{i = 1}^p (u_i^{\dagger})_-.
\]
If $U$ is the $(d + 1)\times p$ matrix whose $i^{\mathrm{th}}$ column
is $u_i$, then we can also write
\[
C^* = P(\transpos{U}F, \mathbf{0}),
\]
where 
\[
 P(\transpos{U}F, \mathbf{0}) = \{v\in \reals^{d+1} \mid
\transpos{U}F \mathbf{v} \leq \mathbf{0}\}.
\]
Consequently, the polar dual of a polyhedral cone w.r.t.
a nondegenerate quadric is a polyhedral cone.
\end{prop}

\proof
The proof is essentially the same as the proof of 
Proposition \ref{polytopdual}. As
\[
C = \mathrm{cone}(u_1, \ldots, u_p) = 
\{\lambda_1 u_1 + \cdots + \lambda_p u_p  \mid \lambda_i \geq 0, \> 
1\leq i \leq p\},
\]
we have
\begin{eqnarray*}
C^* & = & \{v\in \reals^{d+1} \mid (\forall u\in C)
(\varphi(u, v) \leq 0)\} \\
& = &  \{v\in \reals^{d+1} \mid 
\varphi(\lambda_1 u_1 + \cdots + \lambda_p u_p, v) \leq 0,\>
\lambda_i \geq 0,\> 1\leq i \leq p\} \\
& = &  \{v\in \reals^{d+1} \mid 
\lambda_1\varphi(u_1, v) + \cdots +  \lambda_p\varphi(u_p, v)\leq 0,\>
\lambda_i \geq 0,\> 1\leq i \leq p\} \\
& = & \bigcap_{i = 1}^p \{v\in \reals^{d+1} \mid 
\varphi(u_i, v) \leq 0\}\\
& = & \bigcap_{i = 1}^p (u_i^{\dagger})_-.
\end{eqnarray*}
By the equivalence theorem for $\s{H}$-polyhedra and $\s{V}$-polyhedra,
we conclude that $C^*$ is a polyhedral cone.
$\bigsquare$

\medskip
Proposition \ref{conedualp1} allows us to make the following definition:

\begin{defin}
\label{polarquad3}
{\em
Let $Q = V(\Phi(x_1, \ldots, x_{d + 1}))$ be a nondegenerate quadric
with corresponding polar form, $\varphi$.
Given any projective polyhedron,
$P = \mathbb{P}(C)$, where $C$ is a polyhedral cone,
the {\it polar dual (w.r.t. $Q$)\/}, $P^*$, of $P$ is 
the projective polyhedron
\[
P^* = \mathbb{P}(C^*).
\]
}
\end{defin}

\medskip
We also show that projectivities behave well with respect
to polar duality.

\begin{prop}
\label{projectualp1}
Let  $Q = V(\Phi(x_1, \ldots, x_{d + 1}))$ be a nondegenerate quadric
with corresponding polar form, $\varphi$, and matrix,
$F = (f_{i\, j})$.
For every projectivity, $\mapdef{h}{\projr{d}}{\projr{d}}$,
if $h$ is induced by the linear map, $\widehat{h}$, given by the
invertible matrix, $A = (a_{i\, j})$,
for every subset, $X\subseteq \reals^{d+1}$, we have
\[
\widehat{h}(X^*) = (\widehat{h}(X))^*,
\]
where on the left-hand side, $X^*$ is the polar dual of $X$
w.r.t. $Q$ and on the right-hand side, $(\widehat{h}(X))^*$ is the polar dual 
of $\widehat{h}(X)$ w.r.t. the nondegenerate quadric, $h(Q)$,
given by the matrix $\transpos{(A^{-1})}FA^{-1}$.
Consequently, if $X\not= \{0\}$, then
\[
h((\mathbb{P}(X))^*) = (h(\mathbb{P}(X)))^*
\]
and for every projective polyhedron, $P$, we have
\[
h(P^*) = (h(P))^*.
\]
\end{prop}

\proof
As
\[
X^* = \{v\in \reals^{d+1} \mid (\forall u\in X)
(\transpos{\mathbf{u}}F \mathbf{v} \leq 0)\},
\]
we have
\begin{eqnarray*}
\widehat{h}(X^*) & = & \{\widehat{h}(v)\in \reals^{d+1} \mid (\forall u\in X)
(\transpos{\mathbf{u}}F \mathbf{v} \leq 0)\} \\
& = & \{y\in \reals^{d+1} \mid (\forall u\in X)
(\transpos{\mathbf{u}}F A^{-1} \mathbf{y} \leq 0)\} \\
& = & \{y\in \reals^{d+1} \mid (\forall x\in \widehat{h}(X))
(\transpos{\mathbf{x}}\transpos{(A^{-1})}F A^{-1} \mathbf{y} \leq 0)\} \\
& = & (\widehat{h}(X))^*,
\end{eqnarray*}
where $(\widehat{h}(X))^*$ is  the polar dual of 
$\widehat{h}(X)$ w.r.t. the quadric
whose matrix is $\transpos{(A^{-1})}F A^{-1}$, that is,
the polar dual w.r.t. $h(Q)$.

\medskip
The second part of the proposition follows immediately 
by setting $X = C$, where $C$ is the polyhedral cone
defining the projective polyhedron, $P = \mathbb{P}(C)$.
$\bigsquare$

\medskip
We will also need the notion of an affine quadric and polar duality
with respect to an affine quadric. Fortunately, the properties
we need in the affine case are easily derived from the projective case
using the ``trick'' that the affine space, $\eucreal^d$,
can be viewed as the hyperplane, $H_{d+1}\subseteq \reals^{d+1}$, 
of equation, $x_{d+1} = 1$ and that its associated vector space, $\reals^d$,
can be viewed as the hyperplane, $H_{d+1}(0)\subseteq \reals^{d+1}$, of equation
$x_{d+1} = 0$. A point, $a\in \affreal^d$, corresponds to
the vector, $\widehat{a} = \binom{a}{1}\in \reals^{d+1}$, and a vector,
$u\in \reals^d$, corresponds to the vector,
$\widehat{u} = \binom{u}{0}\in \reals^{d+1}$. This way,  the projective
space, $\projr{d} = \mathbb{P}(\reals^{d+1})$, is the
natural {\it projective completion\/} of $\eucreal^d$, which
is isomorphic to the affine patch $U_{d+1}$ where $x_{d+1} \not= 0$.
The hyperplane, $x_{d+1} = 0$, is the ``hyperplane at infinity''
in $\projr{d}$.

\medskip
If we write $x = (x_1, \ldots, x_{d})$,  a polynomial, 
$\Phi(x) = \Phi(x_1, \ldots, x_{d})$, of degree $2$
can be written as
\[
 \Phi(x) = \sum_{i, j = 1}^{d} a_{i\, j} x_i x_j
+ 2\sum_{i = 1}^d b_i x_i + c,
\]
where $A = (a_{i\, j})$ is a symmetric matrix. If we write
$\transpos{b} = (b_1, \ldots, b_d)$, then we have
\[
 \Phi(x) = 
(\transpos{\mathbf{x}}, 1) 
\begin{pmatrix}
A & b \\
\transpos{b} & c
\end{pmatrix}
\begin{pmatrix}
\mathbf{x} \\
1
\end{pmatrix}
= 
\transpos{\widehat{\mathbf{x}}}
\begin{pmatrix}
A & b \\
\transpos{b} & c
\end{pmatrix}
\widehat{\mathbf{x}}.
\]
Therefore, as in the projective case, $\Phi$
is completely determined by a $(d+1)\times (d+1)$ symmetric
matrix, say $F = (f_{i\, j})$, and we have
\[
 \Phi(x) = (\transpos{\mathbf{x}}, 1) F \binom{\mathbf{x}}{1} =
 \transpos{\widehat{\mathbf{x}}} F \widehat{\mathbf{x}}.
\]
We say that $Q\subseteq \reals^d$ is a {\it nondegenerate affine quadric\/} iff
\[
Q = V(\Phi) = \left\{x\in \reals^d \mid 
(\transpos{\mathbf{x}}, 1) F \binom{\mathbf{x}}{1} = 0\right\},
\]
where $F$ is symmetric and invertible.
Given any point $a\in \reals^d$, the {\it polar hyperplane\/},
$a^{\dagger}$, of $a$ w.r.t. $Q$ is defined by
\[
a^{\dagger} = \left\{x\in \reals^d \mid 
(\transpos{\mathbf{a}}, 1) F \binom{\mathbf{x}}{1} = 0\right\}.
\]
From a previous discussion, the equation of the polar hyperplane,
$a^{\dagger}$, is
\[
\sum_{i = 1}^d \frac{\partial \Phi(a)}{\partial x_i}\, (x_i - a_i) = 0.
\]
Given any subset, $X\subseteq \reals^d$, the {\it polar dual\/},
$X^*$, of $X$ is defined by
\[
X^* = \left\{y\in \reals^{d} \mid (\forall x\in X)
 \left( (\transpos{\mathbf{x}}, 1) F \binom{\mathbf{y}}{1} 
\leq  0\right)\right\}.
\]
As noted before, polar duality with respect to the
affine sphere, $S^d\subseteq \reals^{d+1}$, corresponds to the case where
\[
F = 
\begin{pmatrix}
I_{d} & \mathbf{0} \\
\mathbb{O} & -1
\end{pmatrix} 
\] 
and polar duality with respect to the affine paraboloid
$\s{P} \subseteq \reals^{d+1}$, 
 corresponds to the case where
\[
F = 
\begin{pmatrix}
I_{d-1} & \mathbf{0} & \mathbf{0} \\
\mathbb{O}  & 0  & -\frac{1}{2} \\
\mathbb{O}  & -\frac{1}{2} & 0 \\
\end{pmatrix}.
\] 

\medskip
We will need the following version of Proposition
\ref{polytopdual3}:

\begin{prop}
\label{polytopdualaff3}
Let $Q$ be a nondegenerate affine quadric given by the 
$(d + 1)\times (d + 1)$ symmetric matrix, $F$,
let  $\{y_1, \ldots, y_p\}$  be any set of points in
$\eucreal^d$ and 
let  $\{v_1, \ldots, v_q\}$ be any set of nonzero vectors in $\reals^d$.
If $\widehat{Y}$ is the $d\times p$ matrix whose $i^{\mathrm{th}}$ column
is $\widehat{y_i}$ and $V$ is the $d\times q$ matrix whose  $j^{\mathrm{th}}$ 
column is $\widehat{v_j}$, then
\[
 (\mathrm{conv}(\{y_1, \ldots, y_p\}) 
\cup \mathrm{cone}(\{v_1, \ldots, v_q\}))^*
=  P(\transpos{\widehat{Y}}F, \mathbf{0}; \transpos{\widehat{V}}F, \mathbf{0}),
\]
with 
\[
P(\transpos{\widehat{Y}}F, \mathbf{0}; \transpos{\widehat{V}}F, \mathbf{0}) = 
\left\{x\in \reals^d \mid \transpos{\widehat{Y}}F\binom{\mathbf{x}}{1} 
\leq \mathbf{0},\>
\transpos{\widehat{V}}F\binom{\mathbf{x}}{0} \leq \mathbf{0}\right\}.
\]
\end{prop}

\proof
The proof is immediately adpated from that of Proposition
\ref{polytopdual3}.
$\bigsquare$

\medskip
Using Proposition \ref{polytopdualaff3}, we can prove the
following Proposition showing that projective completion
and polar duality commute:

\begin{prop}
\label{commut1}
Let $Q$ be a nondegenerate affine quadric given by the 
$(d + 1)\times (d + 1)$ symmetric, invertible matrix, $F$.
For every polyhedron, $P \subseteq \reals^d$, we have
\[
\widetilde{P^*} = (\widetilde{P})^*,
\]
where on the right-hand side, we use polar duality
w.r.t. the nondegenerate projective quadric,
$\widetilde{Q}$, defined by $F$.
\end{prop}

\proof
By definition, 
we have $\widetilde{P} = \mathbb{P}(C(P))$,
$(\widetilde{P})^* =  \mathbb{P}((C(P))^*)$ and
$\widetilde{P^*} =  \mathbb{P}(C(P^*))$.
Therefore, it suffices to prove that
\[
(C(P))^* = C(P^*).
\]
Now, $P = \mathrm{conv}(Y) + \mathrm{cone}(V)$,
for some finite set of points, $Y$,   and some finite set of
vectors, $V$, and we know that
\[
C(P) = \mathrm{cone}(\widehat{Y} \cup \widehat{V}).
\]  
From Proposition \ref{conedualp1}, 
\[
(C(P))^* = \{v\in \reals^{d+1} \mid \transpos{\widehat{Y}}F \mathbf{v}
\leq \mathbf{0},
\> \transpos{\widehat{V}}F \mathbf{v}\leq \mathbf{0}\}
\]
and by Proposition \ref{polytopdualaff3},
\[
P^* = 
\left\{x\in \reals^d \mid 
\transpos{\widehat{Y}}F\binom{\mathbf{x}}{1} 
\leq \mathbf{0},\>
\transpos{\widehat{V}}F\binom{\mathbf{x}}{0} \leq \mathbf{0}\right\}.
\]
But, by definition of $C(P^*)$, the hyperplanes cutting out
$C(P^*)$ are obtained by homogenizing the equations of the hyperplanes
cutting out $P^*$ and so,
\[
C(P^*) = 
\left\{\binom{\mathbf{x}}{x_{d+1}}\in \reals^{d+1} \mid 
\transpos{\widehat{Y}}F\binom{\mathbf{x}}{x_{d+1}} \leq \mathbf{0},\>
\transpos{\widehat{V}}F\binom{\mathbf{x}}{x_{d+1}} \leq \mathbf{0}\right\}
= (C(P))^*,
\]
as claimed.
$\bigsquare$

\remark
If $Q = V(\Phi(x_1, \ldots, x_{d+1}))$ is
a projective or an affine quadric, it is obvious that
\[
V(\Phi(x_1, \ldots, x_{d+1})) = V(\lambda\Phi(x_1, \ldots, x_{d+1})) 
\]
for every $\lambda\not= 0$.
This raises the following question: If
\[
Q = V(\Phi_1(x_1, \ldots, x_{d+1})) = V(\Phi_2(x_1, \ldots, x_{d+1})),
\]
what is the relationship between $\Phi_1$ and $\Phi_2$?

\medskip
The answer depends crucially on the field over which
projective space or affine space is defined
({\it i.e.\/}, whether $Q \subseteq \mathbb{RP}^d$
or  $Q \subseteq \mathbb{CP}^d$ in the projective case or
whether $Q \subseteq \reals^{d+1}$ or  $Q \subseteq \complex^{d+1}$  
in the affine case).
For example, over $\reals$,
the polynomials
$\Phi_1(x_1, x_2, x_2) = x_1^2 + x_2^2$
and $\Phi_2(x_1, x_2, x_2) = 2x_1^2 + 3x_2^2$ both 
define the point $(0\co 0 \co 1)\in \projr{2}$, since
the only real solution of both $\Phi_1$ and $\Phi_2$ is $(0, 0)$.
However, if $Q$ has some nonsingular point, the following
can be proved (see Samuel \cite{Samuel}, Theorem 46 (Chapter 3)):

\begin{thm}
\label{weaknullsten}
Let $Q = V(\Phi(x_1, \ldots, x_{d+1})$ be a projective or an affine quadric,
over $\mathbb{RP}^d$ or $\reals^{d+1}$.
If $Q$ has a nonsingular point, then for every polynonial, $\Phi'$,
such that $Q =  V(\Phi'(x_1, \ldots, x_{d+1})$, there is some
$\lambda\not= 0$ ($\lambda\in \reals$) so that
$\Phi' = \lambda \Phi$. 
\end{thm} 
 
\medskip
In particular, Theorem \ref{weaknullsten} shows that
the equation of a nondegenerate quadric is unique up to a scalar.

\medskip
Actually, more is true. It turns out that if we allow complex solutions,
that is, if $Q \subseteq \mathbb{CP}^d$ in the projective case
or $Q \subseteq \complex^{d+1}$ in the affine case, then
$Q = V(\Phi_1) = V(\Phi_2)$ always implies $\Phi_2 = \lambda\Phi_1$
for some $\lambda\in \complex$, with $\lambda\not= 0$.
In the real case, the above holds (for some $\lambda\in \reals$,
with $\lambda\not= 0$) unless $Q$ is an affine subspace (resp.
a projective subspace) of dimension at most $d-1$ (resp. of
dimension at most $d - 2$). Even in this case, there is
a bijective affine map, $f$, (resp. a bijective projective map, $h$),
such that $\Phi_2 = \Phi_1\circ f^{-1}$
(resp.  $\Phi_2 = \Phi_1\circ h^{-1}$).
A proof of these facts (and more) can be found in Tisseron \cite{Tisseron}
(Chapter 3).

\medskip
We now have everything we need for a rigorous presentation of the
material of Section \ref{sec13}.
For  a comprehensive treatment of the affine and projective quadrics
and related material, the reader should consult 
Berger (Geometry II) \cite{Berger90b} or Samuel \cite{Samuel}.

\chapter[Basics of Combinatorial Topology]
{Basics of Combinatorial Topology}
\label{chapter3}
\section{Simplicial and Polyhedral Complexes}
\label{sec6}
In order to study and manipulate complex shapes
it is convenient to discretize these shapes and to 
view them as the union of simple building blocks
glued together in a ``clean fashion''. The building blocks
should be simple geometric objects, for example,
points, lines segments, triangles, tehrahedra and more generally
simplices, or even convex polytopes. We will begin 
by using simplices as building blocks.
The material presented in this chapter consists of the most
basic notions of combinatorial topology, going back roughly to  
the 1900-1930 period and it is covered in nearly every algebraic
topology book (certainly the ``classics'').
A classic text (slightly old fashion especially for the
notation and terminology)
is Alexandrov \cite{Alexandrov}, Volume 1 and another more ``modern''
source is Munkres  \cite{Munkresalg}.
An excellent treatment from the point of view of
computational geometry can be found is Boissonnat and
Yvinec \cite{Boissonnat}, especially Chapters 7 and 10.
Another fascinating book covering a lot of the basics
but devoted mostly to three-dimensional
topology and geometry is Thurston \cite{Thurston1}.

\medskip
Recall that a simplex is just the convex hull of a finite number of
affinely independent points. We also need to define faces,
the boundary, and the interior of a simplex.
\index{simplex}%

\begin{defin}
\label{simplexdef}
{\em 
Let $\affs$ be any normed affine space, say $\affs = \eucreal^\mdeg$
with its usual Euclidean norm. Given any $n+1$
affinely independent points $a_0,\ldots,a_n$ in $\affs$, the 
{\it $n$-simplex (or simplex) $\sigma$
defined by $a_0,\ldots,a_n$\/} 
\index{simplex!definition}%
\index{simplex@$n$-simplex}%
is the convex hull
of the points  $a_0,\ldots,a_n$, that is, the set of all
convex combinations $\lambda_0 a_0 + \cdots + \lambda_n a_n$, where 
$\lambda_0  + \cdots + \lambda_n = 1$ and $\lambda_i \geq 0$
for all $i$, $0\leq i \leq n$. We call $n$ the {\it dimension\/} of
\index{dimension!of a simplex}%
the $n$-simplex $\sigma$, and the points $a_0,\ldots,a_n$
are the {\it vertices\/} of $\sigma$.
\index{vertex!of a simplex}%
Given any subset $\{a_{i_{0}},\ldots,a_{i_{k}}\}$ of
$\{a_0,\ldots,a_n\}$ (where $0\leq k \leq n$), the $k$-simplex
generated by  $a_{i_{0}},\ldots,a_{i_{k}}$ is called a {\it $k$-face\/} 
or simply a {\it face\/}
\index{face of a simplex}%
of $\sigma$. A face $s$ of $\sigma$ is a {\it proper face\/}
if $s\not= \sigma$ (we agree that the empty set is a face of any simplex).
\index{proper!face of a simplex}%
For any vertex $a_i$, the face generated by 
$a_0,\ldots,a_{i-1},a_{i+1},\ldots,a_n$ (i.e.,  omitting $a_i$)
is called the {\it face opposite $a_i$\/}. Every face that
is an $(n-1)$-simplex is called a {\it boundary face\/} or {\it facet\/}.
\index{boundary face of a simplex}%
The union of the boundary faces is the {\it boundary of $\sigma$\/},
\index{boundary of a simplex}%
denoted by $\dBd \sigma$, and the complement of $\dBd \sigma$ in $\sigma$
\indsym{\dBd\noexpand\sigma}{boundary of a simplex}%
is the {\it interior\/} 
\index{interior of a simplex}%
$\dInt \sigma = \sigma - \dBd \sigma$ of $\sigma$.
The interior $\dInt \sigma$ of $\sigma$ 
\indsym{\dInt\noexpand\sigma}{interior of a simplex}%
is sometimes called an {\it open simplex\/}.
\index{simplex!open}%
}
\end{defin}

\medskip
It should be noted that for a $0$-simplex consisting of a single
point $\{a_0\}$, $\dBd \{a_0\} = \emptyset$, and
$\dInt \{a_0\} = \{a_0\}$.
Of course, a $0$-simplex is a single point, a $1$-simplex
is the line segment $(a_0, a_1)$, a $2$-simplex is a triangle
$(a_0, a_1, a_2)$ (with its interior),
and a $3$-simplex is a tetrahedron
$(a_0, a_1, a_2, a_3)$ (with its interior).
The inclusion relation between any two faces $\sigma$ and $\tau$
of some simplex, $s$, is written $\sigma \preceq \tau$. 

\medskip
We now state a number of properties of simplices, 
whose proofs are left as an exercise.
Clearly, a point $x$ belongs to the boundary $\dBd \sigma$ of $\sigma$
iff at least one of its barycentric coordinates
\index{barycentric coordinates}%
$(\lambda_0,\ldots,\lambda_n)$ is zero, and a point $x$ belongs to the 
interior $\dInt \sigma$ of $\sigma$ iff 
all of  its barycentric coordinates $(\lambda_0,\ldots,\lambda_n)$
are positive, i.e., $\lambda_i > 0$ for all $i$, $0\leq i\leq n$.
Then, for every $x\in\sigma$, there is a unique face $s$
such that $x\in \dInt s$, the face generated by those points $a_i$ for which
$\lambda_i > 0$, where $(\lambda_0,\ldots,\lambda_n)$ 
are the barycentric coordinates of $x$.

\medskip
A simplex $\sigma$ is convex,  arcwise connected, compact, and closed.
The interior $\dInt \sigma$ of a simplex is convex, 
arcwise connected, open, and
$\sigma$ is the closure of $\dInt \sigma$.

\medskip
We now put simplices together to form more complex shapes,
following  Munkres \cite{Munkresalg}. The intuition behind the
next definition is that the building blocks should
be ``glued cleanly''.

\begin{defin}
\label{complexdef}
{\em 
A {\it simplicial complex in $\eucreal^\mdeg$\/} (for short, a {\it complex\/}
\index{simplicial complex!definition}\index{complex}%
in $\eucreal^\mdeg$) is a set $K$
consisting of a (finite or infinite) set of simplices in $\eucreal^\mdeg$
satisfying the following conditions:
\begin{enumerate}
\item[(1)] 
Every face of a simplex in $K$ also belongs to $K$.
\item[(2)] 
For any two simplices $\sigma_1$ and $\sigma_2$ in $K$,
if $\sigma_1\cap\sigma_2\not=\emptyset$, then $\sigma_1\cap\sigma_2$
is a common face of both $\sigma_1$ and $\sigma_2$.
\end{enumerate}
Every $k$-simplex, $\sigma\in K$,  
is called a {\it $k$-face\/} (or {\it face\/}) of $K$.
A $0$-face $\{v\}$ is called a {\it  vertex\/} and a $1$-face
is called an {\it edge\/}.
\index{vertex}%
\index{dimension!of a complex}%
The {\it dimension\/} of the simplicial complex $K$ is the
maximum of the  dimensions of all simplices in $K$.
If $\mathrm{dim}\, K = d$, then every face of dimension $d$ is
called a {\it cell\/} and every face of dimension $d - 1$
is called a {\it facet\/}.
}
\end{defin}

\medskip
Condition (2) guarantees that the various simplices forming a complex
intersect nicely. It is easily shown that the following condition
is equivalent to condition (2):
\begin{enumerate}
\item[(2$'$)] For any two distinct simplices $\sigma_1, \sigma_2$,
$\dInt \sigma_1 \cap  \dInt \sigma_2 = \emptyset$.
\end{enumerate}

\remarks
\begin{enumerate}
\item
A simplicial complex, $K$, is a combinatorial object, namely, a 
{\it set\/} of simplices satisfying  certain conditions but not 
a subset of $\eucreal^m$. However, every complex, $K$, 
yields a subset of $\eucreal^m$ called the geometric realization
of $K$ and denoted $|K|$. This object will be defined shortly
and should not be confused with the complex.
Figure \ref{Figc1} illustrates this aspect of the definition
of a complex. For clarity, the two triangles ($2$-simplices)
are drawn as disjoint objects even though they share the common edge,
$(v_2, v_3)$ (a $1$-simplex) and similarly for the edges
that meet at some common vertex.
 
\begin{figure}
  \begin{center}
    \begin{pspicture}(0,0)(14,3.5)
\pspolygon*[fillstyle=solid,linecolor=green,linewidth=1.5pt]
(0,1.5)(1.5,0)(1.5,3)
\pspolygon*[fillstyle=solid,linecolor=green,linewidth=1.5pt]
(2.5,0)(4,1.5)(2.5,3)
    \pnode(5,1.5){u1}
    \pnode(6.5,0){u2}
    \pnode(6,1.5){u3}
    \pnode(7.5,3){u4}
    \pnode(8,0){u5}
    \pnode(8,3){u6}
    \pnode(8.5,3){u7}
    \pnode(10,1.5){u8}
    \pnode(9.5,0){u9}
    \pnode(11,1.5){u10}
    \cnode[fillstyle=solid,fillcolor=red](0,1.5){2pt}{w1}
    \cnode[fillstyle=solid,fillcolor=red](1.5,0){2pt}{w2}
    \cnode[fillstyle=solid,fillcolor=red](1.5,3){2pt}{w3}
    \cnode[fillstyle=solid,fillcolor=red](2.5,3){2pt}{w4}
    \cnode[fillstyle=solid,fillcolor=red](2.5,0){2pt}{w5}
    \cnode[fillstyle=solid,fillcolor=red](4,1.5){2pt}{w6}
    \ncline[linewidth=2pt,linecolor=blue]{u1}{u2}
    \ncline[linewidth=2pt,linecolor=blue]{u3}{u4}
    \ncline[linewidth=2pt,linecolor=blue]{u5}{u6}
    \ncline[linewidth=2pt,linecolor=blue]{u7}{u8}
    \ncline[linewidth=2pt,linecolor=blue]{u9}{u10}
    \ncline[linewidth=2pt,linecolor=blue]{w1}{w2}
    \ncline[linewidth=2pt,linecolor=blue]{w2}{w3}
    \ncline[linewidth=2pt,linecolor=blue]{w3}{w1}
    \ncline[linewidth=2pt,linecolor=blue]{w4}{w5}
    \ncline[linewidth=2pt,linecolor=blue]{w5}{w6}
    \ncline[linewidth=2pt,linecolor=blue]{w6}{w4}
    \cnode[fillstyle=solid,fillcolor=red](12,1.5){2.5pt}{v1}
    \cnode[fillstyle=solid,fillcolor=red](13,3){2.5pt}{v2}
    \cnode[fillstyle=solid,fillcolor=red](13,0){2.5pt}{v3}
    \cnode[fillstyle=solid,fillcolor=red](14,1.5){2.5pt}{v4}
    \uput[180](0,1.5){$v_1$}    
    \uput[-90](1.5,0){$v_2$}    
    \uput[90](1.5,3){$v_3$}
    \uput[90](2.5,3){$v_3$}        
    \uput[-90](2.5,0){$v_2$}
    \uput[0](4,1.5){$v_4$}                
    \end{pspicture}
  \end{center}
  \caption{A set of simplices forming a complex}
  \label{Figc1}
\end{figure}
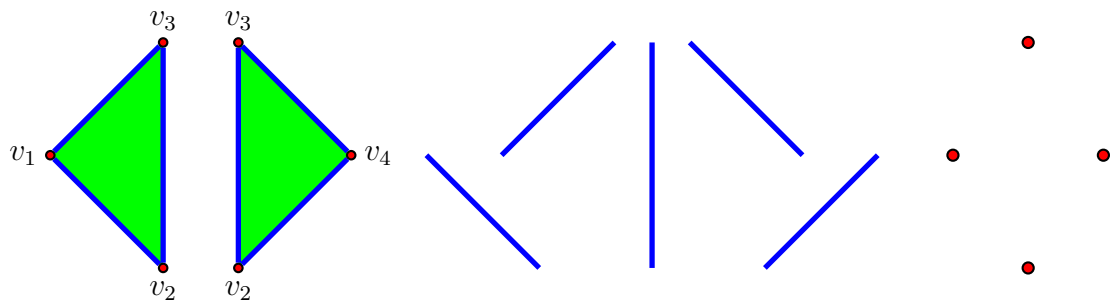

\item
Some authors define a {\it facet\/} of a complex, $K$, of dimension $d$
to be a $d$-simplex
in $K$, as opposed to a $(d - 1)$-simplex, as we did. This practice is not 
consistent with the notion of facet of a polyhedron and this is why
we prefer the terminology {\it cell\/} for the $d$-simplices in $K$.
\item
It is important to note that in order for a complex, $K$,
of dimension $d$ to be realized in $\eucreal^m$, the dimension
of the ``ambient space'', $m$, must be big enough. For example,
there are $2$-complexes that can't be realized in $\eucreal^3$
or even in $\eucreal^4$. There has to be enough room in order
for condition (2) to be satisfied. It is not hard to prove that
$m = 2d + 1$ is always sufficient. Sometimes, $2d$ works, for example
in the case of surfaces (where $d = 2$). 
\end{enumerate}

\begin{figure}
  \begin{center}
    \begin{pspicture}(0,0)(4,4.2)
\pspolygon*[fillstyle=solid,linecolor=lightgray,linewidth=1pt](1,0)(0.5,3)(3,1)
\pspolygon*[fillstyle=solid,linecolor=lightgray,linewidth=1pt](0,1)(4,2)(3,3)
    \cnode[fillstyle=solid,fillcolor=black](1,0){2pt}{v1}
    \cnode[fillstyle=solid,fillcolor=black](0.5,3){2pt}{v2}
    \cnode[fillstyle=solid,fillcolor=black](3,1){2pt}{v3}
    \ncline[linewidth=1pt]{v1}{v2}
    \ncline[linewidth=1pt]{v1}{v3}
    \ncline[linewidth=1pt]{v2}{v3}
    \cnode[fillstyle=solid,fillcolor=black](0,1){2pt}{v4}
    \cnode[fillstyle=solid,fillcolor=black](4,2){2pt}{v5}
    \cnode[fillstyle=solid,fillcolor=black](3,3){2pt}{v6}
    \ncline[linewidth=1pt]{v4}{v5}
    \ncline[linewidth=1pt]{v4}{v6}
    \ncline[linewidth=1pt]{v5}{v6}
    \end{pspicture}
  \hskip 1.5cm
    \begin{pspicture}(0,0)(4,4.2)
\pspolygon*[fillstyle=solid,linecolor=lightgray,linewidth=1pt](0,1)(2,4)(2,1)
\pspolygon*[fillstyle=solid,linecolor=lightgray,linewidth=1pt](2,0)(2,3)(4,1)
    \cnode[fillstyle=solid,fillcolor=black](0,1){2pt}{v1}
    \cnode[fillstyle=solid,fillcolor=black](2,4){2pt}{v2}
    \cnode[fillstyle=solid,fillcolor=black](2,1){2pt}{v3}
    \ncline[linewidth=1pt]{v1}{v2}
    \ncline[linewidth=1pt]{v1}{v3}
    \ncline[linewidth=1pt]{v2}{v3}
    \cnode[fillstyle=solid,fillcolor=black](2,0){2pt}{v4}
    \cnode[fillstyle=solid,fillcolor=black](2,3){2pt}{v5}
    \cnode[fillstyle=solid,fillcolor=black](4,1){2pt}{v6}
    \ncline[linewidth=1pt]{v4}{v5}
    \ncline[linewidth=1pt]{v4}{v6}
    \ncline[linewidth=1pt]{v5}{v6}
    \end{pspicture}
  \hskip 1.5cm
    \begin{pspicture}(0,0)(4,4.2)
\pspolygon*[fillstyle=solid,linecolor=lightgray,linewidth=1pt]
(1,0)(0,3)(4,3)(3,1)
    \pnode(1,0){v1}
    \cnode[fillstyle=solid,fillcolor=black](0,3){2pt}{v2}
    \cnode[fillstyle=solid,fillcolor=black](4,3){2pt}{v3}
    \cnode[fillstyle=solid,fillcolor=black](3,1){2pt}{v4}
    \ncline[linewidth=1pt]{v1}{v4}
    \ncline[linewidth=1pt]{v2}{v3}
    \ncline[linewidth=1pt]{v2}{v4}
    \ncline[linewidth=1pt]{v3}{v4}
    \end{pspicture}
  \end{center}
  \caption{Collections of simplices not forming a complex}
  \label{scomplex1}
\end{figure}
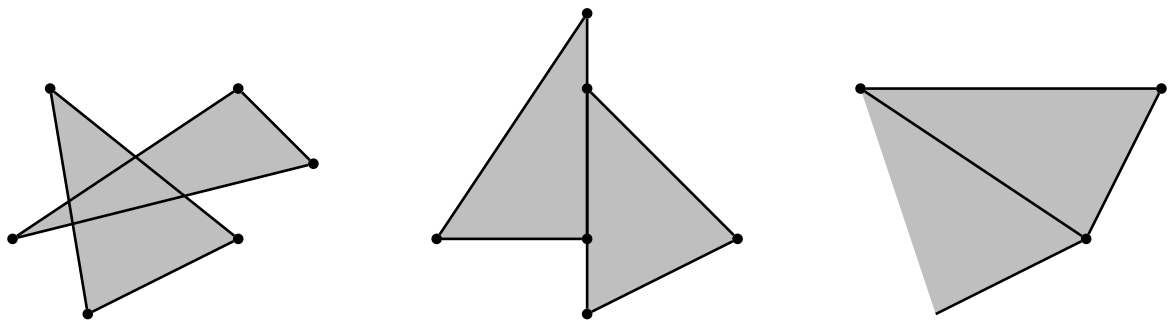

\medskip
Some collections of simplices violating some of the conditions
of Definition \ref{complexdef} are shown in Figure
\ref{scomplex1}. On the left, the intersection of the two $2$-simplices
is neither an edge nor a vertex of either triangle.
In the middle case, two simplices meet along an edge which
is not an edge of either triangle. On the right, there is a
missing edge and a missing vertex.

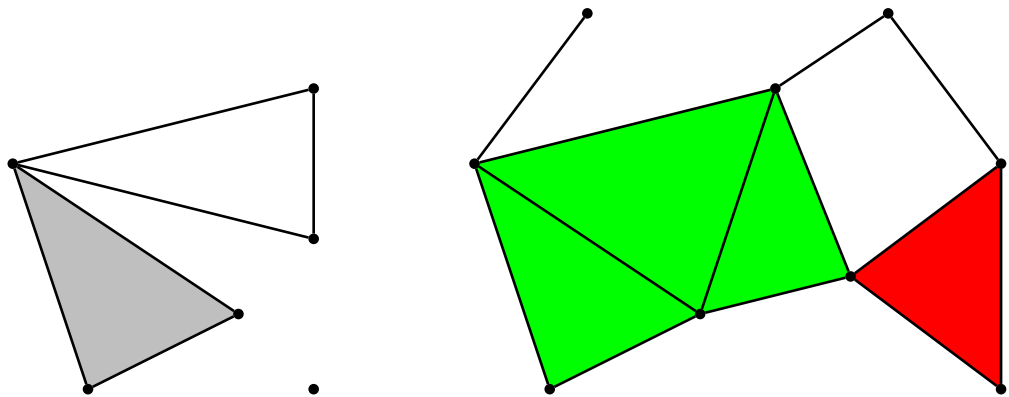
\begin{figure}
  \begin{center}
    \begin{pspicture}(0,0)(4,5.1)
\pspolygon*[fillstyle=solid,linecolor=lightgray,linewidth=1pt](1,0)(0,3)(3,1)
    \cnode[fillstyle=solid,fillcolor=black](1,0){2pt}{v1}
    \cnode[fillstyle=solid,fillcolor=black](0,3){2pt}{v2}
    \cnode[fillstyle=solid,fillcolor=black](3,1){2pt}{v3}
    \ncline[linewidth=1pt]{v1}{v2}
    \ncline[linewidth=1pt]{v1}{v3}
    \ncline[linewidth=1pt]{v2}{v3}
    \cnode[fillstyle=solid,fillcolor=black](4,4){2pt}{v4}
    \cnode[fillstyle=solid,fillcolor=black](4,2){2pt}{v5}
    \cnode[fillstyle=solid,fillcolor=black](4,0){2pt}{v6}
    \ncline[linewidth=1pt]{v2}{v4}
    \ncline[linewidth=1pt]{v4}{v5}
    \ncline[linewidth=1pt]{v5}{v2}
    \end{pspicture}
  \hskip 2cm
    \begin{pspicture}(0,0)(7,5.1)
\pspolygon*[fillstyle=solid,linecolor=green,linewidth=1pt]
(1,0)(0,3)(4,4)(5,1.5)(3,1)
\pspolygon*[fillstyle=solid,linecolor=red,linewidth=1pt]
(5,1.5)(7,3)(7,0)
    \cnode[fillstyle=solid,fillcolor=black](1,0){2pt}{v1}
    \cnode[fillstyle=solid,fillcolor=black](0,3){2pt}{v2}
    \cnode[fillstyle=solid,fillcolor=black](4,4){2pt}{v3}
    \cnode[fillstyle=solid,fillcolor=black](3,1){2pt}{v4}
    \cnode[fillstyle=solid,fillcolor=black](5,1.5){2pt}{v5}
    \cnode[fillstyle=solid,fillcolor=black](7,3){2pt}{v6}
    \cnode[fillstyle=solid,fillcolor=black](7,0){2pt}{v7}
    \cnode[fillstyle=solid,fillcolor=black](5.5,5){2pt}{v8}
    \cnode[fillstyle=solid,fillcolor=black](1.5,5){2pt}{v9}
    \ncline[linewidth=1pt]{v1}{v2}
    \ncline[linewidth=1pt]{v1}{v4}
    \ncline[linewidth=1pt]{v2}{v3}
    \ncline[linewidth=1pt]{v2}{v4}
    \ncline[linewidth=1pt]{v3}{v4}
    \ncline[linewidth=1pt]{v3}{v5}
    \ncline[linewidth=1pt]{v5}{v4}
    \ncline[linewidth=1pt]{v5}{v6}
    \ncline[linewidth=1pt]{v5}{v7}
    \ncline[linewidth=1pt]{v6}{v7}
    \ncline[linewidth=1pt]{v3}{v8}
    \ncline[linewidth=1pt]{v6}{v8}
    \ncline[linewidth=1pt]{v2}{v9}
    \end{pspicture}
  \end{center}
  \caption{Examples of simplicial complexes}
  \label{scomplex2}
\end{figure}

\medskip
Some ``legal'' simplicial complexes are shown in Figure \ref{scomplex2}.

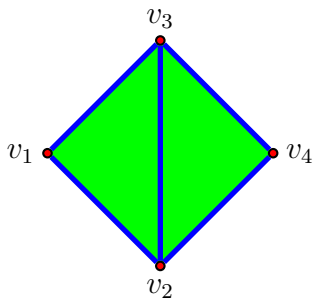
\begin{figure}
  \begin{center}
    \begin{pspicture}(0,0)(3.5,3.5)
\pspolygon*[fillstyle=solid,linecolor=green,linewidth=1.5pt]
(0,1.5)(1.5,0)(1.5,3)
\pspolygon*[fillstyle=solid,linecolor=green,linewidth=1.5pt]
(1.5,0)(3,1.5)(1.5,3)
    \cnode[fillstyle=solid,fillcolor=red](0,1.5){2pt}{w1}
    \cnode[fillstyle=solid,fillcolor=red](1.5,0){2pt}{w2}
    \cnode[fillstyle=solid,fillcolor=red](1.5,3){2pt}{w3}
    \cnode[fillstyle=solid,fillcolor=red](3,1.5){2pt}{w4}
    \ncline[linewidth=2pt,linecolor=blue]{w1}{w2}
    \ncline[linewidth=2pt,linecolor=blue]{w2}{w3}
    \ncline[linewidth=2pt,linecolor=blue]{w3}{w1}
    \ncline[linewidth=2pt,linecolor=blue]{w2}{w4}
    \ncline[linewidth=2pt,linecolor=blue]{w3}{w4}
    \uput[180](0,1.5){$v_1$}    
    \uput[-90](1.5,0){$v_2$}    
    \uput[90](1.5,3){$v_3$}
    \uput[0](3,1.5){$v_4$}                
    \end{pspicture}
  \end{center}
  \caption{The geometric realization of the complex of Figure
\ref{Figc1}}
  \label{Figc2}
\end{figure}

\medskip
The union $|K|$ of all the simplices in $K$ is a subset of $\eucreal^\mdeg$.
We can define a topology on $|K|$ by defining a subset $F$ of $|K|$ 
to be closed
iff $F\cap \sigma$ is closed in $\sigma$ for every face $\sigma\in K$.
It is immediately verified that the axioms of a topological space 
are indeed satisfied. 
\index{topology of the geometric realization}%
The resulting topological space $|K|$ is called
the {\it geometric realization of $K$\/}.
The geometric realization of the complex from Figure \ref{Figc1}
is shown in Figure \ref{Figc2}.
 
Obviously, $|\sigma| = \sigma$ for every simplex, $\sigma$.
Also, note that  distinct complexes may have the
same geometric realization. In fact, all the complexes 
obtained by subdividing the simplices of a given complex
yield the same geometric realization.
\index{geometric realization of $K$}%

%

\medskip
A {\it polytope\/} is the geometric realization
of some simplicial complex. 
\index{polytope!definition}%
A polytope of dimension $1$ is usually
called a {\it polygon\/}, and a polytope of dimension $2$ is usually called
\index{polygon!definition}%
a {\it polyhedron\/}. 
\index{polyhedron!definition}%
When $K$ consists of infinitely many simplices we usually require
that $K$ be {\it locally finite\/}, which means that
every vertex belongs to finitely many faces.
If $K$ is locally finite, then its geometric realization, $|K|$,  is
locally compact. 

\medskip
In the sequel, we will consider only finite simplicial complexes, that is,
complexes $K$ consisting of a finite number of simplices. In this
case, the topology of $|K|$ defined above is identical to
the topology induced from $\eucreal^\mdeg$. 
Also, for any simplex $\sigma$ in $K$, $\dInt \sigma$ coincides
with the interior $\interio{\sigma}$ of $\sigma$ in the topological sense, and
$\dBd \sigma$ coincides with the boundary of $\sigma$ in the topological sense.

\begin{defin}
\label{subcomplex}
{\em
Given any complex, $K_2$, a subset $K_1 \subseteq K_2$ of $K_2$ 
is a {\it subcomplex\/} of $K_2$ iff it is also a complex.
For any complex, $K$, of dimension $d$, for any $i$
with $0 \leq i \leq d$, the subset
\[
K^{(i)} = \{\sigma \in K \mid \mathrm{dim}\, \sigma \leq i\}
\]
is called the {\it $i$-skeleton\/} of $K$. Clearly, $K^{(i)}$
is a subcomplex of $K$. We also let
\[
K^{i} = \{\sigma \in K \mid \mathrm{dim}\, \sigma = i\}.
\]
Observe that $K^0$ is the set of vertices of $K$ and $K^i$ is
not a complex. A simplicial complex, $K_1$ is a
{\it subdivision\/} of a complex $K_2$ iff $|K_1| = |K_2|$
and if every face of $K_1$ is a subset of some face of $K_2$.
A complex $K$ of dimension $d$ is {\it pure\/} (or {\it homogeneous\/})
iff every face of $K$ is a face of some $d$-simplex of $K$
(i.e., some cell of $K$).
A complex is {\it connected\/} iff $|K|$ is connected.
}
\end{defin}

\begin{figure}
  \begin{center}
    \begin{pspicture}(0,0)(4,4.5)
\pspolygon*[fillstyle=solid,linecolor=lightgray,linewidth=1pt](1,0)(0,3)(3,1)
    \cnode[fillstyle=solid,fillcolor=black](1,0){2pt}{v1}
    \cnode[fillstyle=solid,fillcolor=black](0,3){2pt}{v2}
    \cnode[fillstyle=solid,fillcolor=black](3,1){2pt}{v3}
    \ncline[linewidth=1pt]{v1}{v2}
    \ncline[linewidth=1pt]{v1}{v3}
    \ncline[linewidth=1pt]{v2}{v3}
    \cnode[fillstyle=solid,fillcolor=black](4,4){2pt}{v4}
    \cnode[fillstyle=solid,fillcolor=black](4,2){2pt}{v5}
    \cnode[fillstyle=solid,fillcolor=black](4,0){2pt}{v6}
    \ncline[linewidth=1pt]{v2}{v4}
    \ncline[linewidth=1pt]{v4}{v5}
    \ncline[linewidth=1pt]{v5}{v2}
    \uput[90](2,4){(a)}    
    \end{pspicture}
  \hskip 2cm
    \begin{pspicture}(0,0)(7,4.5)
\pspolygon*[fillstyle=solid,linecolor=gray,linewidth=1pt]
(1,0)(0,3)(4,4)(5,1.5)(3,1)
\pspolygon*[fillstyle=solid,linecolor=gray,linewidth=1pt]
(5,1.5)(7,3)(7,0)
    \cnode[fillstyle=solid,fillcolor=black](1,0){2pt}{v1}
    \cnode[fillstyle=solid,fillcolor=black](0,3){2pt}{v2}
    \cnode[fillstyle=solid,fillcolor=black](4,4){2pt}{v3}
    \cnode[fillstyle=solid,fillcolor=black](3,1){2pt}{v4}
    \cnode[fillstyle=solid,fillcolor=black](5,1.5){2pt}{v5}
    \cnode[fillstyle=solid,fillcolor=black](7,3){2pt}{v6}
    \cnode[fillstyle=solid,fillcolor=black](7,0){2pt}{v7}
    \ncline[linewidth=1pt]{v1}{v2}
    \ncline[linewidth=1pt]{v1}{v4}
    \ncline[linewidth=1pt]{v2}{v3}
    \ncline[linewidth=1pt]{v2}{v4}
    \ncline[linewidth=1pt]{v3}{v4}
    \ncline[linewidth=1pt]{v3}{v5}
    \ncline[linewidth=1pt]{v5}{v4}
    \ncline[linewidth=1pt]{v5}{v6}
    \ncline[linewidth=1pt]{v5}{v7}
    \ncline[linewidth=1pt]{v6}{v7}
    \uput[90](2.5,4){(b)}    
    \uput[75](5,1.5){$v$}    
    \end{pspicture}
  \end{center}
  \caption{(a) A complex that is not pure. (b) A pure complex}
  \label{purex1}
\end{figure}
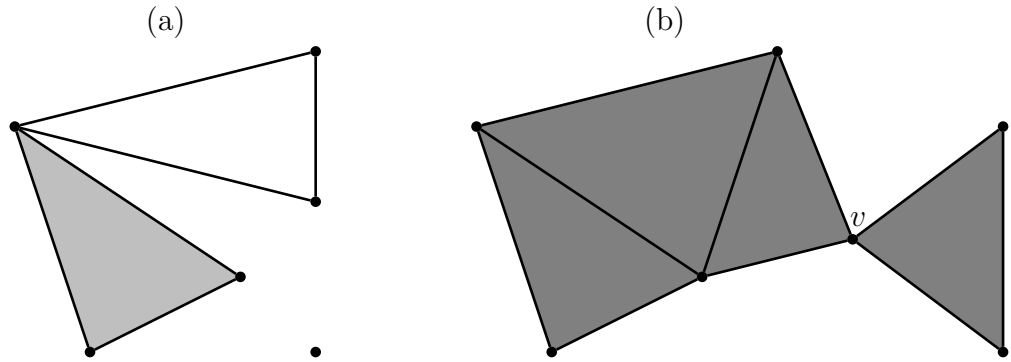

\medskip
It is easy to see that a complex is connected iff its $1$-skeleton
is connected. The intuition behind the notion of a pure complex, $K$,
of dimension $d$ is that a pure complex is the result of gluing
pieces all having the same dimension, namely, $d$-simplices.
For example, in Figure \ref{purex1}, the complex on the left is not pure
but the complex on the right is pure of dimension $2$.

\medskip
Most of the shapes that we will be interested in are well
approximated by pure complexes, in particular, surfaces
or solids. However, pure complexes may still have undesirable
``singularities'' such as the vertex, $v$, in Figure \ref{purex1}(b).
The notion  of link of a vertex provides a technical way
to deal with singularities.

\begin{defin}
\label{linkstar}
{\em
Let $K$ be any complex and let $\sigma$ be any face of $K$.
The {\it star\/}, $\mathrm{St}(\sigma)$ (or if we need to be very precise,
$\mathrm{St}(\sigma, K)$), of $\sigma$ is
the subcomplex of $K$ consisting of all faces, $\tau$,  
containing $\sigma$ 
and of all faces of $\tau$, {\it i.e.\/},
\[
\mathrm{St}(\sigma) = \{s \in K \mid (\exists \tau\in K)
(\sigma \preceq \tau\quad\hbox{and}\quad s \preceq \tau)\}.
\]
The {\it link\/}, $\mathrm{Lk}(\sigma)$ (or 
$\mathrm{Lk}(\sigma, K)$) of $\sigma$ is
the subcomplex of $K$ consisting of all faces in $\mathrm{St}(\sigma)$
that do not intersect $\sigma$, i.e.,
\[
\mathrm{Lk}(\sigma) =
\{\tau \in K \mid \tau\in \mathrm{St}(\sigma) \quad\hbox{and}\quad
\sigma\cap \tau = \emptyset\}.
\]
}
\end{defin}

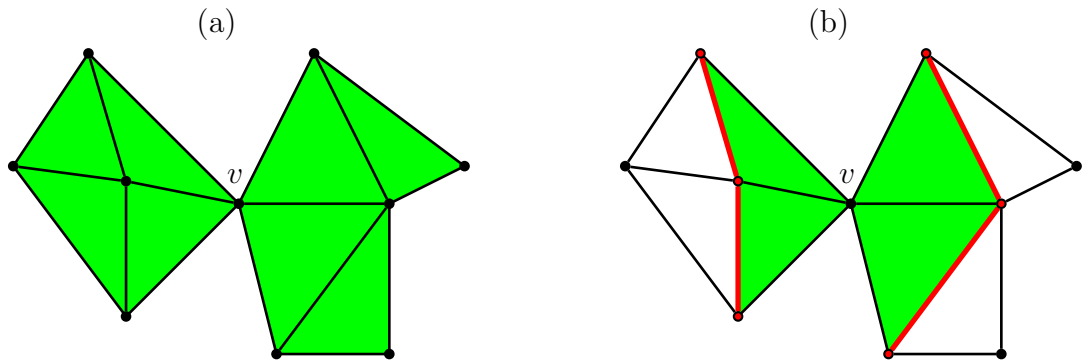
\begin{figure}
  \begin{center}
    \begin{pspicture}(0,0)(6,4.5)
\pspolygon[fillstyle=solid,fillcolor=green,linewidth=1pt]
(3,2)(1.5,0.5)(0,2.5)(1,4)
\pspolygon[fillstyle=solid,fillcolor=green,linewidth=1pt]
(3,2)(3.5,0)(5,0)(5,2)(6,2.5)(4,4)
    \cnode[fillstyle=solid,fillcolor=black](3,2){2pt}{v1}
    \cnode[fillstyle=solid,fillcolor=black](1.5,0.5){2pt}{v2}
    \cnode[fillstyle=solid,fillcolor=black](0,2.5){2pt}{v3}
    \cnode[fillstyle=solid,fillcolor=black](1,4){2pt}{v4}
    \cnode[fillstyle=solid,fillcolor=black](3.5,0){2pt}{v5}
    \cnode[fillstyle=solid,fillcolor=black](5,0){2pt}{v6}
    \cnode[fillstyle=solid,fillcolor=black](5,2){2pt}{v7}
    \cnode[fillstyle=solid,fillcolor=black](6,2.5){2pt}{v8}
    \cnode[fillstyle=solid,fillcolor=black](4,4){2pt}{v9}
    \cnode[fillstyle=solid,fillcolor=black](1.5,2.3){2pt}{v10}
    \ncline[linewidth=1pt]{v1}{v10}
    \ncline[linewidth=1pt]{v2}{v10}
    \ncline[linewidth=1pt]{v3}{v10}
    \ncline[linewidth=1pt]{v4}{v10}
    \ncline[linewidth=1pt]{v5}{v7}
    \ncline[linewidth=1pt]{v7}{v9}
    \ncline[linewidth=1pt]{v7}{v1}
    \uput[90](2.7,4){(a)}    
    \uput[90](2.95,2.1){$v$}    
    \end{pspicture}
 \hskip 2cm
    \begin{pspicture}(0,0)(6,4.5)
\pspolygon[fillstyle=solid,fillcolor=green,linewidth=1pt]
(3,2)(1.5,0.5)(1.5,2.3)(1,4)
\pspolygon[fillstyle=solid,fillcolor=green,linewidth=1pt]
(3,2)(3.5,0)(5,2)(5,2)(4,4)
    \cnode[fillstyle=solid,fillcolor=black](3,2){2pt}{v1}
    \cnode[fillstyle=solid,fillcolor=red](1.5,0.5){2pt}{v2}
    \cnode[fillstyle=solid,fillcolor=black](0,2.5){2pt}{v3}
    \cnode[fillstyle=solid,fillcolor=red](1,4){2pt}{v4}
    \cnode[fillstyle=solid,fillcolor=red](3.5,0){2pt}{v5}
    \cnode[fillstyle=solid,fillcolor=black](5,0){2pt}{v6}
    \cnode[fillstyle=solid,fillcolor=red](5,2){2pt}{v7}
    \cnode[fillstyle=solid,fillcolor=black](6,2.5){2pt}{v8}
    \cnode[fillstyle=solid,fillcolor=red](4,4){2pt}{v9}
    \cnode[fillstyle=solid,fillcolor=red](1.5,2.3){2pt}{v10}
    \ncline[linewidth=1pt]{v1}{v10}
    \ncline[linewidth=1pt]{v2}{v3}
    \ncline[linewidth=1pt]{v3}{v10}
    \ncline[linewidth=1pt]{v4}{v3}
    \ncline[linewidth=1pt]{v5}{v6}
    \ncline[linewidth=1pt]{v6}{v7}
    \ncline[linewidth=1pt]{v7}{v8}
    \ncline[linewidth=1pt]{v8}{v9}
    \ncline[linewidth=1pt]{v1}{v7}
    \ncline[linewidth=2pt,linecolor=red]{v2}{v10}
    \ncline[linewidth=2pt,linecolor=red]{v4}{v10}
    \ncline[linewidth=2pt,linecolor=red]{v5}{v7}
    \ncline[linewidth=2pt,linecolor=red]{v7}{v9}
    \uput[90](2.7,4){(b)}    
    \uput[90](2.95,2.1){$v$}    
    \end{pspicture}
  \end{center}
  \caption{(a) A complex. (b) Star and Link of $v$}
  \label{starlinkex}
\end{figure}

\medskip
To simplify notation, if $\sigma = \{v\}$ is a vertex we write
$\mathrm{St}(v)$ for $\mathrm{St}(\{v\})$ and
$\mathrm{Lk}(v)$ for $\mathrm{Lk}(\{v\})$.
Figure \ref{starlinkex} shows: 
\begin{enumerate}
\item[(a)]
A complex (on the left).
\item[(b)]
The star of the vertex $v$, indicated in gray
and the link of $v$, shown as thicker lines.
\end{enumerate}

\medskip
If $K$ is pure and of dimension $d$, then $\mathrm{St}(\sigma)$
is also pure of dimension $d$ and if $\mathrm{dim}\, \sigma = k$,
then $\mathrm{Lk}(\sigma)$ is pure of dimension $d - k - 1$.

\medskip
For technical reasons,  following Munkres \cite{Munkresalg}, 
besides defining
the complex, $\mathrm{St}(\sigma)$, it is useful to introduce
the {\it open star\/} of $\sigma$,   denoted $\mathrm{st}(\sigma)$,
defined as the subspace of $|K|$ consisting of the
union of the interiors, $\mathrm{Int}(\tau) = \tau - \partial\,\tau$,
of all the faces, $\tau$, containing,  $\sigma$.
According to this definition, the open star of $\sigma$
is not a complex but instead a subset of $|K|$. 

\medskip
Note that
\[
\overline{\mathrm{st}(\sigma)} = |\mathrm{St}(\sigma)|,
\]
that is, the closure of $\mathrm{st}(\sigma)$ is the geometric
realization of the complex $\mathrm{St}(\sigma)$.
Then, \\
$\mathrm{lk}(\sigma) = |\mathrm{Lk}(\sigma)|$
is the union of the simplices in $\mathrm{St}(\sigma)$ that
are disjoint from $\sigma$. If $\sigma$ is a vertex, $v$,
we have
\[
\mathrm{lk}(v) = \overline{\mathrm{st}(v)} - \mathrm{st}(v).
\]
However, beware that if $\sigma$ is not a vertex, then
$\mathrm{lk}(\sigma)$ is properly contained in 
$\overline{\mathrm{st}(\sigma)} - \mathrm{st}(\sigma)$!

\medskip
One of the nice properties of the open star, $\mathrm{st}(\sigma)$, 
of $\sigma$ is that it is open.
To see this, observe that for any point, $a\in |K|$, there is a unique
smallest simplex, $\sigma = (v_0, \ldots, v_k)$, such that
$a\in \mathrm{Int}(\sigma)$, that is, such that
\[
a = \lambda_0 v_0 + \cdots + \lambda_k v_k
\]
with $\lambda_i > 0$ for all $i$, with $0 \leq i \leq k$
(and of course, $\lambda_0 + \cdots + \lambda_k = 1$).
(When $k = 0$, we have $v_0 = a$ and $\lambda_0 = 1$.)
For every arbitrary vertex, $v$, of $K$, we define $t_v(a)$ by
\[
t_v(a) = \cases{
\lambda_i & if $v = v_i$, with $0 \leq i \leq k$, \cr
0 & if $v \notin \{v_0, \ldots, v_k\}$. \cr
}
\]
Using the above notation, observe that
\[
\mathrm{st}(v) = \{a\in |K| \mid t_v(a) > 0\}
\]
and thus, $|K| - \mathrm{st}(v)$ is the union of all
the faces of $K$ that do not contain $v$ as a vertex,
obviously a closed set. Thus, $\mathrm{st}(v)$ is open in $|K|$.
It is also quite clear that $\mathrm{st}(v)$ is path connected.
Moreover, for any $k$-face, $\sigma$, of $K$,
if $\sigma = (v_0, \ldots, v_k)$, then
\[
\mathrm{st}(\sigma) = \{a\in |K| \mid t_{v_i}(a) > 0,\quad 0\leq i \leq k\}, 
\]
that is,
\[
\mathrm{st}(\sigma) =
\mathrm{st}(v_0) \cap \cdots \cap \mathrm{st}(v_k).
\]
Consequently, $\mathrm{st}(\sigma)$ is open and path connected.

\danger
Unfortunately, the ``nice'' equation
\[
\mathrm{St}(\sigma) =
\mathrm{St}(v_0) \cap \cdots \cap \mathrm{St}(v_k)
\]
is false! (and anagolously for $\mathrm{Lk}(\sigma)$.) 
For a counter-example, consider the boundary of a tetrahedron
with one face removed.

\medskip
Recall that in $\eucreal^d$, the {\it (open) unit ball, $B^d$\/},
is defined by
\[
B^d = \{x\in \eucreal^d \mid \norme{x} < 1\},
\]
the {\it closed unit ball, $\overline{B}^d$\/},
is defined by
\[
\overline{B}^d = \{x\in \eucreal^d \mid \norme{x} \leq 1\},
\]
and the {\it $(d-1)$-sphere\/}, $S^{d - 1}$, by
\[
S^{d-1} = \{x\in \eucreal^d \mid \norme{x} = 1\}.
\]
Obviously, $S^{d-1}$ is the boundary of $\overline{B}^d$
(and $B^d$).

\begin{defin}
\label{nonsing}
{\em
Let $K$ be a pure complex of dimension $d$
and let $\sigma$ be any $k$-face of $K$, with $0 \leq k \leq d - 1$.
We say that $\sigma$ is {\it nonsingular\/} iff the geometric
realization, $\mathrm{lk}(\sigma)$, of the link of $\sigma$
is homeomorphic to either $S^{d - k - 1}$ or to $\overline{B}^{d - k - 1}$;
this is written as $\mathrm{lk}(\sigma) \approx S^{d - k - 1}$
or $\mathrm{lk}(\sigma) \approx \overline{B}^{d - k - 1}$,
where $\approx$ means homeomorphic.
}
\end{defin}

\medskip
In Figure \ref{starlinkex}, 
note that the link of $v$ is not homeomorphic to $S^1$ or $B^1$,
so $v$ is singular.

\medskip
It will also be useful to express $\mathrm{St}(v)$
in terms of  $\mathrm{Lk}(v)$, where $v$ is a vertex,
and for this, we define yet another notion of cone.

\begin{defin}
\label{conedef2}
{\em
Given any complex, $K$, in $\eucreal^n$, if $\mathrm{dim}\, K = d < n$,
for any point, $v\in \eucreal^n$, such that
$v$ does not belong to the affine hull of $|K|$, the {\it cone on $K$
with vertex $v$\/}, denoted, $v * K$, 
is the complex consisting of all simplices of the form
$(v, a_0, \ldots, a_k)$ and their faces,
where $(a_0, \ldots, a_k)$
is any $k$-face of $K$. If $K = \emptyset$, we set
$v * K = v$.
}
\end{defin}

\medskip
It is not hard to check that $v * K$ is indeed a complex
of dimension $d + 1$ containing $K$ as a subcomplex.

\remark
Unfortunately, the word ``cone'' is overloaded.
It might have been better to use the locution {\it pyramid\/} instead
of cone as some authors do  (for example, Ziegler). 
However, since we have been following Munkres \cite{Munkresalg},
a standard reference in algebraic topology, we decided to stick
with the terminology used in that book, namely, ``cone''.

\medskip
The following proposition is also easy to prove:

\begin{prop}
\label{comp1}
For any complex, $K$, of dimension $d$
and any vertex, $v\in K$, we have
\[
\mathrm{St}(v) = v * \mathrm{Lk}(v).
\] 
More generally, for any face, $\sigma$, of $K$, we have
\[
\overline{\mathrm{st}(\sigma)} =
|\mathrm{St}(\sigma)|  \approx \sigma \times |v * \mathrm{Lk}(\sigma)|,
\] 
for every $v\in \sigma$ and
\[
\overline{\mathrm{st}(\sigma)} - \mathrm{st}(\sigma) = 
\partial\, \sigma\times |v * \mathrm{Lk}(\sigma)|,
\]
for every $v\in \partial\, \sigma$.
\end{prop}

\begin{figure}
  \begin{center}
    \begin{pspicture}(0,0)(6,6.5)
\pspolygon*[fillstyle=solid,linecolor=green]
(0,2.5)(3,0)(5.5,4)(3,6)(0.5,4)
    \cnode[fillstyle=solid,fillcolor=black](3,0){2pt}{v1}
    \cnode[fillstyle=solid,fillcolor=red](0,2.5){2pt}{v2}
    \cnode[fillstyle=solid,fillcolor=red](2,2){2pt}{v3}
    \cnode[fillstyle=solid,fillcolor=red](4,2.5){2pt}{v4}
    \cnode[fillstyle=solid,fillcolor=red](5.5,4){2pt}{v5}
    \cnode[fillstyle=solid,fillcolor=black](3,6){2pt}{v6}
    \cnode[fillstyle=solid,fillcolor=red](0.5,4){2pt}{v7}
    \ncline[linewidth=1pt]{v1}{v2}
    \ncline[linewidth=1pt]{v1}{v3}
    \ncline[linewidth=1pt]{v1}{v4}
    \ncline[linewidth=1pt]{v1}{v5}
    \ncline[linewidth=2pt,linecolor=red]{v2}{v3}
    \ncline[linewidth=1pt]{v2}{v6}
    \ncline[linewidth=2pt,linecolor=red]{v2}{v7}
    \ncline[linewidth=1pt]{v3}{v6}
    \ncline[linewidth=2pt,linecolor=red]{v3}{v4}
    \ncline[linewidth=1pt]{v4}{v6}
    \ncline[linewidth=2pt,linecolor=red]{v4}{v5}
    \ncline[linewidth=1pt]{v5}{v6}
    \ncline[linewidth=1pt]{v6}{v7}
    \ncline[linewidth=1pt,linestyle=dashed]{v1}{v7}
    \ncline[linewidth=1pt,linestyle=dashed]{v1}{v6}
    \ncline[linewidth=2pt,linestyle=dashed,linecolor=red]{v5}{v7}
    \uput[180](0,2.5){$v_1$}    
    \uput[-30](2,2){$v_2$}    
    \uput[-10](4,2.5){$v_3$}    
    \uput[0](5.5,4){$v_4$}
    \uput[180](0.5,4){$v_5$}
    \uput[90](3,6){$v_6$}
    \uput[-90](3,0){$v_7$}                
    \end{pspicture}
  \end{center}
  \caption{More examples of links and stars}
  \label{starlinkex2}
\end{figure}
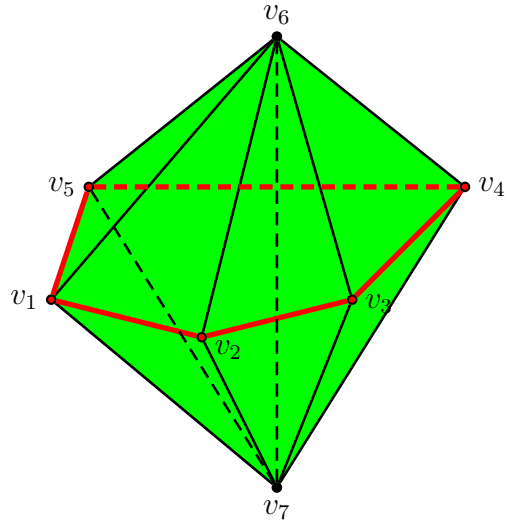

Figure \ref{starlinkex2} shows a $3$-dimensional complex.
The link of the edge $(v_6, v_7)$ is the pentagon 
$P = (v_1, v_2, v_3, v_4, v_5) \approx S^1$. The link of the vertex $v_7$
is the cone $v_6 * P \approx  B^2$. The link of $(v_1, v_2)$
is $(v_6, v_7) \approx B^1$ and the link of $v_1$ is
the union of the triangles $(v_2, v_6, v_7)$ and $(v_5, v_6, v_7)$,
which is homeomorphic to $B^2$.

\medskip
Given a pure complex,
it is necessary to distinguish between two kinds of faces.

\begin{defin}
\label{bounfacetdef}
{\em
Let $K$ be any pure complex of dimension $d$. A $k$-face, $\sigma$, 
of $K$ is a {\it boundary\/}  or {\it external\/}  face
iff it belongs to a single cell 
(i.e., a $d$-simplex) of $K$ and otherwise it is called an {\it internal\/} 
face ($0 \leq k \leq d - 1$).
The {\it boundary\/} of $K$, denoted $\mathrm{bd}(K)$, is the
subcomplex of $K$ consisting of all boundary facets of $K$ together with their
faces. 
}
\end{defin}

\medskip
It is clear by definition that $\mathrm{bd}(K)$ is a pure complex
of dimension $d - 1$. Even if $K$ is connected, $\mathrm{bd}(K)$ is
not connected, in general. For example, if $K$ is a $2$-complex
in the plane, the boundary of $K$ usually consists of several
simple closed polygons (i.e, $1$ dimensional complexes homeomorphic
to the circle, $S^1$).

\begin{prop}
\label{linkbd1}
Let $K$ be any pure complex of dimension $d$. For any $k$-face, $\sigma$, 
of $K$ the boundary complex, $\mathrm{bd}(\mathrm{Lk}(\sigma))$,
is nonempty iff $\sigma$ is a boundary face of $K$ ($0 \leq k \leq d - 2$).
Furthermore, 
$\mathrm{Lk}_{\mathrm{bd}(K)}(\sigma) = \mathrm{bd}(\mathrm{Lk}(\sigma))$
for every face, $\sigma$, of $\mathrm{bd}(K)$, where 
$\mathrm{Lk}_{\mathrm{bd}(K)}(\sigma)$ denotes the link of $\sigma$
in $\mathrm{bd}(K)$.
\end{prop}

\proof
Let $F$ be any facet of $K$ containing $\sigma$.
We may assume that $F = (v_0, \ldots, v_{d-1})$ and
$\sigma = (v_0, \ldots, v_k)$, in which case, 
$F' = (v_{k+1}, \ldots, v_{d-1})$ is a $(d - k - 2)$-face of $K$
and by definition of $\mathrm{Lk}(\sigma)$, we have
$F' \in \mathrm{Lk}(\sigma)$. Now, every cell (i.e., $d$-simplex),
$s$, containing $F$ is of the form $s = \mathrm{conv}(F\cup \{v\})$
for some vertex, $v$,
and $s' = \mathrm{conv}(F'\cup \{v\})$ is a $(d - k - 1)$-face in
$\mathrm{Lk}(\sigma)$ containing $F'$. Consequently,
$F'$ is an external face of $\mathrm{Lk}(\sigma)$ iff
$F$ is an external facet of $K$, establishing the proposition.
The second statement follows immediately from the proof
of the first.
$\bigsquare$

\medskip
Proposition \ref{linkbd1} shows that if every face of $K$ is
nonsingular, then the link of every internal face is a sphere
whereas the link of every external face is a ball.
The following proposition shows that for any pure complex, $K$,
nonsingularity of all the vertices is enough to imply 
that every open star is homeomorphic to $B^{d}$:

\begin{prop}
\label{comp2}
Let $K$ be any pure complex of dimension $d$. 
If every vertex of $K$ is nonsingular, then 
$\mathrm{st}(\sigma) \approx B^d$ for 
every $k$-face, $\sigma$, of $K$  ($1\leq k \leq d - 1$).
\end{prop}

\proof
Let $\sigma$ be any $k$-face of $K$ and assume that $\sigma$
is generated by the vertices $v_0, \ldots, v_k$,
with $1\leq k \leq d - 1$. By hypothesis, $\mathrm{lk}(v_i)$
is homeomorphic to either $S^{d-1}$ or $\overline{B}^{d-1}$. 
Then, it is easy to show that in either case,  we have
\[
|v_i * \mathrm{Lk}(v_i)| \approx \overline{B}^{d},
\]
and  by Proposition \ref{comp1}, we get
\[
|\mathrm{St}(v_i)|  \approx \overline{B}^d.
\]
Consequently, $\mathrm{st}(v_i) \approx B^d$.
Furthermore,  
\[
\mathrm{st}(\sigma) =
\mathrm{st}(v_0) \cap \cdots \cap \mathrm{st}(v_k) \approx B^d
\]
and so,
$\mathrm{st}(\sigma) \approx B^d$, as claimed. 
$\bigsquare$

\medskip
Here are more useful propositions about pure complexes
without singularities.

\begin{prop}
\label{comp3}
Let $K$ be any pure complex of dimension $d$. 
If every vertex of $K$ is nonsingular, then for every point, $a\in |K|$,
there is an open subset, $U \subseteq |K|$, containing $a$ such that
$U \approx B^d$ or $U \approx B^d \cap \mathbb{H}^d$, where
$\mathbb{H}^d = \{(x_1, \ldots, x_{d})\in \reals^{d} \mid x_d \geq  0\}$. 
\end{prop}

\proof
We already know from Proposition \ref{comp2} that
$\mathrm{st}(\sigma) \approx B^d$, for every $\sigma\in K$.
So, if $a\in \sigma$ and $\sigma$ is not a boundary face,
we can take $U= \mathrm{st}(\sigma) \approx B^d$.
If $\sigma$ is a boundary face, then $|\sigma| \subseteq
|\mathrm{bd}(\mathrm{St}(\sigma))|$ and it can be shown that
we can take $U =  B^d \cap \mathbb{H}^d$.
$\bigsquare$

\begin{prop}
\label{comp4}
Let $K$ be any pure complex of dimension $d$. 
If every facet of $K$ is nonsingular, then every facet of $K$,
is contained in at most two cells ($d$-simplices).
\end{prop}

\proof
If $|K| \subseteq \eucreal^d$, then this is an immediate consequence
of the definition of a complex. Otherwise, consider
$\mathrm{lk}(\sigma)$. By hypothesis, either 
$\mathrm{lk}(\sigma) \approx B^0$ or $\mathrm{lk}(\sigma) \approx S^0$.
As $B^0 = \{0\}$, $S^0 = \{-1, 1\}$ and 
$\mathrm{dim}\, \mathrm{Lk}(\sigma) = 0$, we deduce that 
$\mathrm{Lk}(\sigma)$ has either one or two points, which proves
that $\sigma$ belongs to at most two $d$-simplices.
$\bigsquare$

\begin{prop}
\label{comp5}
Let $K$ be any pure and connected complex of dimension $d$. 
If every face of $K$ is nonsingular, then for every pair
of cells ($d$-simplices), $\sigma$ and $\sigma'$, there  is a 
sequence of cells, $\sigma_0, \ldots, \sigma_p$,
with $\sigma_0 = \sigma$ and $\sigma_p = \sigma'$, and such
that $\sigma_i$ and $\sigma_{i+1}$ have a common facet,
for $i = 0, \ldots, p - 1$.
\end{prop}

\proof
We proceed by induction on $d$, using the fact that
the links are connected for $d \geq 2$.
$\bigsquare$

\begin{prop}
\label{comp6}
Let $K$ be any pure complex of dimension $d$. 
If every facet of $K$ is nonsingular, then the boundary,
$\mathrm{bd}(K)$, of $K$ is a pure complex of dimension
$d - 1$ with an empty boundary. Furthermore, if every face of
$K$ is nonsingular, then every face of $\mathrm{bd}(K)$ is
also nonsingular.
\end{prop}

\proof
Left as an exercise.
$\bigsquare$

\medskip
The building blocks of simplicial complexes, namely,
simplicies, are in some sense mathematically ideal.
However, in practice, it may be desirable to use 
a more flexible set of building blocks. We can indeed do this
and use convex polytopes as our  building blocks. 

\begin{defin}
\label{complexdefw}
{\em 
A {\it polyhedral complex in $\eucreal^\mdeg$\/} (for short, a {\it complex\/}
in $\eucreal^\mdeg$) is a set, $K$,
consisting of a (finite or infinite) set of convex
polytopes in $\eucreal^\mdeg$
satisfying the following conditions:
\begin{enumerate}
\item[(1)] 
Every face of a polytope in $K$ also belongs to $K$.
\item[(2)] 
For any two polytopes $\sigma_1$ and $\sigma_2$ in $K$,
if $\sigma_1\cap\sigma_2\not=\emptyset$, then $\sigma_1\cap\sigma_2$
is a common face of both $\sigma_1$ and $\sigma_2$.
\end{enumerate}
Every polytope, $\sigma\in K$, of dimension $k$,  
is called a {\it $k$-face\/} (or {\it face\/}) of $K$.
A $0$-face $\{v\}$ is called a {\it  vertex\/} and a $1$-face
is called an {\it edge\/}.
The {\it dimension\/} of the polyhedral complex $K$ is the
maximum of the  dimensions of all polytopes in $K$.
If $\mathrm{dim}\, K = d$, then every face of dimension $d$ is
called a {\it cell\/} and every face of dimension $d - 1$
is called a {\it facet\/}.
}
\end{defin}

\remark
Since the building blocks of a polyhedral complex are 
convex {\it polytopes\/} it might be more appropriate to use the
term ``polytopal complex'' rather than ``polyhedral complex''
and some authors do that. On the other hand, 
most of the traditional litterature uses the terminology
{\it polyhedral complex\/} so we will stick to it.
There is a notion of complex where the building blocks
are cones but these are called {\it fans\/}.

\medskip
Every convex polytope, $P$, yields two natural polyhedral complexes:
\begin{enumerate}
\item[(i)]
The polyhedral complex, $\s{K}(P)$, consisting of $P$
together with all of its faces. This complex has a single cell,
namely, $P$ itself.
\item[(ii)]
The {\it boundary complex\/}, $\s{K}(\partial P)$,
consisting of all faces of $P$ other than $P$ itself.
The cells of $\s{K}(\partial P)$ are the facets of $P$.
\end{enumerate}

\medskip
The notions of $k$-skeleton and pureness are defined just as
in the simplicial case. The notions of star and link are defined
for polyhedral complexes just as they are defined for
simplicial complexes except that the word ``face''
now means face of a polytope. Now, by Theorem \ref{equivpoly},
every polytope, $\sigma$, is the convex hull of its
vertices. Let $\mathrm{vert}(\sigma)$
denote the set of vertices of $\sigma$. Then,
we have the following crucial observation:
Given any polyhedral complex, $K$, for every point,
$x\in |K|$, there is a {\it unique\/} polytope, $\sigma_x\in K$,
such that $x\in \mathrm{Int}(\sigma_x) = \sigma_x - \partial\, \sigma_x$.
We define a function, $\mapdef{t}{V}{\reals_+}$,
that tests whether $x$ belongs to the interior of
any face (polytope)  of $K$ having $v$ as a vertex as follows: 
For every vertex,
$v$, of $K$,
\[
t_v(x) = \cases{
1 & if $v\in \mathrm{vert}(\sigma_x)$\cr
0 &  if $v\notin \mathrm{vert}(\sigma_x)$,\cr
}
\]
where $\sigma_x$ is the unique face of $K$ such that  
$x\in \mathrm{Int}(\sigma_x)$.

\medskip
Now, just as in the simplicial case, the open star, 
$\mathrm{st}(v)$, of a vertex, $v\in K$, is given by
\[
\mathrm{st}(v) = \{x\in |K| \mid t_v(x) = 1\}
\]
and it is an open subset of $|K|$ (the set $|K| - \mathrm{st}(v)$ is
the union of the polytopes of $K$ that do not contain $v$ as a vertex,
a closed subset of $|K|$).
Also, for any face, $\sigma$, of $K$,
the open star, $\mathrm{st}(\sigma)$, of $\sigma$ is given by
\[
\mathrm{st}(\sigma) = \{x\in |K| \mid t_v(x) = 1,\> \hbox{for all }\>
v \in  \mathrm{vert}(\sigma)\}
= \bigcap_{v\in \mathrm{vert}(\sigma)} \mathrm{st}(v).
\]
Therefore, $\mathrm{st}(\sigma)$ is also open in $|K|$.

\medskip
The next proposition is another result that seems
quite obvious, yet a rigorous proof is more involved
that we might think. This proposition states that a
convex polytope can always be cut up into simplices,
that is, it can be subdivided into a simplicial complex.
In other words, every convex polytope can be triangulated.
This implies that simplicial
complexes are as general as polyhedral complexes.

\medskip
One should be warned that even though, in the plane,
every bounded region (not necessarily convex)  whose boundary
consists of  a finite number of closed polygons 
(polygons homeomorphic to the circle, $S^1$) can be
triangulated, this is no longer true in three dimensions!

\begin{prop}
\label{triangul1}
Every convex $d$-polytope, $P$,  can be subdivided into
a simplicial complex without adding any new vertices,
i.e., every convex polytope can be triangulated.
\end{prop} 

\medskip\noindent
{\it Proof sketch\/}.
It would be tempting to proceed by induction on the dimension, $d$, of
$P$ but we do not know any correct proof of this kind.
Instead, we proceed by induction on the number, $p$, of vertices
of $P$. Since $\mathrm{dim}(P) = d$, we must have $p \geq d + 1$.
The case $p = d + 1$ corresponds to a simplex, so
the base case holds.

\medskip
For $p > d + 1$, we can pick some vertex, $v\in P$,  such that the convex
hull, $Q$,  of the remaining $p - 1$ vertices still has dimension $d$.
Then, by the induction hypothesis, $Q$, has a simplicial
subdivision. Now, we say that a facet, $F$, of $Q$ is {\it visible from
$v$\/} iff $v$ and the interior of $Q$ are strictly separated by 
the supporting hyperplane of $F$. Then, we add the $d$-simplices,
$\mathrm{conv}(F \cup \{v\}) = v * F$, for every facet, $F$, of $Q$
visible from $v$ to those in the triangulation of $Q$. We claim that
the resulting collection of simplices (with their faces)
constitutes a simplicial complex subdividing $P$.
This is the part of the proof that requires a careful
and somewhat tedious case analysis, which we omit.
However, the reader should check that everything really works out!
$\bigsquare$

\medskip
With all this preparation, it is now quite natural to define
combinatorial manifolds.

\section{Combinatorial and Topological Manifolds}
\label{sec7}
The notion of pure complex without singular faces turns out to
be a very good ``discrete'' approximation of the notion of
(topological) manifold because  of its highly
computational nature. This motivates the following definition:

\begin{defin}
\label{compmanif}
{\em
A {\it combinatorial $d$-manifold\/} is any space, $X$, homeomorphic
to the geometric realization, $|K| \subseteq \eucreal^n$, 
of some pure (simplicial or polyhedral) complex, $K$, 
of dimension $d$ 
whose faces are all nonsingular. If the link of every $k$-face of $K$
is homeomorphic to the sphere $S^{d - k - 1}$, we say that
$X$ is a combinatorial manifold {\it without boundary\/}, else it is
a combinatorial manifold {\it with boundary\/}.
}
\end{defin}

\medskip
Other authors use the term {\it triangulation\/}
for what we call a combinatorial manifold. 

\medskip
It is easy to see that the connected components of a 
combinatorial $1$-manifold are either simple closed polygons
or simple chains (``simple'' means that the interiors of distinct edges
are disjoint). 
A combinatorial $2$-manifold which is connected 
is also called a {\it combinatorial
surface\/} (with or without boundary). 
Proposition \ref{comp6} immediately yields the following result:

\begin{prop}
\label{combman1}
If $X$ is a combinatorial $d$-manifold with boundary,
then $\mathrm{bd}(X)$ is a combinatorial $(d-1$)-manifold 
without boundary.
\end{prop}

\medskip
Now, because we are assuming that $X$ sits in some Euclidean space,
$\eucreal^n$, the space $X$ is Hausdorff and second-countable.
(Recall that a topological space is second-countable iff there is
a countable family, $\{U_i\}_{i\geq 0}$, of open sets of $X$ such that
every open subset of $X$ is the union of open sets from this family.)
Since it is desirable to have a good match between manifolds
and combinatorial manifolds, we are led to the definition below.

\medskip
Recall that
\[
\mathbb{H}^d = \{(x_1, \ldots, x_{d})\in \reals^{d} \mid x_d \geq  0\}.
\] 

\begin{defin}
\label{boundman}
{\em
For any $d\geq 1$, a {\it (topological) $d$-manifold
with boundary\/} is
a second-countable, topological Hausdorff space $M$, together
with  an open cover, $(U_i)_{i\in I}$, of open sets in $M$ and
a family, $(\varphi_i)_{i\in I}$, of homeomorphisms,
$\mapdef{\varphi_i}{U_i}{\Omega_i}$, where each $\Omega_i$
is some open subset of $\mathbb{H}^d$ in the subset topology. 
Each pair $(U, \varphi)$
is called a {\it  coordinate system\/}, or
{\it chart\/}, of $M$, 
each homeomorphism $\mapdef{\varphi_i}{U_i}{\Omega_i}$ is called
a {\it coordinate map\/},
and its inverse  $\mapdef{\varphi^{-1}_{i}}{\Omega_i}{U_i}$ is called a 
{\it  parameterization\/} of $U_i$.
The family  $(U_i, \varphi_i)_{i\in I}$ is often called an
{\it atlas\/} for $M$.
A {\it (topological) bordered surface\/} is a connected $2$-manifold
with boundary.
If for every homeomorphism, $\mapdef{\varphi_i}{U_i}{\Omega_i}$,
the open set $\Omega_i\subseteq \mathbb{H}^d$ is actually an open set in
$\reals^d$ (which means that $x_d > 0$ 
for every $(x_1, \ldots, x_d)\in \Omega_i$), then we say that
$M$ is a {\it $d$-manifold\/}.
}
\end{defin}

\medskip
Note that a $d$-manifold is also a $d$-manifold with boundary.

\medskip
If $\mapdef{\varphi_i}{U_i}{\Omega_i}$ is some
homeomorphism onto some open set $\Omega_i$ of $\mathbb{H}^d$ 
in the subset topology, some $p\in U_i$ may be mapped into
$\reals^{d-1}\times\reals_+$, or into the ``boundary''
$\reals^{d-1}\times\{0\}$ of $\mathbb{H}^d$.
Letting $\dBd\mathbb{H}^d = \reals^{d-1}\times\{0\}$, it can be
shown using homology that if
some   coordinate map, $\varphi$, defined on $p$
maps $p$ into $\dBd\mathbb{H}^d$, then
every coordinate map, $\psi$, defined on $p$
maps $p$ into $\dBd\mathbb{H}^d$.

\medskip
Thus, $M$ is the disjoint union of two sets
$\dBd M$ and $\dInt M$, where $\dBd M$ is the
subset consisting of all points $p\in M$ that are mapped
by some (in fact, all)  coordinate map, $\varphi$, defined on $p$
into $\dBd\mathbb{H}^d$, and where 
$\dInt M = M - \dBd M$. The set $\dBd M$ is called
the {\it boundary\/} of $M$, and the set $\dInt M$ is called
the {\it interior\/} of $M$, even though this terminology
clashes with some prior topological definitions.
A good example of a bordered surface is the M\"obius strip.
The boundary of the M\"obius strip is a circle.

\medskip
The boundary $\dBd M$ of $M$
may be empty, but $\dInt M$ is nonempty. Also, it can be shown
using homology that the integer $d$ is unique.
It is clear that $\dInt M$ is open and a $d$-manifold, 
and that $\dBd M$ is closed. If $p\in \dBd M$, and 
$\varphi$ is some coordinate map  defined on $p$,
since $\Omega=\varphi(U)$ is an open subset of $\dBd\mathbb{H}^d$, 
there is some open half ball $B_{o+}^{d}$ centered
at $\varphi(p)$ and contained in  $\Omega$ which intersects
$\dBd\mathbb{H}^d$ along an open ball $B_{o}^{d-1}$,
and if we consider $W = \varphi^{-1}(B_{o+}^{d})$,
we have an open subset of $M$ containing $p$ which is mapped
homeomorphically onto $B_{o+}^{d}$ in such that way
that every point in $W\cap \dBd M$ is mapped
onto the open ball  $B_{o}^{d-1}$. Thus, it is easy to see that
$\dBd M$ is a $(d-1)$-manifold.

\begin{prop}
\label{manp1}
Every combinatorial $d$-manifold is a $d$-manifold with boundary.
\end{prop}

\proof
This is an immediate consequence of Proposition \ref{comp3}.
$\bigsquare$

\medskip
Is the converse of Proposition  \ref{manp1} true?

\medskip
It turns out that answer is yes for $d = 1, 2, 3$ but {\bf no}
for $d \geq 4$. This is not hard to prove for $d = 1$. 
For $d = 2$ and $d = 3$,
this is quite hard to prove; among other things, it is necessary
to prove that triangulations exist and this is very technical.
For $d \geq 4$, not every manifold can be triangulated
(in fact, this is undecidable!).

\medskip
What if we assume that $M$ is a triangulated manifold, which 
means that $M \approx |K|$, for some pure $d$-dimensional
complex, $K$?

\medskip
Surprinsingly, for $d \geq 5$, there are triangulated manifolds 
whose links are not spherical (i.e., not homeomorphic to
$\overline{B}^{d - k - 1}$ or $S^{d - k - 1}$), see
Thurston \cite{Thurston1}. 
 
\medskip
Fortunately, we will only have to deal with $d = 2, 3$!
Another issue that must be addressed is orientability.

\medskip
Assume that we fix a total  ordering of the vertices of a complex,
$K$. Let  $\sigma = (v_0, \ldots, v_k)$ be any simplex.
Recall that every permutation 
(of $\{0, \ldots, k\}$) is a product of
{\it transpositions\/}, where a transposition swaps two distinct
elements, say $i$ and $j$, and leaves every other element fixed.
Furthermore, for any permutation, $\pi$, 
the parity of the number of transpositions needed to obtain
$\pi$ only depends on $\pi$ and it called the {\it signature\/}
of $\pi$. We say that {\it two permutations are equivalent\/} iff
they have the same signature. Consequently, there are two
equivalence classes of permutations: Those of even signature
and those of odd signature.  
Then, an {\it orientation\/} of $\sigma$ is the choice of
one of the two equivalence classes of permutations of its vertices.
If $\sigma$ has been given an 
orientation, then we denote by $-\sigma$ the result of assigning the
other orientation to it (we call it the {\it opposite orientation\/}). 

\medskip
For example, $(0, 1, 2)$ has the two orientation classes:
\[
\{(0, 1, 2), (1, 2, 0), (2, 0, 1)\}\quad\hbox{and}\quad
\{(2, 1, 0), (1, 0, 2), (0, 2, 1)\}.
\]

\begin{defin}
\label{orientdef1}
{\em
Let $X \approx |K|$ be a combinatorial $d$-manifold.
We say that $X$ is {\it orientable\/} if it is possible to 
assign an orientation to all of its cells ($d$-simplices) so that
whenever two cells $\sigma_1$ and $\sigma_2$ have a common
facet, $\sigma$, the two orientations induced by 
$\sigma_1$ and $\sigma_2$ on $\sigma$ are opposite.
A combinatorial $d$-manifold together with a specific
orientation of its cells is called an {\it oriented manifold\/}.
If $X$ is not orientable we say that it is {\it non-orientable\/}.
} 
\end{defin}

\remark
It is possible to define the notion of orientation of a manifold
but this is quite technical and we prefer to avoid digressing 
into this matter. This shows another advantage of combinatorial 
manifolds: The definition of orientability is simple and quite natural.

\medskip
There are non-orientable (combinatorial) 
surfaces, for example, the M\"obius strip
which can be realized in $\eucreal^3$. The  M\"obius strip is
a surface with boundary, its boundary being a circle.
There are also non-orientable (combinatorial) surfaces
such as the Klein bottle or the projective plane but they can only be 
realized in $\eucreal^4$ (in $\eucreal^3$, they must have 
singularities such as self-intersection).
We will only be dealing with orientable manifolds
and, most of the time, surfaces.

\medskip
One of the most important invariants of combinatorial (and topological)
manifolds is their {\it Euler(-Poincar\'e) characteristic\/}. In the next
chapter, we prove a famous formula due to Poincar\'e giving
the Euler characteristic of a convex polytope. For this, we
will introduce a technique of independent interest called {\it shelling\/}.

\chapter[Shellings and the Euler-Poincar\'e Formula]
{Shellings, the Euler-Poincar\'e Formula for Polytopes,
the Dehn-Sommerville Equations and the Upper Bound Theorem}
\label{chap4}
\section[Shellings]
{Shellings}
\label{sec8}
The notion of shellability is motivated by the desire to give
an inductive proof of the Euler-Poincar\'e formula in
any dimension. Historically, this formula was discovered by Euler 
for three dimensional polytopes in 1752 (but it was already known
to Descartes around 1640). 
If $f_0, f_1$ and $f_2$ denote the number
of vertices, edges and triangles of the three dimensional
polytope, $P$, (i.e., the number of $i$-faces of $P$
for $i = 0, 1, 2$), then the {\it Euler formula\/} states that
\[
f_0 - f_1 + f_2 = 2.
\]
The proof of Euler's formula is not very difficult but one
still has to exercise caution. 
Euler's formula was generalized 
to arbitrary $d$-dimensional polytopes by  Schl\"afli (1852) but the first
correct proof was given  by Poincar\'e.
For this, Poincar\'e had to lay the foundations of algebraic topology
and after a first ``proof'' given in 1893 (containing some flaws)
he finally gave the first
correct proof in 1899. If $f_i$ denotes the number of $i$-faces
of the $d$-dimensional polytope, $P$, (with $f_{-1} = 1$
and $f_d = 1$), the {\it Euler-Poincar\'e formula\/} states that:
\[
\sum_{i = 0}^{d - 1} (-1)^i f_i = 1 - (-1)^d,
\]
which can also be written as
\[
\sum_{i = 0}^{d} (-1)^i f_i = 1,
\]
by incorporating  $f_d = 1$ in the first formula or as
\[
\sum_{i = -1}^{d} (-1)^i f_i = 0,
\]
by incorporating both $f_{-1} = 1$ and $f_d = 1$ in the first formula.

\medskip
Earlier inductive ``proofs'' of the above formula were proposed,
notably a proof by Schl\"afli in 1852, but it was later observed
that all these proofs assume that the boundary of every polytope
can be built up inductively in a nice way, what is called
{\it shellability\/}. Actually,  counter-examples
of shellability for various simplicial complexes suggested
that polytopes were perhaps not shellable. However, the fact
that polytopes are shellable was finally proved in 1970 by
Bruggesser and Mani \cite{Bruggesser}
and soon after that (also in 1970)
a striking application
of shellability was made by McMullen \cite{McMullen71}
who gave the first proof of 
the so-called ``upper bound theorem''.

\medskip
As shellability of polytopes is an important tool and as it yields one of
the cleanest inductive proof of the Euler-Poincar\'e formula,
we will sketch its proof in some details. This Chapter
is heavily inspired by Ziegler's excellent treatment
\cite{Ziegler97}, Chapter 8.
We begin with the definition of shellability. It's a bit technical,
so please be patient!

\begin{defin}
\label{shellable}
{\em 
Let $K$ be a pure polyhedral complex of dimension $d$.
A {\it shelling of $K$\/} is a list, $F_1, \ldots, F_s$, of
the cells (i.e., $d$-faces) of $K$ such that either $d = 0$
(and thus, all $F_i$ are points) or the following conditions hold:
\begin{enumerate}
\item[(i)]
The boundary complex, $\s{K}(\partial F_1)$, of the first cell, 
$F_1$, of $K$ has a shelling.
\item[(ii)]
For any $j$, $1<  j \leq s$, the intersection of the cell $F_j$
with the previous cells is nonempty and is an initial segment of
a shelling of the $(d - 1)$-dimensional boundary complex of $F_j$, that is
\[
F_j \cap \left(\bigcup_{i = 1}^{j - 1} F_i  \right) =  G_1 \cup G_2 \cup \cdots
\cup G_r,
\]
for some shelling $G_1, G_2, \ldots, G_r, \ldots, G_t$ of 
$\s{K}(\partial F_j)$, with $1 \leq r \leq t$.
As the intersection should be the  initial segment of a shelling for 
the $(d - 1)$-dimensional complex, $\partial F_j$, it has to be pure
$(d - 1)$-dimensional and connected for $d > 1$.
\end{enumerate}
A polyhedral complex is {\it shellable\/} if it is pure and  has a shelling.
}
\end{defin}

\medskip
Note that shellabiliy is only defined for pure complexes.
Here are some examples of shellable complexes:
\begin{enumerate}
\item[(1)]
Every $0$-dimensional complex, that is, evey set of points, is shellable,
by definition.
\item[(2)]
A $1$-dimensional complex is a graph without loops and parallel edges.
A $1$-dimensional complex is shellable iff it is connected,
which implies that it has no isolated vertices. Any ordering of the edges,
$e_1, \ldots, e_s$, such that $\{e_1, \ldots, e_i\}$
induces a connected subgraph for every $i$ will do. Such an ordering
can be defined inductively, due to the connectivity of the graph.
\item[(3)]
Every simplex is shellable. In fact, any ordering of its facets yields
a shelling. This is easily shown by induction on the dimension, since
the intersection of any two facets $F_i$ and $F_j$ is a facet of both
$F_i$ and $F_j$.
\item[(4)]
The $d$-cubes are shellable. By induction on the dimension, it can be shown
that every ordering of the $2d$ facets $F_1, \ldots, F_{2d}$
such that $F_1$ and $F_{2d}$ are opposite (that is, $F_{2d} = -F_1$)
yields a shelling.
\end{enumerate}

\begin{figure}
  \begin{center}
    \begin{pspicture}(0,0)(3,3)
\pspolygon[fillstyle=solid,fillcolor=lightgray,linewidth=1.5pt]
(0,0)(1.5,0)(1.5,1.5)(0,1.5)
\pspolygon[fillstyle=solid,fillcolor=lightgray,linewidth=1.5pt]
(1.5,1.5)(3,1.5)(3,3)(1.5,3)
    \uput[0](0.5,0.75){$1$}    
    \uput[0](2,2.25){$2$}    
    \end{pspicture}
  \hskip 1.5cm
    \begin{pspicture}(0,0)(3,3)
\pspolygon[fillstyle=solid,fillcolor=lightgray,linewidth=1.5pt]
(0,0)(1,0)(1,1)(0,1)
\pspolygon[fillstyle=solid,fillcolor=lightgray,linewidth=1.5pt]
(0,1)(1,1)(1,2)(0,2)
\pspolygon[fillstyle=solid,fillcolor=lightgray,linewidth=1.5pt]
(0,2)(1,2)(1,3)(0,3)
\pspolygon[fillstyle=solid,fillcolor=lightgray,linewidth=1.5pt]
(1,0)(2,0)(2,1)(1,1)
\pspolygon[fillstyle=solid,fillcolor=lightgray,linewidth=1.5pt]
(1,2)(2,2)(2,3)(1,3)
\pspolygon[fillstyle=solid,fillcolor=lightgray,linewidth=1.5pt]
(2,0)(3,0)(3,1)(2,1)
\pspolygon[fillstyle=solid,fillcolor=lightgray,linewidth=1.5pt]
(2,1)(3,1)(3,2)(2,2)
\pspolygon[fillstyle=solid,fillcolor=lightgray,linewidth=1.5pt]
(2,2)(3,2)(3,3)(2,3)
    \uput[0](0.2,2.5){$1$}    
    \uput[0](1.2,2.5){$2$}
    \uput[0](2.2,2.5){$3$}
    \uput[0](2.2,1.5){$4$}            
    \uput[0](2.2,0.5){$5$} 
    \uput[0](1.2,0.5){$6$}
    \uput[0](0.2,0.5){$7$}
    \uput[0](0.2,1.5){$8$}                                               
    \end{pspicture}
  \hskip 1.5cm
    \begin{pspicture}(0,0)(3,3)
\pspolygon[fillstyle=solid,fillcolor=lightgray,linewidth=1.5pt]
(0,2)(1,2)(1,3)(0,3)
\pspolygon[fillstyle=solid,fillcolor=lightgray,linewidth=1.5pt]
(1,2)(2,2)(2,3)(1,3)
\pspolygon[fillstyle=solid,fillcolor=lightgray,linewidth=1.5pt]
(2,2)(3,2)(3,3)(2,3)
\pspolygon[fillstyle=solid,fillcolor=lightgray,linewidth=1.5pt]
(1,1)(2,1)(2,2)(1,2)
\pspolygon[fillstyle=solid,fillcolor=lightgray,linewidth=1.5pt]
(1,0)(2,0)(2,1)(1,1)
    \uput[0](1.2,0.5){$1$}            
    \uput[0](1.2,1.5){$2$} 
    \uput[0](1.2,2.5){$3$}
    \uput[0](0.2,2.5){$4$}
    \uput[0](2.2,2.5){$5$}                                               
    \end{pspicture}
  \end{center}
  \caption{Non shellable and Shellable $2$-complexes}
  \label{figshell1}
\end{figure}
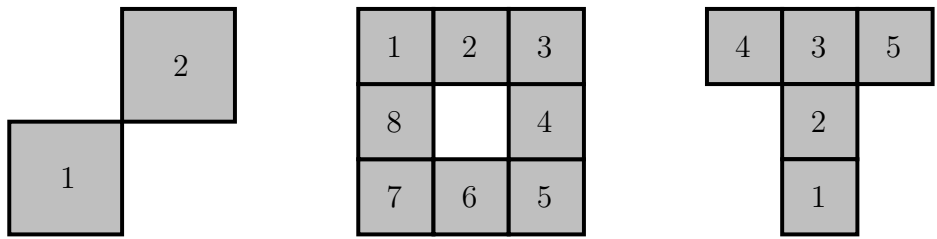

\medskip
However, already for $2$-complexes, problems arise.
For example, in Figure \ref{figshell1}, the left and the middle
$2$-complexes are not shellable but the right complex is shellable.
 
\medskip
The problem with the left complex is that cells $1$ and $2$
intersect at a vertex, which is not $1$-dimensional, and in the middle
complex, the intersection of cell  $8$ with its predecessors  is
not connected. In contrast, the ordering of the right complex
is a shelling. However, observe that the reverse ordering is
not a shelling because cell $4$ has an empty intersection with
cell $5$!

\remarks
\begin{enumerate}
\item
Condition (i) in Definition \ref{shellable} is redundant because,
as we shall prove shortly, every polytope is shellable.
However, if we want to use
this definition for more general complexes,
then condition (i) is necessary.
\item
When $K$ is a simplicial complex, condition (i) is of course redundant,
as every simplex is shellable but condition (ii) can also be simplified to:
\begin{enumerate}
\item[(ii')]
For any $j$, with $1 < j \leq s$, the intersection of $F_j$ with the
previous cells is nonempty and pure $(d - 1)$-dimensional.
This means that for every $i < j$ there is some $l < j$ such that
$F_i\cap F_j \subseteq F_l\cap F_j$ and $F_l\cap F_j$ is a facet of $F_j$.
\end{enumerate}
\end{enumerate}

\medskip
The following proposition yields an important piece of information
about the local structure of shellable simplicial complexes:

\begin{prop}
\label{shellp1}
Let $K$ be a shellable simplicial complex and say 
$F_1, \ldots, F_s$ is a shelling for $K$. Then, for every vertex, $v$,
the restriction of the above sequence to the link, $\mathrm{Lk}(v)$,
and to the star, $\mathrm{St}(v)$, are shellings.
\end{prop}

\medskip
Since the complex, $\s{K}(P)$, associated with a polytope, $P$,
has a single cell, namely $P$ itself, note that by condition (i) in the
definition of a shelling,
$\s{K}(P)$ is shellable iff the complex, $\s{K}(\partial P)$,
is shellable. We will say simply say that ``$P$ is shellable'' instead
of ``$\s{K}(\partial P)$ is shellable''.

\medskip
We have the following useful property of shellings of polytopes
whose proof is left as an exercise (use induction on the dimension):

\begin{prop}
\label{shellingrev}
Given any polytope, $P$, if $F_1, \ldots, F_s$ is a shelling
of $P$, then the reverse sequence $F_s, \ldots, F_1$ is also
a shelling of $P$.
\end{prop}

\danger
Proposition \ref{shellingrev} generally fails for 
complexes that are not polytopes, see the right $2$-complex
in Figure \ref{figshell1}.

\medskip
We will now present the proof that every polytope is
shellable, using a technique invented by Bruggesser and Mani (1970)
known as {\it line shelling\/}  \cite{Bruggesser}.
This is quite a simple
and natural idea if one is willing to ignore the technical details
involved in actually checking that it works. We begin by explaining
this idea in the $2$-dimensional case, a convex polygon,
since it is particularly simple.

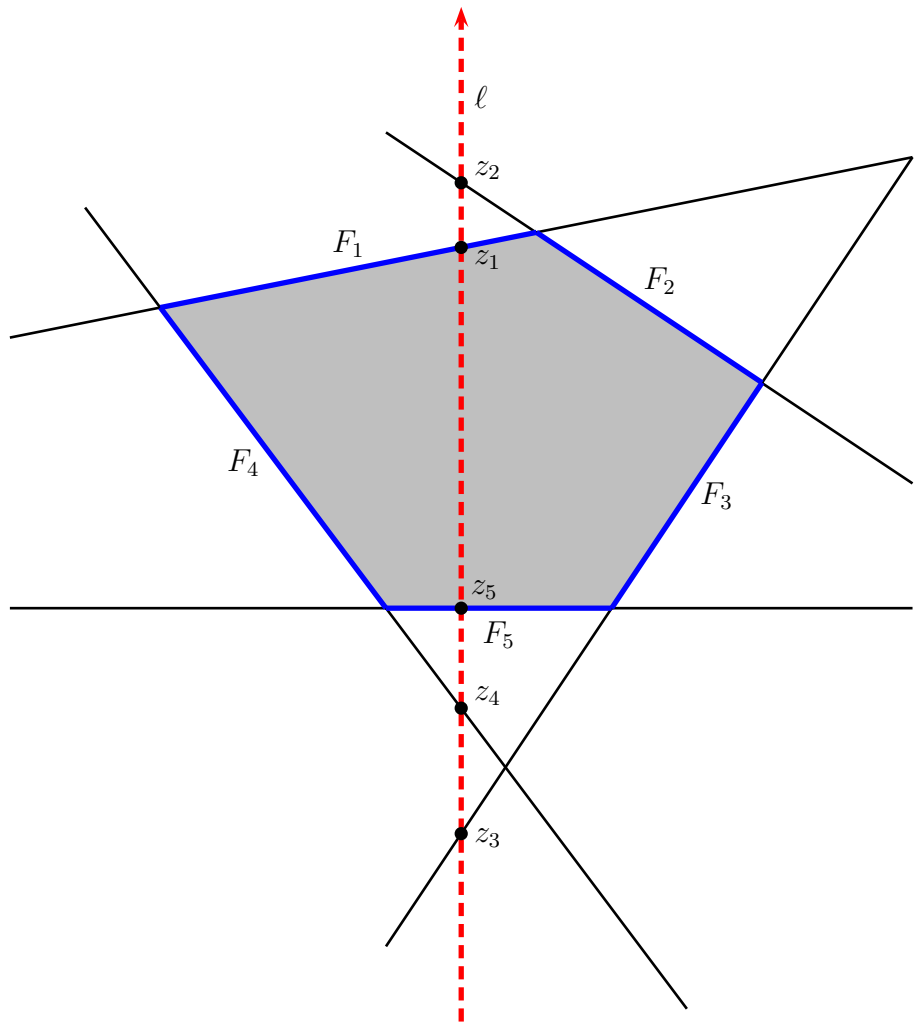
\begin{figure}[H]
  \begin{center}
    \begin{pspicture}(-6.5,-5.5)(7,8)
    \pnode(-6,3.6){v1}
    \pnode(6,6){v2}
    \pnode(-1,6.33){v3}
    \pnode(6,1.66){v4}
    \pnode(6,6){v5}
    \pnode(-1,-4.5){v6}
    \pnode(-6,0){v7}
    \pnode(6,0){v8}
    \pnode(-5,5.33){v9}
    \pnode(3,-5.33){v10}
    \pnode(0,-5.5){v11}
    \pnode(0,8){v12}
    \pnode(-4,4){u1}
    \pnode(1,5){u2}
    \pnode(4,3){u3}
    \pnode(2,0){u4}
    \pnode(-1,0){u5}
    \ncline[linewidth=1pt]{v1}{v2}
    \ncline[linewidth=1pt]{v3}{v4}
    \ncline[linewidth=1pt]{v5}{v6}
    \ncline[linewidth=1pt]{v7}{v8}
    \ncline[linewidth=1pt]{v9}{v10}
    \pspolygon[fillstyle=solid,fillcolor=lightgray,linewidth=1pt]
   (-4,4)(1,5)(4,3)(2,0)(-1,0)
    \ncline[linewidth=2pt,linecolor=red,linestyle=dashed]{->}{v11}{v12}
    \ncline[linewidth=2pt,linecolor=blue]{u1}{u2}
    \ncline[linewidth=2pt,linecolor=blue]{u2}{u3}
    \ncline[linewidth=2pt,linecolor=blue]{u3}{u4}
    \ncline[linewidth=2pt,linecolor=blue]{u4}{u5}
    \ncline[linewidth=2pt,linecolor=blue]{u5}{u1}
     \cnode[fillstyle=solid,fillcolor=black](0,4.8){2.5pt}{z1}
     \cnode[fillstyle=solid,fillcolor=black](0,5.66){2.5pt}{z2}
     \cnode[fillstyle=solid,fillcolor=black](0,-3){2.5pt}{z3}
     \cnode[fillstyle=solid,fillcolor=black](0,-1.33){2.5pt}{z4}
     \cnode[fillstyle=solid,fillcolor=black](0,0){2.5pt}{z5}
    \uput[90](-1.5,4.5){$F_1$}    
    \uput[70](2.5,4){$F_2$}    
    \uput[0](3,1.5){$F_3$}    
    \uput[-90](0.5,0){$F_5$}    
    \uput[190](-2.5,2){$F_4$}    
    \uput[-30](0,4.8){$z_1$}    
    \uput[30](0,5.66){$z_2$} 
    \uput[-10](0,-3){$z_3$}       
    \uput[30](0,-1.33){$z_4$}
    \uput[45](0,0){$z_5$}       
    \uput[0](0,6.8){$\ell$}       
    \end{pspicture}
  \end{center}
  \caption{Shelling a polygon by travelling along a line}
  \label{figshell2}
\end{figure}

\medskip
Consider the $2$-polytope, $P$, shown in Figure \ref{figshell2}
(a polygon) whose faces are labeled $F_1, F_2, F_3, F_4, F_5$.
Pick any line, $\ell$, intersecting the interior of $P$ and 
intersecting the supporting lines
of the facets of $P$ ({\it i.e.\/}, the edges of $P$) in distinct points labeled
$z_1, z_2, z_3, z_4, z_5$ (such a line can always be found, as will be shown
shortly). Orient the line, $\ell$, (say, upward)
and travel on $\ell$ starting from the point of $P$ where
$\ell$ leaves $P$, namely, $z_1$. For a while, only face $F_1$ is visible
but when we reach the intersection, $z_2$, of $\ell$ with the 
supporting line of
$F_2$, the face $F_2$ becomes visible and  $F_1$
becomes invisible as it is now hidden by the supporting line of $F_2$.
So far, we have seen the faces, $F_1$ and $F_2$, {\it in that order\/}.
As we continue traveling along $\ell$, no new face becomes visible
but for a more complicated polygon, other faces, $F_i$,  would become
visible one at a time as we reach  the intersection,
$z_i$, of $\ell$ with the  supporting line of $F_i$
and the order in which these faces become visible corresponds
to the ordering of the $z_i$'s along the line $\ell$.
Then, we  imagine that we travel very fast and when we
reach ``$+\infty$'' in the upward direction on $\ell$, we instantly
come back on $\ell$ from below at ``$-\infty$''. At this
point, we only see the face of $P$
corresponding to the lowest supporting line of faces of $P$, i.e., the line
corresponding to the smallest $z_i$, in our case, $z_3$.
At this stage, the only visible face is $F_3$.
We continue traveling upward on $\ell$ and
we reach $z_3$, the intersection of the supporting line of $F_3$
with $\ell$. At this moment, $F_4$ becomes visible and $F_3$ disappears
as it is now hidden by the supporting line of $F_4$.
Note that $F_5$ is not visible at this stage. Finally, we reach
$z_4$, the intersection of the supporting line of $F_4$ with $\ell$
and at this moment, the last facet, $F_5$, becomes visible
(and $F_4$ becomes invisible, $F_3$ being also invisible).
Our trip stops when we reach $z_5$, the intersection of $F_5$ and $\ell$.
During the second phase of our trip, we saw $F_3, F_4$ and $F_5$
and the entire trip yields the sequence $F_1, F_2, F_3, F_4, F_5$,
which is easily seen to be a shelling of $P$.

\medskip
This is the crux 
of the Bruggesser-Mani method for shelling a polytope: We travel
along a suitably chosen line and record the order in which the
faces become visible during this trip. This is why such shellings
are called {\it line shellings\/}.

\medskip
In order to prove that polytopes are shellable we  need
the notion of points and lines in ``general position''.
Recall from the equivalence of $\s{V}$-polytopes and $\s{H}$-polytopes
that a polytope, $P$,  in $\eucreal^d$ with nonempty interior is cut out by
$t$ irredundant hyperplanes, $H_i$,  and by picking the origin 
in the interior of $P$ the equations of the $H_i$ may be assumed to be
of the form
\[
a_i\cdot z = 1
\]
where $a_i$ and $a_j$ are not proportional for all $i \not= j$,
so that
\[
P = \{z\in \eucreal^d \mid a_i\cdot z \leq 1, \> 1\leq i \leq t\}.
\]

\begin{defin}
\label{genpos}
{\em
Let $P$ be any  polytope in $\eucreal^d$ with nonempty interior 
and assume that $P$ is cut out by the irredudant hyperplanes, $H_i$, of 
equations $a_i\cdot z = 1$, for $i = 1, \ldots, t$.
A point, $c\in \eucreal^d$, is said to be in {\it general position\/}
w.r.t. $P$ is $c$ does not belong to any of the $H_i$, that is, if
$a_i\cdot c \not= 1$ for $i = 1, \ldots, t$.
A line, $\ell$, is  said to be in {\it general position\/} w.r.t. $P$
if $\ell$ is not parallel to any of the $H_i$ and if $\ell$
intersects the $H_i$ in distinct points.
}
\end{defin}

\medskip
The following proposition showing the existence of lines in general position
w.r.t. a polytope illustrates a very useful technique, the
``perturbation method''. The ``trick'' behind this particular
perturbation method is that polynomials (in one variable) have a
finite number of zeros.

\begin{prop}
\label{perturb1}
Let $P$ be any polytope in $\eucreal^d$ with nonempty interior.
For any two points, $x$ and $y$ in $\eucreal^d$, with $x$
outside of $P$; $y$ in the interior of $P$; and $x$ in general position
w.r.t. $P$, for $\lambda\in \reals$ small enough, the line,
$\ell_{\lambda}$, through $x$ and $y_{\lambda}$ with
\[
y_{\lambda} = y + (\lambda, \lambda^2, \ldots, \lambda^d),
\]
intersects $P$ in its interior and is in general position
w.r.t. $P$.
\end{prop}

\proof
Assume that $P$ is defined by $t$ irredundant hyperplanes, $H_i$, where
$H_i$ is given by the equation $a_i\cdot z = 1$ and write
$\Lambda = (\lambda, \lambda^2, \ldots, \lambda^d)$
and $u = y - x$. Then the line $\ell_{\lambda}$ is given by
\[
\ell_{\lambda} = \{x + s(y_{\lambda} - x) \mid s\in \reals\} =
\{x + s(u + \Lambda) \mid s\in \reals\}.
\]
The line, $\ell_{\lambda}$, is not parallel to the hyperplane $H_i$ iff
\[
a_i\cdot (u + \Lambda) \not= 0,\quad i = 1, \ldots, t
\]
and it intersects the $H_i$ in distinct points iff 
there is no $s\in \reals$ such that
\[
a_i\cdot (x + s(u + \Lambda)) = 1
\quad\hbox{and}\quad
a_j\cdot (x + s(u + \Lambda)) = 1
\quad\hbox{for some $i \not= j$}.
\]
Observe that $a_i\cdot (u + \Lambda) = p_i(\lambda)$ is a nonzero polynomial
in $\lambda$ of degree at most $d$. Since a polynomial of degree
$d$ has at most $d$ zeros, if we let  $Z(p_i)$ be the (finite)
set of zeros of $p_i$ we can ensure that
$\ell_{\lambda}$ is not parallel to any of the $H_i$ by picking
$\lambda \notin \bigcup_{i = 1}^t Z(p_i)$ (where
$\bigcup_{i = 1}^t Z(p_i)$ is a finite set).
Now, as $x$ is in general position w.r.t. $P$, we have
$a_i\cdot x \not= 1$, for $i = 1\ldots, t$.
The condition stating that  $\ell_{\lambda}$ intersects the $H_i$ in distinct points
can be written
\[
a_i\cdot x + s a_i\cdot (u + \Lambda) = 1
\quad\hbox{and}\quad
a_j\cdot x + s a_j\cdot (u + \Lambda) = 1
\quad\hbox{for some $i \not= j$},
\]
or
\[
s p_i(\lambda) = \alpha_i
\quad\hbox{and}\quad
s p_j(\lambda) = \alpha_j
\quad\hbox{for some $i \not= j$},
\]
where $\alpha_i = 1 - a_i\cdot x$ and $\alpha_j = 1 - a_j\cdot x$.
As $x$ is in general position w.r.t. $P$, we have $\alpha_i, \alpha_j \not= 0$
and as the $H_i$ are irredundant, the polynomials
$p_i(\lambda) = a_i\cdot (u + \Lambda)$ and
$p_j(\lambda) = a_j\cdot (u + \Lambda)$ are not proportional.
Now, if $\lambda\notin Z(p_i) \cup Z(p_j)$, 
in order for the system 
\begin{eqnarray*}
s p_i(\lambda) & = & \alpha_i \\
s p_j(\lambda) & = & \alpha_j 
\end{eqnarray*}
to have a solution in $s$ we must have
\[
q_{i j}(\lambda) = \alpha_ip_j(\lambda) - \alpha_jp_i(\lambda) = 0,
\]
where $q_{i j}(\lambda)$ is not the zero polynomial since
$p_i(\lambda)$ and $p_j(\lambda)$ are not proportional and
$\alpha_i, \alpha_j \not= 0$.
If we pick $\lambda\notin Z(q_{i j})$, then
$q_{i j}(\lambda) \not= 0$. Therefore, if we pick
\[
\lambda \notin \bigcup_{i = 1}^t Z(p_i) \cup 
\bigcup_{i \not= j}^t Z(q_{i j}),
\]
the line $\ell_{\lambda}$ is in general position w.r.t. $P$.
Finally, we can pick $\lambda$ small enough so that
$y_{\lambda} = y + \Lambda$ is close enough to $y$ so that it is
in the interior of $P$.
$\bigsquare$

\medskip
It should be noted that the perturbation method
involving $\Lambda = (\lambda, \lambda^2, \ldots, \lambda^d)$ 
is quite flexible. For example, by adapting the proof of
Proposition \ref{perturb1} we can prove that for any two 
distinct facets, $F_i$ and $F_j$ of $P$, there is a line
in general position w.r.t. $P$ intersecting $F_i$ and $F_j$.
Start with $x$ outside $P$ and very close to
$F_i$ and $y$ in the interior of $P$ and very close to $F_j$.

\medskip
Finally, before proving the existence of line shellings for polytopes,
we need more terminology. Given any point, $x$, strictly
outside a polytope, $P$,
we say that a facet, $F$, of $P$ is {\it visible from $x$\/}
iff for every $y\in F$ the line through $x$ and $y$ intersects
$F$ only in $y$ (equivalently, $x$ and the interior of $P$ are strictly
separared by the supporting hyperplane of $F$).
We now prove the following fundamental theorem due to
Bruggesser and Mani \cite{Bruggesser} (1970):

\begin{thm} (Existence of Line Shellings for Polytopes)
\label{shellingexist}
Let $P$ be any polytope in $\eucreal^d$ of dimension $d$.
For every point, $x$, outside $P$ and in general position w.r.t. $P$,
there is a shelling of $P$ in which the facets of $P$ that are
visible from $x$ come first.
\end{thm}

\begin{figure}
  \begin{center}
    \begin{pspicture}(0,-0.5)(6.5,8.5)
\pspolygon[fillstyle=solid,fillcolor=white,linewidth=1.5pt]
(0,2.5)(1.5,0.75)(3,0)(2.5,1.7)(0.7,3)
\pspolygon[fillstyle=solid,fillcolor=white,linewidth=1.5pt]
(3,0)(4.9,0.6)(5.8,2)(4.3,3)(2.5,1.7)
\pspolygon[fillstyle=solid,fillcolor=gray,linewidth=1.5pt]
(4.3,3)(5.8,2)(6.5,4)(5,5.7)(3.7,5)
\pspolygon[fillstyle=solid,fillcolor=gray,linewidth=1.5pt]
(2.5,1.7)(4.3,3)(3.7,5)(1.5,5)(0.7,3)
\pspolygon[fillstyle=solid,fillcolor=gray,linewidth=1.5pt]
(0,2.5)(0.7,3)(1.5,5)(1,5.5)(0.4,4.2)
\pspolygon[fillstyle=solid,fillcolor=lightgray,linewidth=1.5pt]
(1.5,5)(3.7,5)(5,5.7)(3.2,6)(1,5.5)
    \cnode[fillstyle=solid,fillcolor=black](0,2.5){2pt}{v1}
    \cnode[fillstyle=solid,fillcolor=black](1.5,0.75){2pt}{v2}
    \cnode[fillstyle=solid,fillcolor=black](3,0){2pt}{v3}
    \cnode[fillstyle=solid,fillcolor=black](4.9,0.6){2pt}{v4}
    \cnode[fillstyle=solid,fillcolor=black](5.8,2){2pt}{v5}
    \cnode[fillstyle=solid,fillcolor=black](6.5,4){2pt}{v6}    
    \cnode[fillstyle=solid,fillcolor=black](5,5.7){2pt}{v7}    
    \cnode[fillstyle=solid,fillcolor=black](3.2,6){2pt}{v8}
    \cnode[fillstyle=solid,fillcolor=black](1,5.5){2pt}{v9}
    \cnode[fillstyle=solid,fillcolor=black](0.4,4.2){2pt}{v10}
    \cnode[fillstyle=solid,fillcolor=black](1.5,5){2pt}{v11}
    \cnode[fillstyle=solid,fillcolor=black](0.7,3){2pt}{v12}
    \cnode[fillstyle=solid,fillcolor=black](2.5,1.7){2pt}{v13}
    \cnode[fillstyle=solid,fillcolor=black](4.3,3){2pt}{v14}
    \cnode[fillstyle=solid,fillcolor=black](3.7,5){2pt}{v15}
    \pnode(1.75,-0.5){w1}
    \pnode(3.3,8.5){w2}
    \ncline[linewidth=2pt,linecolor=red,linestyle=dashed]{->}{w1}{w2}
    \cnode[fillstyle=solid,fillcolor=green](2.75,5.3){2.5pt}{z1}
    \cnode[fillstyle=solid,fillcolor=green](3.07,7.1){2.5pt}{z2}
    \cnode[fillstyle=solid,fillcolor=green](3.18,7.8){2.5pt}{z3}
    \uput[0](1.35,0.1){$\ell$}    
    \uput[190](2.75,5.4){$z_1$}    
    \uput[200](3.05,7.3){$z_2$}   
    \uput[200](3.12,8){$z_3$}   
    \uput[0](2.85,5.5){$F_1$}     
    \uput[0](2.65,3.6){$F_2$}    
    \uput[0](4.75,3.9){$F_3$}    
    \uput[0](0.25,3.85){$F_4$}    
    \end{pspicture}
  \end{center}
  \caption{Shelling a polytope by travelling along a line, $\ell$}
  \label{shellfig3}
\end{figure}
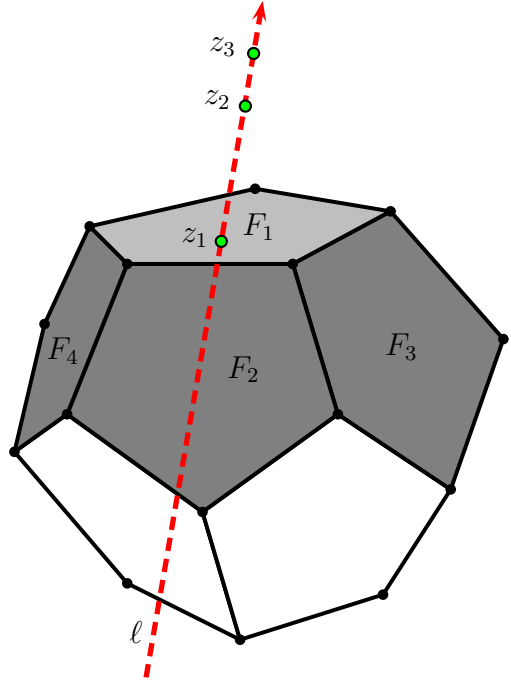

\proof
By Proposition \ref{perturb1}, we can find a line, $\ell$, through
$x$ such that $\ell$ is in general position w.r.t. $P$ and $\ell$ intersects 
the interior of $P$. Pick one of the two faces in which $\ell$
intersects $P$, say $F_1$, let $z_1 = \ell \cap F_1$, 
and orient $\ell$ from the inside of $P$ to $z_1$. As $\ell$
intersects the supporting hyperplanes of the facets of $P$
in distinct points, we get a linearly  ordered list of these
intersection points along $\ell$, 
\[
z_1 , \> z_2 , \> \cdots , \> z_m , \> z_{m+1} , \cdots , \> z_{s},
\]
where $z_{m+1}$ is the smallest element, $z_m$ is the largest
element and where $z_1$ and $z_s$ belong to the faces of $P$
where $\ell$ intersects $P$. Then, as in the example illustrated
by Figure \ref{figshell2}, by travelling ``upward'' along the line $\ell$
starting from $z_1$ we get a total ordering of the facets of
$P$, 
\[
F_1, F_2, \ldots, F_m , F_{m+1}, \ldots , F_s
\]
where $F_i$ is the facet whose supporting hyperplane cuts $\ell$
in $z_i$. 

\medskip
We claim that the above sequence is a shelling of $P$.
This is proved by induction on $d$.
For $d = 1$, $P$ consists a line segment and the theorem clearly holds.

\medskip
Consider the intersection $\partial F_j \cap (F_1\cup \cdots \cup F_{j-1})$.
We need to show that this is an initial segment of a shelling of 
$\partial F_j$. If $j \leq m$, i.e., if $F_j$ become visible
before we reach $\infty$, then the above intersection is exactly
the set of facets of $F_j$ that are visible from 
$z_j = \ell \cap \mathrm{aff}(F_j)$. Therefore, by induction on the dimension,
these facets are shellable and they form an initial segment of a shelling
of the whole boundary $\partial F_j$.

\medskip
If $j \geq m + 1$, that is, after ``passing through $\infty$''
and reentering from $-\infty$, the intersection
$\partial F_j \cap (F_1\cup \cdots \cup F_{j-1})$ is the set of
non-visible facets. By reversing the orientation of the line, $\ell$,
we see that the facets of this intersection are shellable
and we get the reversed ordering of the facets.

\medskip
Finally, when we reach the point $x$ starting from $z_1$,
the facets visible from $x$ form an initial segment of the shelling,
as claimed.
$\bigsquare$

\remark
The trip along the line $\ell$ is often described as a {\it rocket flight\/}
starting from the surface of $P$ viewed as a little planet
(for instance, this is the description given by Ziegler \cite{Ziegler97}
(Chapter 8)). Observe that if we reverse the direction of $\ell$,
we obtain the reversal of the original line shelling. Thus,
the reversal of a line shelling is not only a shelling
but a line shelling as well.

\medskip
We can easily prove the following corollary:

\begin{cor}
\label{shellingexist2}
Given any polytope, $P$, the following facts hold:
\begin{enumerate}
\item[(1)]
For any two facets $F$ and $F'$, 
there is a shelling of $P$ in which $F$ comes first and $F'$
comes last.
\item[(2)]
For any vertex, $v$, of $P$, there is a shelling of $P$ in which
the facets containing $v$ form an initial segment of the shelling.
\end{enumerate}
\end{cor}

\proof
For (1), we use a line in general position and intersecting
$F$ and $F'$ in their interior.
For (2), we pick a point, $x$, beyond $v$ and pick a line in general
position through $x$ intersecting the interior of $P$.
Pick the origin, $O$, in the  interior of $P$.
A point, $x$, is {\it beyond $v$\/} iff $x$ and $O$ lies on different sides
of every hyperplane, $H_i$, supporting a facet of $P$ 
containing $x$ but on the same side of $H_i$ for
every hyperplane, $H_i$, supporting a facet of $P$ {\bf not}
containing $x$. Such a point can be found on a line through
$O$ and $v$, as the reader should check.
$\bigsquare$

\remark
A {\it plane triangulation\/}, $K$, is a pure
two-dimensional complex in the plane such that $|K|$
is homeomorphic to a closed
disk. Edelsbrunner proves that every plane triangulation
has a shelling and from this, that 
$\chi(K) = 1$, where $\chi(K) = f_0 - f_1 + f_2$
is the  Euler-Poincar\'e characteristic of $K$, 
where $f_0$ is the number of vertices, $f_1$ is the number
of edges and $f_2$ is the number of triangles in $K$
(see  Edelsbrunner \cite{Edelsbrunner},
Chapter 3). This result is an immediate consequence
of Corollary \ref{shellingexist2} if one knows about
the stereographic projection map, which will be discussed
in the next Chapter.

\medskip
We now have all the tools needed to prove the famous
Euler-Poincar\'e Formula for Polytopes.

\section[The Euler-Poincar\'e Formula for Polytopes]
{The Euler-Poincar\'e Formula for Polytopes}
\label{sec9}
We begin by defining a very important topological concept,
the Euler-Poincar\'e characteristic of a complex.

\begin{defin}
\label{PoincareEulerdef}
{\em
Let $K$ be a $d$-dimensional complex. For every $i$,
with $0 \leq i \leq d$, we let $f_i$ denote the number
of $i$-faces of $K$ and we let
\[
\mathbf{f}(K)  = (f_0 , \cdots, f_d) \in \natnums^{d+1}
\]
be the {\it $f$-vector\/} associated with $K$ (if necessary
we write $f_i(K)$ instead of $f_i$). 
The {\it Euler-Poincar\'e characteristic\/}, $\chi(K)$, of $K$ is defined by
\[
\chi(K) = f_0 - f_1 + f_2 + \cdots + (-1)^d f_d = \sum_{i = 0}^{d} (-1)^i f_i.
\]
Given any $d$-dimensional polytope, $P$, the {\it $f$-vector\/}
associated with $P$ is the $f$-vector associated with $\s{K}(P)$,
that is, 
\[
\mathbf{f}(P)  = (f_0 , \cdots, f_d) \in \natnums^{d+1},
\]
where $f_i$, is the number of $i$-faces of $P$ ($= $ the number of
$i$-faces of $\s{K}(P)$ and thus, $f_d = 1$), 
and the  {\it Euler-Poincar\'e characteristic\/}, $\chi(P)$, 
of $P$ is defined by
\[
\chi(P) = f_0 - f_1 + f_2 + \cdots + (-1)^d f_d = \sum_{i = 0}^{d} (-1)^i f_i.
\]
Moreover, the {\it $f$-vector\/}
associated with the boundary, $\partial P$, of $P$
is the $f$-vector associated with $\s{K}(\partial P)$,
that is, 
\[
\mathbf{f}(\partial P)  = (f_0 , \cdots, f_{d-1}) \in \natnums^{d}
\]
where $f_i$, is the number of $i$-faces of $\partial P$
(with $0 \leq i \leq d - 1$),
and the  {\it Euler-Poincar\'e characteristic\/}, $\chi(\partial P)$, 
of $\partial P$ is defined by
\[
\chi(\partial P) = f_0 - f_1 + f_2 + \cdots + (-1)^{d-1} f_{d-1} = 
\sum_{i = 0}^{d-1} (-1)^i f_i.
\]
}
\end{defin}

\medskip
Observe that $\chi(P) = \chi(\partial P) + (-1)^d$, since $f_d = 1$. 

\remark
It is convenient to set $f_{-1} = 1$. Then, 
some authors, including Ziegler \cite{Ziegler97} (Chapter 8),
define the {\it reduced Euler-Poincar\'e characteristic\/}, $\chi'(K)$,
of a complex (or a polytope), $K$, as
\[
\chi'(K) = -f_{-1} +  f_0 - f_1 + f_2 + \cdots + (-1)^d f_d = 
\sum_{i = -1}^{d} (-1)^i f_i = -1 + \chi(K),
\]
{\it i.e.\/}, they incorporate $f_{-1} = 1$ into the formula. 

\medskip
A crucial observation for proving the  Euler-Poincar\'e formula
is that the Euler-Poincar\'e characteristic
is additive, which means that if $K_1$ and $K_2$ are any two complexes
such that $K_1\cup K_2$ is also a complex, which implies that
$K_1\cap K_2$ is also a complex (because we must have $F_1\cap F_2
\in K_1\cap K_2$ for every face $F_1$ of $K_1$ and every face
$F_2$ of $K_2$), then
\[
\chi(K_1\cup K_2) = \chi(K_1) + \chi(K_2) - \chi(K_1\cap K_2).
\]
This follows immediately because for any two sets $A$ and $B$
\[
|A \cup B| = |A| + |B| - |A\cap B|.
\]

\medskip
To prove our next theorem we will use complete induction
on $\natnums\times \natnums$ ordered by the lexicographic ordering.
Recall that the lexicographic ordering on $\natnums\times \natnums$ is defined
as follows:
\[
(m, n) < (m', n') \quad\hbox{iff}\quad
\left\{\begin{array}{l}
m = m'\quad{and}\quad n < n' \\
\hbox{or} \\
m < m'.
\end{array}
\right.
\]

\begin{thm} (Euler-Poincar\'e Formula)
\label{PoincareEulerthm1}
For every polytope, $P$, we have
\[
\chi(P) = \sum_{i = 0}^d (-1)^if_i = 1 \qquad
(\hbox{$d \geq 0$}),
\]
and so,
\[
\chi(\partial P) = \sum_{i = 0}^{d-1} (-1)^if_i = 1 - (-1)^d
\qquad
(\hbox{$d \geq 1$}).
\]
\end{thm}

\proof
We prove the following statement:
For every $d$-dimensional polytope, $P$, if $d = 0$ then 
\[
\chi(P) = 1,
\]
else if  $d \geq 1$ then for every
shelling $F_1,\ldots, F_{f_{d-1}}$, of $P$, for every $j$, with
$1\leq j \leq f_{d-1}$, we have
\[
\chi(F_1\cup\cdots\cup F_j) = \cases{
1 & if $1 \leq j < f_{d-1}$ \cr
1 - (-1)^d & if $j = f_{d-1}$. \cr
}
\]
We proceed by complete induction on $(d, j)\geq (0, 1)$.
For $d = 0$ and $j = 1$, the polytope $P$ consists of a single point
and so, $\chi(P) = f_0 = 1$, as claimed.

\medskip
For the induction step, assume that $d\geq 1$.
For $1 = j < f_{d-1}$, since $F_1$ is a polytope of dimension $d - 1$,
by the induction hypothesis, $\chi(F_1) = 1$, as desired.

\medskip
For $1 < j < f_{d - 1}$, we have
\[
\chi(F_1\cup \cdots F_{j-1}\cup F_j) =
\chi\left(\bigcup_{i = 1}^{j-1} F_i\right) + \chi(F_j) -  
\chi\left(\left(\bigcup_{i = 1}^{j-1} F_i \right)\cap F_j\right). 
\]
Since $(d, j - 1) < (d, j)$, by the induction hypothesis,
\[
\chi\left(\bigcup_{i = 1}^{j-1} F_i\right)  = 1
\]
and since $\mathrm{dim}(F_j) = d - 1$, again by the induction
hypothesis,
\[
\chi(F_j)  = 0.
\]
Now, as $F_1, \ldots, F_{f_{d-1}}$ is a shelling and $j < f_{d-1}$, we have
\[
\left(\bigcup_{i = 1}^{j-1} F_i \right)\cap F_j = 
G_1  \cup \cdots \cup G_{r},
\]
for some shelling $G_1,  \ldots, G_r, \ldots, G_t$ of 
$\s{K}(\partial F_j)$, with $r < t = f_{d - 2}(\partial F_j)$.
The fact that $r < f_{d - 2}(\partial F_j)$, {\it i.e.\/}, that
$G_1  \cup \cdots \cup G_{r}$ is not the whole boundary of $F_j$
is a property of line shellings and also follows from
Proposition \ref{shellingrev}.
As $\mathrm{dim}(\partial F_j) = d - 2$,
and $r < f_{d - 2}(\partial F_j)$, by the induction hypothesis, we have
\[
\chi\left(\left(\bigcup_{i = 1}^{j-1} F_i \right)\cap F_j\right) = 
\chi(G_1  \cup \cdots \cup G_{r}) = 1.
\]
Consequently,
\[
\chi(F_1\cup \cdots F_{j-1}\cup F_j) = 1 + 1 - 1 = 1,
\]
as claimed (when $j < f_{d-1}$).

\medskip
If $j = f_{d-1}$, then we have a complete shelling of $\partial F_{f_{d-1}}$,
that is,
\[
\left(\bigcup_{i = 1}^{f_{d-1}-1} F_i \right)\cap F_{f_{d-1}} =
G_1\cup \cdots \cup G_{f_{d - 2}(F_{f_{d-1}})} = 
\partial F_{f_{d-1}}.
\]
As  $\mathrm{dim}(\partial F_j) = d - 2$, by the induction hypothesis,
\[
\chi(\partial F_{f_{d-1}}) = \chi(G_1\cup \cdots \cup G_{f_{d - 2}(F_{f_{d-1}})}) 
= 1 - (-1)^{d- 1}
\]
and it follows that
\[
\chi(F_1\cup \cdots \cup F_{f_{d-1}}) = 1 + 1 - (1 - (-1)^{d-1}) =
1 +  (-1)^{d-1} = 1 -  (-1)^{d},
\]
establishing the induction hypothesis in this last case.
But then,
\[
\chi(\partial P) = \chi(F_1\cup \cdots \cup F_{f_{d-1}}) =  1 -  (-1)^{d}
\]
and 
\[
\chi(P) = \chi(\partial P) +  (-1)^{d} = 1,
\]
proving our theorem.
$\bigsquare$

\remark
Other combinatorial proofs of the
Euler-Poincar\'e formula are given in Gr\"unbaum 
\cite{Grunbaum} (Chapter 8), Boissonnat and Yvinec 
\cite{Boissonnat} (Chapter 7)  and
Ewald  \cite{Ewald} (Chapter 3). Coxeter gives a proof 
very close to Poincar\'e's own proof using
notions of homology theory \cite{Coxeterpoly} (Chapter IX).
We feel that the proof based on shellings is the most direct
and one of the most elegant.
Incidently, the above proof of the 
Euler-Poincar\'e formula is very close to Schl\"afli proof
from 1852 but Schl\"afli did not have shellings
at his disposal so his ``proof'' had a gap.
The Bruggesser-Mani proof that polytopes are shellable
fills this gap!

\section[Dehn-Sommerville Equations for Simplicial Polytopes]
{Dehn-Sommerville Equations for Simplicial \\ Polytopes and
$h$-Vectors}
\label{sec10}
If a $d$-polytope, $P$, 
has the property that its faces
are all simplices,  then it is called a {\it simplicial polytope\/}.
It is easily shown that a polytope is simplicial
iff its facets are simplices, in which case, every
facet has $d$ vertices.
The polar dual of a simplicial polytope is called a {\it simple
polytope\/}. We see immediately that every vertex of a simple
polytope belongs to $d$ facets. 

\medskip
For simplicial (and simple) polytopes
it turns out that other remarkable equations besides the
Euler-Poincar\'e formula hold among the number of $i$-faces.
These equations were discovered by Dehn for $d = 4, 5$ (1905)
and by Sommerville in the general case (1927). Although it is possible 
(and not difficult) to prove
the Dehn-Sommerville equations by ``double counting'',
as in Gr\"unbaum \cite{Grunbaum} (Chapter 9) or Boissonnat and
Yvinec (Chapter 7, but beware, these
are the dual formulae for simple polytopes), 
it turns out that instead of using the
$f$-vector associated with a polytope it is preferable to use
what's known as the $h$-vector because for simplicial polytopes
the $h$-numbers have a natural interpretation in terms of shellings.
Furthermore, the statement of the
Dehn-Sommerville equations in terms of $h$-vectors is 
transparent:
\[
h_{i} = h_{d - i},
\]
and the proof is very simple in terms of shellings.

\medskip
In the rest of this section, we restrict our attention to
simplicial complexes.
In order to motivate $h$-vectors, we begin by examining
more closely the structure of the new faces that are created
during a shelling when the cell $F_j$ is added to the
partial shelling $F_1, \ldots, F_{j-1}$.

\medskip
If $K$ is a simplicial polytope and $V$ is the set of vertices
of $K$, then every $i$-face of $K$ can be identified with an
$(i + 1)$-subset of $V$ (that is, a subset of $V$ of cardinality $i + 1$).

\begin{defin}
\label{restict}
{\em
For any shelling, $F_1, \ldots, F_s$, of a simplicial complex, $K$,
of dimension $d$,
for every $j$, with $1\leq j \leq s$, the {\it restriction\/}, $R_j$, 
of the facet, $F_j$, is the set of ``obligatory'' vertices 
\[
R_j = \{v\in F_j \mid F_j - \{v\} \subseteq F_i,\> \hbox{for some $i$ with}\>
1\leq i < j\}.
\] 
}
\end{defin}

\medskip
The crucial property of the $R_j$ is that the new
faces, $G$, added at step $j$ (when $F_j$ is added to the shelling)
are precisely the faces in the set
\[
I_j = \{G\subseteq V \mid R_j \subseteq G \subseteq F_j\}.
\]
The proof of the above fact is left as an exercise to the reader.

\medskip
But then, we obtain a partition, $\{I_1, \ldots, I_s\}$,
of the set of faces of the simplicial complex (other that $K$ itself).
Note that the empty face is allowed. Now, if we define
\[
h_i = |\{j \mid |R_j| = i, \> 1\leq j \leq s\}|,
\]
for $i = 0, \ldots, d$, then it turns out that
we can recover the $f_k$ in terms of the
$h_i$ as follows:
\[
f_{k-1} = \sum_{j = 1}^s \binom{d - |R_j|}{k - |R_j|} =
\sum_{i = 0}^k h_i \binom{d - i}{k - i}, 
\]
with $1 \leq k \leq d$. 

\medskip
But more is true: The above equations are invertible and the
$h_k$ can be expressed in terms of the $f_i$ as follows:
\[
h_k = \sum_{i = 0}^k (-1)^{k - i} \binom{d- i}{d - k} f_{i-1}, 
\]
with $0 \leq k \leq d$ (remember, $f_{-1} = 1$).

\medskip
Let us explain all this in more detail.
Consider the  example of a connected graph
(a simplicial $1$-dimensional complex) from Ziegler
\cite{Ziegler97} (Section 8.3) shown in Figure \ref{graphfig1}:

\begin{figure}
  \begin{center}
    \begin{pspicture}(0,0)(4.5,2.5)
     \pnode(0,0){u2}
     \pnode(0,2){u1}
     \pnode(2,0){u3}
     \pnode(2,2){u4}
     \pnode(4,2){u5}
     \pnode(4,0){u6}
    \ncline[linewidth=1pt]{u1}{u2}
    \ncline[linewidth=1pt]{u1}{u3}
    \ncline[linewidth=1pt]{u3}{u4}
    \ncline[linewidth=1pt]{u3}{u5}
    \ncline[linewidth=1pt]{u3}{u6}
    \ncline[linewidth=1pt]{u4}{u5}
    \ncline[linewidth=1pt]{u5}{u6}
    \uput[90](0,2){$1$}    
    \uput[-90](0,0){$2$}    
    \uput[-90](2,0){$3$}
    \uput[90](2,2){$4$}
    \uput[90](4,2){$5$}
    \uput[-90](4,0){$6$}                
    \end{pspicture}
  \end{center}
  \caption{A connected $1$-dimensional complex, $G$}
  \label{graphfig1}
\end{figure}
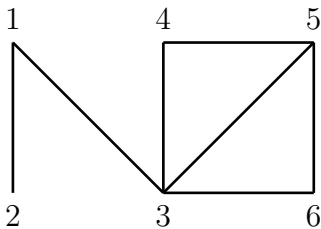

\medskip
A shelling order of its $7$ edges is given by the sequence
\[
12,\> 13,\> 34,\> 35,\> 45, \> 36, \> 56.
\]
The partial order of the faces of $G$ together with the
blocks of the partition $\{I_1, \ldots, I_7\}$ associated with 
the seven edges of $G$ are shown in Figure \ref{graphfig2},
with the blocks $I_j$ shown in boldface:

\begin{figure}
  \begin{center}
    \begin{pspicture}(0,0)(9.5,4)
     \pnode(4.5,0){u0}
     \pnode(0,1.5){u1}
     \pnode(0,3){u12}
     \pnode(1.5,1.5){u2}
     \pnode(1.5,3){u13}
     \pnode(3,1.5){u3}
     \pnode(3,3){u34}
     \pnode(4.5,1.5){u4}
     \pnode(4.5,3){u35}
     \pnode(6,1.5){u5}
     \pnode(6,3){u45}
     \pnode(7.5,1.5){u6}
     \pnode(7.5,3){u36}
     \pnode(9,3){u56}
    \ncline[linewidth=2pt]{u0}{u1}
    \ncline[linewidth=2pt]{u0}{u2}
    \ncline[linewidth=1pt]{u0}{u3}
    \ncline[linewidth=1pt]{u0}{u4} 
    \ncline[linewidth=1pt]{u0}{u5}
    \ncline[linewidth=1pt]{u0}{u6}
    \ncline[linewidth=2pt]{u1}{u12}
    \ncline[linewidth=1pt]{u1}{u13}
    \ncline[linewidth=2pt]{u2}{u12}
    \ncline[linewidth=2pt]{u3}{u13}
    \ncline[linewidth=1pt]{u3}{u34} 
    \ncline[linewidth=1pt]{u3}{u35} 
    \ncline[linewidth=1pt]{u3}{u36} 
    \ncline[linewidth=2pt]{u4}{u34}
    \ncline[linewidth=1pt]{u4}{u45}
    \ncline[linewidth=2pt]{u5}{u35}
    \ncline[linewidth=1pt]{u5}{u45}
    \ncline[linewidth=1pt]{u5}{u56}
    \ncline[linewidth=2pt]{u6}{u36}
    \ncline[linewidth=1pt]{u6}{u56}
    \cnode[fillstyle=solid,fillcolor=black](6,3){2.5pt}{v45}
    \cnode[fillstyle=solid,fillcolor=black](9,3){2.5pt}{v56}
    \uput[-90](4.5,0){$\emptyset$}    
    \uput[180](0,1.5){$1$}    
    \uput[200](1.5,1.5){$2$}    
    \uput[200](3,1.5){$3$}
    \uput[-20](4.5,1.5){$4$}
    \uput[-20](6,1.5){$5$}
    \uput[-10](7.5,1.5){$6$}
    \uput[90](0,3){$12$}
    \uput[90](1.5,3){$13$}
    \uput[90](3,3){$34$}
    \uput[90](4.5,3){$35$}                
    \uput[90](6,3){$45$}
    \uput[90](7.5,3){$36$}
    \uput[90](9,3){$56$}                                                
    \end{pspicture}
  \end{center}
  \caption{the partition associated with a shelling of  $G$}
  \label{graphfig2}
\end{figure}
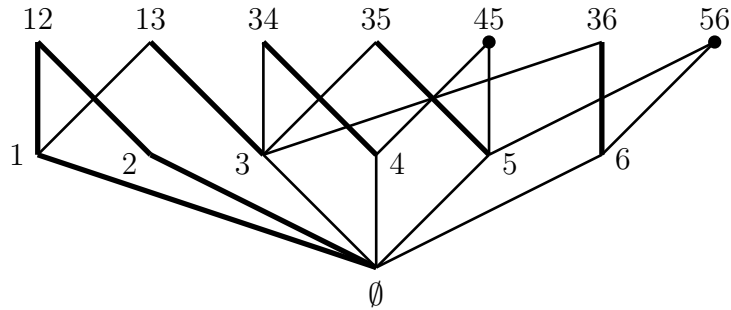

\medskip
The ``minimal'' new faces (corresponding to the $R_j$'s)
added at every stage of the shelling are
\[
\emptyset, \> 3, \> 4, \> 5, \> 45, \> 6, \> 56.
\]
Again, if $h_i$ is the number of blocks, $I_j$,  such that the
corresponding restriction set, $R_j$, has size $i$, that is,
\[
h_i = |\{j \mid |R_j| = i, \> 1\leq j \leq s\}|,
\]
for $i = 0, \ldots, d$, where the simplicial polytope, $K$,  has dimension
$d - 1$, we define the {\it $h$-vector\/} associated with $K$ as
\[
\mathbf{h}(K) = (h_0, \ldots, h_d).
\]
Then, in the above example,  as $R_1 = \{\emptyset\}$,
$R_2 = \{3\}$, $R_3 = \{4\}$, $R_4 = \{5\}$, $R_5 = \{4,5\}$,
$R_6 = \{6\}$ and $R_7 = \{5,6\}$, we get
$h_0 = 1$, $h_1 = 4$ and $h_2 = 2$, that is,
\[
\mathbf{h}(G) = (1, 4, 2).
\]

\medskip
Now, let us show that if $K$ is a shellable simplicial
complex, then the $f$-vector can be recovered from the $h$-vector.
Indeed, if $|R_j| = i$, then each $(k-1)$-face in the block $I_j$
must use all $i$ nodes in $R_j$, so that there are only
$d - i$ nodes available and, among those, $k - i$ must be
chosen. Therefore,
\[
f_{k - 1} = \sum_{j = 1}^s \binom{d - |R_j|}{k - |R_j|}
\]
and, by definition of $h_i$, we get
\begin{equation}
f_{k-1} = \sum_{i = 0}^k h_i\binom{d - i}{k - i}
= h_k + \binom{d - k + 1}{1} h_{k-1} + \cdots + \binom{d - 1}{k - 1} h_1
+ \binom{d}{k} h_0,
\tag{$*$}
\end{equation}
where $1 \leq k \leq d$.
Moreover, the formulae are invertible, that is,
the $h_i$ can be expressed in terms of the $f_k$.
For this, form the two polynomials
\[
f(x) = \sum_{i = 0}^d f_{i - 1} x^{d - i} =
f_{d-1} + f_{d-2} x + \cdots + f_0 x^{d-1} + f_{-1} x^d
\]
with $f_{-1} = 1$ and
\[
h(x) = \sum_{i = 0}^d h_i x^{d - i} =
h_d + h_{d-1} x + \cdots + h_1 x^{d-1} + h_0 x^d.
\]
Then, it is easy to see that
\[
f(x) = \sum_{i = 0}^d h_i(x + 1)^{d - i} = h(x + 1).
\]
Consequently, $h(x) = f(x-1)$ and by comparing the coefficients of
$x^{d - k}$ on both sides of the above equation, we get
\[
h_k = \sum_{i = 0}^k (-1)^{k - i} \binom{d- i}{d - k} f_{i-1}.
\]
In particular, $h_0 = 1$, $h_1 = f_0 - d$, and
\[
h_d = f_{d-1} - f_{d-2} + f_{d- 3} + \cdots + (-1)^{d-1} f_0 + (-1)^d.
\]

\medskip
It is also easy to check that
\[
h_0 + h_1 + \cdots + h_d = f_{d-1}.
\]
Now, we just showed that if $K$ is shellable, then
its $f$-vector and its $h$-vector are related as above.
But even if $K$ is not shellable, the above suggests  defining
the $h$-vector from the $f$-vector as above. Thus, we make the
definition:

\begin{defin}
\label{hvecdef}
{\em 
For any $(d-1)$-dimensional simplicial complex, $K$, the
{\it $h$-vector\/} associated with $K$ is the vector
\[
\mathbf{h}(K) = (h_0, \ldots, h_d) \in \integs^{d+1},
\]
given by
\[
h_k = \sum_{i = 0}^k (-1)^{k - i} \binom{d- i}{d - k} f_{i-1}.
\]
}
\end{defin}

\medskip
Note that if $K$ is shellable, then the interpretation of 
$h_i$ as the number of cells, $F_j$,  such that the
corresponding restriction set, $R_j$, has size $i$ shows
that $h_i \geq 0$. However, for an arbitrary simplicial complex,
some of the $h_i$ can be strictly negative. Such an example
is given in Ziegler \cite{Ziegler97} (Section 8.3).

\medskip
We summarize below most of what we just showed:

\begin{prop}
\label{hvecp1}
Let $K$ be a $(d-1)$-dimensional pure simplicial complex.
If $K$ is shellable, then its $h$-vector is nonnegative
and $h_i$ counts the number of cells in a shelling
whose restriction set has size $i$.  Moreover, the $h_i$
do not depend on the particular shelling of $K$. 
\end{prop}

\medskip
There is a way of computing the $h$-vector of a
pure simplicial complex from its $f$-vector
reminiscent of the Pascal triangle (except that
negative entries can turn up). Again, the reader is referred
to  Ziegler \cite{Ziegler97} (Section 8.3).

\medskip
We are now ready to prove the Dehn-Sommerville equations.
For $d = 3$, these are easily obtained by double counting.
Indeed, for a simplicial polytope, every edge belongs to two facets
and every facet has three edges. It follows that
\[
2f_1 = 3f_2.
\]
Together with Euler's formula
\[
f_0 - f_1 + f_2 = 2,
\]
we see that 
\[
f_1 = 3f_0 - 6
\quad\hbox{and}\quad
f_2 = 2f_0 - 4,
\]
namely, that the number of vertices of a simplicial $3$-polytope
determines its number of edges and faces, these being linear
functions of the number of vertices. For arbitrary dimension $d$,
we have

\begin{thm} (Dehn-Sommerville Equations)
\label{DehnSommerville}
If $K$ is any simplicial $d$-polytope, then the components
of the $h$-vector satisfy
\[
h_k = h_{d - k}\qquad k = 0, 1, \ldots, d.
\]
Equivalently
\[
f_{k-1} = \sum_{i = k}^d (-1)^{d - i}\binom{i}{k} f_{i-1}
\qquad k = 0, \ldots, d.
\]
Furthermore, the equation $h_0 = h_d$ is equivalent to the
Euler-Poincar\'e formula.
\end{thm}

\proof 
We present a short and elegant proof due to McMullen.
Recall from Proposition \ref{shellingrev} that the reversal, 
$F_s, \ldots, F_1$, of a shelling, $F_1, \ldots, F_s$,  
of a polytope is also a shelling. From this, we see that for every $F_j$,
the restriction set of $F_j$ in the reversed shelling is equal to
$R_j - F_j$, the complement of the restriction set of $F_j$
in the original shelling. Therefore, if $|R_j| = k$, then
$F_j$ contributes ``$1$'' to $h_k$ in the original shelling iff
it contributes ``$1$'' to $h_{d - k}$ in the reversed shelling
(where $|R_j - F_j| = d - k$). It follows that the value of $h_k$
computed in the original shelling is the same as the value of $h_{d - k}$
computed in the reversed shelling. However, by Proposition \ref{hvecp1},
the $h$-vector is independent of the shelling and hence,
$h_k = h_{d - k}$.

\medskip
Define the polynomials $F(x)$ and $H(x)$ by
\[
F(x) = \sum_{i = 0}^d f_{i-1}x^i; \quad
H(x) = (1 - x)^dF\left(\frac{x}{1-x}\right).
\]
Note that $H(x) = \sum_{i = 0}^d f_{i - 1}x^i(1 - x)^{d - i}$ and an easy 
computation shows that the coefficient of $x^k$ is equal to
\[
 \sum_{i = 0}^k (-1)^{k - i} \binom{d- i}{d - k} f_{i-1} = h_k.
\]
Now, the equations $h_{k} = h_{d - k}$ are equivalent to
\[
H(x)  = x^dH(x^{-1}),
\]
that is,
\[
F(x - 1) = (-1)^dF(-x).
\]
As 
\[
F(x-1) = \sum_{i = 0}^d f_{i-1} (x - 1)^i = 
\sum_{i = 0}^d f_{i-1} \sum_{j = 0}^i \binom{i}{i - j} x^{i - j}(-1)^j,
\]
we see that the coefficient of $x^k$ in $F(x-1)$ (obtained when $i - j = k$,
that is, $j = i - k$) is
\[
\sum_{i = 0}^d   (-1)^{i-k} \binom{i}{k} f_{i-1}
= \sum_{i = k}^d  (-1)^{i-k} \binom{i}{k} f_{i-1}.
\]
On the other hand, the coefficient of $x^k$ in $(-1)^dF(-x)$ is
$(-1)^{d + k}f_{k-1}$. By equating the coefficients of $x^k$, we get
\[
(-1)^{d + k}f_{k-1} = \sum_{i = k}^d (-1)^{i-k} \binom{i}{k}  f_{i-1},
\]
which, by multiplying both sides by $(-1)^{d + k}$, is equivalent to
\[
f_{k-1} = \sum_{i = k}^d (-1)^{d+i} \binom{i}{k}  f_{i-1}
= \sum_{i = k}^d (-1)^{d-i} \binom{i}{k}  f_{i-1},
\]
as claimed.
Finally, as we already know that
\[
h_d = f_{d-1} - f_{d-2} + f_{d- 3} + \cdots + (-1)^{d-1} f_0 + (-1)^d
\]
and $h_0 = 1$, by multiplying both sides of the equation
$h_d = h_0 = 1$ by $(-1)^{d -1}$ and moving $(-1)^d(-1)^{d-1} = -1$
to the right hand side, we get the Euler-Poincar\'e formula.
$\bigsquare$

\medskip
Clearly, the Dehn-Sommerville equations, $h_k = h_{d - k}$,
are linearly independent for \\
$0 \leq k < \lfloor\frac{d+1}{2}\rfloor$.
For example, for $d = 3$, we have the two independent equations
\[
h_0 = h_3, \> h_1 = h_2,
\]
and for $d = 4$, we also have two independent equations
\[
h_0 = h_4,\> h_1 = h_3,
\]
since $h_2 = h_2$ is trivial. When $d = 3$, we know that
$h_1 = h_2$ is equivalent to $2f_1 = 3f_2$ and when $d = 4$,
if one unravels $h_1 = h_3$ in terms of
the $f_i$' one finds
\[
2f_2 = 4f_3,
\]
that is $f_2 = 2f_3$. More generally, it is easy to check that
\[
2f_{d-2} = df_{d-1}
\]
for all $d$.
For $d = 5$, we find three independent equations
\[
h_0 = h_5,\> h_1 = h_4,\> h_2 = h_3,
\]
and so on.

\medskip
It can be shown that for general $d$-polytopes, the 
Euler-Poincar\'e formula is the only equation satisfied
by all $h$-vectors and for simplicial $d$-polytopes,
the $\lfloor\frac{d+1}{2}\rfloor$ Dehn-Sommerville equations,
$h_k = h_{d - k}$,  are the only equations
satisfied by all $h$-vectors (see Gr\"unbaum \cite{Grunbaum},
Chapter 9). 

\remark
Readers familiar with homology and cohomology may suspect that
the Dehn-Sommerville equations are a consequence
of a type of Poincar\'e duality. Stanley proved that this is
indeed the case. It turns out that the $h_i$ are the
dimensions of cohomology groups of a
certain {\it toric variety\/} associated with the polytope.
For more on this topic, see Stanley \cite{Stanley} 
(Chapters II and III) and
Fulton \cite{Fultontoric} (Section 5.6).

\medskip
As we saw for $3$-dimensional simplicial polytopes, 
the number of vertices, $n = f_0$,  determines the number of edges
and the number of faces, and these are linear in $f_0$.
For $d \geq 4$, this is no longer true and the number
of facets is no longer linear in $n$ but in fact quadratic.
It is then natural
to ask which $d$-polytopes with a prescribed number of vertices
have the maximum number of $k$-faces. This question
which remained an open problem for some twenty years
was eventually settled by McMullen in 1970 \cite{McMullen71}.
We will present this result (without proof) in the
next section.

\section[The Upper Bound Theorem]
{The Upper Bound Theorem and Cyclic Polytopes}
\label{sec11}
Given a $d$-polytope with $n$ vertices, 
what is an upper bound on the number of its $i$-faces?
This  question is not only important from a theoretical point of view
but also from a computational point of view because
of its implications for algorithms in combinatorial
optimization and in  computational geometry.  

\medskip
The answer to the above problem is that there is a class of polytopes
called {\it cyclic polytopes\/} such that the cyclic $d$-polytope,
$C_d(n)$,  has the maximum number of $i$-faces among all
$d$-polytopes with  $n$ vertices.
This result stated by Motzkin in 1957 became known as
the {\it upper bound conjecture\/} until it 
was  proved by
McMullen in 1970, using shellings \cite{McMullen71}
(just after Bruggesser and Mani's
proof that polytopes are shellable).  It is now known as the
{\it upper bound theorem\/}.
Another proof of the upper bound theorem was given later by Alon
and Kalai \cite{AlonKalai} (1985).
A version of this proof can also be found in Ewald \cite{Ewald} (Chapter 3).

\medskip  
McMullen's proof is not really very difficult but
it is still quite involved so we will only state some
propositions needed for its proof. We urge the reader to 
read Ziegler's account of this beautiful proof \cite{Ziegler97}
(Chapter 8). We begin with cyclic polytopes.

\medskip
First, consider the cases $d = 2$ and $d = 3$.
When $d = 2$, our polytope is a polygon in which case
$n = f_0 = f_1$.  Thus, this case is trivial.

\medskip
For $d = 3$, we claim that $2f_1 \geq 3f_2$. Indeed, every edge belongs to 
exactly two faces so
if we add up the number of sides for all faces, we get $2f_1$.
Since every face has at least three sides, we get
$2f_1 \geq 3f_2$. Then, using Euler's relation, it is easy to show that
\[
f_1 \leq 6n - 3\quad f_2 \leq 2n - 4
\]
and we know that equality is achieved for simplicial polytopes.

\medskip
Let us now consider the general case. 
The rational curve, $\mapdef{c}{\reals}{\reals^d}$, given parametrically by
\[
c(t) = (t, t^2, \ldots, t^d)
\]
is at the heart of the story. This curve if often called
the {\it moment curve\/} or {\it rational normal curve\/}
of degree $d$.
For $d = 3$, it is known as the {\it twisted cubic\/}.
Here is the definition of the cyclic polytope, $C_d(n)$.

\begin{defin}
\label{cyclicpolydf}
{\em
For any sequence,  $t_1 < \ldots < t_n$, of  distinct real number,
$t_i\in \reals$, with $n > d$, the convex hull,
\[
C_d(n) = \mathrm{conv}(c(t_1), \ldots, c(t_n))
\]
of the $n$ points, $c(t_1), \ldots, c(t_n)$, on the moment curve of
degree $d$ is called a {\it cyclic polytope\/}.
}
\end{defin}

\medskip
The first interesting fact about the cyclic  polytope is that
it is simplicial.

\begin{prop}
\label{cyclicpoly1}
Every $d + 1$ of the points $c(t_1), \ldots, c(t_n)$ are affinely
independent. Consequently, $C_d(n)$ is a simplicial polytope
and the $c(t_i)$ are vertices.
\end{prop}

\proof
We may assume that $n = d + 1$. Say $c(t_1), \ldots, c(t_n)$ 
belong to a hyperplane, $H$, given by
\[
\alpha_1 x_1 + \cdots + \alpha_d x_d = \beta.
\]
(Of course, not all the $\alpha_i$ are zero.)
Then, we have the polynomial, $H(t)$, given by
\[
H(t) = -\beta + \alpha_1 t + \alpha_2 t^2 + \cdots + \alpha_d t^d, 
\]
of degree at most $d$ and as each $c(t_i)$ belong to $H$,
we see that each $c(t_i)$ is a zero of $H(t)$. However, 
there are $d + 1$ distinct $c(t_i)$, so $H(t)$ would
have $d + 1$ distinct roots. As $H(t)$ has degree at most $d$,
it must be the zero polynomial, a contradiction.
Returing to the original $n > d + 1$, we just proved
every $d + 1$ of the points $c(t_1), \ldots, c(t_n)$ are affinely
independent. Then, every proper face of $C_d(n)$ has at most
$d$ independent vertices, which means that it is a simplex.
$\bigsquare$

\medskip
The following proposition already shows that the cyclic
polytope, $C_d(n)$,  has $\binom{n}{k}$ $(k-1)$-faces if
$1 \leq k \leq \lfloor\frac{d}{2}\rfloor$. 

\begin{prop}
\label{cyclicpoly2}
For any $k$ with $2 \leq 2k \leq d$, every subset of $k$ vertices
of $C_d(n)$ is a $(k - 1)$-face of $C_d(n)$. Hence
\[
f_k(C_d(n)) = \binom{n}{k+1}\qquad \hbox{if}\quad
0 \leq k < \left\lfloor\frac{d}{2}\right\rfloor.
\]
\end{prop}

\proof
Consider any sequence $t_{i_1} < t_{i_2} < \cdots < t_{i_k}$.
We will prove that there is a hyperplane separating
$F = \mathrm{conv}(\{c(t_{i_1}),  \ldots,  c(t_{i_k})\})$ and $C_d(n)$.
Consider the polynomial
\[
p(t) = \prod_{j = 1}^k (t - t_{i_{j}})^2
\]
and write 
\[
p(t) = a_0 + a_1 t + \cdots + a_{2k} t^{2k}.
\]
Consider the vector
\[
a = (a_1, a_2, \ldots, a_{2k}, 0, \ldots, 0)\in \reals^d 
\]
and the hyperplane, $H$, given by
\[
H = \{x\in \reals^d \mid x\cdot a = -a_0\}.
\]
Then, for each $j$ with $1 \leq j \leq k$, we have
\[
c(t_{i_j})\cdot a = a_1t_{i_j} + \cdots + a_{2k}t_{i_j}^{2k} =
p(t_{i_j}) - a_0 = -a_0,
\]
and so, $c(t_{i_j})\in H$.
On the other hand, for any other point, $c(t_i)$, distinct
from any of the $c(t_{i_j})$, we have
\[
c(t_i) \cdot a = -a_0 + p(t_i) = -a_0 + \prod_{j = 1}^k (t_i - t_{i_{j}})^2
> -a_0,
\]
proving that $c(t_i) \in H_+$. But then, $H$ is a supporting hyperplane
of $F$ for $C_d(n)$ and $F$ is a $(k - 1)$-face.
$\bigsquare$

\medskip
Observe that Proposition \ref{cyclicpoly2} shows that
any subset of $\lfloor\frac{d}{2}\rfloor$ vertices of $C_d(n)$
forms a face of $C_d(n)$. When a $d$-polytope
has this property it is called a {\it neighborly polytope\/}.
Therefore, cyclic polytopes are neighborly. 
Proposition  \ref{cyclicpoly2} also shows a phenomenon
that only manifests itself in dimension at least $4$:
For $d \geq 4$, the polytope $C_d(n)$ has $n$ pairwise
adjacent vertices. For $n >> d$, this is counter-intuitive.

\medskip
Finally, the combinatorial structure of cyclic polytopes
is completely determined as follows:

\begin{prop} (Gale evenness condition, Gale (1963)).
\label{cyclicpoly3}
Let $n$ and $d$ be integers with $2 \leq d < n$. For any sequence
$t_1 < t_2 < \cdots < t_n$, consider the cyclic polytope
\[
C_d(n) = \mathrm{conv}(c(t_1), \ldots, c(t_n)).
\]
A subset, $S\subseteq \{t_1, \ldots, t_n\}$ with $|S| = d$
determines a facet
of $C_d(n)$ iff for all $i < j$ not in $S$, then the number of $k\in S$
between $i$ and $j$ is even:
\[
|\{k \in S \mid i < k < j\}| \equiv 0\> (\mathrm{mod}\> 2)
\quad\hbox{for}\quad i, j \notin S
\]
\end{prop}

\proof
Write $S = \{s_1, \ldots, s_d\} \subseteq \{t_1, \ldots, t_n\}$.
Consider the polyomial
\[
q(t) = \prod_{i = 1}^d (t - s_i) = \sum_{j = 0}^d b_j t^j,
\]
let $b = (b_1, \ldots, b_d)$, and let $H$ be
the hyperplane given by
\[
H = \{x\in \reals^d \mid x\cdot b = -b_0\}.  
\]
Then, for each $i$, with $1\leq i \leq d$, we have
\[
c(s_i)\cdot b = \sum_{j = 1}^d b_j s_i^j = q(s_i) - b_0 = -b_0,
\]
so that $c(s_i) \in H$. For all other $t\not= s_i$,
\[
q(t) = c(t)\cdot b + b_0 \not= 0,
\]
that is, $c(t)\notin H$.
Therefore, $F = \{c(s_1), \ldots, c(s_d)\}$ is a facet of $C_d(n)$
iff $\{c(t_1), \ldots, c(t_n)\} - F$ lies in one of the two open half-spaces
determined by $H$.  This is equivalent to $q(t)$ changing its sign an even
number of times while, increasing $t$, we pass through 
the vertices in $F$. Therefore, the proposition is proved.
$\bigsquare$

\medskip
In particular, Proposition \ref{cyclicpoly3} shows that the 
combinatorial structure
of $C_d(n)$ does not depend on the specific choice of the sequence
$t_1 < \cdots < t_n$. This justifies our notation $C_d(n)$.

\medskip
Here is the celebrated upper bound theorem first proved
by McMullen \cite{McMullen71}.

\begin{thm} (Upper Bound Theorem, McMullen (1970))
\label{upperbth1}
Let $P$ be any $d$-polytope with $n$ vertices. Then, for every $k$,
with $1 \leq k \leq d$, the polytope $P$ has at most as many
$(k-1)$-faces as the cyclic polytope, $C_d(n)$, that is
\[
f_{k - 1}(P) \leq f_{k - 1}(C_d(n)).
\]
Moreover, equality for some $k$ with  $\lfloor\frac{d}{2}\rfloor
\leq k \leq d$ implies that $P$ is neighborly.
\end{thm}

\medskip
The first step in the proof of Theorem \ref{upperbth1}
is to prove that among all $d$-polytopes with a given number, $n$, 
of vertices, the maximum number of $i$-faces is achieved
by simplicial $d$-polytopes. 

\begin{prop}
\label{upperbp1}
Given any $d$-polytope, $P$, with $n$-vertices, it is possible to
form a simplicial polytope, $P'$, by perturbing the vertices
of $P$ such that $P'$ also has $n$ vertices and
\[
f_{k-1}(P) \leq f_{k-1}(P')
\qquad\hbox{for}\quad 1 \leq k \leq d.
\]
Furthermore, equality for $k >  \lfloor\frac{d}{2}\rfloor$
can occur only if $P$ is simplicial.
\end{prop}

\medskip\noindent
{\it Sketch of proof\/}. First, we apply Proposition \ref{triangul1}
to triangulate the facets of $P$ without adding any vertices.
Then, we can perturb the vertices to obtain a simplicial polytope, $P'$,
with at least as many facets (and thus, faces) as $P$.
$\bigsquare$

\medskip
Proposition \ref{upperbp1} allows us to restict our attention
to simplicial polytopes. Now, it is obvious that
\[
f_{k - 1} \leq \binom{n}{k}
\]
for any polytope $P$ (simplicial or not) and we also know
that equality holds if $k \leq  \lfloor\frac{d}{2}\rfloor$
for neighborly polytopes such as the cyclic polytopes. 
For  $k >  \lfloor\frac{d}{2}\rfloor$, it turns out that equality can only 
be achieved for simplices. 

\medskip
However, for a {\it simplicial\/} polytope, the Dehn-Sommerville
equations $h_k = h_{d - k}$ together with the equations $(*)$ giving
$f_k$ in terms of the $h_i$'s show that
$f_0, f_1, \ldots, f_{\lfloor\frac{d}{2}\rfloor}$ already determine
the whole $f$-vector. Thus, it is possible to express the
$f_{k-1}$ in terms of $h_0, h_1, \ldots, h_{\lfloor\frac{d}{2}\rfloor}$
for $k \geq \lfloor\frac{d}{2}\rfloor$. It turns out that we get
\[
f_{k-1} = \sideset{}{^*}\sum_{i = 0}^{\lfloor\frac{d}{2}\rfloor}
\left(\binom{d - i}{k - i} + \binom{i}{k - d + i}
\right) h_i,
\]
where the meaning of the superscript $*$ is that when $d$ is even
we only take half of the last term for $i = \frac{d}{2}$
and when $d$ is odd we take the whole last term for $i = \frac{d - 1}{2}$
(for details, see Ziegler \cite{Ziegler97}, Chapter 8).
As a consequence if we can show that the neighborly polytopes
maximize not only  $f_{k-1}$ but also  $h_{k-1}$ when 
$k \leq \lfloor\frac{d}{2}\rfloor$, the upper bound theorem will 
be proved. Indeed, McMullen proved the following theorem
which is ``more than enough'' to yield the desired result
(\cite{McMullen71}):

\begin{thm} (McMullen (1970))
\label{upperbp2}
For every simplicial $d$-polytope with $f_0 = n$ vertices, we have
\[
h_k(P) \leq \binom{n - d - 1  + k}{k}
\qquad\hbox{for}\quad 0 \leq k \leq d.
\]
Furthermore, equality holds for all $l$ and all $k$
with $0 \leq k \leq l$ iff $l \leq \lfloor\frac{d}{2}\rfloor$
and $P$ is $l$-neighborly.
(a polytope is $l$-neighborly iff any subset of $l$ or less
vertices determine a face of $P$.)
\end{thm}

\medskip
The proof of Theorem \ref{upperbp2} is too involved
to be given here, which is unfortunate, since 
it is really beautiful. It makes a clever
use of shellings and a careful analysis of the $h$-numbers of 
links of vertices. Again, the reader is referred to 
Ziegler \cite{Ziegler97}, Chapter 8.

\medskip
Since cyclic $d$-polytopes are neighborly (which means that
they are $\lfloor\frac{d}{2}\rfloor$-neighborly),
Theorem \ref{upperbth1} follows from Proposition \ref{upperbp1},
and Theorem \ref{upperbp2}.

\begin{cor}
\label{upperbp3}
For every simplicial neighborly $d$-polytope with $n$ vertices,
we have
\[
f_{k - 1} =  \sideset{}{^*}\sum_{i = 0}^{\lfloor\frac{d}{2}\rfloor}
\left(\binom{d - i}{k - i} + \binom{i}{k - d + i}
\right) \binom{n - d - 1  + i}{i}
\qquad\hbox{for}\quad 1 \leq k \leq d.
\]
This gives the maximum number of $(k - 1)$-faces
for any $d$-polytope with $n$-vertices, for all $k$
with $1 \leq k \leq d$.  In particular, the number
of facets of the cyclic polytope, $C_d(n)$, is
\[
f_{d-1} = \sideset{}{^*}\sum_{i = 0}^{\lfloor\frac{d}{2}\rfloor} 
2\binom{n - d - 1 + i}{i}
\]
and, more explicitly,
\[
f_{d-1} = \binom{n - \lfloor\frac{d+1}{2}\rfloor}{n - d} +
\binom{n - \lfloor\frac{d+2}{2}\rfloor}{n - d}.
\]
\end{cor}

\medskip
Corollary \ref{upperbp3} implies that the number of facets
of any $d$-polytope is $O(n^{\lfloor\frac{d}{2}\rfloor})$.
An unfortunate consequence of this upper bound is that the
complexity of any convex hull algorithms for $n$ points
in $\eucreal^d$ is  $O(n^{\lfloor\frac{d}{2}\rfloor})$.

\medskip
The  $O(n^{\lfloor\frac{d}{2}\rfloor})$ upper bound can be obtained
more directly using a pretty argument using shellings due to 
R. Seidel \cite{Seidel95}.
Consider any shelling of any 
simplicial $d$-polytope, $P$.
For every facet, $F_j$, of a shelling either the restriction
set $R_j$ or its complement $F_j - R_j$ has at most 
$\lfloor\frac{d}{2}\rfloor$ elements. So, either in the shelling or
in the reversed shelling, the restriction set of $F_j$ has at most
$\lfloor\frac{d}{2}\rfloor$ elements. Moreover, the restriction
sets are all distinct, by construction. Thus, the number of facets
is at most twice the number of $k$-faces of $P$ with
$k \leq  \lfloor\frac{d}{2}\rfloor$. It follows that
\[
f_{d-1} \leq 2\sum_{i = 0}^{\lfloor\frac{d}{2}\rfloor} \binom{n}{i}
\]
and this rough estimate  yields a  $O(n^{\lfloor\frac{d}{2}\rfloor})$
bound.

\remark
There is also a {\it lower bound theorem\/} 
due to Barnette (1971, 1973) which gives
a lower bound on the $f$-vectors all $d$-polytopes
with $n$ vertices. In this case, there is an analog
of the cyclic polytopes called {\it stacked polytopes\/}.
These polytopes, $P_d(n)$, are simplicial polytopes obtained from
a simplex by building shallow pyramids over the facets of 
the simplex. Then, it turns out that  if $d \geq 2$, then
\[
f_{k} \geq \cases{
\binom{d}{k}n - \binom{d + 1}{k + 1}k &
if $0 \leq k \leq d - 2$ \cr
(d - 1)n - (d + 1)(d - 2) & if
 $k =   d - 1$. \cr
}
\]

\medskip
There has been a lot of progress on the combinatorics
of $f$-vectors and $h$-vectors since 1971, especially by R. Stanley,
G. Kalai and L. Billera and K. Lee, among others.
We recommend two excellent surveys: 
\begin{enumerate}
\item 
Bayer and Lee
\cite{BayerLee} summarizes progress in this
area up to 1993.
\item
Billera and Bj\"orner \cite{BilleraBjorner} is a more advanced
survey which reports on results up to 1997.
\end{enumerate}
In fact, many of the chapters in Goodman and O'Rourke 
\cite{DiscGeomHand} should be of interest to the reader.

\medskip
Generalizations of the Upper Bound Theorem using 
sophisticated techniques (face rings) due to Stanley can be found
in Stanley \cite{Stanley} 
(Chapters II) and connections with toric varieties
can be found in Stanley \cite{Stanley} (Chapters III) and
Fulton \cite{Fultontoric}.

\def\dInt{{\rm Int}\>}
\def\dBd{\partial}
\chapter[Dirichlet--Voronoi Diagrams]
{Dirichlet--Voronoi Diagrams and Delaunay Triangulations}
\label{chap22}
\section{Dirichlet--Voronoi Diagrams}
\label{sec221}
In this chapter we present very briefly the concepts of a Voronoi diagram
and of a Delaunay triangulation. 
\index{Voronoi diagram}\index{Delaunay triangulation}%
\index{Dirichlet--Voronoi diagram}%
These are important tools in computational
geometry, and Delaunay triangulations are important in problems
where it is necessary to fit $3$D data using surface splines.
It is usually useful to compute a good mesh for the projection
\index{mesh}%
of this set of data points onto the $xy$-plane, and 
a Delaunay triangulation is a good candidate.
Our presentation will be rather sketchy. We are primarily interested
in defining these concepts and stating their most important
properties. For a comprehensive exposition of
Voronoi diagrams,  Delaunay triangulations, and more topics in 
computational geometry,
our readers may consult O'Rourke \cite{ORourke}, 
\index{computational geometry}
Preparata and Shamos \cite{Preparata}, Boissonnat and Yvinec \cite{Boissonnat}, 
de Berg, Van Kreveld, Overmars, and Schwarzkopf \cite{Berg},
or Risler \cite{Risler92}.
The survey by Graham and Yao \cite{Graham} contains a very gentle and lucid
introduction to computational geometry.
Some practical applications
of Voronoi diagrams and Delaunay triangulations
are briefly discussed in Section \ref{sec225}.

\medskip
Let $\affs$ be a Euclidean space of finite dimension, 
that is, an affine space $\affs$
whose underlying vector space $\affvec{\affs}$ is equipped with an inner
product (and has finite dimension).
For concreteness, one may safely assume that
$\affs = \eucreal^m$, although what follows applies to
any Euclidean space of finite dimension.
\index{Euclidean space}%
Given a set  $P=\{p_1,\ldots,p_\ndeg\}$ of $\ndeg$ points in $\affs$,
it is often useful to find a partition of the space $\affs$ into regions
each containing a single point of $P$ and having some
nice properties. It is also often useful to find triangulations
of the convex hull of $P$  having some 
nice properties. We shall see that this can be done and that the two problems
are closely related. In order to solve the first problem, we
need to introduce bisector lines and bisector planes.
\index{bisector!line}\index{bisector!plane}%

\begin{figure}
  \begin{center}
    \begin{pspicture}(-1,-1.5)(5,5)
    \psline[linewidth=1pt]{->}(-1,0)(5,0)
    \psline[linewidth=1pt]{->}(0,-1)(0,5)
    \psline[linewidth=1pt](-1,5)(5,-1)
    \psline[linewidth=1pt](1,1)(3,3)
    \psdots[dotstyle=*,dotscale=1.5](1,1)
    \psdots[dotstyle=*,dotscale=1.5](3,3)
    \uput[45](1,3){$L$}
    \uput[135](1,1){$a$}
    \uput[135](3,3){$b$}
    \end{pspicture}
  \end{center}
  \caption{The bisector line $L$ of $a$ and $b$}
\label{bissectfig}
\end{figure}
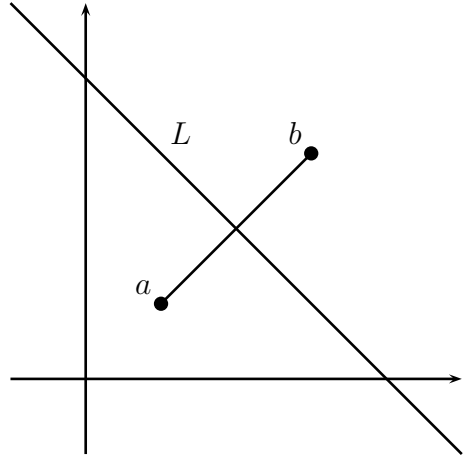

\medskip
For simplicity, let us first assume that $\affs$ is a plane i.e.,
has dimension $2$.
Given any two distinct points $a, b\in \affs$, the line orthogonal to the line segment
$(a, b)$ and passing through the midpoint of this segment is the locus
of all points having equal distance to $a$ and $b$. It is called the
{\it bisector line of $a$ and $b$\/}. The bisector line
of two points is illustrated in Figure \ref{bissectfig}.

If $h = \frac{1}{2}\, a + \frac{1}{2}\, b$
is the midpoint of the line segment $(a, b)$, letting $m$ be an arbitrary
point on the bisector line, the equation of this line can be found
by writing that $\libvecbo{h}{m}$ is orthogonal to $\libvecbo{a}{b}$.
In any orthogonal frame, letting $m = (x, y)$, $a = (a_1, a_2)$, $b = (b_1, b_2)$,
the equation of this line is
\[(b_1 - a_1)(x - (a_1 + b_1)/2) + (b_2 - a_2)(y - (a_2 + b_2)/2) = 0,\]
which can also be written as
\[(b_1 - a_1)x  + (b_2 - a_2)y  = 
(b_1^2 + b_2^2)/2 - (a_1^2+  a_2^2)/2.\]
The closed half-plane $H(a, b)$ containing $a$ and with boundary the bisector line
\indsym{H(a,b)}{closed half-plane}%
is the locus of all points such that
\[(b_1 - a_1)x  + (b_2 - a_2)y  \leq 
(b_1^2 + b_2^2)/2 - (a_1^2+  a_2^2)/2,\]
and the  closed half-plane $H(b, a)$ containing $b$ and with boundary the bisector line
is the locus of all points such that
\[(b_1 - a_1)x  + (b_2 - a_2)y  \geq
(b_1^2 + b_2^2)/2 - (a_1^2+  a_2^2)/2.\]
The closed half-plane $H(a, b)$ is the set of all points whose distance to $a$
is less that or equal to the distance to $b$, and vice versa for $H(b, a)$.
Thus, points in the closed half-plane $H(a, b)$ are closer to $a$ than they
are to $b$.

\medskip
We now consider a problem called the
{\it post office problem\/} by Graham and Yao \cite{Graham}.
\index{post office problem}%
Given any set  $P=\{p_1,\ldots,p_\ndeg\}$ of $\ndeg$ points in the plane
(considered as {\it post offices\/} or {\it sites\/}), for any arbitrary point
$x$, find out which post office is closest to $x$.
Since $x$ can be arbitrary, it seems desirable to precompute the
sets $V(p_i)$ consisting of all points 
that are closer to  $p_i$ than to any other point $p_j\not= p_i$.
Indeed, if the sets $V(p_i)$ are known,
the answer is any post office $p_i$ such that
$x\in V(p_i)$. Thus, it remains to compute the sets $V(p_i)$.
For this, if $x$ is closer to   $p_i$ than to any other point $p_j\not= p_i$,
then $x$ is on the same side as $p_i$ with respect to the bisector line
of $p_i$ and $p_j$ for every $j\not= i$, and thus
\[V(p_i) = \bigcap_{j\not= i} H(p_i, p_j).\]

\medskip
If $\affs$ has dimension $3$, the locus of
all points having equal distance to $a$ and $b$ is a plane.
It is called the {\it bisector plane of $a$ and $b$\/}.
\index{bisector!plane}%
The equation of this plane is also found 
by writing that $\libvecbo{h}{m}$ is orthogonal to $\libvecbo{a}{b}$.
The equation of this plane is
\[
\eqaligneno{
(b_1 - a_1)(&x - (a_1 + b_1)/2) + (b_2 - a_2)(y - (a_2 + b_2)/2) \cr
&+ (b_3 - a_3)(z - (a_3 + b_3)/2) = 0,\cr
}
\]
which can also be written as
\[(b_1 - a_1)x  + (b_2 - a_2)y +  (b_3 - a_3)z =
(b_1^2 + b_2^2 + b_3^2)/2 - (a_1^2+  a_2^2 + a_3^2)/2.\]
The closed half-space $H(a, b)$ containing $a$ and with boundary the bisector plane
is the locus of all points such that
\[(b_1 - a_1)x  + (b_2 - a_2)y +  (b_3 - a_3)z \leq
(b_1^2 + b_2^2 + b_3^2)/2 - (a_1^2+  a_2^2 + a_3^2)/2,\]
and the  closed half-space $H(b, a)$ containing $b$ and with boundary the bisector plane
is the locus of all points such that
\[(b_1 - a_1)x  + (b_2 - a_2)y +  (b_3 - a_3)z \geq
(b_1^2 + b_2^2 + b_3^2)/2 - (a_1^2+  a_2^2 + a_3^2)/2.\]
The closed half-space $H(a, b)$ is the set of all points whose distance to $a$
is less that or equal to the distance to $b$, and vice versa for $H(b, a)$.
Again, points in the closed half-space $H(a, b)$ are closer to $a$ than they
are to $b$.

\medskip
Given any set  $P=\{p_1,\ldots,p_\ndeg\}$ of $\ndeg$ points in $\affs$ (of dimension
$\mdeg= 2, 3$), it is often useful to find for every point $p_i$
the region consisting of all points
that are closer to  $p_i$ than to any other point $p_j\not= p_i$, that is,
the set
\[V(p_i) = \{x\in \affs\ |\ d(x, p_i) \leq d(x, p_j),\> \hbox{for all $j\not= i$}\},\]
where $d(x,y) = (\libvecbo{x}{y}\cdot\libvecbo{x}{y})^{1/2}$,
the Euclidean distance  associated with the inner product
$\cdot$ on $\affs$. From the definition of the bisector line (or plane), 
it is immediate that
\[V(p_i) = \bigcap_{j\not= i} H(p_i, p_j).\]

\medskip
Families of sets of the form $V(p_i)$ were investigated by Dirichlet 
\cite{Dirichlet} (1850) and
Voronoi \cite{Voronoi} (1908). 
\index{Dirichlet}\index{Voronoi}%
Voronoi diagrams also arise in crystallography 
(Gilbert \cite{Gilbert}). Other applications, including
facility location and path planning,
are discussed in O'Rourke \cite{ORourke}.
For simplicity, we also denote the set $V(p_i)$ by $V_i$, and
we introduce the following definition.

\begin{defin}
\label{Voronoi}
{\em
Let $\affs$ be a Euclidean space of dimension $\mdeg\geq 1$. Given any
set $P=\{p_1$, $\ldots$, $p_\ndeg\}$ of $\ndeg$ points in $\affs$, the
{\it Dirichlet--Voronoi diagram\/} 
\index{Dirichlet--Voronoi diagram!definition}%
$\s{V}\mathit{or}(P)$ of $P=\{p_1,\ldots,p_\ndeg\}$ is the family of subsets
of $\affs$ 
consisting of the sets $V_i =  \bigcap_{j\not= i} H(p_i, p_j)$ and of all of their
intersections.
}
\end{defin}

\medskip
Dirichlet--Voronoi diagrams are also  called {\it Voronoi diagrams\/},
\index{Voronoi diagram}%
{\it Voronoi tessellations\/}, or 
{\it Thiessen polygons\/}. 
\index{Thiessen polygons}%
Following common usage, we
will use the terminology {\it Voronoi diagram\/}.
As intersections of convex sets (closed half-planes or closed half-spaces),
the {\it Voronoi regions\/} $V(p_i)$ are convex sets. 
\index{Voronoi region}%
\indsym{V(p_i)}{Voronoi region}%
\index{convex!set}%
In dimension two,
the boundaries of these regions are convex polygons, and in dimension
three, the boundaries are convex polyhedra.

\medskip
Whether a region $V(p_i)$ is bounded or not depends on the location of $p_i$.
If $p_i$ belongs to the boundary of the convex hull of the set $P$,
then $V(p_i)$ is unbounded, and otherwise bounded.
\index{convex!hull}%
In dimension two, the convex hull is a convex polygon, and in dimension
three, the convex hull is a convex polyhedron.
\index{convex!polygon}%
\index{convex!polyhedron}%
As we will see later,
there is an intimate relationship between convex hulls and
Voronoi diagrams.

\medskip
Generally, if $\affs$ is a Euclidean space of dimension $\mdeg$,
given any two distinct points $a, b\in \affs$, the locus of all points
having equal distance to $a$ and $b$ is a hyperplane.
It is called the {\it bisector hyperplane of $a$ and $b$\/}. 
\index{bisector!hyperplane}%
The equation of this hyperplane is still found 
by writing that $\libvecbo{h}{m}$ is orthogonal to $\libvecbo{a}{b}$.
The equation of this hyperplane  is
\[(b_1 - a_1)(x_1 - (a_1 + b_1)/2) + \cdots + 
(b_\mdeg - a_\mdeg)(x_\mdeg - (a_\mdeg + b_\mdeg)/2) = 0,\]
which can also be written as
\[(b_1 - a_1)x_1  + \cdots +  (b_\mdeg - a_\mdeg)x_\mdeg =
(b_1^2 + \cdots + b_\mdeg^2)/2 - (a_1^2+  \cdots + a_\mdeg^2)/2.\]
The closed half-space $H(a, b)$ containing $a$ and with boundary the bisector hyperplane
is the locus of all points such that
\[(b_1 - a_1)x_1  + \cdots +  (b_\mdeg - a_\mdeg)x_\mdeg \leq
(b_1^2 + \cdots + b_\mdeg^2)/2 - (a_1^2+  \cdots + a_\mdeg^2)/2,\]
and the  closed half-space $H(b, a)$ containing $b$ and with boundary the bisector hyperplane
is the locus of all points such that
\[(b_1 - a_1)x_1  + \cdots +  (b_\mdeg - a_\mdeg)x_\mdeg \geq
(b_1^2 + \cdots + b_\mdeg^2)/2 - (a_1^2+  \cdots + a_\mdeg^2)/2.\]
The closed half-space $H(a, b)$ is the set of all points whose distance to $a$
is less than or equal to the distance to $b$, and vice versa for $H(b, a)$.

\begin{figure}[H]
\centerline{
\psfig{figure=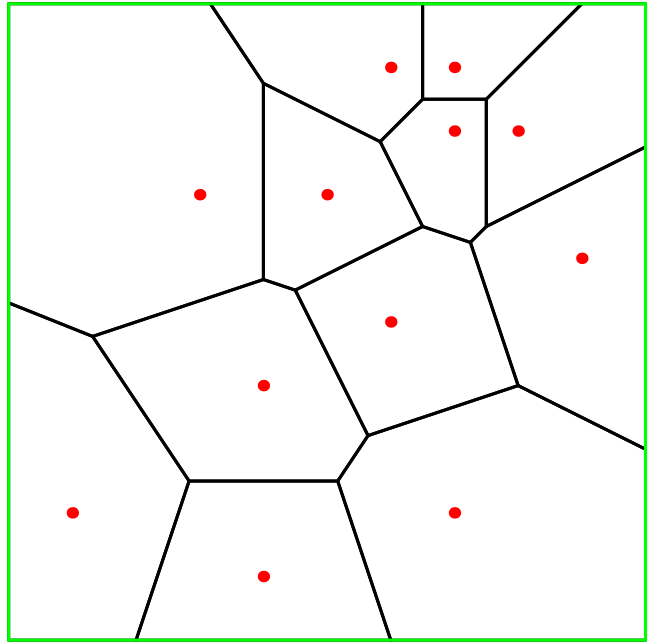,width=4truein}}
\caption{A Voronoi diagram}
\label{vorpic1}
\end{figure}

\medskip
Figure \ref{vorpic1} shows the Voronoi diagram of a set of twelve points.
\medskip
In the general case where $\affs$ has dimension $\mdeg$,
the definition of the Voronoi diagram $\s{V}\mathit{or}(P)$ of $P$ is  the same as Definition
\ref{Voronoi}, 
\indsym{\s{V}\mathit{or}(P)}{Voronoi diagram of $P$}%
except that  $H(p_i, p_j)$ is the  closed half-space
containing $p_i$ and having the bisector hyperplane of $a$ and $b$
as boundary. Also, observe that the convex hull of $P$ is a convex polytope.

\medskip
We will now state a lemma listing the main properties of Voronoi 
diagrams. It turns out that certain degenerate situations can be avoided
if we assume that if $P$ is a set of points
in an affine space of dimension $\mdeg$, then no $\mdeg + 2$ points from $P$
belong to the same $(\mdeg - 1)$-sphere. We will say that the points
of $P$ are in {\it general position\/}.
\index{points!in general position}%
Thus when $\mdeg = 2$, no $4$ points in $P$ are cocyclic, and when $\mdeg = 3$, 
no $5$ points in $P$ are on the same sphere.

\begin{lemma}
\label{voronoilem}
Given a set $P = \{p_1,\ldots,p_n\}$ of $n$ points in some Euclidean space $\affs$
of dimension $\mdeg$ (say $\eucreal^{m}$), 
if the points in $P$ are in general position and not in a common hyperplane 
then the 
Voronoi diagram of $P$ satisfies the following conditions:
\begin{enumerate}
\item[(1)]
Each region $V_i$ is convex and contains $p_i$ in its interior.
\item[(2)]
Each vertex of $V_i$ belongs to $\mdeg + 1$ regions $V_j$ and to
$\mdeg + 1$ edges.
\item[(3)]
The region $V_i$ is unbounded iff $p_i$ belongs to the
boundary of the convex hull of $P$.
\item[(4)]
If $p$ is a vertex that belongs to the regions $V_1, \ldots, V_{\mdeg + 1}$,
then $p$ is the center of the $(\mdeg - 1)$-sphere $S(p)$ determined by
$p_1, \ldots, p_{\mdeg + 1}$. Furthermore,  no point in $P$
is inside the sphere $S(p)$ (i.e., in the open ball associated with
the sphere $S(p)$).
\item[(5)]
If $p_j$ is a nearest neighbor of $p_i$, then one of the faces of $V_i$
is contained in the bisector hyperplane of $(p_i,p_j)$.
\item[(6)]
\[\bigcup_{i=1}^{n} V_i = \affs,\quad\hbox{and}\quad \interio{V_i}\cap \interio{V_j} = \emptyset,
\quad\hbox{for all $i, j$, with $i\not= j$},\]
where $\interio{V_i}$ denotes the interior of $V_i$.
\end{enumerate}
\end{lemma}

\medskip
\proof We prove only some of the statements, leaving the others
as an exercise (or see Risler \cite{Risler92}).

\medskip
(1) Since $V_i =\bigcap_{j\not= i} H(p_i, p_j)$ and each half-space $H(p_i, p_j)$ is
convex, as an intersection of convex sets, $V_i$ is convex. Also, since $p_i$
belongs to the interior of each $H(p_i, p_j)$, the point $p_i$ belongs to the
interior of $V_i$.

\medskip
(2)
Let $F_{i, j}$ denote $V_i\cap V_j$. Any vertex $p$ of the Vononoi diagram of $P$
must belong to $r$ faces $F_{i, j}$. 
Now, given a vector space $E$ and any two subspaces $M$ and $N$ of $E$,
recall that we have the {\it Grassmann relation\/} 
\[\dimm{(M)} + \dimm{(N)} = \dimm{(M + N)} + \dim{(M\cap N)}.\]  
Then since $p$ belongs to the intersection of the hyperplanes that form
the boundaries of the $V_i$, and since a hyperplane has dimension $\mdeg -1$,
by the Grassmann relation, we must have $r\geq \mdeg$.
For simplicity of notation, let us denote these faces by
$F_{1, 2}, F_{2, 3}, \ldots, F_{r, r+1}$. 
Since $F_{i,j} = V_i\cap V_j$, we have
\[F_{i, j} = \{p\ |\ d(p, p_i) = d(p, p_j) \leq d(p, p_k),\ 
\hbox{for all $k \not= i, j$}\},\]
and since $p \in F_{1, 2}\cap F_{2, 3}\cap \cdots\cap F_{r, r+1}$,
we have
\[d(p,p_1) = \cdots = d(p, p_{r+1}) < d(p, p_k)\ 
\hbox{for all $k\notin \{1, \ldots, r+1\}$}.\]
This means that $p$ is the center of a sphere passing through
$p_1, \ldots, p_{r+1}$ and containing no other point in $P$.
By the assumption that points in $P$ are in general position, we must have $r\leq \mdeg$,
and thus $r = \mdeg$.
Thus, $p$ belongs to $V_1\cap \cdots\cap V_{\mdeg + 1}$, but to no other $V_j$ with
$j\notin \{1, \ldots, \mdeg + 1\}$. Furthermore, every edge of the Voronoi diagram
containing $p$ is the intersection of $\mdeg$ of the regions $V_1, \ldots, V_{\mdeg + 1}$, 
and so there are $\mdeg + 1$ of them.
$\bigsquare$
\endproof

\begin{figure}[H]
\centerline{
\psfig{figure=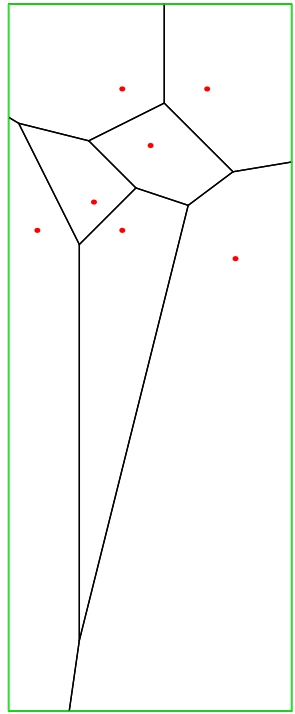,height=4truein}}
\caption{Another Voronoi diagram
\label{vorpic2}
}
\end{figure}

\medskip
For simplicity, let us again consider the case where $\affs$ is a plane.
It should be noted that certain Voronoi regions,
although closed, may extend very far.
Figure \ref{vorpic2}  shows such an example.

\medskip
It is also possible for certain unbounded regions to have parallel edges.

\medskip
There are a number of methods for computing Voronoi diagrams.
A fairly simple (although not very efficient) method is to compute
each Voronoi region $V(p_i)$ by intersecting the half-planes
$H(p_i, p_j)$. One way to do this is to construct successive
convex polygons that converge to the boundary of the region.
At every step we intersect the current convex polygon with
the bisector line of $p_i$ and $p_j$. There are at most two intersection
points. We also need a starting polygon, and for this we can pick
a square containing all the points. A naive implementation
will run in $O(n^3)$. However, the intersection of half-planes
can be done in $O(n\log n)$, using the fact that the vertices of a convex
polygon can be sorted. Thus, the above method runs in $O(n^2\log n)$.
Actually, there are faster methods (see Preparata and Shamos \cite{Preparata}
or O'Rourke \cite{ORourke}), and it is possible to design algorithms
running in $O(n\log n)$. The most direct method to obtain 
fast algorithms is to use the ``lifting method'' discussed in Section \ref{sec224}, 
whereby the original set of points is lifted onto a paraboloid,
and to use fast algorithms for finding a convex hull.

\medskip
A very interesting (undirected) graph can be obtained from the Voronoi diagram as follows:
The vertices of this graph are the points $p_i$ (each corresponding to
a unique region of $\s{V}\mathit{or}(P)$), and there is an edge between $p_i$ and $p_j$ iff
the regions $V_i$ and $V_j$ share an edge. The resulting graph is called
a {\it Delaunay triangulation\/} of the convex hull of $P$, after Delaunay, who
invented this concept in 1934. Such triangulations have remarkable properties.
\index{Delaunay triangulation}%

\begin{figure}
\centerline{
\psfig{figure=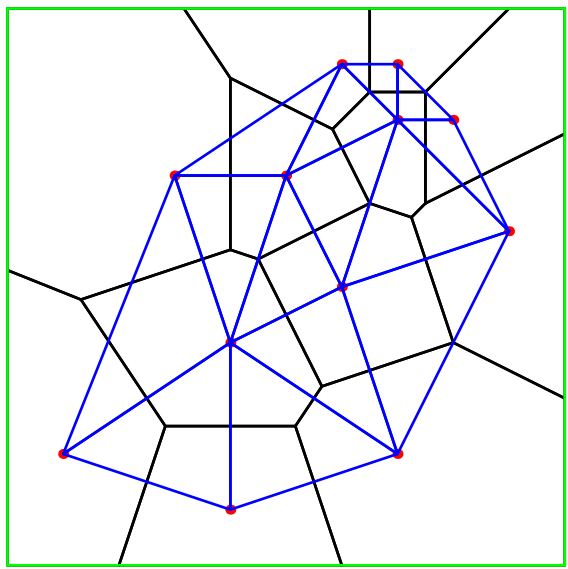,width=3.5truein}}
\caption{Delaunay triangulation associated with a Voronoi diagram}
\label{vorpic3}
\end{figure}

\medskip
Figure \ref{vorpic3} shows the Delaunay triangulation
associated with the earlier Voronoi diagram of a set of twelve points.

\begin{figure}
\centerline{
\psfig{figure=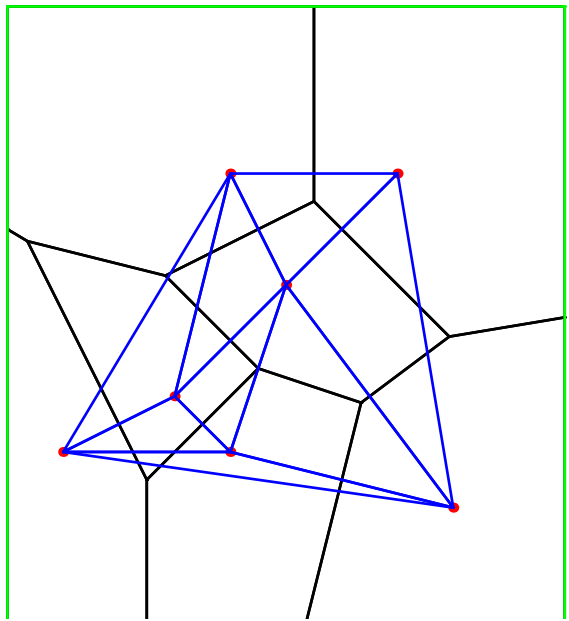,width=3.5truein}}
\caption{Another Delaunay triangulation associated with a Voronoi diagram}
\label{vorpic4}
\end{figure}

\medskip
One has to be careful to make sure that all the Voronoi vertices have been computed
before computing a Delaunay triangulation, since otherwise, some
edges could be missed. 
In Figure \ref{vorpic4} illustrating such a situation,
if the lowest Voronoi vertex had not been computed 
(not shown on the diagram!),
the lowest edge of the Delaunay triangulation would be missing.

\medskip
The concept of a triangulation can be
generalized to dimension $3$, or even to any dimension $\mdeg$.
\section{Triangulations}
\label{sec222}
The concept of a triangulation relies on the notion of
pure simplicial complex defined in Chapter \ref{chapter3}.
The reader should review Definition \ref{complexdef}
and Definition \ref{subcomplex}.

\begin{defin}
\label{triangdef}
{\em
Given a subset, $S\subseteq \eucreal^\mdeg$ (where $\mdeg \geq 1$),
a {\it triangulation of $S$\/} 
\index{triangulation!definition}%
is a pure (finite) simplicial complex, $K$, of dimension $m$ such that
$S = |K|$, that is, $S$ is equal to the geometric realization of $K$.
}
\end{defin}

\medskip
Given a finite set $P$ of $n$ points in the plane, and given a triangulation
of the convex hull of $P$ having $P$ as its set of vertices,
observe that the boundary of $P$ is
a convex polygon. Similarly, given a finite set $P$ of points in $3$-space,
and given a triangulation
of the convex hull of $P$ having $P$ as its set of vertices,
observe that the boundary of $P$ is
a convex polyhedron. It is interesting to know how many
triangulations exist for a set of $n$ points (in the plane or in $3$-space),
and it is also interesting to know  the number of edges and faces
in terms of the number of vertices in $P$.
These questions can be settled using the Euler--Poincar\'e characteristic.
\index{Euler--Poincar\'e characteristic}%
We say that a polygon in the plane is a {\it simple polygon\/}
\index{simple polygon}%
iff it is a connected closed polygon such that
no two edges intersect (except at a common vertex). 

\begin{lemma}
\label{Euler}
\
\begin{enumerate}
\item[(1)] For any triangulation of a region of the plane whose boundary is a simple
polygon, letting $v$ be the number of vertices, $e$ the number of edges, and
$f$ the number of triangles, we have the ``Euler formula''
\[v - e + f = 1.\]
\index{Euler's formula}%
\item[(2)]
For any region, $S$, in $\eucreal^3$ homeomorphic to a closed ball
and for any triangulation of $S$, letting
$v$ be the number of vertices, $e$ the number of edges, 
$f$ the number of triangles, and $t$ the number of tetrahedra,
we have the ``Euler formula''
\[v - e + f - t = 1.\]
\index{Euler's formula}%
\item[(3)]
Furthermore, for any triangulation of the combinatorial surface,
$B(S)$, that is the boundary of $S$, 
letting $v'$ be the number of vertices, $e'$ the number of edges, and
$f'$ the number of triangles, we have the ``Euler formula''
\[v' - e' + f' = 2.\]
\end{enumerate}
\index{Euler's formula}%
\end{lemma}

\proof 
All the statements are immediate consequences of Theorem
\ref{PoincareEulerthm1}. For example, part (1) is obtained by mapping the
triangulation onto a sphere using  inverse stereographic projection,
say from the North pole. Then, we get a polytope on the sphere
with an extra facet corresponding to the ``outside'' 
of the triangulation. We have to deduct this facet from the
Euler characteristic of the polytope  and this is why we get $1$ instead of $2$.
$\bigsquare$
\endproof

\medskip
It is now easy to see that in case (1),  the number of edges and faces
is a linear function of the number of vertices and boundary edges, and that
in case (3),   the number of edges and faces
is a linear function of the number of vertices.
Indeed, in the case of a planar triangulation, each face has $3$ edges, 
and if there are $e_b$ edges in the boundary and $e_i$ edges not
in the boundary, each nonboundary edge
is shared by two faces, and thus $3f = e_b + 2e_i$.
Since $v - e_b - e_i + f = 1$, we get
\[
\eqaligneno{
v - e_b - e_i + e_b/3 + 2e_i/3  &= 1,\cr
2e_b/3 + e_i/3 &= v - 1,\cr
}
\]
and thus $e_i = 3v - 3 - 2e_b$.  Since $f = e_b/3 + 2e_i/3$,
we have $f = 2v -2 - e_b$. 

\medskip
Similarly, since
$v' - e' + f' = 2$ and $3f' = 2e'$, we easily get
$e = 3v - 6$ and $f = 2v - 4$.
Thus, given a set $P$ of $n$ points, the number of triangles (and edges)
for any triangulation of the convex hull of $P$  using the $n$ points
in $P$ for its vertices is fixed.

\medskip
Case (2) is  trickier, but it can be shown that
\[v - 3 \leq t \leq (v-1)(v-2)/2.\]
Thus, there can be different numbers of tetrahedra for different 
triangulations of the convex hull of $P$.

\medskip
\remark
The numbers of the form $v - e + f$ and
$v - e + f - t$ are called {\it Euler--Poincar\'e characteristics\/}.
\index{Euler--Poincar\'e characteristic}%
They are topological invariants, in the sense that they are the same
for all triangulations of a given polytope. 
\index{topological invariants}%
This is a fundamental
fact of algebraic topology.
\endremark

\medskip
We shall now investigate triangulations induced by Voronoi diagrams.
\section{Delaunay Triangulations}
\label{sec223}
Given a set $P = \{p_1,\ldots,p_n\}$ of $n$ points in the plane
and the Voronoi diagram $\s{V}\mathit{or}(P)$ for $P$, we explained in Section 
\ref{sec221} how to define an (undirected) graph:
The vertices of this graph are the points $p_i$ (each corresponding to
a unique region of $\s{V}\mathit{or}(P)$), and there is an edge between $p_i$ and $p_j$ iff
the regions $V_i$ and $V_j$ share an edge. The resulting graph turns out
to be a triangulation of the convex hull of $P$ having $P$ as its
set of vertices. Such a complex can be defined in general.
For any set $P = \{p_1,\ldots,p_n\}$ of $n$ points in $\eucreal^\mdeg$,
we say that a triangulation of the convex hull of $P$ is
{\it associated with $P$\/} if its set of vertices is the set $P$.

\begin{defin}
\label{Delaudef}
{\em
Let  $P = \{p_1,\ldots,p_n\}$ be a set of $n$ points in $\eucreal^\mdeg$,
and let $\s{V}\mathit{or}(P)$ be the Voronoi diagram of $P$. We define
a complex $\s{D}\mathit{el}(P)$ as follows.
The complex $\s{D}\mathit{el}(P)$  contains the $k$-simplex
$\{p_1,\ldots, p_{k+1}\}$ iff 
$V_1\cap \cdots \cap V_{k+1} \not= \emptyset$, where $0\leq k \leq \mdeg$.
The complex $\s{D}\mathit{el}(P)$ is called the {\it Delaunay triangulation of 
the convex hull of $P$\/}.
\index{Delaunay triangulation!definition}%
}
\end{defin}

\medskip
Thus, $\{p_i, p_j\}$ is an edge iff $V_i\cap V_j \not= \emptyset$,
$\{p_i, p_j, p_h\}$ is a triangle  iff $V_i\cap V_j \cap V_h\not= \emptyset$,
$\{p_i, p_j, p_h, p_k\}$ is a tetrahedron 
iff $V_i\cap V_j \cap V_h\cap V_k\not= \emptyset$,
etc.

\begin{figure}
\centerline{
\psfig{figure=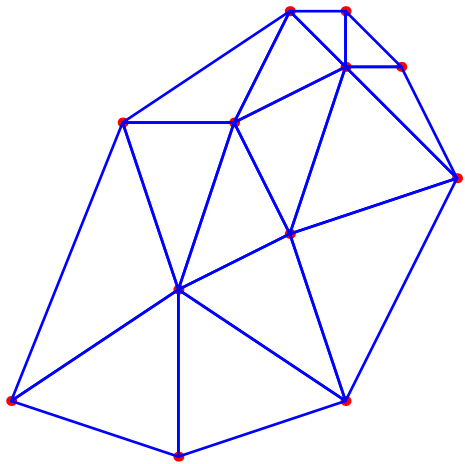,width=3.5truein}}
\caption{A Delaunay triangulation}
\label{delaupic1}
\end{figure}

\medskip
For simplicity, we often write $\s{D}\mathit{el}$ instead of $\s{D}\mathit{el}(P)$.
A Delaunay triangulation for a set of twelve points is shown in
Figure \ref{delaupic1}.
\indsym{\s{D}\mathit{el}(P)}{Delaunay triangulation of $P$}%

\medskip
Actually, it is not obvious that  $\s{D}\mathit{el}(P)$ is a triangulation
of the convex hull of $P$, but this can be shown, as well as the
properties listed in the following lemma.

\begin{lemma}
\label{Delaunaylem}
Let  $P = \{p_1,\ldots,p_n\}$ be a set of $n$ points in $\eucreal^\mdeg$,
and assume that they are in general position.
Then the Delaunay triangulation of the convex hull of $P$ is
indeed a triangulation associated with $P$, and it satisfies the following
properties:
\begin{enumerate}
\item[(1)]
The boundary of $\s{D}\mathit{el}(P)$ is the convex hull of $P$.
\item[(2)]
A triangulation $T$ associated with $P$ is the Delaunay triangulation 
$\s{D}\mathit{el}(P)$ iff every $(\mdeg - 1)$-sphere $S(\sigma)$ circumscribed about
an $\mdeg$-simplex $\sigma$ 
of $T$ contains no other point from $P$
(i.e., the open ball associated with $S(\sigma)$ contains no point from $P$).
\end{enumerate}
\end{lemma}

\medskip
The proof can be found in Risler \cite{Risler92} and O'Rourke \cite{ORourke}.
In the case of a planar set $P$, it can also be shown that the Delaunay
triangulation has the property that it maximizes the minimum angle
of the triangles involved in any triangulation of $P$.
However, this does not characterize the Delaunay triangulation.
Given a connected graph in the plane, it can also be shown that any minimal
spanning tree is contained in the Delaunay triangulation of the
convex hull of the set of vertices of the graph  (O'Rourke \cite{ORourke}).
\index{spanning tree}%

\medskip
We will now explore briefly the connection between Delaunay
triangulations and convex hulls.
\section{Delaunay Triangulations and Convex Hulls}
\label{sec224}
In this section we show that there is an intimate relationship between
convex hulls and Delaunay triangulations. We will see that
given a set $P$ of points in the Euclidean space $\eucreal^\mdeg$
 of dimension $\mdeg$,
we can ``lift'' these points onto a paraboloid living
in the  space $\eucreal^{\mdeg + 1}$ of dimension $\mdeg+1$,
and that the Delaunay triangulation of $P$ is the projection
of the downward-facing faces of the convex hull of the set of lifted
points. This remarkable connection was first discovered by
Edelsbrunner and Seidel \cite{Eldel}.
For simplicity, we consider the case of a set $P$ of points in the plane
$\eucreal^2$, and we assume that they are in general position.

\medskip
Consider the paraboloid of revolution of equation
$z = x^2 + y^2$. 
\index{paraboloid of revolution}%
A point $p = (x, y)$ in the plane is lifted to
the point $l(p) = (X, Y, Z)$ in $\eucreal^3$, where
$X = x$, $Y = y$, and $ Z = x^2 + y^2$.

\medskip
The first crucial observation is that a circle in the plane is lifted
into a plane curve (an ellipse). Indeed, if such a circle $C$ is defined
by the equation 
\[x^2 + y^2 + ax + by + c = 0,\]
since $X = x$, $Y = y$, and $ Z = x^2 + y^2$, by eliminating $x^2 + y^2$
we get 
\[Z = -ax - by - c,\]
and thus $X, Y, Z$ satisfy the linear equation
\[aX + bY + Z + c = 0,\]
which is the equation of a plane. 
Thus, the intersection of the cylinder of revolution consisting of the
lines parallel to the $z$-axis and passing through a point of 
the circle $C$ with the paraboloid $z = x^2 + y^2$ is
a planar curve (an ellipse).

\medskip
We can compute the convex hull of the set of lifted points.
Let us focus on the downward-facing faces of this convex hull.
Let $(l(p_1), l(p_2), l(p_3))$ be such a face.
The points $p_1, p_2, p_3$ belong to the set $P$.
We claim
that no other point from $P$ is inside the circle $C$.
Indeed, a point $p$ inside the circle $C$ would lift to a point
$l(p)$ on the paraboloid. Since no four points are cocyclic,
one of the four points $p_1, p_2, p_3, p$ is further
from $O$ than the others; say this point is $p_3$. Then,
the face $(l(p_1), l(p_2), l(p))$ would be below the face
$(l(p_1), l(p_2), l(p_3))$, contradicting the fact that
$(l(p_1), l(p_2), l(p_3))$ is one of the downward-facing
faces of the convex hull of
$P$. But then, by property (2) of Lemma \ref{Delaunaylem},
the triangle $(p_1, p_2, p_3)$ would belong to the Delaunay triangulation
of  $P$.

\begin{figure}
\centerline{
\psfig{figure=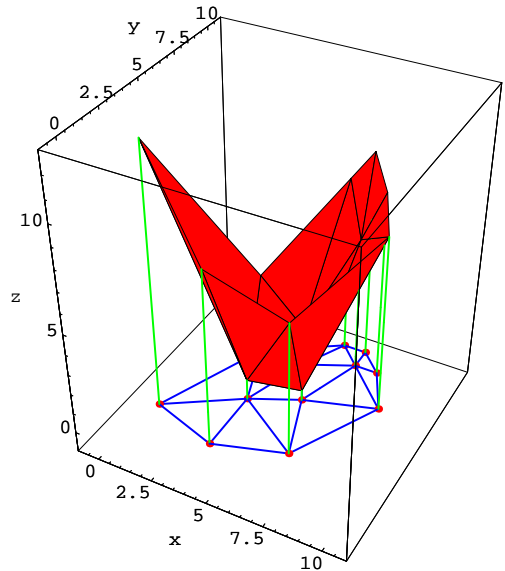,width=3truein}}
\caption{A Delaunay triangulation and its lifting to a paraboloid}
\label{parapic1}
\end{figure}

\medskip
Therefore, we have shown that 
{\it the projection of the part of the convex hull
of the lifted set $l(P)$ consisting of the downward-facing faces is the
Delaunay triangulation of $P$\/}.
\index{convex hull!and Delaunay triangulation}%
Figure \ref{parapic1} shows the lifting of the
Delaunay triangulation shown earlier.

\begin{figure}
\centerline{
\psfig{figure=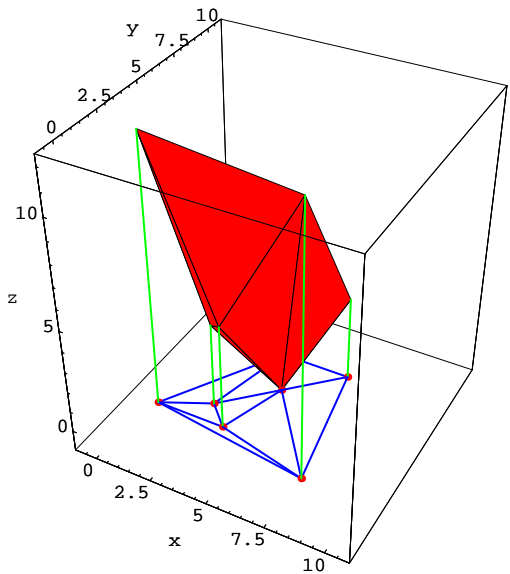,width=3truein}}
\caption{Another Delaunay triangulation and its lifting to a paraboloid}
\label{parapic2}
\end{figure}

Another example of the lifting of a Delaunay triangulation is shown 
in Figure \ref{parapic2}.

\medskip
The fact that a Delaunay triangulation can be obtained by projecting
a lower convex hull  
can be used to find efficient algorithms for
computing a Delaunay triangulation. It also holds for higher
dimensions.

\medskip
The Voronoi diagram itself can also be obtained from the lifted
set $l(P)$. However, this time, we need to consider tangent
planes to the paraboloid at the lifted points.
It is fairly obvious that the tangent plane at the lifted
point $(a, b, a^2 + b^2)$ is
\[z = 2ax + 2by - (a^2 + b^2).\]
Given two distinct lifted points $(a_1, b_1, a_1^2 + b_1^2)$ and
$(a_2, b_2, a_2^2 + b_2^2)$, the intersection of the tangent planes
at these points is a line belonging to the plane of equation
\[(b_1 - a_1)x  + (b_2 - a_2)y  = 
(b_1^2 + b_2^2)/2 - (a_1^2+  a_2^2)/2.\]
Now, if we project this plane onto the $xy$-plane, we see that
the above is precisely the equation of 
the bisector line of the two points
$(a_1, b_1)$ and $(a_2, b_2)$.
Therefore, {\it if we look at the paraboloid from $z = +\infty$
(with the paraboloid transparent), the projection of the tangent
planes at the lifted points is the Voronoi diagram\/}!
\index{convex hull!and Voronoi diagram}%

\medskip
It should be noted that the ``duality'' between the Delaunay
triangulation, which is the projection of the convex hull of 
the lifted set $l(P)$ viewed from $z = -\infty$, and the 
Voronoi diagram, which is the projection of the tangent
planes at the lifted set $l(P)$ viewed from $z = +\infty$,
is reminiscent of the polar duality with respect to a quadric.
\index{duality!Delaunay triangulations, Voronoi diagrams}%
This duality will be thoroughly investigated  in Section \ref{sec13}.

\medskip
The reader interested in algorithms for finding
Voronoi diagrams and Delaunay triangulations is referred to
O'Rourke \cite{ORourke}, Preparata and Shamos \cite{Preparata},
Boissonnat and Yvinec \cite{Boissonnat}, 
de Berg, Van Kreveld, Overmars, and Schwarzkopf \cite{Berg},
and Risler \cite{Risler92}.

\section[Stereographic Projection and the Space of Spheres]
{Stereographic Projection  and the  Space of \\
Generalized Spheres}
\label{sec12}
Brown appears to be the first person who
observed  that Voronoi diagrams and convex hulls
are related {\it via\/}  inversion with respect to a
sphere   \cite{Brown}.

\medskip
In fact, more generally, it turns out that   
Voronoi diagrams,
Delaunay Triangulations and their properties can also be nicely explained
using  stereographic projection  and its inverse,
although a rigorous justification of why this
``works'' is not as simple as it might appear.

\medskip
The advantage of stereographic projection
over the lifting onto a paraboloid is that the \\
($d$-)sphere is compact.
Since the stereographic projection and its inverse map $(d-1)$-spheres
to $(d-1)$-spheres (or hyperplanes), all the crucial properties
of Delaunay triangulations are preserved.
The purpose of this section is to establish the properties
of stereographic projection (and its inverse) that will
be needed in Section \ref{sec13}.

\medskip
Recall that the {\it $d$-sphere\/}, $S^d \subseteq \eucreal^{d+1}$,
is given by
\[
S^d = \{(x_1, \ldots, x_{d+1})\in \eucreal^{d+1} \mid
x_1^2 + \cdots + x_d^2  + x_{d+1}^2 = 1\}.
\]
It will be convenient to write a point,
$(x_1, \ldots, x_{d+1})\in \eucreal^{d+1}$, as 
$z = (x, x_{d+1})$, with  \\
$x = (x_1, \ldots, x_{d})$.
We denote $N = (0, \ldots, 0, 1)$ (with $d$ zeros) as
$(\mathbf{0}, 1)$ and call it the {\it north pole\/} and
$S = (0, \ldots, 0, -1)$ (with $d$ zeros) as
$(\mathbf{0}, -1)$ and call it the {\it south pole\/}.
We also write $\norme{z} = (x_1^2+  \cdots +  x_{d+1}^2)^{\frac{1}{2}}
= (\norme{x}^2 + x_{d+1}^2)^{\frac{1}{2}}$
(with $\norme{x} = (x_1^2+  \cdots +  x_{d}^2)^{\frac{1}{2}}$).
With these notations,
\[
S^d = \{(x, x_{d+1}) \in \eucreal^{d+1} \mid \norme{x}^2 + x_{d+1}^2 = 1\}.
\]

\medskip
The {\it stereographic projection from the north pole,
$\mapdef{\sigma_N}{(S^d - \{N\})}{\eucreal^d}$\/}, is the
restriction to $S^d$ of the central projection from $N$
onto the hyperplane, 
$H_{d+1}(0) \cong \eucreal^d$, of equation $x_{d+1} = 0$; that is,
$M \mapsto \sigma_N(M)$ where $\sigma_N(M)$ is the intersection
of the line, $\lag N, M\rag$,  through $N$  and $M$ with $H_{d+1}(0)$. 
Since the line through $N$ and $M = (x, x_{d+1})$ is given parametrically by
\[
\lag N, M\rag  =
\{(1 - \lambda)(\mathbf{0}, 1) + \lambda (x, x_{d+1}) \mid \lambda\in \reals\}, 
\]
the intersection, $\sigma_N(M)$, of this line with the hyperplane
$x_{d+1} = 0$ corresponds to the value of $\lambda$ such that
\[
(1 - \lambda) + \lambda x_{d + 1} = 0,
\]
that is, 
\[
\lambda = \frac{1}{1 - x_{d+1}}.
\]
Therefore, the coordinates of $\sigma_N(M)$, with $M = (x, x_{d+1})$, are 
given by
\[
\sigma_N(x, x_{d+1}) = \left(\frac{x}{1 - x_{d+1}}, 0  \right).
\]
Let us find the inverse, $\tau_N = \sigma_N^{-1}(P)$, 
of any $P \in H_{d+1}(0)$. This time, $\tau_N(P)$ is the intersection
of the line, $\lag N, P\rag$, through $P\in H_{d+1}(0)$ and $N$
with the sphere, $S^d$. 
Since the line through $N$ and $P = (x, 0)$ is given parametrically by
\[
\lag N, P\rag  =
\{(1 - \lambda)(\mathbf{0}, 1) + \lambda (x, 0) \mid \lambda\in \reals\}, 
\]
the intersection, $\tau_N(P)$, of this line with the sphere $S^d$
corresponds to the nonzero value of $\lambda$ such that
\[
\lambda^2 \norme{x}^2 + (1 - \lambda)^2 = 1,
\]
that is
\[
\lambda(\lambda (\norme{x}^2 + 1) - 2) = 0. 
\]
Thus, we get
\[
\lambda = \frac{2}{\norme{x}^2 + 1},
\]
from which we get
\begin{eqnarray*}
\tau_N(x, 0) & = & \left(\frac{2x}{\norme{x}^2 + 1}, 
1 - \frac{2}{\norme{x}^2 + 1} \right) \\
 & = & \left(\frac{2x}{\norme{x}^2 + 1}, 
\frac{\norme{x}^2 - 1}{\norme{x}^2 + 1} \right).
\end{eqnarray*}

\medskip
We leave it as an exercise to the reader to verify that
$\tau_N\circ \sigma_N = \id$ and $\sigma_N\circ \tau_N = \id$.
We can also define the 
{\it stereographic projection from the south pole,
$\mapdef{\sigma_S}{(S^d - \{S\})}{\eucreal^d}$\/}, and its inverse,
$\tau_S$. Again, the computations are left as a simple exercise to the reader. 
The above computations are summarized in the following
definition:

\begin{defin}
\label{stereoproj}
{\em 
The {\it stereographic projection from the north pole,
$\mapdef{\sigma_N}{(S^d - \{N\})}{\eucreal^d}$\/}, is the
map given by
\[
\sigma_N(x, x_{d+1}) = \left(\frac{x}{1 - x_{d+1}}, 0  \right)
\qquad (x_{d+1} \not= 1).
\]
The inverse of $\sigma_N$, denoted $\tau_N$ and called
{\it inverse stereographic projection from the north pole\/}
is given by
\[
\tau_N(x, 0) 
  =  \left(\frac{2x}{\norme{x}^2 + 1}, 
\frac{\norme{x}^2 - 1}{\norme{x}^2 + 1} \right).
\]
}
\end{defin}

\remark
An {\it inversion of center $C$ and power $\rho > 0$\/}
is a geometric transformation, \\
$\mapdef{f}{(\eucreal^{d+1} - \{C\})}{\eucreal^{d+1}}$, defined
so that for any $M\not= C$, the points
$C$, $M$ and $f(M)$ are collinear and
\[
\smnorme{\libvecbo{C}{M}}\smnorme{\libvecbo{C}{f(M)}} = \rho.
\]
Equivalently, $f(M)$ is given by 
\[
f(M) = C + \frac{\rho}{\smnorme{\libvecbo{C}{M}}^2}\,\libvecbo{C}{M}.
\]
Clearly, $f\circ f = \id$ on $\eucreal^{d+1} - \{C\}$,
so $f$ is invertible and the reader will check that
if we pick the center of inversion to be the north pole
and if we set $\rho = 2$, then 
the coordinates of $f(M)$ are given by
\begin{eqnarray*}
y_i & = & \frac{2x_i}{x_1^2 + \cdots + x_d^2 + x_{d+1}^2 -2x_{d+1} + 1}
\qquad\qquad 1\leq i \leq d \\
y_{d+1} & = & \frac{x_1^2 + \cdots + x_d^2 + x_{d+1}^2 - 1}
{x_1^2 + \cdots + x_d^2 + x_{d+1}^2 -2x_{d+1} + 1},
\end{eqnarray*}
where $(x_1, \ldots, x_{d+1})$ are the coordinates of $M$.
In particular, if we restrict our inversion to the unit sphere,
$S^d$, as $x_1^2 + \cdots + x_d^2 + x_{d+1}^2 = 1$,
we get
\begin{eqnarray*}
y_i & = & \frac{x_i}{1 - x_{d+1}}
\qquad\qquad 1\leq i \leq d \\
y_{d+1} & = & 0,
\end{eqnarray*}
which means that our inversion restricted to $S^d$ is
simply the stereographic projection, $\sigma_N$
(and the inverse of our inversion restricted to the
hyperplane, $x_{d+1} = 0$, is the inverse stereographic
projection, $\tau_N$).

\medskip
We will now show that the image of any $(d-1)$-sphere, $S$, on
$S^d$ not passing through the north pole, that is, the intersection,
$S = S^d\cap H$, of $S^d$ 
with any hyperplane, $H$,  not passing through $N$  is
a $(d - 1)$-sphere. Here, we are assuming that 
$S$ has positive radius, that is, $H$ is {\it not\/} tangent to $S^d$.

\medskip
Assume that $H$ is given by
\[
a_1x_1 + \cdots + a_dx_d + a_{d+1} x_{d+1} + b = 0.
\]
Since $N\notin H$, we must have $a_{d+1} + b \not= 0$.
For any $(x, x_{d+1}) \in S^d$, write $\sigma_N(x, x_{d+1}) = (X, 0)$.
Since
\[
X = \frac{x}{1 - x_{d+1}},
\]
we get $x = X(1 - x_{d+1})$ and 
using the fact that $(x, x_{d+1})$ also belongs to $H$ we
will express $x_{d+1}$ in terms of $X$ and then find
an equation for $X$ which will show that $X$ belongs to a 
$(d - 1)$-sphere. Indeed, $(x, x_{d+1})\in H$ implies that
\[
\sum_{i = 1}^d a_i X_i(1 - x_{d+1}) + a_{d+1}x_{d+1} + b = 0,
\]
that is,
\[
\sum_{i = 1}^d a_i X_i +  (a_{d+1} - \sum_{j = 1}^d a_j X_j) x_{d+1} +  b = 0.
\]
If $\sum_{j = 1}^d a_j X_j = a_{d+1}$, then $a_{d+1} + b = 0$,
which is impossible. Therefore, we get
\[
x_{d+1} = \frac{-b - \sum_{i = 1}^d a_i X_i}{a_{d+1}  - \sum_{i = 1}^d a_i X_i}
\]
and so,
\[
1 - x_{d+1} = \frac{a_{d+1} + b}{a_{d+1}  - \sum_{i = 1}^d a_i X_i}.
\]
Plugging $x = X(1 - x_{d+1})$ in the equation, $\norme{x}^2 + x_{d+1}^d = 1$, 
of $S^d$, we get
\[
(1 - x_{d+1})^2\norme{X}^2 + x_{d+1}^2 = 1,
\]
and replacing $x_{d+1}$ and $1 - x_{d+1}$ by their expression in terms of
$X$, we get
\[
(a_{d+1} + b)^2\norme{X}^2 + (-b - \sum_{i = 1}^d a_i X_i)^2 
 =  (a_{d+1}  - \sum_{i = 1}^d a_i X_i)^2 
\]
that is, 
\begin{eqnarray*}
(a_{d+1} + b)^2\norme{X}^2 
& = &  (a_{d+1}  - \sum_{i = 1}^d a_i X_i)^2  - (b + \sum_{i = 1}^d a_i X_i)^2  \\
& = &  (a_{d+1} + b)(a_{d+1} - b  - 2\sum_{i = 1}^d a_i X_i)^2  \\
\end{eqnarray*}
which yields
\[
(a_{d+1} + b)^2\norme{X}^2 + 2(a_{d+1} + b)(\sum_{i = 1}^d a_i X_i) = 
(a_{d+1} + b)(a_{d+1} - b),
\]
that is,
\[
\norme{X}^2 + 2\sum_{i = 1}^d \frac{a_i}{a_{d+1} + b} X_i  
- \frac{a_{d+1} - b}{a_{d+1} + b} = 0,
\]
which is indeed the equation of a $(d-1)$-sphere in $\eucreal^d$.
Therefore, when $N\notin H$, the image of $S = S^d\cap H$
by $\sigma_N$ is a $(d-1)$-sphere in $H_{d+1}(0) = \eucreal^d$.

\medskip
If the hyperplane, $H$, contains the north pole, then
$a_{d+1} + b = 0$, in which case, for every $(x, x_{d+1})\in S^d\cap H$,
we have
\[
\sum_{i = 1}^d a_i x_i + a_{d+1}x_{d+1} - a_{d+1} = 0,
\]
that is,
\[
\sum_{i = 1}^d a_i x_i - a_{d+1}(1 - x_{d+1})  = 0,
\]
and except for the north pole, we have
\[
\sum_{i = 1}^d a_i \frac{x_i}{1- x_{d+1}} - a_{d+1}  = 0,
\]
which shows that
\[
\sum_{i = 1}^d a_i X_i - a_{d+1} = 0,
\]
the intersection of the hyperplanes $H$ and $H_{d+1}(0)$
Therefore,  the image of $S^d\cap H$ by $\sigma_N$
is the hyperplane in $\eucreal^d$ which is the intersection 
of $H$ with $H_{d+1}(0)$.

\medskip
We will also prove that $\tau_N$ maps $(d-1)$-spheres 
in $H_{d+1}(0)$ to $(d-1)$-spheres on $S^d$ not passing through the north
pole. Assume that $X\in \eucreal^d$ belongs to the $(d-1)$-sphere
of equation
\[
\sum_{i = 1}^d X_i^2 + \sum_{j = 1}^d a_j X_j +  b = 0.
\]
For any $(X, 0)\in H_{d+1}(0)$, we know that $(x, x_{d+1}) = \tau_N(X, 0)$
is given by
\[
(x, x_{d+1})  =  \left(\frac{2X}{\norme{X}^2 + 1}, 
\frac{\norme{X}^2 - 1}{\norme{X}^2 + 1} \right).
\]
Using the equation of the $(d-1)$-sphere, we get
\[
x = \frac{2X}{-b + 1 - \sum_{j = 1}^d a_j X_j}
\]
and
\[
x_{d+1} = \frac{-b - 1 - \sum_{j = 1}^d a_j X_j}{-b + 1 - \sum_{j = 1}^d a_j X_j}.
\]
Then, we get
\[
\sum_{i = 1}^d a_ix_i = \frac{2\sum_{j = 1}^d a_jX_j}
{-b + 1 - \sum_{j = 1}^d a_j X_j},
\]
which yields
\[
(-b + 1)(\sum_{i = 1}^d a_ix_i) - (\sum_{i = 1}^d a_ix_i)(\sum_{j = 1}^d a_jX_j)
= 2\sum_{j = 1}^d a_jX_j.
\]
From the above, we get
\[
\sum_{i = 1}^d a_iX_i = \frac{(-b + 1)(\sum_{i = 1}^d a_ix_i)}
{\sum_{i = 1}^d a_ix_i + 2}.
\]
Plugging this expression in the formula for $x_{d+1}$ above, we get
\[
x_{d+1} = \frac{-b - 1 - \sum_{i = 1}^d a_ix_i}{-b + 1},
\]
which yields 
\[
\sum_{i = 1}^d a_ix_i + (-b + 1)x_{d+1} + (b + 1) = 0,
\]
the equation of a hyperplane, $H$,  not passing through the north pole.
Therefore, the image of a $(d-1)$-sphere in $H_{d+1}(0)$ is indeed
the intersection, $H\cap S^d$, of $S^d$ with a hyperplane
not passing through $N$, that is, a $(d-1)$-sphere on $S^d$.

\medskip
Given any hyperplane, $H'$,  in $H_{d+1}(0) = \eucreal^d$, say of equation
\[
\sum_{i = 1}^d a_i X_i +  b = 0,
\]
the image of $H'$ under $\tau_N$ is a $(d - 1)$-sphere on $S^d$,
the intersection of $S^d$ with the hyperplane, $H$, passing through $N$
and determined as follows:
For any $(X, 0)\in H_{d+1}(0)$, if $\tau_N(X, 0) = (x, x_{d+1})$, then
\[
X = \frac{x}{1 - x_{d+1}}
\]
and so, $(x, x_{d+1})$ satisfies the equation
\[
\sum_{i = 1}^d a_i x_i +  b(1 - x_{d+1}) = 0,
\]
that is,
\[
\sum_{i = 1}^d a_i x_i - b x_{d+1} +  b = 0,
\]
which is indeed the equation of a hyperplane, $H$,
passing through $N$.
We summarize all this in the following proposition:

\begin{prop}
\label{stereop1}
The stereographic projection, $\mapdef{\sigma_N}{(S^d - \{N\})}{\eucreal^d}$,
induces a bijection, $\sigma_N$, 
between the set of  $(d - 1)$-spheres on $S^d$
and the union of the set of  $(d - 1)$-spheres in $\eucreal^d$
with the set of hyperplanes in $\eucreal^d$;
every $(d - 1)$-sphere on $S^d$ not passing through the
north pole is mapped to a $(d - 1)$-sphere in $\eucreal^d$
and every  $(d - 1)$-sphere on $S^d$ passing through the north pole
is mapped to a hyperplane in $\eucreal^d$.
In fact,  $\sigma_N$  maps the hyperplane 
\[
a_1x_1 + \cdots + a_dx_d + a_{d+1} x_{d+1} + b = 0
\]
not passing through the north pole ($a_{d+1} + b \not= 0$) to the
$(d - 1)$-sphere
\[
\sum_{i = 1}^d X_i^2 + 2\sum_{i = 1}^d \frac{a_i}{a_{d+1} + b} X_i  
- \frac{a_{d+1} - b}{a_{d+1} + b} = 0
\]
and the hyperplane
\[
\sum_{i = 1}^d a_i x_i + a_{d+1}x_{d+1} - a_{d+1} = 0
\]
through the north pole to the hyperplane
\[
\sum_{i = 1}^d a_i X_i - a_{d+1} = 0;
\]
the map $\tau_N = \sigma^{-1}_N$ 
maps the $(d - 1)$-sphere
\[
\sum_{i = 1}^d X_i^2 + \sum_{j = 1}^d a_j X_j +  b = 0
\]
to the hyperplane
\[
\sum_{i = 1}^d a_ix_i + (-b + 1)x_{d+1} + (b + 1) = 0
\]
not passing through the north pole and the hyperplane
\[
\sum_{i = 1}^d a_i X_i +  b = 0
\]
to the hyperplane
\[
\sum_{i = 1}^d a_i x_i - b x_{d+1} +  b = 0
\]
through the north pole.
\end{prop}

\medskip
Proposition \ref{stereop1} raises a natural question:
What do the hyperplanes, $H$, in $\eucreal^{d+1}$ that do not intersect
$S^d$ correspond to, if they correspond to anything at all?

\medskip
The first thing to observe is that the geometric definition
of the stereographic projection and its inverse make it clear that
the hyperplanes associated with $(d-1)$-spheres in $\eucreal^d$
(by $\tau_N$) do intersect $S^d$.
Now, when we write the equation of a $(d-1)$-sphere, $S$, say
\[
\sum_{i = 1}^d X_i^2 + \sum_{i = 1}^d a_iX_i + b = 0
\]
we are implicitly assuming a condition on the $a_i$'s  and $b$
that ensures that $S$ is not the empty sphere, that is,
that its radius, $R$, is positive (or zero). By ``completing the square'',
the above equation can be rewritten as
\[
\sum_{i = 1}^d \left(X_i + \frac{a_i}{2} \right)^2 =  
 \frac{1}{4}\sum_{i = 1}^d a_i^2 - b, 
\]
and so the radius, $R$, of our sphere is given by
\[
R^2 = \frac{1}{4}\sum_{i = 1}^d a_i^2 - b
\]
whereas its center is the point, $c = -\frac{1}{2}(a_1, \ldots, a_d)$. 
Thus, our sphere is a ``real'' sphere of positive radius iff
\[
 \sum_{i = 1}^d a_i^2 > 4b
\]
or a single point,  $c = -\frac{1}{2}(a_1, \ldots, a_d)$,
iff  $\sum_{i = 1}^d a_i^2 = 4b$. 

\medskip
What happens when 
\[
\sum_{i = 1}^d a_i^2 < 4 b\mathrm{?}
\]
In this case, if we allow ``complex points'', that is, if
we consider solutions of our equation
\[
\sum_{i = 1}^d X_i^2 + \sum_{i = 1}^d a_iX_i + b = 0
\]
over $\complex^d$, then we get a ``complex'' sphere of (pure)
imaginary radius, $ \frac{i}{2}\sqrt{4b -\sum_{i = 1}^d a_i^2}$.
The funny thing is that our computations carry over unchanged
and the image of the complex sphere, $S$, 
is still the hyperplane, $H$, given
\[
\sum_{i = 1}^d a_ix_i + (-b + 1)x_{d+1} + (b + 1) = 0.
\]
However, this time, even though $H$ does not have any ``real'' intersection
points with $S^d$, we can show that it
does intersect the ``complex sphere'',
\[
S^d = \{(z_1, \ldots, z_{d+1})\in \complex^{d+1} \mid
z_1^2 + \cdots + z_{d+1}^2 = 1\}
\]
in a nonempty set of points in $\complex^{d+1}$.

\medskip
It follows from all this that $\sigma_N$
and $\tau_N$ establish a bijection between
the set of all hyperplanes in $\eucreal^{d+1}$ minus
the hyperplane, $H_{d+1}$ (of equation $x_{d+1} = 1$),
tangent to $S^d$ at the north pole,  with the
union of four sets:
\begin{enumerate}
\item[(1)]
The set of all (real) $(d - 1)$-spheres of positive radius;
\item[(2)]
The set of all (complex) $(d - 1)$-spheres of imaginary radius;
\item[(3)]
The set of all hyperplanes in $\eucreal^d$;
\item[(4)]
The set of all points of $\eucreal^d$ (viewed as spheres of radius $0$).
\end{enumerate}
Moreover, set (1) corresponds to the hyperplanes that intersect
the interior of $S^d$ and do not pass through the north pole;
set (2) corresponds to the hyperplanes that do not intersect
$S^d$; set (3) corresponds to the hyperplanes that pass through
the north pole minus the tangent hyperplane at the north pole;
and set (4)  corresponds to the hyperplanes that are tangent
to $S^d$, minus the tangent hyperplane at the north pole.

\medskip
It is convenient to add the ``point at infinity'', $\infty$, to
$\eucreal^d$, because then the above bijection
can be extended to map the tangent hyperplane at the north pole
to $\infty$. The union of these four sets (with $\infty$ added)
is called the {\it set of generalized spheres\/}, sometimes,
denoted $\s{S}(\eucreal^d)$. This is a fairly complicated 
space. For one thing, topologically, $\s{S}(\eucreal^d)$
is homeomorphic to the projective space $\projr{d+1}$ with one point removed
(the point corresponding to the ``hyperplane at infinity''),
and this is not a simple space. We can get a slightly more concrete
```picture'' of  $\s{S}(\eucreal^d)$ by looking at the polars
of the hyperplanes w.r.t. $S^d$. Then, the ``real'' spheres
correspond to the points strictly outside $S^d$ which do
not belong to the tangent hyperplane at the norh pole;
the complex spheres correspond to the points in the interior of $S^d$;
the points of $\eucreal^d\cup\{\infty\}$ correspond to the points
on $S^d$;  the hyperplanes in $\eucreal^d$ correspond
to the points in the tangent hyperplane at the norh pole
expect for the north pole. Unfortunately, the poles of hyperplanes
through the origin are undefined. This can be fixed by
embedding $\eucreal^{d+1}$ in its projective completion,
$\projr{d+1}$, but we will not go into this.

\medskip
There are other ways of dealing rigorously with the
set of generalized spheres. One method 
described by Boissonnat \cite{Boissonnat}
is to use the embedding
where the sphere, $S$,  of equation
\[
\sum_{i = 1}^d X_i^2 -2\sum_{i = 1}^d a_iX_i + b = 0
\]
is mapped to the point 
\[
\varphi(S)  = (a_1, \ldots, a_d, b)\in \eucreal^{d+1}.
\]
Now, by a previous computation we know that
\[
b = \sum_{i = 1}^d a_i^2 - R^2,
\] 
where $c = (a_1, \ldots, a_d)$ is the center of $S$ and
$R$ is its radius. The quantity $\sum_{i = 1}^d a_i^2 - R^2$
is known as the {\it power\/} of the origin w.r.t. $S$.  
In general, the {\it power\/} of a point, $X\in \eucreal^d$,
is defined as $\rho(X) = \norme{\libvecbo{c}{X}}^2 - R^2$,
which, after a moment of thought, is just
\[
\rho(X) = \sum_{i = 1}^d X_i^2 -2\sum_{i = 1}^d a_iX_i + b.
\]
Now, since points correspond to spheres of radius $0$,
we see that the image of the point,
$X = (X_1, \ldots, X_d)$, is 
\[
l(X) = (X_1, \ldots, X_d, \sum_{i = 1}^d X_i^2). 
\] 
Thus, in this model, points of $\eucreal^{d}$ are lifted to the
hyperboloid, $\s{P} \subseteq \eucreal^{d+1}$,  of equation
\[
x_{d+1} = \sum_{i = 1}^d x_i^2.
\]
Actually, this method does not deal with hyperplanes
but it is possible to do so. The trick is to consider
equations of a slightly more general form that
capture both spheres and hyperplanes, namely, equations of the form
\[
c\sum_{i = 1}^d X_i^2 + \sum_{i = 1}^d a_iX_i + b = 0.
\]
Indeed, when $c = 0$, we do get a hyperplane! 
Now, to carry out this method we really need to consider
equations up to a nonzero scalars, that is, we consider
the projective space, $\mathbb{P}(\widehat{S}(\eucreal^d))$, 
associated with the vector space,
$\widehat{S}(\eucreal^d)$,
consisting of the above equations. Then, it turns out that 
the quantity
\[
\varrho(a, b, c) = \frac{1}{4}(\sum_{i = 1}^d a_i^2 - 4bc)
\]
(with $a = (a_1, \ldots, a_d)$)
defines a quadratic form on $\widehat{S}(\eucreal^d)$
whose corresponding bilinear form,
\[
\rho((a, b, c), (a', b', c')) = 
\frac{1}{4}(\sum_{i = 1}^d a_ia_i' - 2bc' - 2b'c),
\]
has a natural interpretation 
(with $a = (a_1, \ldots, a_d)$ and $a' = (a_1', \ldots, a_d')$).
Indeed, orthogonality with respect to $\rho$ (that is,
when $\rho((a, b, c), (a', b', c')) = 0$) says that
the corresponding spheres defined by $(a, b, c)$ and $(a', b', c')$
are orthogonal, that the 
corresponding hyperplanes defined by $(a, b, 0)$ and $(a', b', 0)$
are orthogonal, etc. The reader who wants to read more about
this approach should consult Berger (Volume II) 
\cite{Berger90b}.

\medskip
There is a simple relationship between
the lifting onto a hyperboloid and the lifting onto
$S^d$ using the inverse stereographic projection map
because the sphere and the paraboloid are projectively equivalent,
as we showed for $S^2$ in Section \ref{sec5d}. 

\medskip
Recall that the hyperboloid, $\s{P}$, in $\eucreal^{d+1}$ is given by 
the equation
\[
x_{d+1} = \sum_{i = 1}^d x_i^2
\]
and of course, the sphere $S^d$ is given by 
\[
 \sum_{i = 1}^{d+1} x_i^2 = 1.
\]
Consider the ``projective transformation'', $\Theta$,
of $\eucreal^{d+1}$ given by
\begin{eqnarray*}
z_i & = & \frac{x_i}{1 - x_{d+1}},\qquad 1\leq i \leq d \\
z_{d+1} & = & \frac{x_{d+1} + 1 }{1 - x_{d+1}}.
\end{eqnarray*}
Observe that $\Theta$ is undefined on the hyperplane, $H_{d+1}$,
tangent to $S^d$ at the north pole and that its 
first $d$ component are identical to those of the stereographic projection!
Then, we immediately find that
\begin{eqnarray*}
x_i & = & \frac{2z_i}{1 + z_{d+1}},\qquad 1\leq i \leq d \\
x_{d+1} & = & \frac{z_{d+1} - 1 }{1 + z_{d+1}}.
\end{eqnarray*}
Consequently, $\Theta$ is a bijection between $\eucreal^{d+1} - H_{d+1}$
and  $\eucreal^{d+1} - H_{d+1}(-1)$, where $H_{d+1}(-1)$ is the
hyperplane of equation  $x_{d+1} = -1$.

\medskip
The fact that $\Theta$ is undefined on the hyperplane, $H_{d+1}$ is
not a problem as far as mapping the sphere to the paraboloid
because  the north pole is the only point that does have not an image.
However,  later on when we consider the Voronoi polyhedron, $\s{V}(P)$,
of a lifted set of points, $P$,  
we will have more serious problems because in general, such
a polyhedron intersects both hyperplanes $H_{d+1}$ and  $H_{d+1}(-1)$.
This means that $\Theta$ will not be well-defined on the whole
of $\s{V}(P)$ nor will it be surjective on its image. To remedy this difficulty,
we will work with  projective completions. Basically, this amounts
to chasing denominators and homogenizing equations but we also have to be 
careful in dealing with convexity and this is where the
projective polyhedra (studied in Section \ref{sec5e}) will come handy.

\medskip 
So, let us consider the projective sphere, $S^d\subseteq \projr{d+1}$, 
given by the 
equation
\[
 \sum_{i = 1}^{d+1} x_i^2 = x_{d+2}^2
\]
and the paraboloid, $\s{P}\subseteq \projr{d+1}$, given by the  equation
\[
x_{d+1}x_{d+2}  = \sum_{i = 1}^d x_i^2.
\]
Let $\mapdef{\theta}{\projr{d+1}}{\projr{d+1}}$ be the projectivity
induced by the linear map, 
$\mapdef{\widehat{\theta}}{\reals^{d+2}}{\reals^{d+2}}$,
given by
\begin{eqnarray*}
z_i & = & x_i, \qquad\qquad\qquad 1 \leq i \leq d \\
z_{d+1} & = & x_{d+1} + x_{d+2} \\
z_{d+2} & = & x_{d+2} - x_{d+1},
\end{eqnarray*}
whose inverse  is given by
\begin{eqnarray*}
x_i & = & z_i, \qquad\qquad\qquad 1 \leq i \leq d \\
x_{d+1} & = & \frac{z_{d+1} - z_{d+2}}{2} \\
x_{d+2} & = & \frac{z_{d+1} + z_{d+2}}{2}.
\end{eqnarray*}
If we plug these formulae in the equation
of $S^d$, we get
\[
4(\sum_{i = 1}^d z_i^2) + (z_{d+1} - z_{d+2})^2 = (z_{d+1} + z_{d+2})^2,
\]
which simplifies to
\[
z_{d+1}z_{d+2} = \sum_{i = 1}^d z_i^2.
\] 
Therefore, $\theta(S^d) = \s{P}$, that is, $\theta$ maps the sphere
to the hyperboloid. Observe that the north pole, 
$N = (0\co \cdots\co 0\co 1\co 1)$,
is mapped to the point at infinity, $(0\co \cdots\co 0\co 1\co 0)$.

\medskip
The map $\Theta$ is the restriction of $\theta$ to the affine patch,
$U_{d+1}$,  and as such,  it 
can be fruitfully described as the composition
of $\widehat{\theta}$ with a suitable  projection
onto $\eucreal^{d+1}$. For this, as we have done before, 
we identify $\eucreal^{d+1}$ with the hyperplane, $H_{d+2}
\subseteq \eucreal^{d+2}$, of equation $x_{d+2} = 1$ (using the injection,
$\mapdef{i_{d+2}}{\eucreal^{d+1}}{\eucreal^{d+2}}$, where
$\mapdef{i_{j}}{\eucreal^{d+1}}{\eucreal^{d+2}}$ is the injection given by
\[
(x_1, \ldots, x_{d+1})  \mapsto (x_1, \ldots, x_{j-1}, 1, x_{j+1}, 
\ldots, x_{d+1})
\]
for any $(x_1, \ldots, x_{d+1})\in \eucreal^{d+1}$).
For each $i$, with $1\leq i \leq d+2$,
let $\mapdef{\pi_{i}}{(\eucreal^{d+2} - H_{i}(0))}{\eucreal^{d+1}}$ be
the projection of center $0\in \eucreal^{d+2}$ onto the hyperplane,
$H_{i} \subseteq \eucreal^{d+2}$, of equation $x_i = 1$
($H_{i} \cong \eucreal^{d+1}$ and  $H_{i}(0)\subseteq \eucreal^{d+2}$ 
is the hyperplane of equation $x_{i} = 0$) 
given by
\[
\pi_{i}(x_1, \ldots, x_{d+2}) = \left(\frac{x_1}{x_{i}}, \ldots,
\frac{x_{i-1}}{x_{i}}, \frac{x_{i+1}}{x_{i}}, \ldots, 
 \frac{x_{d+2}}{x_{i}}\right) \qquad (x_{i} \not= 0).
\]
Geometrically, for any $x\notin H_{i}(0)$, the image, $\pi_{i}(x)$,
of $x$ is the intersection of the line through the origin and $x$ with
the hyperplane, $H_{i}\subseteq \eucreal^{d+2}$ of equation $x_{i} = 1$. 
Observe that the map,  
$\mapdef{\pi_{i}}{(\eucreal^{d+2} - H_{d+2}(0))}{\eucreal^{d+1}}$,
is an ``affine'' version of the bijection, \\
$\mapdef{\varphi_{i}}{U_{i}}{\reals^{d+1}}$, of Section \ref{sec5d}.
Then, we have
\[
 \Theta =  \pi_{d+2} \circ \widehat{\theta}\circ i_{d+2}.
\]
If we identify $H_{d+2}$ and $\eucreal^{d+1}$, we may
write  with a slight abuse of notation,
$\Theta = \pi_{d+2} \circ \widehat{\theta}$.

\medskip
Besides $\theta$, we need to define a few more  maps
in order to establish the connection between
the Delaunay complex on $S^d$ and the Delaunay complex on $\s{P}$.
We use the convention of denoting the extension to projective spaces 
of a map, $f$, defined between Euclidean spaces,
by $\widetilde{f}$.

\medskip
The Euclidean orthogonal projection, $\mapdef{p_i}{\reals^{d+1}}{\reals^d}$,
is given by
\[
p_i(x_1, \ldots, x_{d+1}) = (x_1, \ldots, x_{i-1}, x_{i+1}, \ldots, x_{d+1})
\]
and   $\mapdef{\widetilde{p}_{i}}{\projr{d+1}}{\projr{d}}$
denotes the  projection
from $\projr{d+1}$ onto $\projr{d}$ given by
\[
\widetilde{p}_{i}(x_1\co \cdots \co x_{d+2}) = 
(x_1 \co \cdots\co x_{i-1}\co x_{i+1} 
\co \cdots \co x_{d+2}),
\]
which is undefined at the point $(0\co \cdots \co 1 \co 0 \co \cdots \co 0)$,
where the ``$1$'' is in the $i^{\mathrm{th}}$ slot. 
The map $\mapdef{\widetilde{\pi}_N}{(\projr{d+1} - \{N\})}{\projr{d}}$ 
is the central projection
from the north pole onto $\projr{d}$ given by
\[
\widetilde{\pi}_N(x_1\co \cdots\co x_{d+1}\co x_{d+2}) = \left(x_1\co \cdots\co
x_d\co x_{d+2} - x_{d+1}\right).
\]
A geometric interpretation of $\widetilde{\pi}_N$ will be needed 
later in certain proofs. If we identify $\projr{d}$ with the
hyperplane, $H_{d+1}(0) \subseteq \projr{d+1}$, of equation $x_{d+1} = 0$, 
then we claim that for any, 
$x\not = N$, the point $\widetilde{\pi}_N(x)$ is the intersection
of the line through $N$ and $x$ with the hyperplane, $H_{d+1}(0)$.
Indeed, parametrically, the line, $\lag N, x\rag$,  through 
$N = (0\co \cdots \co 0 \co 1\co 1)$ and $x$ is given by
\[
\lag N, x\rag = \{(\mu x_1\co \cdots \co \mu x_d \co \lambda + \mu x_{d+1}\co
 \lambda + \mu x_{d+2}) \mid \lambda, \mu\in \reals,\> \lambda \not= 0
\quad\hbox{or}\quad \mu\not= 0\}.
\]
The line $\lag N, x\rag$ intersects the hyperplane $x_{d+1} = 0$ iff
\[
\lambda + \mu x_{d+1} = 0,
\]
so we can pick $\lambda = - x_{d+1}$ and $\mu = 1$,
which yields the intersection point,
\[
(x_1\co \cdots \co x_d \co 0 \co x_{d+2} - x_{d+1}),
\]
as claimed.

\medskip
We also have the projective versions of $\sigma_N$ and $\tau_N$, 
denoted $\mapdef{\widetilde{\sigma}_N}{(S^d - \{N\})}{\projr{d}}$ and 
$\mapdef{\widetilde{\tau}_N}{\projr{d}}{S^d \subseteq \projr{d+1}}$, 
given by:
\[
\widetilde{\sigma}_N(x_1\co \cdots\co x_{d+2}) = 
(x_1\co \cdots \co x_d\co x_{d+2} - x_{d+1})
\]
and
\[
\widetilde{\tau}_N(x_1\co \cdots\co x_{d+1}) = 
\left(2x_1x_{d+1}\co \cdots \co 2x_dx_{d+1}\co
\sum_{i = 1}^d x_i^2 - x_{d+1}^2 \co  \sum_{i = 1}^d x_i^2 + x_{d+1}^2\right).
\]
It is an easy exercise to check that the image of $S^d - \{N\}$
by $\widetilde{\sigma}_N$ is $U_{d+1}$ and that $\widetilde{\sigma}_N$ and 
$\widetilde{\tau}_N\res U_{d+1}$
are mutual inverses.
Observe that $\widetilde{\sigma}_N = \widetilde{\pi}_N\res S^d$, 
the restriction of the projection,
$\widetilde{\pi}_N$, to the sphere, $S^d$.
The lifting,
$\mapdef{\widetilde{l}}{\eucreal^d}{\s{P}} \subseteq \projr{d+1}$, is given by
\[
\widetilde{l}(x_1, \ldots, x_d) = 
\left(x_1\co \cdots\co x_d\co \sum_{i = 1}^d x_i^2\co 1\right) 
\] 
and the embedding, $\mapdef{\psi_{d+1}}{\eucreal^d}{\projr{d}}$,  
(the map  $\psi_{d+1}$ defined in
Section \ref{sec5d}) is given by
\[
\psi_{d+1}(x_1, \ldots, x_d) = (x_1\co\cdots \co x_d\co 1).
\] 
Then, we easily check 

\begin{prop}
\label{Thetaprop1}
The maps, $\theta, \widetilde{\pi}_N, \widetilde{\tau}_N, 
\widetilde{p}_{d+1}, \widetilde{l}$ and $\psi_{d+1}$ defined
before satisfy the equations
\begin{eqnarray*}
\widetilde{l} & = & \theta\circ \widetilde{\tau}_N\circ \psi_{d+1} \\
\widetilde{\pi}_N & =  & \widetilde{p}_{d+1} \circ \theta \\
\widetilde{\tau}_N\circ \psi_{d+1} & = & \psi_{d+2}\circ \tau_N \\
\widetilde{l} & = & \psi_{d+2}\circ l \\
l & = & \Theta\circ \tau_N.
\end{eqnarray*}
\end{prop}

\proof
Let us check the first equation leaving the others as an exercise.
Recall that $\theta$ is given by
\[
\theta(x_1\co \cdots\co x_{d+2}) = 
(x_1\co \cdots\co x_{d}\co x_{d+1} + x_{d+2}\co  x_{d+2} - x_{d+1}).
\]
Then, as
\[
\widetilde{\tau}_N\circ \psi_{d+1}(x_1, \ldots, x_{d}) =
\left(2x_1\co \cdots \co 2x_d\co
\sum_{i = 1}^d x_i^2 - 1 \co  \sum_{i = 1}^d x_i^2 + 1\right),
\]
we get
\begin{eqnarray*}
 \theta\circ \widetilde{\tau}_N\circ \psi_{d+1}(x_1, \ldots, x_{d}) 
& = &
\left(2x_1\co \cdots \co 2x_d\co
2\sum_{i = 1}^d x_i^2  \co  2\right) \\
& = &
\left(x_1\co \cdots \co x_d\co
\sum_{i = 1}^d x_i^2  \co  1\right) = \widetilde{l}(x_1, \ldots, x_{d}),
\end{eqnarray*}
as claimed. 
$\bigsquare$

\medskip
We will also need some properties of the  projection $\pi_{d+2}$
and of $\Theta$ and for this, let
\[
\mathbb{H}^d_+ =  \{(x_1, \ldots, x_d) \in \eucreal^d \mid x_d > 0\} 
\quad\hbox{and}\quad
\mathbb{H}^d_- =  \{(x_1, \ldots, x_d) \in \eucreal^d \mid x_d < 0\}.
\]

\begin{prop}
\label{Thetaprop2}
The projection, $\pi_{d+2}$, has the following properties:
\begin{enumerate}
\item[(1)]
For every hyperplane, $H$, through the origin, 
$\pi_{d+2}(H)$ is a hyperplane in $H_{d+2}$.
\item[(2)]
Given any set of points,
$\{a_1, \ldots, a_n\}\subseteq \eucreal^{d+2}$, 
if $\{a_1, \ldots, a_n\}$ is contained in the open half-space
above the hyperplane $x_{d+2} = 0$ or 
$\{a_1, \ldots, a_n\}$ is contained in the open half-space
below the hyperplane $x_{d+2} = 0$, 
then the image by $\pi_{d+2}$
of the convex hull of the $a_i$'s is the convex hull of the
images of these points, that is,
\[
\pi_{d+2}(\mathrm{conv}(\{a_1, \ldots, a_n\})) = 
\mathrm{conv}(\{\pi_{d+2}(a_1), \ldots, \pi_{d+2}(a_n)\}).
\]
\item[(3)]
Given any set of points,
$\{a_1, \ldots, a_n\}\subseteq \eucreal^{d+1}$, 
if $\{a_1, \ldots, a_n\}$ is  contained in the open half-space
above the hyperplane $H_{d+1}$
or $\{a_1, \ldots, a_n\}$ is contained in the open  half-space 
below $H_{d+1}$, then 
\[
\Theta(\mathrm{conv}(\{a_1, \ldots, a_n\})) = 
\mathrm{conv}(\{\Theta(a_1), \ldots, \Theta(a_n)\}).
\]
\item[(4)]
For any  set $S\subseteq \eucreal^{d+1}$, if $\mathrm{conv}(S)$
does not intersect $H_{d+1}$, then 
\[
\Theta(\mathrm{conv}(S)) = \mathrm{conv}(\Theta(S)). 
\]
\end{enumerate}
\end{prop}

\proof
(1)
The image, $\pi_{d+2}(H)$, of a hyperplane, $H$, through the origin
is the intersection of $H$ with
$H_{d+2}$, which is a hyperplane in $H_{d+2}$.

\medskip
(2)
This seems fairly clear geometrically but the result fails
for arbitrary sets of points so to be on the safe side
we give an algebraic proof. We will prove the following two facts
by induction on $n\geq 1$:
\begin{enumerate}
\item[(1)]
For  all $\lambda_1, \ldots, \lambda_n\in \reals$ with
$\lambda_1 + \cdots + \lambda_n = 1$ and $\lambda_i \geq 0$, for all 
$a_1, \ldots, a_n\in \mathbb{H}_+^{d+2}$ (resp. $\in \mathbb{H}_-^{d+2}$)
there exist some $\mu_1,\ldots, \mu_n\in \reals$ with
$\mu_1 + \cdots + \mu_n = 1$ and $\mu_i\geq 0$, so that
\[
\pi_{d+2}(\lambda_1 a_1 + \cdots + \lambda_n a_n) =
\mu_1\pi_{d+2}(a_1) + \cdots + \mu_n\pi_{d+2}(a_n).
\]
\item[(2)]
For  all $\mu_1, \ldots, \mu_n\in \reals$ with
$\mu_1 + \cdots + \mu_n = 1$ and $\mu_i\geq 0$, for all
$a_1, \ldots, a_n\in \mathbb{H}_+^{d+2}$ (resp. $\in \mathbb{H}_-^{d+2}$)
there exist some $\lambda_1,\ldots, \lambda_n\in \reals$ with
$\lambda_1 + \cdots + \lambda_n = 1$ and $\lambda_i \geq 0$, so that
\[
\pi_{d+2}(\lambda_1 a_1 + \cdots + \lambda_n a_n) =
\mu_1\pi_{d+2}(a_1) + \cdots + \mu_n\pi_{d+2}(a_n).
\]
\end{enumerate}

\medskip
(1) The base case is clear. Let us assume for the moment that
we proved (1) for $n = 2$ and consider the induction step for
$n \geq 2$. Since $\lambda_1 + \cdots + \lambda_{n+1} = 1$
and $n\geq 2$, there is some $i$ such that $\lambda_i \not= 1$, 
and without loss of generality, say $\lambda_1 \not= 1$. Then, we can write
\[
\lambda_1 a_1 + \cdots + \lambda_{n+1} a_{n+1} = 
\lambda_1 a_1 +  (1 - \lambda_1)
\left(
\frac{\lambda_2}{1 - \lambda_1}  a_2 + \cdots + 
\frac{\lambda_{n+1}}{1 - \lambda_1} a_{n+1} 
\right)
\]
and since  $\lambda_1 + \lambda_2 + \cdots + \lambda_{n+1} = 1$,
we have
\[
\frac{\lambda_2}{1 - \lambda_1} + \cdots + \frac{\lambda_{n+1}}{1 - \lambda_1}
= 1.
\]
By the induction hypothesis, for $n = 2$, there exist $\alpha_1$
with $0 \leq \alpha_1 \leq 1$, 
such that
\begin{multline*}
\pi_{d+2}(\lambda_1 a_1 + \cdots + \lambda_{n+1} a_{n+1})  = 
\pi_{d+2}\left(
\lambda_1 a_1 +  (1 - \lambda_1)
\left(
\frac{\lambda_2}{1 - \lambda_1} a_2 + \cdots + 
\frac{\lambda_{n+1}}{1 - \lambda_1} a_{n+1}
\right) \right) \\
 =  
(1 - \alpha_1) \pi_{d+2}(a_1) + \alpha_1 \pi_{d+2}\left(
\frac{\lambda_2}{1 - \lambda_1} a_2 + \cdots + 
\frac{\lambda_{n+1}}{1 - \lambda_1} a_{n+1}
\right) 
\end{multline*}
Again, by the induction hypothesis (for $n$), there exist
$\beta_2, \ldots, \beta_{n+1}$ with 
$\beta_2 + \cdots +  \beta_{n+1} = 1$ and $\beta_i \geq 0$, so that
\[
\pi_{d+2}\left(
\frac{\lambda_2}{1 - \lambda_1} a_2 + \cdots + 
\frac{\lambda_{n+1}}{1 - \lambda_1} a_{n+1}
\right)  = 
\beta_2 \pi_{d+2}(a_2) + \cdots + \beta_{n+1} \pi_{d+2}(a_{n+1}), 
\]
so we get
\begin{eqnarray*}
\pi_{d+2}(\lambda_1 a_1 + \cdots + \lambda_{n+1} a_{n+1}) &  = & 
(1 - \alpha_1) \pi_{d+2}(a_1) + \alpha_1 
(\beta_2 \pi_{d+2}(a_2) + \cdots + \beta_{n+1} \pi_{d+2}(a_{n+1})) \\
& = & 
(1 - \alpha_1) \pi_{d+2}(a_1) + 
\alpha_1 \beta_2 \pi_{d+2}(a_2) + \cdots + 
\alpha_1 \beta_{n+1} \pi_{d+2}(a_{n+1}) 
\end{eqnarray*}
and clearly, $1 - \alpha_1 + \alpha_1 \beta_2 + \cdots + \alpha_1 \beta_{n+1}
= 1$ as $\beta_2 + \cdots +  \beta_{n+1} = 1$;
$1 - \alpha_1 \geq 0$; and  $\alpha_1 \beta_i \geq 0$, as
$0 \leq \alpha_1 \leq 1$ and  $\beta_i \geq 0$.
This establishes the induction step and thus, all is left is to prove
the case $n = 2$.

\medskip
(2) 
The base case $n = 1$ is also clear. As in (1), let us assume for a
moment that (2) is proved for $n = 2$ and consider the induction step.
The proof is quite similar to that of (1) but this time,
we may assume that $\mu_1\not = 1$ and we write
\begin{multline*}
\mu_1 \pi_{d+2}(a_1) + \cdots + \mu_{n+1} \pi_{d+2}(a_{n+1})  \\
= 
\mu_1 \pi_{d+2}(a_1) +  (1 - \mu_1)
\left(
\frac{\mu_2}{1 - \mu_1} \pi_{d+2}(a_2) \cdots + 
\frac{\mu_{n+1}}{1 - \mu_1} \pi_{d+2}(a_{n+1}) 
\right).
\end{multline*}
By the induction hypothesis, there are some $\alpha_2, \ldots, \alpha_{n+1}$
with  $\alpha_2 +  \cdots + \alpha_{n+1} = 1$ and $\alpha_i \geq 0$
such that
\[
\pi_{d+2}(\alpha_2 a_2 + \cdots + \alpha_{n+1} a_{n+1}) =
\frac{\mu_2}{1 - \mu_1} \pi_{d+2}(a_2) + \cdots + 
\frac{\mu_{n+1}}{1 - \mu_1} \pi_{d+2}(a_{n+1}) .
 \]
By the induction hypothesis for $n = 2$, there is some
$\beta_1$ with $0 \leq \beta_1 \leq 1$, so that
\[
\pi_{d+2}((1 - \beta_1) a_1 + 
\beta_1(\alpha_2 a_2 + \cdots + \alpha_{n+1} a_{n+1})) = 
\mu_1 \pi_{d+2}(a_1) +  (1 - \mu_1)
\pi_{d+2}(\alpha_2 a_2 + \cdots + \alpha_{n+1} a_{n+1}), 
\]
which establishes the induction hypothesis.
Therefore, all that remains is to prove (1) and (2) for $n = 2$.

\medskip
As $\pi_{d+2}$ is given by
\[
\pi_{d+2}(x_1, \ldots, x_{d+2}) = \left(\frac{x_1}{x_{d+2}}, \ldots, 
\frac{x_{d+1}}{x_{d+2}}\right) \qquad (x_{d+2} \not= 0)
\]
it is enough to treat the case when $d = 0$, that is,
\[
\pi_2(a, b) = \frac{a}{b}.
\]

\medskip
To prove (1) it is enough to show that for any $\lambda$,
with $0 \leq \lambda \leq 1$, if $b_1b_2 > 0$ then
\[
\frac{a_1}{b_1} \leq 
\frac{(1 - \lambda)a_1 + \lambda a_2}{(1 - \lambda)b_1 + \lambda b_2}
\leq \frac{a_2}{b_2} \quad \hbox{if}\quad 
\frac{a_1}{b_1} \leq \frac{a_2}{b_2} 
\]
and 
\[
\frac{a_2}{b_2} \leq 
\frac{(1 - \lambda)a_1 + \lambda a_2}{(1 - \lambda)b_1 + \lambda b_2}
\leq \frac{a_1}{b_1} \quad \hbox{if}\quad 
\frac{a_2}{b_2} \leq \frac{a_1}{b_1}, 
\]
where, of course $(1 - \lambda)b_1 + \lambda b_2 \not= 0$.
For this, we compute (leaving some steps as an exercise)
\[
\frac{(1 - \lambda)a_1 + \lambda a_2}{(1 - \lambda)b_1 + \lambda b_2} 
- \frac{a_1}{b_1} =
\frac{\lambda(a_2b_1 - a_1b_2)}{((1 - \lambda)b_1 + \lambda b_2)b_1}
\]
and 
\[
\frac{(1 - \lambda)a_1 + \lambda a_2}{(1 - \lambda)b_1 + \lambda b_2} 
- \frac{a_2}{b_2} =
-\frac{(1 - \lambda)(a_2b_1 - a_1b_2)}{((1 - \lambda)b_1 + \lambda b_2)b_2}.
\]
Now, as $b_1b_2 > 0$, that is, $b_1$ and $b_2$ have the same sign and as
$0 \leq \lambda \leq 1$, we have both
$((1 - \lambda)b_1 + \lambda b_2)b_1 > 0$ and 
$((1 - \lambda)b_1 + \lambda b_2)b_2 > 0$. Then, if
$a_2b_1 - a_1b_2 \geq  0$, that is
$\frac{a_1}{b_1} \leq \frac{a_2}{b_2}$ (since  $b_1b_2 > 0$),
the first two inequalities holds and if 
$a_2b_1 - a_1b_2 \leq  0$, that is
$\frac{a_2}{b_2} \leq \frac{a_1}{b_1}$ (since  $b_1b_2 > 0$),
the last two inequalities holds. This proves (1).

\medskip
In order to prove (2), given any $\mu$, with $0 \leq \mu \leq 1$,
if $b_1b_2 > 0$, we show that we can find $\lambda$
with $0 \leq \lambda \leq 1$, so that
\[
(1 - \mu) \frac{a_1}{b_1} + \mu \frac{a_2}{b_2} =
\frac{(1 - \lambda)a_1 + \lambda a_2}{(1 - \lambda)b_1 + \lambda b_2}. 
\]
If we let 
\[
\alpha = (1 - \mu) \frac{a_1}{b_1} + \mu \frac{a_2}{b_2},
\]
we find that $\lambda$ is given by the equation
\[
\lambda(a_2 - a_1 + \alpha (b_1 - b_2)) = \alpha b_1 - a_1.
\]
After some (tedious) computations (check for yourself!) we find:
\begin{eqnarray*}
a_2 - a_1 + \alpha (b_1 - b_2) & = & 
\frac{((1 - \mu)b_2 + \mu b_1)(a_2b_1 - a_1b_2)}{b_1b_2} \\
\alpha b_1 - a_1 & = & \frac{\mu b_1(a_2b_1 - a_1b_2)}{b_1b_2}.
\end{eqnarray*}
If $a_2b_1 - a_1b_2 = 0$, then  $\frac{a_1}{b_1} = \frac{a_2}{b_2}$
and $\lambda = 0$ works. If  $a_2b_1 - a_1b_2 \not= 0$, then
\[
\lambda = \frac{\mu b_1}{(1 - \mu)b_2 + \mu b_1} =
\frac{\mu}{(1 - \mu)\frac{b_2}{b_1} + \mu}.
\]
Since $b_1b_2 > 0$, we have $\frac{b_2}{b_1} > 0$, and since
$0 \leq \mu \leq 1$, we conclude that $0 \leq \lambda \leq 1$,
which proves (2).


\medskip
(3)
Since
\[
\Theta =  \pi_{d+2} \circ \widehat{\theta}\circ i_{d+2},
\]
as $i_{d+2}$ and $\widehat{\theta}$ are linear,
they preserve convex hulls, so by (2), 
we simply have to show that either
$\widehat{\theta}\circ i_{d+2}(\{a_1, \ldots, a_n\})$
is strictly below the hyperplane, $x_{d+2} = 0$, or
strictly above it.
But, 
\[
\widehat{\theta}(x_1, \ldots, x_{d+2})_{d+2} = 
x_{d+2} - x_{d+1}
\]
and $i_{d+2}(x_1, \ldots, x_{d+1}) = (x_1, \ldots, x_{d+1}, 1)$, so
\[
(\widehat{\theta}\circ i_{d+2})(x_1, \ldots, x_{d+1})_{d+2} = 
1 - x_{d+1},
\]
and this quantity is positive iff $x_{d+1} < 1$, negative
iff  $x_{d+1} > 1$; that  is, either
all the points $a_i$ are strictly below the hyperplane
$H_{d+1}$ or all strictly above it.

\medskip
(4) This follows immediately from (4) as $\mathrm{conv}(S)$
consists of all finite convex combinations of points in $S$.
$\bigsquare$

\danger
If a set, $\{a_1, \ldots, a_n\}\subseteq \eucreal^{d+2}$, contains
points on {\it both sides\/} of the hyperplane, $x_{d+2} = 0$, then 
$\pi_{d+2}(\mathrm{conv}(\{a_1, \ldots, a_n\}))$ is {\bf not}
necessarily convex (find such an example!).

\section[Stereographic Projection and Delaunay Polytopes]
{Stereographic Projection, Delaunay Polytopes and Voronoi Polyhedra}
\label{sec13}
We saw in an earlier section that lifting a set 
of points, $P \subseteq \eucreal^d$,
to the paraboloid, $\s{P}$,
{\it via\/} the lifting function, $l$, was fruitful
to better understand Voronoi diagrams and Delaunay triangulations.
As far as we know, 
Edelsbrunner and Seidel \cite{Eldel} were the first
to find the relationship between Voronoi diagrams
and the polar dual of the convex hull of 
a lifted set of points onto a paraboloid.
This connection is described in Note 3.1 of Section 3 in
\cite{Eldel}. The connection between the Delaunay
triangulation and the convex hull of the lifted set of
points is described in Note 3.2 of the same paper.
Polar duality is not  mentioned and seems
to enter the scene only with
Boissonnat and Yvinec \cite{Boissonnat}.

\medskip
It turns out that instead of using a paraboloid we can use
a sphere and instead of the lifting function $l$ we can use the
composition of $\psi_{d+1}$ with the
inverse stereographic projection, $\widetilde{\tau}_N$.
Then, to get back down
to $\eucreal^d$, we use the composition of 
the projection, $\widetilde{\pi}_N$, with $\pi_{d+1}$,
instead of the orthogonal projection, $p_{d+1}$.

\medskip
However, we have to be a bit careful because $\Theta$
does map all convex polyhedra to convex polyhedra. Indeed,
$\Theta$ is the composition of $\pi_{d+2}$ with some linear maps,
but $\pi_{d+2}$ does not behave well with 
respect to arbitrary convex sets.
In particular, $\Theta$ is not well-defined on any face
that intersects the hyperplane $H_{d+1}$ (of equation $x_{d+1} = 1$).
Fortunately, we can circumvent these difficulties by using 
the concept of a projective polyhedron introduced in
Chapter \ref{chap2b}.

\medskip
As we said in the previous section, the correspondence between Voronoi
diagrams and convex hulls {\it via\/} inversion 
was first observed by Brown \cite{Brown}.
Brown takes a set of points, $S$, for simplicity
assumed to be in the plane, first lifts these points
to the unit sphere $S^2$ using inverse stereographic
projection (which is equivalent to an inversion 
of power $2$ centered at the north pole),
getting $\tau_N(S)$, and then takes
the convex hull, $\s{D}(S) = \mathrm{conv}(\tau_N(S))$, of the lifted set.
Now, in order to obtain the Voronoi diagram of $S$, 
apply our inversion (of power $2$ centered at the north pole)
to each of the faces of $\mathrm{conv}(\tau_N(S))$,
obtaining spheres passing through the center of $S^2$
and then intersect these spheres with the plane
containing $S$, obtaining circles.  The centers of some
of these circles are
the Voronoi vertices. Finally,  a simple
criterion can be used to  retain the
``nearest Voronoi points'' and to connect up
these vertices.

\medskip
Note that Brown's  method is {\it not\/} the method that uses the
polar dual of the polyhedron  $\s{D}(S) = \mathrm{conv}(\tau_N(S))$,
as we might have expected from the lifting method
using a paraboloid. In fact, it is more natural
to get the {\it Delaunay triangulation of $S$\/} from Brown's method,
by applying the  stereographic projection
(from the north pole) to $\s{D}(S)$,
as we will prove below. As $\s{D}(S)$ is strictly below
the plane $z = 1$, there are no problems.
Now, in order to get the Voronoi diagram, we take the polar dual,
$\s{D}(S)^*$, of $\s{D}(S)$ and then apply the 
central projection w.r.t. the north pole.
This is where problems arise, as some faces of $\s{D}(S)^*$
may intersect the hyperplane $H_{d+1}$ and this is why
we have recourse to projective geometry.

\medskip
First, we show that $\theta$ has a good behavior with respect to tangent
spaces. Recall from Section \ref{sec5e}
that for any point, 
$a = (a_1\co \cdots\co a_{d+2})\in \projr{d+1}$,
the tangent hyperplane, $T_a S^d$, to the sphere $S^d$ at $a$ is given by
the equation 
\[
\sum_{i = 1}^{d+1} a_ix_i - a_{d+2} x_{d+2}  = 0.
\]
Similarly, the tangent hyperplane, $T_a \s{P}$, to the paraboloid
$\s{P}$ at $a$ is given by the equation 
\[
2\sum_{i = 1}^d a_i x_i  - a_{d+2} x_{d+1} - a_{d+1}x_{d+2} = 0.
\]
If we lift a point $a\in \eucreal^d$ to $S^d$ by 
$\widetilde{\tau}_N\circ \psi_{d+1}$
and to $\s{P}$ by $\widetilde{l}$, it turns out that the image of
the tangent hyperplane to $S^d$ at 
$\widetilde{\tau}_N\circ \psi_{d+1}(a)$ by $\theta$
is the tangent hyperplane to $\s{P}$ at $\widetilde{l}(a)$.

\begin{prop}
\label{Thetatang}
The map $\theta$ has the following properties:
\begin{enumerate}
\item[(1)]
For any point, $a = (a_1, \ldots, a_d) \in \eucreal^d$, we have
\[
\theta(T_{\widetilde{\tau}_N\circ \psi_{d+1}(a)} S^d) = T_{\widetilde{l}(a)} \s{P},
\]
that is, $\theta$ preserves tangent hyperplanes.
\item[(2)]
For every $(d - 1)$-sphere, $S\subseteq \eucreal^d$, we have
\[
\theta(\widetilde{\tau}_N\circ \psi_{d+1}(S)) = \widetilde{l}(S),
\]
that is, $\theta$ preserves lifted $(d-1)$-spheres.
\end{enumerate}
\end{prop}

\proof
(1)
By Proposition \ref{Thetaprop1}, we know that
\[
\widetilde{l} = \theta\circ \widetilde{\tau}_N\circ \psi_{d+1}
\]
and we proved in Section \ref{sec5e} that projectivities
preserve tangent spaces. Thus,
\[
\theta(T_{\widetilde{\tau}_N\circ \psi_{d+1}(a)} S^d) = 
T_{\theta\circ\widetilde{\tau}_N\circ \psi_{d+1}(a)} \theta(S^d)
= T_{\widetilde{l}(a)} \s{P},
\]
as claimed.

\medskip
(2)
This follows immediately from the equation
$\widetilde{l} = \theta\circ \widetilde{\tau}_N\circ \psi_{d+1}$.
$\bigsquare$

\medskip
Given any two distinct points, 
$a = (a_1, \ldots, a_d)$ and  $b = (b_1, \ldots, b_d)$
in $\eucreal^d$, recall that the bisector hyperplane, $H_{a, b}$,
of $a$ and $b$ is given by
\[
(b_1 - a_1)x_1  + \cdots +  (b_d - a_d)x_d =
(b_1^2 + \cdots + b_d^2)/2 - (a_1^2+  \cdots + a_d^2)/2.
\]
We have the following useful proposition:

\begin{prop}
\label{pimap1}
Given any two distinct
points, $a = (a_1, \ldots, a_d)$ and  $b = (b_1, \ldots, b_d)$
in $\eucreal^d$, the image under the projection, $\widetilde{\pi}_N$,
of the intersection,
$T_{\widetilde{\tau}_N\circ \psi_{d+1}(a)} S^d\cap 
T_{\widetilde{\tau}_N\circ \psi_{d+1}(b)} S^d$, of the tangent 
hyperplanes at the
lifted points $\widetilde{\tau}_N\circ \psi_{d+1}(a)$ and 
$\widetilde{\tau}_N\circ \psi_{d+1}(b)$ on the 
sphere $S^d\subseteq \projr{d+1}$
is the embedding of the
bisector hyperplane, $H_{a, b}$, of $a$ and $b$, into $\projr{d}$,
that is,
\[
\widetilde{\pi}_N(T_{\widetilde{\tau}_N\circ \psi_{d+1}(a)} S^d
\cap T_{\widetilde{\tau}_N\circ \psi_{d+1}(b)} S^d) = \psi_{d+1}(H_{a, b}).
\]
\end{prop}

\proof
In view of the geometric interpretation of $\widetilde{\pi}_N$ given earlier,
we need to find the equation of the hyperplane, $H$, passing through
the intersection of the tangent hyperplanes, 
$T_{\widetilde{\tau}_N\circ \psi_{d+1}(a)}$ 
and $T_{\widetilde{\tau}_N\circ \psi_{d+1}(b)}$ 
and passing through the north pole and then, it is geometrically
obvious that 
\[
\widetilde{\pi}_N(T_{\widetilde{\tau}_N\circ \psi_{d+1}(a)} S^d
\cap T_{\widetilde{\tau}_N\circ \psi_{d+1}(b)} S^d) = 
H\cap H_{d+1}(0),
\]
where $H_{d+1}(0)$ is the hyperplane (in $\projr{d+1}$)
of equation $x_{d+1} = 0$.
Recall that $T_{\widetilde{\tau}_N\circ \psi_{d+1}(a)} S^d$ and 
$T_{\widetilde{\tau}_N\circ \psi_{d+1}(b)} S^d$ are given by
\[
E_1 = 2 \sum_{i = 1}^d a_i x_i + (\sum_{i = 1}^d a_i^2 - 1)x_{d+1} 
- (\sum_{i = 1}^d a_i^2 + 1) x_{d+2} = 0
\]
and
\[
E_2 = 2 \sum_{i = 1}^d b_i x_i + (\sum_{i = 1}^d b_i^2 - 1)x_{d+1} 
- (\sum_{i = 1}^d b_i^2 + 1) x_{d+2} = 0.
\]
The hyperplanes passing through 
$T_{\widetilde{\tau}_N\circ \psi_{d+1}(a)} S^d
\cap T_{\widetilde{\tau}_N\circ \psi_{d+1}(b)} S^d$ are
given by an equation of the form
\[
\lambda E_1 + \mu E_2 = 0,
\]
with $\lambda, \mu \in \reals$. Furthermore, in order to
contain the north pole, this equation must vanish for
$x = (0\co \cdots\co 0\co 1\co 1)$. But, observe that 
setting $\lambda = -1$
and $\mu = 1$ gives a solution since the corresponding equation is
\[
2 \sum_{i = 1}^d (b_i - a_i) x_i + 
(\sum_{i = 1}^d b_i^2 - \sum_{i = 1}^d a_i^2)x_{d+1} 
- (\sum_{i = 1}^d b_i^2 - \sum_{i = 1}^d a_i^2) x_{d+2} = 0
\]
and it vanishes on $(0\co \cdots\co 0\co 1\co 1)$. But then, the 
intersection of $H$ with the hyperplane $x_{d+1} = 0$ is given by
\[
2 \sum_{i = 1}^d (b_i - a_i) x_i  
- (\sum_{i = 1}^d b_i^2 - \sum_{i = 1}^d a_i^2)x_{d+2} = 0,
\]
which is equivalent to the equation  of $\psi_{d+1}(H_{a, b})$
(except that $x_{d+2}$ is replaced by $x_{d+1}$).
$\bigsquare$

\medskip
In order to define precisely Delaunay complexes as projections of
objects obtained by deleting some faces from a projective polyhedron
we need to define the notion of ``projective (polyhedral) complex''.
However, this is easily done by defining the notion  of
cell complex where the cells are polyhedral cones.
Such objects are known as {\it fans\/}. The definition below
is basically Definition \ref{complexdefw} in which the cells are cones
as opposed to polytopes.

\begin{defin}
\label{fandef}
{\em
A {\it fan in $\eucreal^\mdeg$\/}  is a set, $K$,
consisting of a (finite or infinite) set of 
polyhedral cones in $\eucreal^\mdeg$
satisfying the following conditions:
\begin{enumerate}
\item[(1)] 
Every face of a cone in $K$ also belongs to $K$.
\item[(2)] 
For any two cones $\sigma_1$ and $\sigma_2$ in $K$,
if $\sigma_1\cap\sigma_2\not=\emptyset$, then $\sigma_1\cap\sigma_2$
is a common face of both $\sigma_1$ and $\sigma_2$.
\end{enumerate}
Every cone, $\sigma\in K$, of dimension $k$,  
is called a {\it $k$-face\/} (or {\it face\/}) of $K$.
A $0$-face $\{v\}$ is called a {\it  vertex\/} and a $1$-face
is called an {\it edge\/}.
The {\it dimension\/} of the fan $K$ is the
maximum of the  dimensions of all cones in $K$.
If $\mathrm{dim}\, K = d$, then every face of dimension $d$ is
called a {\it cell\/} and every face of dimension $d - 1$
is called a {\it facet\/}.

\medskip
A {\it projective (polyhedral) complex\/}, $\s{K}\subseteq \projr{d}$, 
is a set of projective polyhedra of the form, \\
$\{\mathbb{P}(C) \mid C\in K\}$, where $K\subseteq \reals^{d+1}$ is a fan. 
}
\end{defin}

\medskip
Given a projective  complex, the notions
of face, vertex, edge, cell, facet, are dedined in the obvious way. 

\medskip
If $K\subseteq \reals^d$ is a polyhedral complex, then it is easy to 
check that the set \\ 
$\{C(\sigma) \mid \sigma \in K\}\subseteq \reals^{d+1}$ is a fan 
and we get the projective complex
\[
\widetilde{K} = \{\mathbb{P}(C(\sigma)) \mid \sigma \in K\}
\subseteq \projr{d}.
\]
The projective complex, $\widetilde{K}$,
is called  the {\it projective completion\/} of $K$.
Also, it is easy to check that if $\mapdef{f}{P}{P'}$
is an injective affine map between two polyhedra $P$ and $P'$, 
then $f$ extends uniquely to a projectivity, 
$\mapdef{\widetilde{f}}{\widetilde{P}}{\widetilde{P'}}$,
between the projective completions of $P$ and $P'$.

\medskip
We now have all the facts needed to show that
Delaunay triangulations and Voronoi diagrams can be defined
in terms of the lifting, $\widetilde{\tau}_N\circ \psi_{d+1}$, 
and the projection,
$\widetilde{\pi}_N$, and to establish their duality {\it via\/} polar duality
with respect to $S^d$. 

\begin{defin}
\label{Delaunaystereo}
{\em 
Given any set of points, $P = \{p_1, \ldots, p_n\} \subseteq 
\eucreal^d$, the polytope, 
$\s{D}(P) \subseteq \reals^{d+1}$,
called the {\it  Delaunay polytope\/} associated with $P$
is the convex hull of the union of the lifting of
the points of $P$ onto the sphere $S^d$ ({\it via\/} inverse
stereographic projection) with the north pole, that is,
$\s{D}(P) = \mathrm{conv}(\tau_N(P)\cup \{N\})$.
The  {\it projective Delaunay polytope\/}, 
$\widetilde{\s{D}}(P) \subseteq \projr{d+1}$,  associated with $P$ 
is the projective completion of $\s{D}(P)$. 
The polyhedral complex,  $\s{C}(P) \subseteq \reals^{d+1}$,   
called the {\it lifted Delaunay complex of $P$\/} 
is the complex obtained from $\s{D}(P)$ by deleting
the facets containing the north pole (and their faces)
and $\widetilde{\s{C}}(P) \subseteq \projr{d+1}$ is the projective   
completion of $\s{C}(P)$. 
The polyhedral complex, 
$\s{D}\mathit{el}(P) = \pi_{d+1}\circ\widetilde{\pi}_N(\widetilde{\s{C}}(P))
\subseteq \eucreal^d$, is the {\it Delaunay complex of $P$\/}
or {\it Delaunay triangulation of $P$\/}.
}
\end{defin}

\medskip
The above is not the ``standard'' definition of the Delaunay
triangulation of $P$ but it is equivalent to the definition, say
given in Boissonnat and Yvinec \cite{Boissonnat},
as we will prove shortly.
It also has certain advantages over lifting onto a paraboloid,
as we will explain.

\medskip
It it possible and useful to define $\s{D}\mathit{el}(P)$ more directly 
in terms of $\s{C}(P)$.
The projection, $\mapdef{\widetilde{\pi_N}}{(\projr{d+1} -\{N\})}{\projr{d}}$,
comes from the linear map,
$\mapdef{\widehat{\pi}_N}{\reals^{d+2}}{\reals^{d+1}}$, given by
\[
\widehat{\pi}_N(x_1, \ldots, x_{d+1}, x_{d+2}) = 
(x_1, \ldots, x_d, x_{d+2} - x_{d+1}).
\]
Consequently, as $\widetilde{\s{C}}(P) = \widetilde{\s{C}(P)}
= \mathbb{P}(C(\s{C}(P)))$, we immediately check that
\[
\s{D}\mathit{el}(P) = \pi_{d+1}\circ\widetilde{\pi}_N(\widetilde{\s{C}}(P))
= \pi_{d+1}\circ \widehat{\pi}_N(C(\s{C}(P)))
= \pi_{d+1}\circ \widehat{\pi}_N(\mathrm{cone}(\widehat{\s{C}(P)})),
\]
where $\widehat{\s{C}(P)} = \{\widehat{u} \mid u \in \s{C}(P)\}$ and
$\widehat{u} = (u, 1)$.

\medskip
This suggests defining the map,
$\mapdef{\pi_N}{(\reals^{d+1} - H_{d+1})}{\reals^d}$, 
by 
\[
\pi_N = \pi_{d+1}\circ \widehat{\pi}_N \circ i_{d+2},
\]
which is explicity given by
\[
\pi_N(x_1, \ldots, x_d, x_{d+1}) = 
\frac{1}{1 - x_{d+1}}(x_1, \ldots, x_d).
\]
Then, as $\s{C}(P)$ is strictly below the hyperplane $H_{d+1}$,
we have
\[
\s{D}\mathit{el}(P) = \pi_{d+1}\circ\widetilde{\pi}_N(\widetilde{\s{C}}(P))
= \pi_N(\s{C}(P)).
\]

\medskip
First, note that $\s{D}\mathit{el}(P) = 
\pi_{d+1}\circ\widetilde{\pi}_N(\widetilde{\s{C}}(P))$ is indeed
a polyhedral complex whose geometric realization is the convex hull,
$\mathrm{conv}(P)$, of $P$. Indeed, by Proposition \ref{Thetaprop2},
the images of the facets of $\s{C}(P)$ are 
polytopes and when any two such polytopes meet,
they meet along a common face. Furthermore, if 
$\mathrm{dim}(\mathrm{conv}(P)) = m$,
then $\s{D}\mathit{el}(P)$ is pure $m$-dimensional. 
First,  $\s{D}\mathit{el}(P)$ contains
at least one $m$-dimensional cell. If  $\s{D}\mathit{el}(P)$
was not pure, as the complex is connected there would be
some cell, $\sigma$, of dimension $s < m$
meeting some other cell, $\tau$,
of dimension $m$ along a common face of dimension at most $s$ and
because $\sigma$ is not contained in any face of dimension $m$,
no facet of $\tau$ containing $\sigma\cap \tau$ can be adjacent to
any cell of dimension $m$ and so, $\s{D}\mathit{el}(P)$ would not be convex,
a contradiction.

\medskip
For any polytope, $P \subseteq \eucreal^d$, 
given any point, $x$, not in $P$, recall that
a facet, $F$, of $P$ is {\it visible from $x$\/} iff for every
point, $y\in F$, the line through $x$ and $y$ intersects $F$ only in $y$.
If $\mathrm{dim}(P) = d$, this is equivalent to saying that $x$ and
the interior of $P$ are strictly separated by the supporting hyperplane of $F$.
Note that if  $\mathrm{dim}(P) < d$, it possible that
every facet of $P$ is visible from $x$.

\medskip
Now, assume that $P \subseteq \eucreal^d$ is a polytope with nonempty interior. 
We say that a facet,
$F$, of $P$ is a {\it lower-facing facet\/} of $P$ iff
the unit normal to the supporting hyperplane
of $F$ pointing towards the interior of $P$ has non-negative
$x_{d+1}$-coordinate. A facet, $F$, that is not lower-facing is called an
{\it upper-facing facet\/} (Note that in this case
the $x_{d+1}$ coordinate of
the unit normal to the supporting hyperplane
of $F$ pointing towards the interior of $P$ is strictly negative).

\medskip
Here is a convenient way to characterize lower-facing
facets.

\begin{prop}
\label{lowerfacep1}
Given any polytope, $P \subseteq \eucreal^d$, with nonempty interior,
for any point, $c$, on the $O x_{d}$-axis, if $c$ lies strictly above
all the intersection points of the $O x_{d}$-axis with 
the supporting hyperplanes of all the upper-facing facets of $F$, then
the lower-facing facets of $P$ are exactly the facets not visible from $c$.
\end{prop}

\proof
Note that the intersection points of the $O x_{d}$-axis with 
the supporting hyperplanes of all the upper-facing facets of $P$
are strictly above the intersection points of the $O x_{d}$-axis with 
the supporting hyperplanes of all the lower-facing facets.
Suppose $F$ is visible from $c$. 
Then, $F$ must not be lower-facing
as otherwise,  for any $y\in F$, the line
through $c$ and $y$ has to intersect some upper-facing facet
and $F$ is not be visible from $c$, a contradiction. 

\medskip
Now, as $P$ is the intersection of the closed half-spaces
determined by the supporting hyperplanes of its facets, by
the definition of an upper-facing facet, any point, $c$,
on the $O x_{d}$-axis that lies strictly above the 
the intersection points of the $O x_{d}$-axis with 
the supporting hyperplanes of all the upper-facing facets of $F$
has the property that $c$ and the interior of $P$ are strictly separated
by all these supporting hyperplanes. Therefore, all
the upper-facing facets of $P$ are visible from $c$.
It follows that the facets visible from $c$ are exactly the upper-facing
facets, as claimed.
$\bigsquare$

\medskip
We will also need the following fact when $\mathrm{dim}(P) < d$.

\begin{prop}
\label{lowerfacep2}
Given any polytope, $P \subseteq \eucreal^d$, 
there is a point, $c$, on the $O x_{d}$-axis, 
such that for all points, $x$, on the $O x_{d}$-axis and above $c$, 
the set of facets of $\mathrm{conv}(P \cup \{x\})$ not
containing $x$ is identical. Moreover, the set of facets of $P$
not visible from $x$ is the set of facets of 
$\mathrm{conv}(P \cup \{x\})$ that do not contain $x$. 
\end{prop}

\proof
If $\mathrm{dim}(P) = d$ then pick any $c$ on 
the $O x_{d}$-axis above the intersection points of the $O x_{d}$-axis with 
the supporting hyperplanes of all the upper-facing facets of $F$.
Then, $c$ is in general position w.r.t. $P$ in the sense
that $c$ and any $d$ vertices of $P$ do not lie in a common hyperplane.
Now, our result follows by lemma 8.3.1 of Boissonnat and
Yvinec \cite{Boissonnat}.  If  $\mathrm{dim}(P) < d$,
consider the affine hull of $P$ with the $O x_{d+1}$ -axis
and use the same argument.
$\bigsquare$

\begin{cor}
\label{lowerfacep3}
Given any polytope, $P \subseteq \eucreal^d$, with nonempty interior,
there is a point, $c$, on the $O x_{d}$-axis, so that
for all $x$ on the $O x_{d}$-axis and above $c$, 
the lower-facing facets of $P$ are exactly the facets 
of $\mathrm{conv}(P \cup \{x\})$ that do not contain $x$. 
\end{cor}

\medskip
As usual, let $e_{d+1} = (0, \ldots, 0, 1)\in \reals^{d+1}$.

\begin{thm}
\label{Delaunayp1}
Given any set of points, $P = \{p_1, \ldots, p_n\} \subseteq 
\eucreal^d$, let $\s{D}'(P)$  denote the polyhedron
$\mathrm{conv}(l(P)) + \mathrm{cone}(e_{d+1})$ 
and let $\widetilde{\s{D}}'(P)$ be the projective
completion of $\s{D}'(P)$. Also, let $\s{C}'(P)$
be the polyhedral complex consisting of the bounded facets
of the polytope  $\s{D}'(P)$ and
let $\widetilde{\s{C}}'(P)$ be the
projective completion of $\s{C}'(P)$.
Then
\[
\theta(\widetilde{\s{D}}(P)) = \widetilde{\s{D}}'(P)
\quad\hbox{and}\quad
\theta(\widetilde{\s{C}}(P)) = \widetilde{\s{C}}'(P).
\]
Furthermore, if $\s{D}\mathit{el}'(P) = 
\pi_{d+1}\circ \widetilde{p}_{d+1}(\widetilde{\s{C}}'(P)) 
= p_{d+1}(\s{C}'(P))$ is 
the ``standard'' Delaunay complex of $P$, that is, the
orthogonal projection of $\s{C}'(P)$ onto $\eucreal^d$, then
\[
\s{D}\mathit{el}(P) = \s{D}\mathit{el}'(P).
\]
Therefore, the two notions of a Delaunay complex agree.
If $\mathrm{dim}(\mathrm{conv}(P)) = d$, then the bounded facets of 
$\mathrm{conv}(l(P)) + \mathrm{cone}(e_{d+1})$  are precisely
the lower-facing facets of $\mathrm{conv}(l(P))$.
\end{thm}

\proof
Recall that
\[
\s{D}(P) = \mathrm{conv}(\tau_N(P)\cup \{N\})
\]
and $\widetilde{\s{D}}(P) = \mathbb{P}(C(\s{D}(P)))$ is the projective
completion of $\s{D}(P)$. If we write $\widehat{\tau_N(P)}$ for \\
$\{\widehat{\tau_N(p_i)} \mid p_i \in P\}$, then
\[
C(\s{D}(P)) = \mathrm{cone}(\widehat{\tau_N(P)}\cup \{\widehat{N}\}).
\]
By definition, we have
\[
\theta(\widetilde{\s{D}}) = \mathbb{P}(\widehat{\theta}(C(\s{D}))).
\]
Now, as $\widehat{\theta}$ is linear,
\[
\widehat{\theta}(C(\s{D})) = \widehat{\theta}
(\mathrm{cone}(\widehat{\tau_N(P)}\cup \{\widehat{N}\}))
= \mathrm{cone}(\widehat{\theta}(\widehat{\tau_N(P)})\cup 
\{\widehat{\theta}(\widehat{N})\}).
\]
We claim that
\begin{eqnarray*}
\mathrm{cone}(\widehat{\theta}(\widehat{\tau_N(P)})\cup 
\{\widehat{\theta}(\widehat{N})\}) 
& = & \mathrm{cone}(\widehat{l(P)} \cup \{(0, \ldots, 0, 1, 1)\})  \\
& = & C(\s{D}'(P)),
\end{eqnarray*}
where
\[
\s{D}'(P) = \mathrm{conv}(l(P)) + \mathrm{cone}(e_{d+1}).
\]
Indeed,
\[
\widehat{\theta}(x_1, \ldots, x_{d+2}) = 
(x_1, \ldots, x_{d}, x_{d+1} + x_{d+2},  x_{d+2} - x_{d+1}),
\]
and for any $p_i = (x_1, \ldots, x_d)\in P$,
\begin{eqnarray*}
\widehat{\tau_N(p_i)} 
& = &
\left(
\frac{2x_1}{\sum_{i = 1}^d x_i^2 + 1}, \ldots, 
\frac{2x_d}{\sum_{i = 1}^d x_i^2 + 1},  
\frac{\sum_{i = 1}^d x_i^2 - 1}{\sum_{i = 1}^d x_i^2 + 1},  1
\right) \\
& = &
\frac{1}{\sum_{i = 1}^d x_i^2 + 1}
\left(
2x_1, \ldots, 2x_d, 
\sum_{i = 1}^d x_i^2 - 1, \sum_{i = 1}^d x_i^2 + 1
\right),
\end{eqnarray*}
so we get
\begin{eqnarray*}
\widehat{\theta}(\widehat{\tau_N(p_i)}) 
& = & 
\frac{2}{\sum_{i = 1}^d x_i^2 + 1}
\left(x_1, \ldots, x_d, \sum_{i = 1}^d x_i^2, 1\right) \\
& = &  \frac{2}{\sum_{i = 1}^d x_i^2 + 1}\, \widehat{l(p_i)}.
\end{eqnarray*}
Also, we have
\[
\widehat{\theta}(\widehat{N}) = 
\widehat{\theta}(0, \ldots, 0, 1, 1) = 
(0, \ldots, 0, 2, 0) = 2\widehat{e_{d+1}},
\]
and by definition of $\mathrm{cone}(-)$
(scalar factors are irrelevant), we get
\[
\mathrm{cone}(\widehat{\theta}(\widehat{\tau_N(P)})\cup 
\{\widehat{\theta}(\widehat{N})\}) = 
\mathrm{cone}(\widehat{l(P)} \cup \{(0, \ldots, 0, 1, 1)\})  =
C(\s{D}'(P)),
\]
with 
$\s{D}'(P) = \mathrm{conv}(l(P)) + \mathrm{cone}(e_{d+1})$,
as claimed. 
This proves that
\[
\theta(\widetilde{\s{D}}(P)) = \widetilde{\s{D}}'(P).
\]

\medskip
Now, it is clear that the facets of $\mathrm{conv}(\tau_N(P)\cup \{N\})$
that do not contain $N$ are mapped to the bounded facets of 
$\mathrm{conv}(l(P)) + \mathrm{cone}(e_{d+1})$,
since $N$ goes the point at infinity, so
\[
\theta(\widetilde{\s{C}}(P)) = \widetilde{\s{C}}'(P).
\]
As $\widetilde{\pi}_N =  \widetilde{p}_{d+1} \circ \theta$ by 
Proposition \ref{Thetaprop1},
we get
\[
\s{D}\mathit{el}'(P) = \pi_{d+1} \circ 
\widetilde{p}_{d+1}(\widetilde{\s{C}}'(P)) = 
 \pi_{d+1} \circ (\widetilde{p}_{d+1} \circ \theta)(\widetilde{\s{C}}(P)) = 
\pi_{d+1} \circ \widetilde{\pi}_N(\widetilde{\s{C}}(P)) = \s{D}\mathit{el}(P), 
\]
as claimed. Finally, if  $\mathrm{dim}(\mathrm{conv}(P)) = d$, 
then, by Corollary \ref{lowerfacep3}, we can pick a point, $c$,
on the $O x_{d+1}$-axis, so that
the facets of 
$\mathrm{conv}(l(P)\cup \{c\})$ that do not contain $c$ are precisely
the lower-facing facets of $\mathrm{conv}(l(P))$.
However, it is also clear that the facets of $\mathrm{conv}(l(P)\cup \{c\})$
that contain $c$ tend to the unbounded facets of 
$\s{D}'(P) = \mathrm{conv}(l(P)) + \mathrm{cone}(e_{d+1})$
when $c$ goes to $+\infty$.
$\bigsquare$

\medskip
We can also characterize when the Delaunay complex, $\s{D}\mathit{el}(P)$,
is simplicial. Recall  that we say that a set of points, 
$P\subseteq \eucreal^d$, is in {\it general position\/} iff
no $d + 2$ of the points in $P$ belong to a common
$(d - 1)$-sphere.

\begin{prop}
\label{Delaunayp2}
Given any set of points, $P = \{p_1, \ldots, p_n\} \subseteq 
\eucreal^d$, if $P$ is in general position, 
then the Delaunay complex, $\s{D}\mathit{el}(P)$, is
a pure simplicial complex.
\end{prop}

\proof
Let $\mathrm{dim}(\mathrm{conv}(P)) = r$. Then, $\tau_N(P)$
is contained in a $(r - 1)$-sphere of $S^d$, so we may assume
that $r = d$. Suppose $\s{D}\mathit{el}(P)$ has some facet, $F$, which
is not a $d$-simplex. If so, $F$ is the convex hull of
at least $d + 2$ points, $p_1, \ldots, p_k$
of $P$ and since $F = \pi_N(\widehat{F})$,  
for some facet, $\widehat{F}$,
of $\s{C}(P)$, we deduce that $\tau_N(p_1), \ldots, \tau_N(p_k)$
belong to the supporting hyperplane, $H$, of $\widehat{F}$.
Now, if $H$ passes through the north pole, then we know
that  $p_1, \ldots, p_k$ belong to some hyperplane of $\eucreal^d$,
which is impossible since $p_1, \ldots, p_k$ are the vertices of
a facet of dimension $d$. Thus, $H$ does not pass through $N$
and so, $p_1, \ldots, p_k$ belong to some $(d-1)$-sphere in $\eucreal^d$.
As $k \geq d+2$, this contradicts the assumption that the points in
$P$ are in general position.
$\bigsquare$

\remark
Even when the points in $P$ are in general position, the 
Delaunay polytope, $\s{D}(P)$, may not be a simplicial
polytope. For example, if $d + 1$ points belong to a hyperplane
in $\eucreal^d$, then
the lifted points belong to a hyperplane passing through
the north pole and these $d+1$  lifted points together with $N$
may form a non-simplicial facet. For example, consider the
polytope obtained by lifting  our original $d+1$ points on
a hyperplane, $H$, plus
one more point not in the the hyperplane $H$.

\medskip
We can also characterize the Voronoi diagram of $P$
in terms of the polar dual of  $\s{D}(P)$. 
Unfortunately, we can't simply take the polar dual, 
$\s{D}(P)^*$, of $\s{D}(P)$ and project it
using $\pi_N$ because some of the facets of $\s{D}(P)^*$
may intersect the hyperplane, $H_{d+1}$, and $\pi_N$ is undefined
on $H_{d+1}$. However, using projective completions,
we can indeed recover the Voronoi diagram of $P$.

\begin{defin}
\label{voroidef2}
{\em
Given any set of points, $P = \{p_1, \ldots, p_n\} \subseteq 
\eucreal^d$, 
the {\it  Voronoi polyhedron\/} associated with $P$
is the polar dual (w.r.t. $S^d\subseteq \reals^{d+1}$),
$\s{V}(P) = (\s{D}(P))^* \subseteq \reals^{d+1}$,
of the Delaunay polytope,
$\s{D}(P) = \mathrm{conv}(\tau_N(P)\cup \{N\})$.
The  {\it projective Voronoi polytope\/}, 
$\widetilde{\s{V}}(P) \subseteq \projr{d+1}$,  associated with $P$ 
is the projective completion of $\s{V}(P)$. 
The polyhedral complex, 
$\s{V}\mathit{or}(P) = \pi_{d+1}\circ\widetilde{\pi}_N(\widetilde{\s{V}}(P))
\subseteq \eucreal^d$, is the {\it Voronoi complex of $P$\/}
or {\it Voronoi diagram  of $P$\/}.
}
\end{defin}

\medskip
Given any set of points, $P = \{p_1, \ldots, p_n\} \subseteq 
\eucreal^d$, let $\s{V}'(P) = (\s{D}'(P))^*$ be the polar dual 
(w.r.t. $\s{P}\subseteq \reals^{d+1}$) of the ``standard''
Delaunay polyhedron defined in Theorem \ref{Delaunayp1}
and let 
$\widetilde{\s{V}}'(P) = \widetilde{\s{V}'(P)}\subseteq \projr{d}$ 
be its projective completion. It is not hard to check that
\[
p_{d+1}(\s{V}'(P)) = \pi_{d+1}\circ \widetilde{p}_{d+1}(\widetilde{\s{V}}'(P))
\]
is the ``standard'' Voronoi diagram, denoted $\s{V}\mathit{or}'(P)$.

\begin{thm}
\label{Voronoipol1}
Given any set of points, $P = \{p_1, \ldots, p_n\} \subseteq 
\eucreal^d$, we have
\[
\theta(\widetilde{\s{V}}(P)) = \widetilde{\s{V}}'(P)
\]
and
\[
\s{V}\mathit{or}(P) = \s{V}\mathit{or}'(P).
\]
Therefore, the two notions of Voronoi diagrams agree.
\end{thm}

\proof
By definition,
\[
\widetilde{\s{V}}(P) = \widetilde{\s{V}(P)} = 
\widetilde{(\s{D}(P))^*}
\]
and by Proposition \ref{commut1},
\[
\widetilde{(\s{D}(P))^*} = \left(\widetilde{\s{D}(P)}\right)^* =
(\widetilde{\s{D}}(P))^*,
\]
so
\[
\widetilde{\s{V}}(P) = (\widetilde{\s{D}}(P))^*.
\]
By Proposition \ref{projectualp1},
\[
\theta(\widetilde{\s{V}}(P)) = 
\theta((\widetilde{\s{D}}(P))^*) =
(\theta(\widetilde{\s{D}}(P)))^*
\]
and by Theorem \ref{Delaunayp1},
\[
\theta(\widetilde{\s{D}}(P)) = \widetilde{\s{D}}'(P),
\]
so we get
\[
\theta(\widetilde{\s{V}}(P)) = (\widetilde{\s{D}}'(P))^*.
\]
But, by  Proposition \ref{commut1} again,
\[
(\widetilde{\s{D}}'(P))^* = \left(\widetilde{\s{D}'(P)}\right)^*
= \widetilde{(\s{D}'(P))^*} = 
\widetilde{\s{V}'(P)} = \widetilde{\s{V}}'(P).
\]
Therefore,
\[
\theta(\widetilde{\s{V}}(P)) = \widetilde{\s{V}}'(P),
\]
as claimed.

\medskip
As $\widetilde{\pi}_N =  \widetilde{p}_{d+1} \circ \theta$ by 
Proposition \ref{Thetaprop1}, we get
\begin{eqnarray*}
\s{V}\mathit{or}'(P) & = &
\pi_{d+1}\circ\widetilde{p}_{d+1}(\widetilde{\s{V}}'(P))\\
& = & \pi_{d+1}\circ\widetilde{p}_{d+1}\circ \theta(\widetilde{\s{V}}(P))\\
& = & \pi_{d+1}\circ\widetilde{\pi}_{N}(\widetilde{\s{V}}(P))\\
& = & \s{V}\mathit{or}(P),
\end{eqnarray*}
finishing the proof.
$\bigsquare$

\medskip
We can also prove the proposition below which
shows directly that $\s{V}\mathit{or}(P)$ is 
the Voronoi diagram of $P$.
Recall that that $\widetilde{\s{V}}(P)$ is the
projective completion of $\s{V}(P)$. We observed in Section \ref{sec5e}
(see page \pageref{faceref}) that in the patch $U_{d+1}$,
there is a bijection between the faces of
$\widetilde{\s{V}}(P)$ and the faces of $\s{V}(P)$.
Furthermore, the projective completion, $\widetilde{H}$,
of every hyperplane, $H\subseteq \reals^d$, is also a hyperplane
and it is easy to see that if $H$ is tangent to  $\s{V}(P)$,
then $\widetilde{H}$ is tangent to $\widetilde{\s{V}}(P)$.

\begin{prop}
\label{Voronoiclx1}
Given any set of points, $P = \{p_1, \ldots, p_n\} \subseteq 
\eucreal^d$,
for every $p\in P$, if $F$ is the facet of $\s{V}(P)$
that contains $\tau_N(p)$, if $H$ is the tangent hyperplane at 
$\tau_N(p)$ to $S^d$ and if $F$ is cut out by the
hyperplanes $H, H_1, \ldots, H_{k_p}$, in the sense that
\[
F = (H\cap H_1)_- \cap \cdots \cap (H\cap H_{k_p})_-,
\] 
where $(H\cap H_i)_-$ denotes the closed half-space in $H$
containing $\tau_N(p)$ 
determined by $H\cap H_i$, then
\[
V(p) = \pi_{d+1}\circ \widetilde{\pi}_N(
\widetilde{H}\cap \widetilde{H}_1)_- \cap \cdots \cap 
\widetilde{\pi}_N(\widetilde{H}\cap \widetilde{H}_{k_p})_-
\]
is the Voronoi region of $p$
(where  $\pi_{d+1}\circ\widetilde{\pi}_N(\widetilde{H}\cap \widetilde{H}_i)_-$ 
is the closed half-space
containing $p$).
If $P$ is in general position, then $\s{V}(P)$
is a simple polyhedron (every vertex belongs to $d + 1$ facets).
\end{prop}

\proof
Recall that by Proposition \ref{Thetaprop1}, 
\[
\widetilde{\tau}_N\circ \psi_{d+1}  =  \psi_{d+2}\circ \tau_N.
\]
Each $H_i = T_{\tau_N(p_i)} S^d$ is the tangent hyperplane to $S^d$
at $\tau_N(p_i)$, for some $p_i\in P$.
Now, by definition of the projective completion,  the embedding,
$\s{V}(P)\longrightarrow \widetilde{\s{V}}(P)$,
is given by $a \mapsto \psi_{d+2}(a)$.
Thus, every point, $p\in P$, is mapped to the point 
$\psi_{d+2}(\tau_N(p)) = \widetilde{\tau_N}(\psi_{d+1}(p))$ and
we also have 
$\widetilde{H}_i =  T_{\widetilde{\tau}_N\circ \psi_{d+1}(p_i)} S^d$ and
$\widetilde{H} =  T_{\widetilde{\tau}_N\circ \psi_{d+1}(p)} S^d$.
By Proposition \ref{pimap1},
\[
\widetilde{\pi}_N (T_{\widetilde{\tau}_N\circ \psi_{d+1}(p)} S^d 
\cap T_{\widetilde{\tau}_N\circ \psi_{d+1}(p_i)} S^d) = 
\psi_{d+1}(H_{p, p_i})
\]
is the embedding of the
bisector hyperplane of $p$ and $p_i$ in $\projr{d}$,
so the first part holds. 

\medskip
Now, assume that some vertex,
$v\in \s{V}(P) = \s{D}(P)^*$,  belongs to $k \geq d + 2$
facets of $\s{V}(P)$. 
By polar duality, this means that the facet, $F$, 
dual of $v$ has  $k \geq d + 2$ vertices
$\tau_N(p_1), \ldots, \tau_N(p_k)$ of $\s{D}(P)$. We claim that
$\tau_N(p_1), \ldots, \tau_N(p_k)$ must belong to some hyperplane
passing through the north pole. Otherwise, 
$\tau_N(p_1), \ldots, \tau_N(p_k)$
would belong to a hyperplane not passing through the north pole
and so they would belong to 
a $(d - 1)$ sphere of $S^d$ and thus,
$p_1, \ldots, p_k$ would belong to a $(d - 1)$-sphere even though
$k \geq d + 2$, contradicting that $P$ is in general position.
But then, by polar duality, $v$ would be a point at infinity,
a contradiction.
$\bigsquare$

\medskip
Note that when $P$ is in general position, even though
the polytope, $\s{D}(P)$, may not be simplicial, 
its dual,  $\s{V}(P) =  \s{D}(P)^*$,
is a simple {\it polyhedron\/}. What is happening is that
$\s{V}(P)$ has unbounded faces which have ``vertices
at infinity'' that do not count! In fact,  
the faces of  $\s{D}(P)$ that fail to be simplicial
are those that are contained in some hyperplane through the north
pole. By polar duality, these faces correspond to
a vertex at infinity. Also, if  $m = \mathrm{dim}(\mathrm{conv}(P)) < d$,
then $\s{V}(P)$ may not have any vertices!

\medskip
We conclude our  presentation of
Voronoi diagrams and Delaunay triangulations
with a short section on applications.

\section[Applications]
{Applications of Voronoi Diagrams and Delaunay Triangulations}
\label{sec225}
The examples below are taken from O'Rourke \cite{ORourke}.
Other examples can be found in
Preparata and Shamos \cite{Preparata},
Boissonnat and Yvinec \cite{Boissonnat}, and
de Berg, Van Kreveld, Overmars, and Schwarzkopf \cite{Berg}.
\index{applications!of Voronoi diagrams}%
\index{applications!of Delaunay triangulations}%

\medskip
The first example is the {\it nearest neighbors\/} problem.
There are actually two subproblems: {\it Nearest neighbor queries\/}
and  {\it all nearest neighbors\/}.
\index{nearest neighbors problems}%

\medskip
The nearest neighbor queries problem
is as follows. Given a set $P$ of points
and a query point $q$, find the nearest neighbor(s) of $q$ in $P$.
This problem can be solved by computing the Voronoi diagram of $P$
and determining in which Voronoi region $q$ falls.
This last problem, called {\it point location\/}, has been heavily studied
(see O'Rourke \cite{ORourke}).
\index{point location problem}%
The all neighbors problem is as follows:
Given a set $P$ of points, find the nearest neighbor(s) to all
points in $P$. This problem can be solved by building a graph,
the {\it nearest neighbor graph\/}, for short {\it nng\/}.
The nodes of this undirected 
graph are the points in $P$, and there is an arc from 
$p$ to $q$ iff $p$ is a nearest neighbor of $q$ or vice versa.
Then it can be shown that this graph is contained in the Delaunay
triangulation of $P$.

\medskip
The second example  is the {\it largest empty circle\/}.
\index{largest empty circle}%
Some practical applications of this problem
are to locate a new store (to avoid competition),
or to locate a nuclear plant as far as possible 
from a set of towns.
More precisely, the problem is as follows.
Given a set $P$ of points, find a largest empty circle
whose center is in the (closed) convex hull 
of $P$, empty in that it contains no points from $P$ inside it,
and largest in the sense that there is no other circle with 
strictly larger radius. The Voronoi diagram of $P$ can be used
to solve this problem. It can be shown that if the center $p$
of a largest empty circle is strictly inside the convex hull of $P$,
then $p$ coincides with a Voronoi vertex.
However, not every Voronoi vertex is a good candidate.
It can also be shown that if  the center $p$
of a largest empty circle lies on the boundary of the convex hull of $P$,
then $p$ lies on  a Voronoi edge.

\medskip
The third example is the {\it minimum spanning tree\/}.
\index{minimum spanning tree}%
Given a graph $G$, a minimum spanning tree of $G$ is a subgraph of $G$
that is a tree, contains every vertex of the graph $G$,
and minimizes the sum of the lengths of the tree edges.
It can be shown that a minimum spanning tree is a subgraph
of the Delaunay triangulation of the vertices of the graph.
This can be used to improve algorithms for finding
minimum spanning trees, for example Kruskal's algorithm
(see O'Rourke \cite{ORourke}).

\medskip
We conclude by mentioning that Voronoi diagrams 
have applications to {\it motion planning\/}.
\index{motion!planning}%
For example, consider the problem of moving a disk on
a plane while avoiding a set of polygonal obstacles.
If we ``extend'' the obstacles by the diameter of the disk,
the problem reduces to finding a collision--free path between two points
in the extended obstacle space.
One needs to generalize the notion of a Voronoi diagram.
Indeed, we need to define the distance to an object,
and medial curves (consisting of points equidistant to
two objects) may no longer be straight lines.
A collision--free path with maximal clearance from the obstacles
can be found by moving along the edges of the generalized
Voronoi diagram.
This is an active area of research in robotics.
For more on this topic, see  O'Rourke \cite{ORourke}.

\bigskip
{\it Acknowledgement}.
I wish to thank Marcelo Siqueira for suggesting many improvements
and for catching many bugs with his ``eagle eye''.

\bibliography{../basicmath/algeom}

\begin{thebibliography}{10}

\bibitem{Alexandrov}
P.S. Alexandrov.
\newblock {\em Combinatorial Topology}.
\newblock Dover, first edition, 1998.
\newblock Three volumes bound as one.

\bibitem{AlonKalai}
Noga Alon and Gil Kalai.
\newblock A simple proof of the upper-bound theorem.
\newblock {\em European J. Comb.}, 6:211--214, 1985.

\bibitem{Barvinok}
Alexander Barvinok.
\newblock {\em A Course in Convexity}.
\newblock GSM, Vol. 54. AMS, first edition, 2002.

\bibitem{BayerLee}
Margaret~M. Bayer and Carl~W. Lee.
\newblock Combinatorial aspects of convex polytopes.
\newblock In P.M. Gruber and J.M. Wills, editors, {\em Handbook of Convex
  Geometry}, pages 485--534. Elsevier Science, 1993.

\bibitem{Berg}
M.~Berg, M.~Van~Kreveld, M.~Overmars, and O.~Schwarzkopf.
\newblock {\em Computational Geometry. Algorithms and Applications}.
\newblock Springer, first edition, 1997.

\bibitem{Berger90b}
Marcel Berger.
\newblock {\em G\'eom\'etrie 2}.
\newblock Nathan, 1990.
\newblock English edition: Geometry 2, Universitext, Springer Verlag.

\bibitem{BilleraBjorner}
J.~Billera, Louis and Anders Bj\"orner.
\newblock Faces numbers of polytopes and complexes.
\newblock In J.E. Goodman and Joe O'Rourke, editors, {\em Handbook of Discrete
  and Computational Geometry}, pages 291--310. CRC Press, 1997.

\bibitem{Boissonnat}
J.-D. Boissonnat and M.~Yvinec.
\newblock {\em G\'eom\'etrie Algorithmique}.
\newblock Ediscience International, first edition, 1995.

\bibitem{BourbakiEVT}
Nicolas Bourbaki.
\newblock {\em Espaces Vectoriels Topologiques}.
\newblock El\'ements de Math\'ematiques. Masson, 1981.

\bibitem{Brannan}
David~A. Brannan, Matthew~F. Esplen, and Jeremy~J. Gray.
\newblock {\em Geometry}.
\newblock Cambridge University Press, first edition, 1999.

\bibitem{Brown}
K.Q. Brown.
\newblock Voronoi diagrams from convex hulls.
\newblock {\em Inform. Process. Lett.}, 9:223--228, 1979.

\bibitem{Bruggesser}
Heinz Bruggesser and Peter Mani.
\newblock Shellable decompositions of cells and spheres.
\newblock {\em Math. Scand.}, 29:197--205, 1971.

\bibitem{Coxeterpoly}
H.S.M. Coxeter.
\newblock {\em Regular Polytopes}.
\newblock Dover, third edition, 1973.

\bibitem{Cromwell}
Peter Cromwell.
\newblock {\em Polyhedra}.
\newblock Cambridge University Press, first edition, 1994.

\bibitem{Dirichlet}
G.L. Dirichlet.
\newblock {\"Uber} die reduktion der positiven quadratischen formen mid drei
  unbestimmten ganzen zahlen.
\newblock {\em Journal f\"ur die reine und angewandte Mathematik}, 40:209--227,
  1850.

\bibitem{Eldel}
H.~Edelsbrunner and R.~Seidel.
\newblock Voronoi diagrams and arrangements.
\newblock {\em Discrete Computational Geometry}, 1:25--44, 1986.

\bibitem{Edelsbrunner}
Herbert Edelsbrunner.
\newblock {\em Geometry and Topology for Mesh Generation}.
\newblock Cambridge University Press, first edition, 2001.

\bibitem{Ewald}
G\"unter Ewald.
\newblock {\em Combinatorial Convexity and Algebraic Geometry}.
\newblock GTM No. 168. Springer Verlag, first edition, 1996.

\bibitem{Fultontoric}
William Fulton.
\newblock {\em Introduction to Toric Varieties}.
\newblock Annals of Mathematical Studies, No. 131. Princeton University Press,
  1997.

\bibitem{Gallbook2}
Jean~H. Gallier.
\newblock {\em {Geometric Methods and Applications, For Computer Science and
  Engineering}}.
\newblock TAM, Vol. 38. Springer, first edition, 2000.

\bibitem{Gilbert}
E.N. Gilbert.
\newblock Random subdivisions of space into crystals.
\newblock {\em Annals of Math. Stat.}, 33:958--972, 1962.

\bibitem{DiscGeomHand}
Jacob~E. Goodman and Joseph O'Rourke.
\newblock {\em Handbook of Discrete and Computational Geometry}.
\newblock CRC Press, first edition, 1997.

\bibitem{Graham}
R.~Graham and F.~Yao.
\newblock A whirlwind tour of computational geometry.
\newblock {\em American Mathematical Monthly}, 97(8):687--701, 1990.

\bibitem{Grunbaum}
Branko Gr\"unbaum.
\newblock {\em Convex Polytopes}.
\newblock GTM No. 221. Springer Verlag, second edition, 2003.

\bibitem{Lang96}
Serge Lang.
\newblock {\em Real and Functional Analysis}.
\newblock GTM 142. Springer Verlag, third edition, 1996.

\bibitem{Lax}
Peter Lax.
\newblock {\em Functional Analysis}.
\newblock Wiley, first edition, 2002.

\bibitem{McMullen71}
Peter McMullen.
\newblock The maximum number of faces of a convex polytope.
\newblock {\em Mathematika}, 17:179--184, 1970.

\bibitem{Munkresalg}
James~R. Munkres.
\newblock {\em Elements of Algebraic Topology}.
\newblock Addison-Wesley, first edition, 1984.

\bibitem{ORourke}
Joseph O'Rourke.
\newblock {\em Computational Geometry in C}.
\newblock Cambridge University Press, second edition, 1998.

\bibitem{Preparata}
F.P. Preparata and M.I. Shamos.
\newblock {\em Computational Geometry: An Introduction}.
\newblock Springer Verlag, first edition, 1988.

\bibitem{Risler92}
J.-J. Risler.
\newblock {\em Mathematical Methods for CAD}.
\newblock Masson, first edition, 1992.

\bibitem{Rockafellar}
R.~Tyrrell Rockafellar.
\newblock {\em Convex Analysis}.
\newblock Princeton Landmarks in Mathematics. Princeton University Press, 1970.

\bibitem{Samuel}
Pierre Samuel.
\newblock {\em Projective Geometry}.
\newblock Undergraduate Texts in Mathematics. Springer Verlag, first edition,
  1988.

\bibitem{Seidel95}
Raimund Seidel.
\newblock The upper-bound theorem for polytopes: an easy proof of its
  asymptotic version.
\newblock {\em Comput. Geometry: Theory and Applications}, 5:115--116, 1995.

\bibitem{Stanley}
Richard~P. Stanley.
\newblock {\em Combinatorics and Commutative Algebra}.
\newblock Progress in Mathematics, No. 41. Birkh\"auser, second edition, 1996.

\bibitem{Stolfi}
J.~Stolfi.
\newblock {\em Oriented Projective Geometry}.
\newblock Academic Press, first edition, 1991.

\bibitem{Sturmfels96}
Bernd Sturmfels.
\newblock {\em Gr\"obner Bases and Convex Polytopes}.
\newblock ULS, Vol. 8. AMS, first edition, 1996.

\bibitem{Thomas}
Rekha~R. Thomas.
\newblock {\em Lectures in Geometric Combinatorics}.
\newblock STML, Vol. 33. AMS, first edition, 2006.

\bibitem{Thurston1}
P.~Thurston, William.
\newblock {\em Three-Dimensional Geometry and Topology, Vol. 1}.
\newblock Princeton University Press, first edition, 1997.
\newblock Edited by Silvio Levy.

\bibitem{Tisseron}
Claude Tisseron.
\newblock {\em G\'eom\'etries affines, projectives, et euclidiennes}.
\newblock Hermann, first edition, 1994.

\bibitem{Valentine}
Frederick~A. Valentine.
\newblock {\em Convex Sets}.
\newblock McGraw--Hill, first edition, 1964.

\bibitem{Voronoi}
M.G. Voronoi.
\newblock Nouvelles applications des param\`etres continus \`a la th\'eorie des
  formes quadratiques.
\newblock {\em J. Reine u. Agnew. Math.}, 134:198--287, 1908.

\bibitem{Ziegler97}
Gunter Ziegler.
\newblock {\em Lectures on Polytopes}.
\newblock GTM No. 152. Springer Verlag, first edition, 1997.

\end{thebibliography}
\bibliographystyle{plain} 
\end{document}